\documentclass[11pt]{article}
\usepackage[a4paper,margin=22mm]{geometry}
\usepackage{amsmath,amssymb,amsthm,mathtools,booktabs,array,graphicx,tikz}
\usepackage{longtable}
\usetikzlibrary{arrows.meta}
\usepackage[hidelinks,bookmarksnumbered=true,bookmarksdepth=2,hypertexnames=false]{hyperref}
\hypersetup{pdftitle={Quantitative stability for the polygonal Szegő–Weinberger inequality},pdfauthor={Hedong Wang}}
\newtheorem{theorem}{Theorem}[section]
\newtheorem{lemma}[theorem]{Lemma}
\newtheorem{corollary}[theorem]{Corollary}
\newtheorem{proposition}[theorem]{Proposition}
\theoremstyle{remark}
\newtheorem{remark}[theorem]{Remark}

\newcommand{\sinc}{\operatorname{sinc}}
\let\mathscr\mathcal
\newsavebox{\matbox}
\newcommand{\fitmath}[2][.98\linewidth]{\sbox{\matbox}{$\displaystyle #2$}%
\ifdim\wd\matbox>#1\relax\resizebox{#1}{!}{\usebox{\matbox}}\else\usebox{\matbox}\fi}
\allowdisplaybreaks[2]
\title{Quantitative stability for the polygonal\\Szeg\H{o}--Weinberger inequality}
\author{Hedong Wang}
\date{29 November 2026}
\begin{document}
\maketitle
\begin{abstract}
Let $P$ be a simple planar polygon with at most $n$ sides, and let $R_n$
be a regular $n$-gon of the same area. We measure the distance from $P$
to the regular $n$-gon by minimizing $|P\mathbin{\triangle}R_n|/|P|$
over all positions and orientations of $R_n$. We show that $|P|(\mu_1(R_n)-\mu_1(P))$ is at least an
explicit positive constant times the square of this distance for every
$n\ge3$. The exponent two is optimal for each fixed $n\ge4$. For
triangles the optimal exponent is one. For $n\ge10$, we first show that
a sufficiently small area-normalized deficit forces $P$ to have $n$ strictly
convex vertices and bounds the sector angles and vertex radii about a
suitable center. Under these bounds our stability constant is of order
$n^{-1}$. For all polygons with at most $n$ sides we obtain a constant of
order $n^{-5}$.
\end{abstract}
\section{Introduction}\label{sec:intro}
For a bounded polygonal domain $P$ in the plane, write $\mu_1(P)$ for the
first positive Neumann eigenvalue and $F(P)$ for its product with the area. Thus
\[
 \mu_1(P)=\inf\left\{\frac{\int_P |\nabla u|^2}{\int_P u^2}:
 u\in H^1(P),\ u\not\equiv0,\ \int_Pu=0\right\},\qquad
 F(P)=|P|\mu_1(P).
\]
Under a dilation by $s>0$, area is multiplied by $s^2$ and the eigenvalue by $s^{-2}$. Thus $F$ depends only on the similarity class. Write $R_n$ for any regular $n$-gon.

Szeg\H{o} proved that the disk maximizes the first positive Neumann eigenvalue
among simply connected planar domains of prescribed area \cite{Szego}.
Weinberger extended the comparison to domains of prescribed volume in
$\mathbb R^N$, with the ball as maximizer \cite{Weinberger}.
Brasco and Pratelli showed that the normalized deficit from the ball controls
the square of the Fraenkel asymmetry, the least relative symmetric difference
from an equal-volume ball, and that the power two is sharp \cite{BP}.

For each $n\ge4$, Laugesen asked whether the regular $n$-gon maximizes
$|P|\mu_1(P)$ among $n$-gons \cite[Problem~(1), p.~405]{ABLW}.
Laugesen and Siudeja later stated this as the polygonal
Szeg\H{o}--Weinberger conjecture \cite[Conjecture~6.36]{LSH}.

For the polygonal problem define the equal-area asymmetry by
\[
 \mathcal A_n(P)=\inf_{R}\frac{|P\mathbin{\triangle}R|}{|P|},
\]
where $R$ ranges over all regular $n$-gons of area $|P|$.
\begin{theorem}[Quantitative stability]\label{thm:stability}
For every integer $n\ge3$ and every simple polygon $P$ of positive area with at most $n$ sides,
\[
 F(R_n)-F(P)\ge c_n\mathcal A_n(P)^2,\qquad
 c_n=\begin{cases}
 F(R_3)/144,& n=3,\\
 1/60000,& n=4,\\
 1/(10^6 n^2),&5\le n\le9,\\
 1/(60n^5),& n\ge10 .
 \end{cases}
\]
\end{theorem}
In particular, the bound holds with $1/(8000n^5)$ for every $n\ge3$. For triangles we also prove the linear estimate
$D_3\ge F(R_3)\mathcal A_3/200$, where $D_n=F(R_n)-F(P)$.
Proposition~\ref{prop:optimal-power} determines the optimal power for each fixed $n$.
The underlying extremal result is the following.
\begin{theorem}\label{thm:main}
For each integer $n\ge3$ and each positive-area simple planar polygon $P$ with at most $n$ sides,
\[
                   F(P)\le F(R_n).
\]
Equality holds if and only if $P$ is similar to $R_n$.
\end{theorem}
For the first Dirichlet eigenvalue, Indrei proved a sharp stability
inequality for triangles \cite{Indrei}, and Cheng, Gui, Hu and Li
used quantitative Faber--Krahn stability to localize minimizing
convex polygons with many sides \cite{CGHL}. For the present Neumann
problem, consider a polygon with at most $n$ sides, where $n\ge5$.
If the deficit is sufficiently small, the polygon has $n$ strictly
convex vertices. Joining these vertices to a suitable interior point
divides the polygon into triangles. We call this collection of triangles
a fan. The deficit also bounds the angles at the center of the fan and
the distances from the center to the vertices. We refer to the permitted
ranges as the geometric window. Within these ranges we estimate the
deficit in terms of the deviations of the angles and radii from their
values on the regular polygon (Section~\ref{sec:stability}).

We note that Luo and Xie \cite{LX} proved the extremal comparison stated
in Theorem~\ref{thm:main}. The author independently developed the extremal argument
in an earlier draft.
Their paper was submitted to arXiv on 19 November 2026 and appeared in the
math.AP announcement of 22 November (US Eastern time), listed under
23 November.
By 20 November (US Eastern time), an unpublished draft stated the extremal
comparison for every side count and used computational certificates for
several estimates. Both proofs
transport, by the affine map of each fan sector, the difference between a
regular-polygon eigenfield and a radial Bessel field \cite[Section~4.2]{LX}.
Section~\ref{IV:sec:background} defines the projected field and
Lemma~\ref{IV:lem:projected-fan-identities} proves the identities used here.
Afterwards the computer-assisted steps were replaced throughout by displayed
estimates, and the bounds of
Proposition~\ref{prop:quantitative-window-entry} and the stability estimates
of Theorem~\ref{thm:stability} were obtained; none of this relies on \cite{LX}.

Our proof of the extremal comparison uses the triangle theorem of
Laugesen and Siudeja \cite{LS}. Endo and Osting \cite{EO} proved the
comparison for convex quadrilaterals. We give an affine proof of that
case in Section~\ref{I:sec:all-real-quads} and treat quadrilaterals with
a flat or reflex corner in Section~\ref{I:sec:lowertri}.

To bound the eigenvalue of a polygon with fewer than $n$ strictly
convex vertices, we use the author's shifted circle-mean estimate for signed
turning angles. Lemma~\ref{IV:lem:signed-circle} applies it with the center
chosen so that the Neumann trial field has mean zero.
Theorem~\ref{IV:thm:corner-gap} gives the resulting spectral gap.
\subsection*{Structure of the proof}
We first bound the eigenvalue when the polygon has fewer than $n$
strictly convex vertices. For a strictly convex $n$-gon we then choose
a trial function according to the number of sides, as shown below.
\begin{figure}[ht]\centering
\begin{tikzpicture}[x=1cm,y=1cm,block/.style={draw,rounded corners=2pt,align=center,inner sep=4pt,font=\small},dep/.style={-{Stealth[length=1.6mm]},semithick}]
\node[block] (c) at (0,1.4) {Fewer than $n$ strictly convex vertices\\Sections~\ref{I:sec:five-six-regular-data}--\ref{I:sec:five-six-thick}, Theorem~\ref{IV:thm:corner-gap}};
\node[block] (g) at (7.2,1.4) {Uniform sector estimates\\$0<h\le\pi/6$};
\node[block] (s) at (-2.1,-.2) {Affine and centroid trials\\$3\le n\le9$};
\node[block] (d) at (3,-.2) {Separate sector estimates\\$n=10,11$};
\node[block] (t) at (8.2,-.2) {Projected disk field\\$n\ge12$};
\node[block] (m) at (3,-1.6) {Extremal comparison, Theorem~\ref{thm:main}};
\node[block] (p) at (3,-2.6) {Angle and radius bounds for small deficits, Section~\ref{sec:stability}};
\node[block] (z) at (3,-3.6) {Quantitative stability, Theorem~\ref{thm:stability}};
\draw[dep] (m)--(p);\draw[dep] (p)--(z);
\draw[dep] (c)--(s);\draw[dep] (c)--(d);\draw[dep] (g)--(t);
\draw[dep] (c.south east)--(t.north west);\draw[dep] (s)--(m);\draw[dep] (d)--(m);\draw[dep] (t)--(m);
\end{tikzpicture}
\caption{Sector comparisons and their quantitative refinement.}\label{fig:deps}
\end{figure}
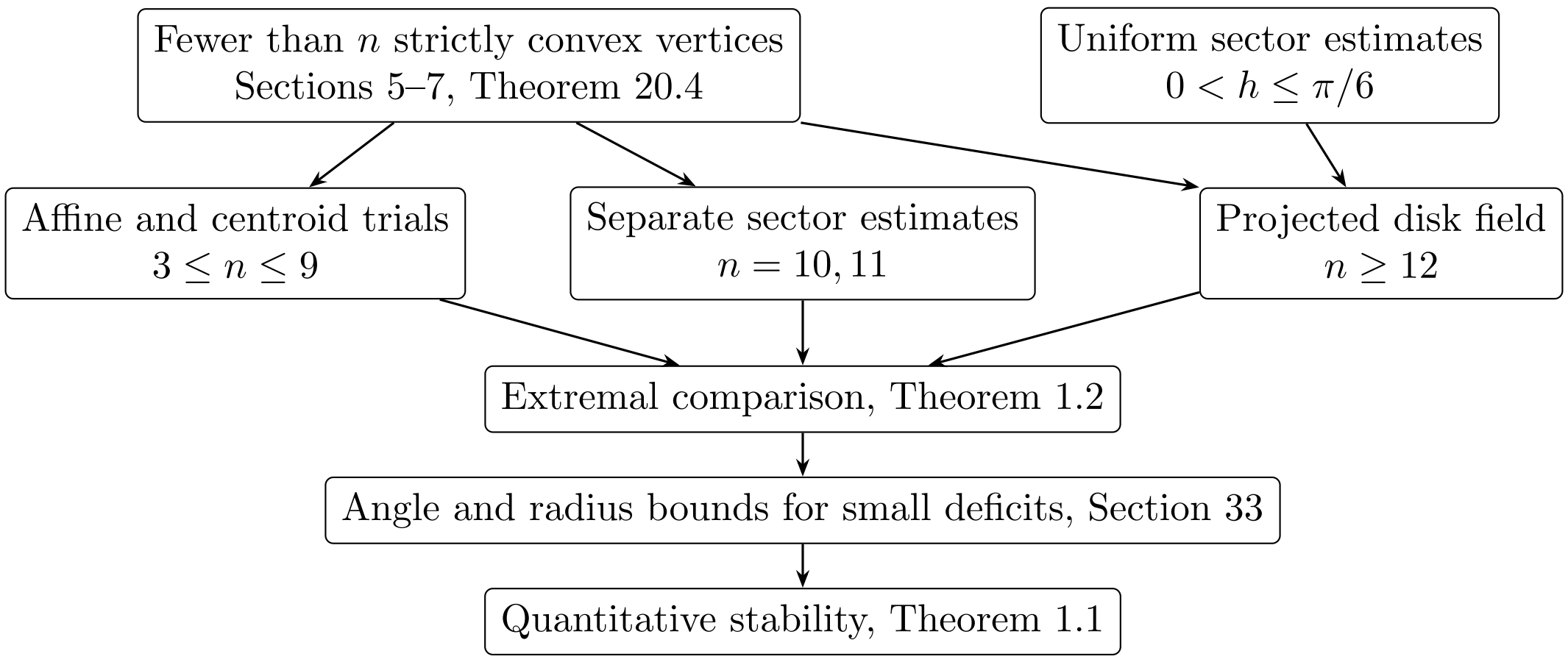
\paragraph{Angles and scales.}
When the center is the centroid, we write $a_i$ for half the angle of
sector $i$ and set $b=\pi/n$ and $x_i=\tan(a_i)/\tan b$. The angles
give $\sum_i a_i=nb$, so $\sum_i\arctan(x_i\tan b)=nb$. For the
projected disk field we use the ratio $q_i=a_i/b$ of the sector
half-angle to the regular half-angle, which gives $\sum_iq_i=n$.
We use the term normalized aperture for $x_i$ in the centroid construction and for $q_i$ in the projected-field construction. Both are one on a regular sector. The skew of a sector
measures the difference between the distances from the center to its
two vertices.
We normalize $R_n$ to have circumradius one for the centroid trials
and area $\pi$ for the conformal radius estimates. We state the change of scale whenever it affects a Rayleigh quotient.
\paragraph{Why the estimate improves as $n$ increases.}
Write $h=2\pi/n$ and let $q$ be the sector angle divided by $h$.
With both endpoint radii fixed at one, let $\Delta$ be the difference
between the tangent line at the regular angle and the angular part of
the trial functional. Let $C$ be the increase in the functional allowed
by varying the two radii. For $n\ge12$, Sections~\ref{mid-angular}--\ref{mid-radial}
give the following bounds, with $F=q(1+q)^2$ and $q\ne1$:
\[
 \frac{C}{h^4(q-1)^2}\le .0123\,hF,\qquad
 \frac{\Delta}{h^4(q-1)^2}\ge .002+.0123\,hF .
\]
After division by $h^4(q-1)^2$, the angular gain is at least $.002$
and the upper bound on the radial cost tends to zero with $h$, uniformly
over the allowed angles. At $q=1$ the denominator vanishes. We then use the second
derivatives in the two radii to prove strict inequality unless both radii
are one. For $n=10,11$ we estimate the same functional separately on
each side of $q=1$. Appendix~\ref{sec:verification} lists the numerical
bounds used in the proof.

\medskip
\begin{center}
\begin{tabular}{@{}lr@{}}\toprule
Proof step & Page\\\midrule
Trial field on the centroid fan and its sector bounds & \pageref{I:sec:fan}\\
Moments of the trial field for $n=7,8,9$ & \pageref{seven-moments}\\
Ordering of the Fourier coefficients for $n=7,8,9$ & \pageref{seven-fourier}\\
Sector comparison for $n\ge10$ & \pageref{prop:common-projected-comparison}\\
Bound for fewer than $n$ strictly convex vertices & \pageref{IV:thm:corner-gap}\\
Global quantitative stability & \pageref{sec:stability}\\
Optimal exponent of the asymmetry & \pageref{prop:optimal-power}\\
Numerical assertions and margins & \pageref{sec:verification}\\\bottomrule
\end{tabular}
\end{center}

\part{Triangles, quadrilaterals, pentagons and hexagons}\label{part:small}
\section{Triangles and quadrilaterals with a flat or reflex corner}\label{I:sec:lowertri}
Write $F(P)=|P|\mu_1(P)$. For a triangle of perimeter $p$, a
theorem of Laugesen and Siudeja gives
$\mu_1(P)p^2\le16\pi^2$ \cite[Theorem 3.1, formula (3.3)]{LS}.
The triangle isoperimetric inequality $12\sqrt3\,|P|\le p^2$
follows from Heron's formula and the arithmetic--geometric mean
inequality on the three differences between the semiperimeter
and the sides. Hence
\[
 F(P)\le\frac{4\pi^2}{3\sqrt3}=F(R_3)<\pi^2=F(R_4).
\]
The triangle isoperimetric inequality is strict unless the triangle is equilateral, so the first inequality has the same equality condition.
The last equality is the first positive eigenvalue of the unit square,
obtained by separation of variables. A positive-area simple polygon
with at most three sides is a triangle. A flat vertex on a simple
quadrilateral may be deleted, giving the same triangular domain.

\paragraph{A reflex corner.}
Let $A=|P|$, let $\mu=\mu_1(P)$, and let
$\Sigma$ be the covariance matrix of a uniform point in $P$.
Testing all centered linear coordinates gives $\mu\Sigma\preceq I$.
If the two covariance eigenvalues are $\sigma_1\le\sigma_2$,
then $\mu\sigma_1\sigma_2\le\sigma_1$, so
$\mu\det(\Sigma)I\preceq\Sigma$. For any nonconstant affine
coordinate $X$ and cubic $f$ of positive variance the centered
test $f(X)$ gives
\[
 F(P)^2\le\frac{A^2}{\det\Sigma}\,
 \frac{\operatorname{Var}(X)\,\mathbb E[f'(X)^2]}
      {\operatorname{Var}(f(X))}.
\]
The right side is invariant under affine changes of coordinates
when $X$ is transported with the domain.

Write the vertices in boundary order as $O,B,D,C$, with $D$ reflex.
It lies inside the triangle $OBC$. In affine coordinates
$O=(0,0), B=(1,0), C=(0,1)$ we can therefore write
$D=d(1-s,s)$ with $0<d<1$ and, after interchanging $B,C$,
$0<s\le1/2$. The triangles $ODB,ODC$ have area fractions $s,1-s$
and $A=d/2$. A point chosen uniformly from the quadrilateral has the form $LD+MZ$, where $Z$
takes the values $B,C$ with these probabilities, and $(L,M)$
is independently uniform on the standard unit simplex. Write
$v=B-C$, $z_0=\mathbb E Z$. Using $\mathbb EL=\mathbb EM=1/3$,
$\mathbb E L^2=\mathbb E M^2=1/6$, $\mathbb E LM=1/12$ gives
\[
 18\Sigma=(D-z_0/2)(D-z_0/2)^T+3z_0z_0^T/4
             +3s(1-s)vv^T .
\]
Taking its determinant, using $\det(D,z_0)=d(1-2s)$,
$\det(D,v)=-d$, $\det(z_0,v)=-1$, yields
\[
 108\,\det\Sigma/A^2=1+4s(1-s)(1-d)/d^2 .
\]

Consider two affine coordinates. The first has values $(s,0,1)$
on $(O,B,C)$ and value $s$ at $D$. Its density on $[0,s]$ is
$2t/s$, and on $[s,1]$ it is $2(1-t)/(1-s)$ by triangular slicing.
The second, $1-x-y$, has vertex values $1,0,1-d$ in each triangle,
so also has triangular density, with mode $1-d$.
Choose the coordinate with the smaller mode $b=\min(s,1-d)$ and
put $z=b(1-b)$. Since $s(1-s)\ge z$ and
$(1-d)/d^2\ge b/(1-b)^2\ge z$,
\[
 108\,\det\Sigma/A^2\ge1+4z^2,\qquad 0<z\le1/4 .
\]
Integration of the triangular density gives
$\mathbb E X^j=2(1+b+\cdots+b^j)/[(j+1)(j+2)]$.
For $f(t)=3t^2-2t^3$ this yields
\[
 \operatorname{Var}(X)=(1-z)/18,\quad
 \mathbb E[f'(X)^2]=6(1+z+2z^2)/5,\quad
 \operatorname{Var}(f(X))=3(19-2z-25z^2-24z^3)/700 .
\]
Consequently $F(P)^2\le1680(1+z^2-2z^3)/D(z)$, where
\[
 D(z)=19-2z+51z^2-32z^3-100z^4-96z^5\ge19-2z+35z^2 .
\]
The last estimate uses $32z+100z^2+96z^3\le16$. The polynomial
$89(19-2z+35z^2)-1680(1+z^2)=11-178z+1435z^2$ is positive
by its negative discriminant. Thus $F(P)<\sqrt{89}<\pi^2$,
using $\pi>31/10$. This excludes equality at every reflex
quadrilateral.

\section{The square and all convex quadrilaterals}
\label{I:sec:all-real-quads}

\begin{theorem}\label{I:thm:all-real-quads}
Every strictly convex quadrilateral $P$ satisfies $F(P)\le\pi^2$. Equality holds precisely for squares.
\end{theorem}
Endo and Osting proved this comparison for convex quadrilaterals using a certified finite covering \cite{EO}. Their vertex coordinates separate the affine and bilinear
parameters. We use that separation, three elementary test functions
and a covariance estimate. We treat all terminating decimals below as exact rational numbers.

\begin{proof}
Write $p=\pi^2$ and use $3.14159<\pi<3.14160$, in particular $9.85<p<10$. Write $\operatorname{sinc}(v)=\sin(v)/v$ with value one at zero. Apply an affine map to the given strictly convex
quadrilateral so that its diagonals are the vectors $(1,1)$ and $(-1,1)$. After a translation, the resulting
quadrilateral $Q$ is the image of $[-1/2,1/2]^2$ under
\[
 x=u(1-2bv),\qquad y=v(1+2cu),\qquad J=1-2bv+2cu .
\]
Indeed the four vertices in cyclic order are those of the original
square plus alternating vectors $(-b,c)/2$ and $(b,-c)/2$.
The quadrilateral is strictly convex precisely when $|b|+|c|<1$, as the Jacobian at each corner shows. The map is one-to-one: for fixed $x$, $u=x/(1-2bv)$ and
$y$ increases in $v$ with derivative $J/(1-2bv)>0$ throughout the
admissible interval. We may therefore use this change of variables on $Q$, and integrating the Jacobian gives
$|Q|=1$. Axis reflections and an interchange allow $b\ge c\ge0$.
Put $r=b^2+c^2$. The two signed offsets of the diagonal intersection
from the midpoints, measured in half-diagonal units, have magnitudes
$|b+c|,|b-c|$: subtract the midpoints and resolve the difference
along the two diagonal vectors.

Scale the original quadrilateral to area one. It has the form
$TQ$ up to a translation with $|\det T|=1$. Write
\[
 M=T^{-1}T^{-\top}=
 \begin{pmatrix}m+u_0&v_0\\v_0&m-u_0\end{pmatrix},\qquad
 t=(u_0^2+v_0^2)^{1/2},\qquad m=\sqrt{1+t^2},\quad z=m-1 .
\]
In the rest of the proof, brackets denote averages with respect to the uniform probability measure on $Q$.
If $F(TQ)\ge p$, the symmetric bilinear form
\[
 B(f,g)=\langle\nabla f^{\top}M\nabla g\rangle-p\operatorname{Cov}(f,g)
\]
is positive semidefinite. Denote the form at $M=I$ by $B_I$.

\paragraph{An affine estimate when the diagonals meet far from a midpoint.}
In affine coordinates another model for the same domain has vertices
$(s\pm1,0),(0,d\pm1)$, where $|s|=b+c$ and $|d|=b-c$ after possibly
exchanging diagonals. Its area is two and
\[
 \Sigma=\frac1{18}\begin{pmatrix}3+s^2&-2sd\\-2sd&3+d^2\end{pmatrix}.
\]
Here $\Sigma$ is its covariance. For instance, a uniform point is
$(LU,HV)$, where $(L,H)$ is uniform on the probability simplex
$L,H\ge0,L+H\le1$. Independently $U=s\pm1$ with respective
probabilities $(1\pm s)/2$, and $V=d\pm1$ likewise. Integration
over the four triangles gives the formula. Apply the affine covariance inequality from the argument for
a reflex quadrilateral in Section~\ref{I:sec:lowertri}, which holds for any domain. In the axial model choose $Z=(x-s+1)/2$; its triangular
density on $[0,1]$ has mode $a=(1-s)/2$. Its $j$th moment is
$2(1+a+\cdots+a^j)/[(j+1)(j+2)]$. For $g(Z)=3Z^2-2Z^3$ the
last factor is
\[
 \frac{140}{9}\frac{1+q^2-2q^3}{19-2q-25q^2-24q^3},
 \qquad q=a(1-a).
\]
If $b+c\ge .85$, then $q=(1-s^2)/4<.07$ and $A^2/\det\Sigma=432/[3+s^2+(1-s^2)d^2]\le432/3.72$.
Replacing the numerator by $1+.07^2$ and the denominator by its
value at $.07$ bounds this factor by $.835$. Thus $F^2<(3600/31)(.835)<97<p^2$.
We may therefore assume $b+c<.85$.

\paragraph{Bounding the affine distortion of a competitor.}
The coordinate covariances of $Q$ are
\[
 \Sigma_Q=I/12+\frac1{36}
 \begin{pmatrix}3b^2-c^2&-2bc\\-2bc&3c^2-b^2\end{pmatrix}.
\]
They follow by integrating polynomials with weight $J$; the coordinate
means are $c/6,-b/6$. The least eigenvalue is at least
$1/12-r/36$, since the two eigenvalues of the displayed correction
matrix before division by $36$ are $r\pm2\sqrt{r^2-3b^2c^2}$.
Positive semidefiniteness on linear functions gives $M\succeq p\Sigma_Q$.
As $\det M=1$, its largest eigenvalue $m+t$ is bounded by
$k=12/[9.85(1-r/3)]$. In particular
\begin{equation}\label{I:eq:quad-t-affine}
 t\le (k-k^{-1})/2<.201+.405r=:T(r),\qquad 0\le r\le S=.7225 .
\end{equation}
The middle expression is convex in $r$ (a reciprocal plus a linear
function), so the last estimate follows at $0,S$.

Set $U=\sin(\pi x)$, $V=\sin(\pi y)$, $H=xy$.
One has $\langle\cos(2\pi x)\rangle=\langle\cos(2\pi y)\rangle=0$:
in the first identity the $2cu$ term drops out and integration in
$u$ gives $\sin(\pi(1-2bv))/\pi$, odd in $v$.
Write $I=\langle U\rangle$, $J_0=\langle V\rangle$.
Direct integration gives
\[
 I=2c\int_{-1/2}^{1/2}u\sin(\pi u)\operatorname{sinc}(\pi b u)\,du,
 \qquad -J_0=2b\int_{-1/2}^{1/2}u\sin(\pi u)\operatorname{sinc}(\pi c u)\,du .
\]
In particular $0\le I\le4c/p$ and $0\le-J_0\le4b/p$.
Since $\operatorname{sinc}(h)\ge1-h^2/6$ here, also
$-I J_0\ge(16bc/p^2)(1-pr/24)$.

The diagonal projections determine the remaining two sine moments.
For $d=b-c$, the three projection levels of $x+y$ are
$-1-d/2,d/2,1-d/2$, with the two middle vertices projecting
to the same point. By linearity of the lengths of the parallel
chords on either side, the density is triangular of height one.
The law of $x-y$ has the same description with $d$ replaced by
$-s$, $s=b+c$. Integrating its two linear pieces gives
$\langle\cos(\pi(x+y))\rangle=4A(d)/p$ and
$\langle\cos(\pi(x-y))\rangle=4A(s)/p$, where
\[
 A(h)=\frac{\cos(\pi h/2)}{1-h^2}
     =\prod_{j\ge1}\left(1-\frac{h^2}{(2j+1)^2}\right).
\]
For $0\le u\le v\le1$ the product gives
$0\le A(u)-A(v)\le\alpha(v^2-u^2)$ by telescoping factors bounded
by one. Here $\alpha=\sum_{j\ge1}(2j+1)^{-2}=p/8-1<.24$.
The addition formulas therefore give
\[
 C:=\langle\cos(\pi x)\cos(\pi y)\rangle
 =2[A(s)+A(d)]/p\ge(4/p)(1-\alpha r),
 \quad 0\le-\langle UV\rangle\le(8\alpha/p)bc .
\]
Since $8\alpha/p=1-8/p<.20$ and $32(1-pr/24)/p^2>.22$,
we have $0\le-\langle UV\rangle\le-2IJ_0$.
Thus the symmetric matrix with diagonal $I^2,J_0^2$ and
off-diagonal $I J_0-\langle UV\rangle$ is positive semidefinite
and has norm at most $I^2+J_0^2\le16r/p^2$.

Subtracting this matrix from the matrix of $B/p$ on $U,V$ leaves
the diagonal $(z\pm u_0)/2$ and off-diagonal $v_0 C$.
Its least eigenvalue is at most $z/2-Ct$, since $C\le4/p<1/2$.
Hence positivity of $B$ and $z\le t^2/2$ give
$t(C-t/4)\le16r/p^2$. As $4/p>.405$, $r<3/4$ and
$t<T(S)<1/2$, one has $C-t/4>.405(1-.24\cdot3/4)-1/8>.207$.
We conclude that
\begin{equation}\label{I:eq:quad-t-small}
                         t\le 4r/5 \qquad(0\le r\le S).
\end{equation}

\paragraph{The mixed term produced by a quadratic test.}
In $B_I(U,H)$ the centered mass contributes $pI\langle H\rangle
=-pbcI/12$. For the other terms put $A=1-2bv$ and
$K(A)=A\int_{-1/2}^{1/2}u^2\cos(\pi uA)\,du$. Then
\[
 B_I(U,H)=\pi\int_{-1/2}^{1/2}v\sin(\pi bv)\,dv
       +4\pi c^2\int_{-1/2}^{1/2}v K'(1-2bv)\,dv-pbcI/12 .
\]
Indeed $yJ$ has even part $v(A+4c^2u^2)$ in $u$, while
$\cos(\pi x)-\pi x\sin(\pi x)=\partial_A[A\cos(\pi uA)]$.
The first part integrates in $u$ to $vA\cos(\pi A/2)$,
whose even part in $v$ is $v\sin(\pi bv)$.

Writing $l=\pi|u|$, subtraction of the two derivatives gives
\begin{align*}
 K'(1-h)-K'(1+h)
 &=2\int_{-1/2}^{1/2}u^2[(\sin l+l\cos l)\sin(lh)
                      +lh\sin l\cos(lh)]\,du\\
 &\ge2h\int_{-1/2}^{1/2}u^2\sin^2 l\,du\ge h/10 .
\end{align*}
Here $0\le h\le1$, $\sin(lh)\ge h\sin l$ by concavity, and
$\sin l\ge2|u|$. Pairing $v,-v$ with $h=2bv$ thus gives at least
$\pi bc^2/30$ from the $K'$ integral with its prefactor.
The cubic lower bound for sine in the first integral gives
$b(p/12-p^2b^2/480)$. As $I\le4c/p$, it follows that
\[
 B_I(U,H)\ge b(.82-.21b^2-.23c^2).
\]
Interchanging the axes with a reflection gives the analogous bound
for $-B_I(V,H)$ with $b,c$ exchanged.

Two first moments suffice to change the metric. In $\pi\langle y\cos(\pi x)\rangle$,
the part of $yJ$ free of $c$ integrates in $u$ to a multiple of
$v\cos(\pi b v)$. The remaining even part gives
$4\pi c^2\int u^2v\sin(\pi u)\sin(2\pi buv)\,du\,dv$,
which lies between zero and $(p/48)bc^2$ since
$|\sin(2\pi buv)|\le2\pi b|uv|$.
Likewise $\pi\langle x\cos(\pi x)\rangle$ equals $2\pi c$ times
the integral of $u^2[\cos l\cos(lh)-h\sin l\sin(lh)]$, now with
$h=2b|v|$. Both terms inside the brackets are nonnegative before
subtraction. The first contributes at most $(1-8/p)c$ by integrating
$u^2\cos l$, and the second at most $(p/48)b^2c$ by $h\sin(lh)\le lh^2$.
Thus these moments have absolute bounds $.21bc^2$ and $.21c$.
Exchanging axes gives the corresponding bounds $.21cb^2$ and $.21b$,
with $\langle x\cos(\pi y)\rangle\le0$.

For $r>0$ put $w=2bc/r$ and $P_0=[bB(U,H)-cB(V,H)]/r$.
The Euclidean bounds combine to $.82-.21r-.01rw^2$.
In the change of metric the $z$ contribution is nonnegative by the
two signs just proved. The $u_0$ term is a difference of nonnegative
quantities, each of size at most $.21t b^2c^2/r$ after multiplication
by $|u_0|$. The $v_0$ term costs at most $.21tw$. Consequently
\begin{equation}\label{I:eq:quad-mixed-lower}
 P_0\ge .82-.21r-.01rw^2-.21t(w+rw^2/4).
\end{equation}

\paragraph{Contradiction by Cauchy--Schwarz.}
Cauchy--Schwarz for the positive semidefinite form $B$, followed by
the Euclidean Cauchy--Schwarz inequality in $(b,-c)$ yields
\begin{equation}\label{I:eq:quad-schur}
 P_0^2\le B(H,H)\,[B(U,U)+B(V,V)]/r
          \le B(H,H)[5t^2/r+1.625].
\end{equation}
The trace in the middle is $p z+p(\langle U\rangle^2+
\langle V\rangle^2)$. Polynomial integration with weight $J$ gives
$\langle x^2\rangle=(1+b^2)/12$, $\langle y^2\rangle=(1+c^2)/12$,
and $\operatorname{Var}(H)=1/144+r/80+b^2c^2/180$.
Substitution in $B(H,H)$ yields
\[
 B(H,H)\le .099-.039r+\tfrac1{12}
 [z(2+r)+u_0(c^2-b^2)-2v_0bc] .
\]
Here $1/6-p/144<.099$ and $1/12-p/80<-.039$.
The last two terms inside the bracket have sum at most $tr$.

The constraints $r(1+w)=(b+c)^2<S$ and $w\le1$ give $rw<S/2$.
By \eqref{I:eq:quad-t-small}, the losses containing $w$ in
\eqref{I:eq:quad-mixed-lower} sum to at most
$(.01+.168+.042S/2)rw<.21rw$. Thus $P_0>.82-.21S>2/3$.
For the right side of \eqref{I:eq:quad-schur} use $t<T(S)<.494$
in the preceding bound for $B(H,H)$. This gives
\[
 B(H,H)<.099+.494^2(2+S)/24+S(.494/12-.039)<.129 .
\]
Finally, multiplying $t\le .201+.405r$ by $t/r$, we find
$t^2/r\le .201(4/5)+.405(.494)<.361$.
Now \eqref{I:eq:quad-schur} gives
$P_0^2<.129(1.625+5\cdot.361)<(2/3)^2$, a contradiction.

If $r=0$, the global sine estimate \eqref{I:eq:quad-t-small} gives
$t=0$ and $M=I$. The normalized domain is a square. Its value $p$
follows by separation of variables. Together with Section~\ref{I:sec:lowertri}, this proves the theorem for all competitors with at most four sides.

\end{proof}

\section{A continuous trial field and its sector bounds}\label{I:sec:fan}
The variational comparison in this section holds for any $n\ge3$. We will use it here for $n=5,6$. 
By scale invariance take $R_n$ of circumradius one and $|P|=|R_n|$. A polygon is called competitive if $F(P)\ge F(R_n)$. The geometric estimates below show that a competitive polygon has $n$ strictly convex corners and hence is a strictly convex $n$-gon. Place its area centroid at the origin. Write $a_i$ for half the angle between the radii of the two vertices of edge $i$ and put
\[\fitmath{b=\pi/n,\quad k=\tan b,\quad t=k^2,\quad
 x_i=\tan(a_i)/k,\quad z_i=(r_{i+1}-r_i)/(r_{i+1}+r_i),\quad u_i=z_i^2.}\]
We call $x_i$ the aperture and $z_i$ the skew of sector $i$; the regular polygon has $x_i=1$ and $z_i=0$ for every $i$. The angles satisfy $\sum_i a_i=\pi$. For any one field and its fixed tangent slope $\gamma$ we have
\[\fitmath{T_i=1+\gamma(a_i-b),\qquad \sum_i T_i=n .}\]

Let $I$ denote area integration on $0<X<1$, $|Y|<X$, a triangle of area one. Take the bisector of a regular sector as the first coordinate axis. The field then has components $(P_0,kQ_0)$ in coordinates proportional to $(X,kY)$. Its entries are respectively even and odd in $Y$ and $P_0(X,X)=Q_0(X,X)=h(X)$. Define
\[\fitmath{c_0=I(P_0),\quad M_0=I(P_0^2+tQ_0^2),\quad
 E_0=I(P_{0X}^2+P_{0Y}^2/t+tQ_{0X}^2+Q_{0Y}^2),\quad R=(1+t)E_0/M_0 .}\]
Suppress the sector index and write $y=kx$, $D=1-u$, $j=xD$ and $\sigma=(r_++r_-)\cos(a)/2$. On this sector set
\[\fitmath{{\bf x}=\sigma(X+zY,y(Y+zX)),\qquad {\bf W}=\sigma(p+zr,y(q+z\nu)).}\]
Here $p,\nu$ are even and $q,r$ odd in $Y$. All equal $h$ on the upper ray. The traces agree on both sides of each shared ray, where the trial field equals $h(X)$ times the vertex vector. Use $I_c(f)=I(f^2)-(I(f)-c_0)^2$ and $J(f,g,w)=I[(f-ug)^2+uw^2]$. Define
\[\fitmath{\begin{aligned}
 M_0m&=I_c(p)+uI(r^2)+y^2[I(q^2)+uI_c(\nu)],\\
 J_1&=J(p_X,r_Y,r_X-p_Y),&J_2&=J(p_Y,r_X,r_Y-p_X),\\
 J_3&=J(q_X,\nu_Y,\nu_X-q_Y),&J_4&=J(q_Y,\nu_X,\nu_Y-q_X),\\
 V&=y^2J_1+J_2+y^4J_3+y^2J_4,& e&=V/(y^2D^2E_0),\\
 H_\alpha&=\frac{j(1+\alpha e)^2}{4m}.&&
\end{aligned}}\]
On an unmodified sector take $p=\nu=P_0$, $q=r=Q_0$. We call a pair of first Neumann eigenfunctions of $R_n$ a \emph{first eigenfield} if, written as a complex-valued function $u$ of a complex coordinate $z$ whose real axis bisects a sector, it satisfies $u(\omega z)=\omega u(z)$ for $\omega=e^{2\pi i/n}$ and $u(\bar z)=\overline{u(z)}$. By Lemma~\ref{IV:lem:regular-character} it is unique up to a real factor. To obtain a first eigenfield, project the polynomial pair onto the first eigenspace. For the polynomial fields used below this projection is nonzero: their Rayleigh quotients are below the next eigenvalue in this rotation class (see Sections~\ref{I:sec:five-six-regular-data} and~\ref{seven-spectral}). Dihedral symmetry ensures the same trace property for that pair.

Let $d_i=(\sigma_i/\cos b)^2$. Then $\sum_i j_i d_i=n$ by the area constraint. The vector $c_0$ times the midpoint of the boundary edge is $3c_0/2$ times the centroid of the triangle. The area-weighted sum of these vectors is zero. By Jensen's inequality, subtracting the mean over the polygon reduces the mass by at most the sum, weighted by the sector areas, of the squared differences between the sector means and these vectors. Parity then gives the lower bound $m_i$ for the mass. Differentiating the coordinate formula and integrating gives the energy $e_i$. The mass bounds for the fields used below are given in Sections~\ref{I:sec:five-six-thick}--\ref{paper-small-sum} and~\ref{seven-central}--\ref{sec:seven-narrow}, with the shear estimate in Section~\ref{seven-shear}. Provided $m_i>0$ on every sector, we obtain
\[\fitmath{\frac{\mu_1(P)}{R}
 \le\frac{\sum_i j_i e_i d_i}{\sum_i j_i m_i d_i^2}
 \le\frac1n\sum_i H_{1,i}.}\]
Indeed, set $v_i=d_i\sqrt{j_im_i}$, $A_i=\sqrt{j_i/m_i}$, $B_i=e_iA_i$. Then $n$ times the numerator equals $(A\cdot v)(B\cdot v)$, bounded above by $\|A+B\|^2\|v\|^2/4$. The denominator is $\|v\|^2$.

On a narrow sector, fix the aperture and vary $u$. If $c_0>0$ and $I(p),I(\nu)>c_0/2$, then $I_c(f)=\operatorname{Var}(f)+2c_0I(f)-c_0^2$ shows that $m$ is positive, increasing and affine. Hence $Dm$ is positive and concave. Since $V$ is quadratic with $V(1)\ge0$, $V/D$ is affine plus $V(1)/D$. Consequently $D(1+\alpha e)$ is positive and convex for $\alpha\ge0$. The convex function $(w,v)\mapsto w^2/v$ is increasing in its first argument and decreasing in its second for $w\ge0$, $v>0$. Its composition with the preceding two functions is therefore convex, as required for $H_\alpha$. For a fixed aperture it therefore suffices to check the two endpoints of the skew interval.

\section{Eigenvalue bounds for regular polygons and projection onto the first eigenspace}\label{I:sec:five-six-regular-data}
Normalize the circumradius to one and put $c=\cos(\pi/n)$, $t=\tan^2(\pi/n)$,
$c^2=1/(1+t)$. The dimensionless area and boundary inner products of the
polynomial fields, and the normalized quotients $R$ and $J/m_{\rm phys}$, belong to
$\mathbb Q(\sqrt5)$ for $n=5$ and to $\mathbb Q$ for $n=6$.
In the reference sector $0<X<1$, $-X<Y<X$, the physical
Cartesian coordinate is $c(X+i\sqrt tY)$. Thus the physical area
element is $c^2\sqrt t\,dX\,dY$.

Write a field as a pair $(P,Q)$ representing $P+i\sqrt t Q$.
Represent its scalar entries in homogeneous coordinates $X^dv^j$
where $v=Y/X$. If $e_j=(1+(-1)^j)/(j+1)$, the scalar
multiplication tables for the area inner product $I$ and chord inner
product $[\ ,\ ]$ are
\[\fitmath{I(X^dv^i,X^fv^j)=e_{i+j}/(d+f+2),\qquad [v^i,v^j]=e_{i+j}.}\]
Their pair extensions use $\langle(P,Q),(A,B)\rangle
=\langle P,A\rangle+t\langle Q,B\rangle$. The partial derivative
operations $D_x,D_y$ in $(X,Y)$ on a scalar basis element give
$(d-j)X^{d-1}v^j$ and $jX^{d-1}v^{j-1}$. Define
$E(U)=I(D_xU,D_xU)+I(D_yU,D_yU)/t$.

In this section only, $k$ denotes the radial-series coefficient. In the complex coordinate $z=X+i\sqrt tY$ the polynomial is
\[\fitmath{U=\sum_m a_m Z_m\sum_{l=0}^L\frac{(-k|z|^2)^l}{d_{ml}},
 \qquad Z_m=(X+i\epsilon_m\sqrt tY)^m,}\]
where $d_{ml}=l!(m+l)!/m!$, including $d_{m0}=1$.
The sign $\epsilon_m=+1$ or $-1$ is chosen so
$\epsilon_m m\equiv1\pmod n$. Hence $U(\omega z)=\omega U(z)$ for $\omega=e^{2\pi i/n}$, and $U(\bar z)=\overline{U(z)}$.
At five sides the sector field uses $t=5-2\sqrt5$, $L=2$, $k=5/7$. At six sides take
$t=1/3,L=3,\sigma=4043/1000,k=\sigma c^2/4=12129/16000$.
The table gives all coefficients for the sector fields. At six sides the last column also gives the denominator for the flux correction $T$.
\[\fitmath{\begin{array}{cc|rr|rrrr}
n&m&\epsilon_m&a_m&d_{m1}&d_{m2}&d_{m3}&d_{m4}\\\hline
5&1 & 1&1            &2 &12 &-& -\\
5&4 &-1&-1/20        &5 &60 &-& -\\
5&6 & 1&1/150        &7 &112&-& -\\\hline
6&1 & 1&1            &2 &12 &144&2880\\
6&5 &-1&-1/34        &6 &84 &2016&72576\\
6&7 & 1&1/170        &8 &144&4320&190080\\
6&11&-1&1/1130       &12&312&13104&786240\\
6&13& 1&-399/940160  &14&420&20160&1370880
\end{array}}\]
In particular the hexagonal $U$ has degree nineteen.

We can compute the area and boundary integrals together by writing the fields as complex monomials.
The term indexed by $(m,l)$ is just $c_{ml}z^a\bar z^b$, with
$c_{ml}=a_m(-k)^l/d_{ml}$, $a=l+\max(\epsilon_m m,0)$ and
$b=l+\max(-\epsilon_m m,0)$. For any $A,B\ge0$, let $p=\max(A,B)$
and $r=\min(A,B)$. Set $u_0=0,u_1=1$ and
$u_{j+2}=2u_{j+1}-(1+t)u_j$. Then
\[
 K_{AB}=\sum_{j=0}^r(-1)^j\binom rj 2^{r-j}
          \frac{2u_{p+j+1}}{p+j+1},\qquad
 I(z^A\bar z^B)=\frac{K_{AB}}{A+B+2}.
\]
Indeed $K_{AB}$ is the integral on the chord $X=1$, where
$\bar z=2-z$ and the primitive of $z^j$ gives
$2\operatorname{Im}(1+i\sqrt t)^{j+1}/((j+1)\sqrt t)$.
Conjugation preserves this real integral, and radial integration gives
the denominator $A+B+2$. Consequently two complex monomials with
exponents $(a,b),(c,d)$ contribute the single scalar
$K_{a+d,b+c}/(a+b+c+d+2)$ to their real area inner product.
For the energy use $2I(|U_z|^2+|U_{\bar z}|^2)$, differentiating
monomials by lowering the corresponding exponent. For example the
single $z$ term has mass $(3+t)/6$ and energy $2$. These rules also
give the boundary inner products using $K_{AB}$ alone.

To recover each scalar area product directly from the coefficient table,
for
$0\le j\le m,0\le h\le l$ a single entry of order $(m,l)$
contributes to degree $(d,i)=(m+2l,j+2h)$ the scalar
\[\fitmath{a_m\frac{(-k)^l m!}{l!(m+l)!}\binom mj\binom lh
       (-1)^{\lfloor j/2\rfloor}t^{\lfloor j/2\rfloor+h}.}\]
For even $j$ add it to $P$. For odd $j$ multiply by $\epsilon_m$
and add it to $Q$. Summing $0\le l\le L$ gives $U=(P,Q)$,
whereas taking $l=L+1$ gives $T$. Then $M=I(U,U)$,
$R=E(U)/(c^2M)$ and $m_{\rm phys}=c^2\sqrt t\,M$.
Indeed $\Delta(z^a\overline z^b)=4ab z^{a-1}\overline
z^{b-1}/c^2$ for the physical Laplacian. Successive radial coefficients
have ratio $-k/[(l+1)(m+l+1)]$, so $-\Delta(U+T)=\sigma U$.

The flux calculation below applies to any field $U$ constructed from these radial series. The truncation order may differ among harmonics. Let
$T$ denote the sum of the first omitted terms and put $\sigma=4k/c^2$. We apply it to the hexagonal
sector field and to the pentagonal trial of Lemma~\ref{lem:degree-six-pentagon-floor}.
To bound the flux residual, let
$Z=D_x(U+T)|_{X=1}$. Take a pair $V_0=(F,G)$ of primitives of
$Z$ with $F(0)=0$ and $G(1)=-F(1)/t$. Then the physical stream
trace $\sqrt t V_0/\sigma$ has the boundary derivative needed to cancel
the normal part of $-\nabla(U+T)/\sigma$. The endpoint relation
for $F,G$ makes it equivariant across adjacent sectors.
At six sides use the physical harmonic fields represented, up to arbitrary
constant scale factors, by
$H_j=i\sqrt t(X+i\epsilon_j\sqrt t Y)^j$, $j=1,n-1$.
Set $H_j^\partial=H_j|_{X=1}$,
\[\fitmath{\begin{aligned}
 D_j&=[H_j^\partial,D_xH_j|_{X=1}],&
 a'_j&=[V_0,D_xH_j|_{X=1}]/D_j,\\
 S&=\sum_j(a'_j)^2D_j,&V&=V_0-\sum_j a'_jH_j^\partial .
\end{aligned}}\]
The two harmonic gradients are orthogonal pointwise by the
Cauchy--Riemann equations for one analytic and one antianalytic
monomial. By Green's formula and the choices of $a'_j$, both are
orthogonal to every equivariant extension of $V$.
Use $X^\beta V$ for the remaining extension, with $\beta=28$ at six
sides. An empty harmonic sum is also allowed; then $S=0$ and $V=V_0$.
Write
\[\fitmath{A=[V,V],\quad B_f=[V',V']+t[vV',vV'],\quad C_f=[vV,V'],
 \qquad D_T=[T|_{X=1},Z].}\]
Define the physical flux $\mathcal F$ as $-\nabla(U+T)/\sigma$ plus the rotated gradient
of the extended stream function
$\psi=(\sqrt t/\sigma)(\sum a'_j H_j+X^\beta V)$, where
$\nabla^\perp=(\partial_y,-\partial_x)$ in physical coordinates. Consequently $\operatorname{div}\mathcal F=U$
and $\mathcal F$ has vanishing outer normal trace. Pairing with $\nabla U$ gives $-m_{\rm phys}$.
Thus replacing $1/R$ by $1/\sigma$ in its squared residual with
the trial gradient gives an upper bound $J$ for the optimal
residual. Integration of $X^{2\beta-1}$ and Green's formula give
\[\fitmath{\frac{J}{m_{\rm phys}}
 =\frac{t(S+\beta A/2-C_f)+B_f/(2\beta)+E(T)-2D_T}
 {\sigma^2c^2M}.}\]
For example, the contribution of $X^\beta V$ before stream
rescaling is
$\sqrt t(\beta A/2-C_f)+B_f/(2\beta\sqrt t)$,
since its coordinate derivatives are
$X^{\beta-1}(\beta V-vV')$ and $X^{\beta-1}V'$.
The last cross term is the boundary term
$\int\nabla T:\nabla^\perp\psi$ after integration by parts.
For the central moments of the original sector fields use
\[\fitmath{\begin{gathered}
 B=I(P_Y,P_Y)/(tE),\quad C=t I(Q_X,Q_X)/E,\quad v_p=tI(Q,Q)/M,\\
 \tau=[I(P_X-Q_Y,P_X-Q_Y)+I(P_Y-tQ_X,P_Y-tQ_X)/t]/E,
\end{gathered}}\]
where $E=E(U)$.

In the tables below, each entry is a lower bound. Adding one unit in its last decimal place gives the upper bound. The residual in this first table
belongs to the hexagonal field:
\[
\begin{array}{c|rrrrrrr}
 &10^6 S&10^6 A&10^6 B_f&10^6 C_f&10^6 E(T)&10^6 D_T&10^6J/m_{\rm phys}\\\hline
 n=6&.9277&6.054&1587.1&59.85&.3514&.1220&9.589
\end{array}
\]
The formulas above give the sector fractions directly. To test a bound $r$
on the Rayleigh quotient, take the sign of $E-c^2rM$.
\[
\begin{array}{c|rr}
 \text{sector quantity}&5&6\\\hline
 M&.33863&.31527\\ E(U)&.96529&.95607\\ R&4.355244&4.043368\\
 B&.002787&.004723\\ C&.001590&.002901\\
 v_p&.16470&.10609\\ \tau&.148630&.14522\\
 1+2(C-B-v_p)&.6681983&.7841602
\end{array}
\]
For instance $E(T)$ is obtained by
applying the displayed area pairing to the differentiated coefficient
list restricted to $l=L+1$. In the pentagonal column use
$t=5-2\sqrt5$ and reduce every product by $(\sqrt5)^2=5$;
the sign of a difference $a+b\sqrt5$ with $ab<0$ is determined
by $a^2-5b^2$ and the sign of $a$. At six sides use $t=1/3$.

We say that a complex-valued function $u$ on $R_n$ has \emph{rotation character one}
if $u(\omega z)=\omega u(z)$ for $\omega=e^{2\pi i/n}$. Let $G$ denote a lower bound for the
second Neumann eigenvalue of $R_n$ among such functions.
By Lemma~\ref{IV:lem:disk-angular-floors}, with
$K=(\sqrt{4+t}+\sqrt t)^2/4$, this eigenvalue is at least $28/K$ at six sides, and
$28/K>15.84$ by squaring positive quantities. At five sides
the map with constant Jacobian described immediately after that lemma gives
the larger lower bound $G=55/3$.
By Lemma~\ref{IV:lem:regular-character} the first eigenspace is two-dimensional
and the dihedral group acts on it as on $\mathbb R^2$. Hence the first eigenfields
form a real line. The fractions below are invariant under multiplication by a
nonzero real scalar. To obtain the flux estimates, let $K_N=(-\Delta_N)^{-1}$
on the functions of rotation character one and divide every integral on the
polygon by $n$. Integration by parts and the minimizing choice $1/R$
give
$\langle U,K_N U\rangle\le\|\mathcal F\|^2
 \le m_{\rm phys}/R+J$ for the constructed flux $\mathcal F$.
The spectral measure of $U$ is supported at $\mu$ and in $[G,\infty)$.
Integrating $(\lambda-G)(\lambda-\mu)/\lambda\ge0$ against this
measure gives the first inequality below.
The component of $\mathcal F+\nabla U/R$ in the gradient space is
$\nabla(U/R-K_NU)$, whose norm is therefore at most $\sqrt J$.
On the higher spectral subspace $I-RK_N\ge1-R/G$,
and the mass norm there is at most $1/\sqrt G$ times the energy
norm. Thus the sines $s_E$ and $s_M$ of the energy and mass angles
from $U$ to its first spectral projection satisfy
\[\fitmath{\begin{gathered}
 \mu^{-1}\le R^{-1}+G(J/m_{\rm phys})/(G-R),\\
 s_E\le \sqrt{RJ/m_{\rm phys}}/(1-R/G),\qquad
 s_M\le \sqrt{R/G}\,s_E
\end{gathered}}\]
with $\mu=\mu_1(R_n)$.

\begin{lemma}[A sixth-degree field on the regular pentagon]
\label{lem:degree-six-pentagon-floor}
For the circumradius-one regular pentagon, $\mu_1(R_5)>4.35415$.
\end{lemma}
\begin{proof}
Use the mass, energy and flux pairing just derived, with $s=|z|^2$,
$k=32/45$ and $\sigma=4k/c^2$. Take the sixth-degree field
\[
 \begin{aligned}
 W&=z(1-ks/2+k^2s^2/12)-\bar z^4(1-ks/5)/20+z^6/160,\\
 T&=-k^3zs^3/144-k^2\bar z^4s^2/1200-kz^6s/1120 .
 \end{aligned}
\]
The three angular orders are $1,-4,6$, all congruent to one modulo
five. Each last radial term of $W$ gets one additional power in $T$, so
$-\Delta(W+T)=\sigma W$. On the outer edge of the sector put
$Z=D_x(W+T)|_{X=1}$ and extend the primitive $V_0$ itself by $X^{14}V_0$.
Thus the harmonic sum is empty. With $M=I(W,W)$ the Rayleigh quotient
and the residual in the preceding inequalities are
\[
 R=4.355+D/M,\qquad
 J/m_{\rm phys}=N/[16k^2(1+t)M],
\]
where $D=(1+t)E(W)-4.355M$ and
$N=t(7A-C_f)+B_f/28+E(T)-2D_T$. The sparse form of $W$ allows the
three scalar substitutions below. Put $v=t-66/125$ and for each row
write $P=\sum_i p_i v^i$. The entries are rational bounds.
\[
\begin{array}{c|c|c|c|c}
P&\deg&10^6p_0&p_1&\sum_{i\ge2}|p_i|<\\\hline
M&6&(339290,339300)&(.063,.064)&.054\\
D&6&(129.5,129.6)&(.294,.295)&.23\\
N&8&(184.1,184.3)&(.00055,.00057)&.068
\end{array}
\]
For direct evaluation list the six terms of $W$ and three of $T$ as
$z^a\bar z^b$, using the monomial pairing above. If
$Z=(\sum c_jy^j,\sum d_jy^j)$ in the edge coordinate $y=Y/X$,
the two entries of $V_0$ are
$F=\sum c_jy^{j+1}/(j+1)$ and
$\sum d_j(y^{j+1}-1)/(j+1)-F(1)/t$. Use
$[y^i,y^j]=(1+(-1)^{i+j})/(i+j+1)$; all factors $1/t$ from the
constant in the second entry cancel in $N$. For the last column one can
use the following separate upper bounds, in increasing degree from two:
\[
\begin{array}{c|rrrrrrr}
 M&.0387&.0125&.00211&.00017&.000001&0&0\\
 D&.157&.0412&.0264&.00485&.000067&0&0\\
 N&.0085&.0281&.0209&.0081&.0016&.00015&.000007
\end{array}
\]
Zero entries are exact.

Squaring gives $-.00014<v<0$. The rows give
$.339<M<.340$,
$D>.0001295-.295(.00014)-.23(.00014)^2$,
$D<.0001296+.23(.00014)^2$ and
$N<.0001843+.068(.00014)^2$. Thus
$4.355259<R<4.355383$ and $J/m_{\rm phys}<.000043981$.
Since $G=55/3$,
\[
 \frac1{4.355259}+\frac{G\cdot.000043981}{G-4.355383}
 <\frac1{4.35415}.
\]
\end{proof}

For the hexagonal field the tabulated bounds on the residual and Rayleigh quotient
give $\mu_1(R_6)>4.04315$, $s_E<.008362$ and $s_M<.004226$. For the
pentagonal field we need only the upper bound $R<4.355245$.
Indeed, expand the mass-normalized field in eigenvectors. Since all
higher eigenvalues are at least $G$, the fraction of the mass in the higher modes is at most
$(R-\mu)/(G-\mu)$, and the corresponding energy fraction is
$1-(\mu/R)(1-s_M^2)$. Lemma~\ref{lem:degree-six-pentagon-floor} therefore gives
\[
 s_M^2\le\frac{R-\mu}{G-\mu}<(.0089)^2,\qquad
 s_E^2\le\frac{G(R-\mu)}{R(G-\mu)}<(.0182)^2 .
\]
The spectral ratios $R/\mu$ are below $1.000252,1.000053$ for the pentagon and hexagon, respectively. The resulting errors over the entire fan, $n(1-\mu/R)$, are below
$.001258,.000313$.

The angle estimate applies on one sector because equivariance makes all sector contributions to energy, mass and their
polarizations equal. For $a=B,C,\tau/2$ we project the unit gradient orthogonally onto a coordinate direction in each of the first two cases, and onto the two normalized signed differences of coordinates in the third. For $a=v_p$
use the second coordinate of the mass-normalized field instead. If the sine of the angle is at most $s$, the new squared fraction lies between
$(\sqrt{a(1-s^2)}-s\sqrt{1-a})^2$ and
$(\sqrt{a(1-s^2)}+s\sqrt{1-a})^2$ for the above ranges.
The quantity being squared in the first expression is positive because $s^2<a$.
The second is increasing in the angle here since $a+s^2<1$.
The projected first field consequently satisfies
\[\fitmath{\begin{array}{c|cccc}
 n&B&C&v_p&\tau\\\hline
 5&(.00099,.00549)&(.00034,.00375)&(.15582,.17379)&(.1270,.17175)\\
 6&(.00364,.00594)&(.00207,.00388)&(.10350,.10872)&(.13666,.15403).
\end{array}}\]
Both the polynomial and its projection therefore satisfy the bounds used below for thick sectors.

\section{Bounds on the sector geometry for pentagons and hexagons}
A polygon is called competitive when $F(P)\ge F(R_n)$.
For a convex polygon
let $G$ denote its area centroid, $2a$ the angle of a side at $G$,
$z=(r_+-r_-)/(r_++r_-)$, and $x=\tan a/\tan(\pi/n)$.

\subsection{The pentagonal aperture and skew}\label{sec:pent-entry}
Write $A=|P|$, $F=A\mu_1(P)$ and let $\Sigma$ be the covariance matrix
of uniform area measure. The one-dimensional marginals of a convex
polygon, rescaled to $[0,1]$, are mixtures of unit tents with random
apex $T$: the nonnegative mixture weight at an interior break $t$ is
$t(1-t)$ times its downward slope jump divided by two; at each
endpoint the weight is half the endpoint value. With $u=\mathbb E T$,
$B=\operatorname{Var}T$, integration gives
$\mathbb EX^j=2\sum_{k=0}^j\mathbb ET^k/[(j+1)(j+2)]$ and hence
\[\fitmath{15\{(12/5)(\operatorname{Var}X)^2-\mathbb E(X-\mathbb E X)^4\}
 = B^2+\mathbb E[(T-u)^2 T(1-T)]\ge0 .}\]
Fix a direction in the plane. Let $X$ now be the centered coordinate in that direction, $V=\mathbb E X^2$,
and $\kappa=\mathbb EX^4/V^2\le12/5$. The trials $X$ and
$X^3-3VX-\mathbb EX^3$ are orthogonal in energy, with energies
$1,9V^2(\kappa-1)$. Projecting their mass onto $X$ and optimizing
their two coefficients gives
\[\fitmath{1/\mu_1(P)\ge V\{1+(3-\kappa)^2/[9(\kappa-1)]\}\ge36V/35.}\]
For a unit tent the conditional variance $(1-t+t^2)/18$ is at least
$(1+t)^2/72$, so total variance gives $V\ge(\mathbb EX)^2/8$ in
the interval coordinates. Also $\mathbb ET^2\le u$ yields
$\mathbb E X^2\le(1+2u)/6\le3(\mathbb E X)^2/2$. Reflect and
rescale the interval. In every supporting-line direction the centroid
distance $d$ consequently satisfies $2V\le d^2\le 8V$.
If $\sigma_-\le\sigma_+$ are the eigenvalues of $\Sigma$, rays from
$G$ therefore have lengths between $\sqrt{2\sigma_-}$ and
$\sqrt{8\sigma_+}$. A determinant-one linear change followed by the
polar moment inequality for a disk of area $A$ gives
$\sqrt{\det\Sigma}\ge A/(4\pi)$.
The ratio of these two bounding radii is thus at most
$2\sigma_+/\sqrt{\det\Sigma}\le K=70\pi/(9F)$.

Let $\delta$ be the angle of the side normal from its ray bisector.
Both endpoint rays make angle at most $\arccos(1/K)$ with the normal,
so $a+|\delta|\le\arccos(1/K)$ and
$|z|=\tan a\,|\tan\delta|$ by the sine rule.
Put $y=\tan a$ and $S=\sqrt{K^2-1}$. It follows that
$0<y\le S$, $|z|\le y(S-y)/(1+Sy)$ and
$|z|\le(K-1)/(K+1)$.
These bounds on aperture and skew require only $F>10.3518$.
The lower bound for the regular pentagon gives this inequality for a competitive polygon, and
Proposition~\ref{prop:quantitative-window-entry} gives it for small positive deficits.
Use $\pi<355/113$ and $.72654<k=\tan(\pi/5)<.72655$
(the latter follows by squaring, since $k^2=5-2\sqrt5$).
Then $S<2.13813$, $K<2.360422$, $x<2.13813/.72654<2.943$,
and $|z|<.405$. For $x\le.65$, the function
$y(S-y)/(1+Sy)$ increases in $S$.
With $1.553<2.13813k<1.555$ and $k^2>.5278$ this gives
\[\fitmath{|z|\le x(1.555-.5278x)/(1+1.553x).}\]
Thus every sector of a competitive pentagon either has aperture in $[.65,3]$ and absolute skew at most $.41$, or satisfies the preceding bounds for a narrow sector.

\subsection{Bounds on hexagonal sectors from three pairs of harmonics}\label{sec:hex-entry}
Normalize area to $\pi$. Let $f$ map the unit disk conformally onto the
polygon, centered so the first disk field $g(r)e^{i\theta}$ has zero
pullback mean, and rotated so $\alpha=f'(0)>0$. Put $O=f(0)$,
$r_*=\operatorname{dist}(O,\partial P)$, $h=|G-O|/r_*$,
$v=f-O-\alpha z$, and $T=\alpha^{-2}-1$.
Write $f'/\alpha=1+\sum_{j\ge1}\gamma_jz^j$,
$w_j=|\gamma_j|^2/(j+1)$, so $T=\sum w_j$.
Let $\lambda$ be the first positive disk eigenvalue and
$M_j=2j\int_0^1g(r)^2r^{2j-1}dr$ for the normalization $g(1)=1$.
In Proposition~\ref{IV:prop:harmonic-budget-65}, the centered harmonic estimate gives
$\sum_jD_j w_j\le H$ with $D_j=1/\lambda-H-(1-M_{j+1})$;
$H=(\lambda-1)/\mu_*-(1-1/\lambda)$ for any lower bound
$\mu_*\le\mu_1(P)$ in area-$\pi$ units. It suffices that $H(P)=(1-1/\lambda)(\pi\lambda/F(P)-1)\le H_0$.
In that case one may take $\mu_*=\lambda/(1+\lambda H_0/(\lambda-1))$.
For a competitive hexagon this threshold follows by converting the regular hexagon bound $\mu_1(R_6)>4.04315$ to the area normalization used here. The conversion factor is
$3\sqrt3/(2\pi)$. With $p_+=3.38996,p_-=3.3899$,
$\pi<355/113$, $\sqrt3>1.7320508$ we obtain
\[\fitmath{H\le (1-1/p_+)\{2\pi p_+/(3\sqrt3\cdot4.04315)-1\}
 <.009763< H_0=.009764 .}\]
Lemma~\ref{IV:lem:two-harmonic-secants} gives $|\gamma_i|\le H_0g_i$ at $i=1,2$,
$g_1=6.812+37H_0$, $g_2=9.13+37.5H_0$. By Lemma~\ref{IV:lem:uniform-harmonic-correlations},
$D_i\ge d_i$ where $d_1=(5/p_+-1)/3-H_0$,
$d_j=1/p_+-H_0-4.92/[(2j+3)(2j+4)]$ for $j\ge2$.
For $J(D)=\sum_{i=1,2}(1-d_i/D)_+g_i^2/(i+1)$ this gives
$T\le T_0=H_0/d_3+H_0^2J(.295)<.04452$, since $D_i\ge d_3$ for $i\ge3$, $d_3<.295$, and $J$ is increasing.
Thus the harmonic multiplier $\mu_* /(1+T)$ exceeds
$p_-/[(1+p_+H_0/(p_--1))(1+T_0)]>3.201$.
The estimate following Lemma~\ref{IV:lem:two-harmonic-secants} gives
$X=\sum_{i\ge1}|\gamma_i|^2/(i+2)
 \le H_0/(1/3.39-H_0)+7H_0^2<.03491$.

Take $s_0=3.20$ in Proposition~\ref{IV:prop:harmonic-budget-65}, and retain only the three pairs
of signed frequencies $1\pm k$, $k=3,4,5$. Choose
\[\fitmath{(t_3,t_4,t_5)=(2/19,1/11,1/18),\qquad
 (b_3,b_4,b_5)=(3/68,2/75,5/408),
 \qquad b_k=\frac{k t_k}{(k^2-1)(1-t_k)}.}\]
Set $t_j=0$ for $j\ge6$ and
$g_j^*=[(j^2-1)^{-1}+(j+2)^{-2}]/2$.
In either signed harmonic, denote its squared mass pairing by
$|p+e|^2$, where $p$ is the contribution from the constant Fourier
coefficient of $f'/\alpha$. The inequality
$|p+e|^2\ge t_k|p|^2-t_k|e|^2/(1-t_k)$ and $g(r)\ge r$ give
the diagonal gain $s_0t_kg_k^*$ on $w_k$. Cauchy--Schwarz on
the remaining terms gives a loss of at most $s_0 Xb_k(1-(k-2)/j)$ on each
$w_j$ with $j>k$; these are exactly the two terms obtained in
the proof of Proposition~\ref{IV:prop:harmonic-budget-65}. There is no contribution for $j<k$.
Write $A(j)=1.23/[(j+3/2)(j+2)]$; Lemma~\ref{IV:lem:uniform-harmonic-correlations} gives $1-M_{j+1}\le A(j)$.
The resulting diagonal coefficient falls short of $1/\lambda$ by at most
\[\fitmath{H_0+A(j)+.112\!\sum_{\substack{k\in\{3,4,5\}\\k<j}}
           b_k(1-(k-2)/j)-3.20t_jg_j^*.}\]
At $j=3,4,5,6$ substitution gives upper bounds
$.03665,.03701,.03703,.03706$, respectively. At $j\ge7$ one
bound suffices: $H_0+A(7)+.112(b_3+b_4+b_5)<.03515$.
Hence every diagonal from the third on exceeds $1/\lambda-\varepsilon$
with $\varepsilon=.03706$. Summing the resulting coefficients and
using $C=6^5\log(1+T)-4.808(36)-12.4$ gives
\[\fitmath{\begin{aligned}
 T&\le H_0/(1/p_+-\varepsilon)+H_0^2J(1/p_--\varepsilon)<.03957,\\
 C&<6^5(.03957-.03957^2/2+.03957^3/3)-4.808(36)-12.4
 <116.3.
\end{aligned}}\]
For $c_1=\gamma_1/2$, $c_2=\gamma_2-\gamma_1^2/2$,
the two leading coefficient bounds give $|c_1|<.03502$, $|c_2|<.09518$,
$\Delta=[1-4(|c_1|+|c_2|)]/3>.15973$.
Use $e^u\le1+u+u^2/[2(1-u/3)]$ at $u=2(.03502)+.09518$.
The boundary estimate in Lemma~\ref{IV:lem:convex-boundary-two-heads} gives
\[\fitmath{\|v\|_\infty/\alpha<
 \sqrt{[5(.694)/.15973-1.5](.03957)}
      +1.17967/[64(.15973)]+2^{-15}<1.011 .}\]
At $r=.9$, set $N=1.662-r^2-r^4/2-r^6/3$.
The coefficient bound for $f$ on this circle gives
\[\fitmath{r_*/\alpha\ge r-H_0g_1r^2/2-H_0g_2r^3/3-\sqrt{.03957N}>.7319 .}\]
Indeed on the radius in any supporting normal direction take the
corresponding real component of $f$. Bound its first two
nonlinear coefficients by $H_0g_1/2,H_0g_2/3$, then apply
Cauchy--Schwarz to the remaining Taylor terms using
$-\log(1-r^2)<1.662$. The support line bounds this component above.
With $\alpha=(1+T)^{-1/2}$, the centroid bound of Lemma~\ref{IV:lem:physical} and $5/16<.502$ give
\[\fitmath{h<T[.502\sqrt{1+T}+1.011]/[(1+T)(.7319)]<.0793}\]
by monotonicity of $T/(1+T)$ and $T/\sqrt{1+T}$.

\paragraph{Moving the center of the hexagonal fan to the centroid.}
The estimate from the three pairs of harmonics gives $T<.03957$, $C<117$ and
$h=|G-O|/\operatorname{dist}(O,\partial P)<.080$ at the conformal center
$O=f(0)$. To bound the aperture, use $q=6a_O/\pi$, $y=(6-q)/5$ and
\[\fitmath{A(q)=28.8(q^4+5y^4-6)+(12.4/6)(q^6+5y^6-6)
                      +(q^8+5y^8)/6 .}\]
Dropping the nonnegative relative entropy and the remaining nonnegative skew terms in \eqref{IV:eq:single-aperture-budget} gives $A(q)+6^5(q/6)W(q/6,\delta_O/\pi)\le C$,
where $W\ge0$ is the skew term of Lemma~\ref{IV:lem:skew-fan}.
Strict convexity gives the minimum of $A$ at one.
Direct rational substitution gives $A(.125)>119.51$ and
$A(1.575)>123.49$. Hence $.125<q<1.575$. Set
$a_0=A(0)-117$, $b=-A'(0)$ and $m=6(28.8)(36)/126$.
The quartic alone gives $A''\ge2m$, by
$(q+5y)^2\le126(q^2+y^2/5)$. Thus
$(C-A(q))/q^2\le b^2/(4a_0)-m=:K_6<397.21$
($a_0>30.83674$, $b<234.69835$).
For $v=(\delta_O/\pi)^2$, the rational term in the sector inequality that bounds the skew is at least $6^3q^2v[32/(1-4v)+41/25]$.
With $c=32+41/25$, the bracket is at least $c/(1-128v/c)$ by convexity.
Consequently $v\le K_6/(6^3c+128K_6/c)$, so $|\delta_O|<.669$.
The first inequality of Corollary~\ref{IV:cor:aperture-skew-product} gives
\[\fitmath{|z_O|\le z_0=
 \sqrt{117\pi^4/[32(6^5)(1-4(.669)^2/\pi^2)]}<.23654 .}\]
Set $p=.92$, $b_0=\pi/6$, $d=.669$. Integrating the changes of angle and relative skew along $OG$ as in Section~\ref{seven-transport} yields
\[\fitmath{\begin{aligned}
 \frac{p\tan(.125b_0)}{\tan b_0}<x
 &<\frac{\tan(1.575b_0)}{p\tan b_0},\\
 |z|&\le\frac{z_0+(1-p)/(1+p)}{1+z_0(1-p)/(1+p)},\\
 |z|/x&\le\frac{\tan b_0}{2p}\{(1+p^2)\tan d+(1-p^2)\sec d\}.
\end{aligned}}\]
These imply $.104<x<2.037$, $|z|<.276$, $|z|/x<.520$.
One may bound $\pi$ by $3.14159,3.14160$ and use the alternating
sine and cosine polynomials through degrees nine and ten on the
three arguments, all less than $.825$. These bounds put each sector within the range of either the thick-sector or the narrow-sector estimate below.

\subsection{Excluding a missing convex corner}\label{sec:five-six-missing}
\paragraph{Up to four convex corners.}
Assume a polygon with at most four strictly convex corners has
$F\ge56\pi/17$. Center the conformal disk against
$g(r)e^{i\theta}$ with $g(r)=(3r-r^3)/2$. Write
$f'/f'(0)=1+\sum_{j\ge1}\gamma_j z^j$,
$H(s)=\sum_{j\ge1}|\gamma_j|^2e^{-js}$, $\tau=17/24$,
$W(x)=x(3-x)^2/4-\tau$. The Dirichlet energy of the pair divided
by $\pi$ is $7/3$ and its disk mean squared mass is $11/16$.
Using the two components as test functions gives
\[\fitmath{\int_0^1W(x)H(-\log x)\,dx\le 1/48 .}\]
The following radial comparison gives a strict bound. Put $q=e^{-4s}$, $Y=1+H$, $v=\log Y$, $z=-v'$.
Lemma~\ref{IV:lem:signed-circle}, gives $H''/4+bH'+b^2Y\ge(1+b)^2q$ for every $b\ge0$.
Taking $b=(2q-H')/[2(Y-q)]$ gives
$v''\ge(z+2)^2q/(Y-q)$, so $z$ increases with $q$. For $z\le4$ put
$J(z)=3z/(2+z)-\log(1+z/2)$.
Then $[J(z)+\log(1-q/Y)]'\le0$ in $s$, with limit zero at infinity.
This gives $Y/q\ge1/(1-3/e^2)>37/22$ up to $z=4$, and thereafter
$Y/q$ is nondecreasing in $q$.

For the inner positive-weight range only $q\le(24/25)^4<17/20$ is
needed. Let $h(q)=q/4+9q^2/[64(1-2q/3)]$,
$Z=4qh'/(1+h)$. On this shorter range $Z<4$ since
$1+h-qh'=1-(3q/8)^2/(1-2q/3)^2>0$.
Write $D=J(Z)+\log(1-q/(1+h))$. Differentiation shows the sign of
$D'$ in $q$ is the opposite of that of
$320q^2-1395q+1008$, which is decreasing and positive at $17/20$. Thus $D\le D(0)=0$.
Whenever $Y<1+h$, the first integral and monotonicity of $J$ up to
four imply $z\ge Z$ (also immediate if $z>4$).
There $\partial_q\log[Y/(1+h)]\ge0$, so $Y\ge1+h$ starting from
equality at $q=0$.

For the possible negative mass at $0\le x\le a=3/5$, one needs
only the first Fourier coefficient. Indeed $W$ has a single zero
in this interval and $\int_0^a x^2W(x)dx=129/250000>0$, so
$\int_0^a x^jW(x)dx\ge0$ for all $j\ge2$.
For completeness put $d_j=1-3(j+5)/[2(j+2)(j+3)(j+4)]-\tau$,
$u_j=|\gamma_j|\sqrt{d_j/(j+1)}$. The trial functions give
$\sum_{j\ge1}u_j^2\le1/48$. The centered first moment of $g$ has
coefficients $B_j=(2j+7)/[2(j+2)(j+3)]$. Hence
$B_0|\gamma_1|\le\sum_{j\ge1}c_j u_j u_{j+1}$, where
$c_j=B_j\sqrt{(j+1)(j+2)/(d_jd_{j+1})}$.
Here $c_1^2=24300/833$, $c_2^2=26136/1379$, $c_3^2=3136/197$,
and $c_j\le1/\sqrt{d_4d_5}<4$ for $j\ge4$.
Apply $2u_ju_{j+1}\le (t_{j+1}/t_j)u_j^2+
(t_j/t_{j+1})u_{j+1}^2$ with $t_1=2/3$, $t_j=1$ otherwise:
the first three row sums are at most $33/8,(11/3+22/5)/2,
(22/5+4)/2$ and all others at most $4$. Consequently
$|\gamma_1|\le(21/5)/(48B_0)=3/20$.

On $a\le x\le b=24/25$, $q=x^4$ and $W\ge0$. Since
$1/(1-y)\ge1+y+4y^3$ for $0\le y<1$, we have
$H\ge L(q)=q/4+9q^2/64+3q^3/32+q^5/6$.
On $b\le x\le1$ use $H\ge37x^4/22-1$. All three integrals are
obtained from the single primitive
\[\fitmath{U_j(t)=t^{j+1}\left(\frac{9t}{4(j+2)}
 -\frac{3t^2}{2(j+3)}+\frac{t^3}{4(j+4)}-\frac{17}{24(j+1)}\right).}\]
In particular $U_1(a)=-.010212$,
$\int_a^b W L(x^4)>.0146859$ and
$\int_b^1W(37x^4/22-1)>.0064384$. Their sum including the possible
first-coefficient loss exceeds
$.0146859+.0064384-.010212(3/20)^2>1/48+.000061$.
This contradicts the upper bound from the two test functions. Finally $56\pi/17<10.35$, whereas
Lemma~\ref{lem:degree-six-pentagon-floor} gives $F(R_5)>10.3526$.

\paragraph{Up to five convex corners.} Theorem~\ref{IV:thm:corner-gap} excludes a competitive hexagon with a missing convex vertex. The integrated circle mean there gives $F(P)<10.4959$, while the conformal estimate of Lemma~\ref{small-schur-regular} gives $F(R_6)>10.4967$. Comparing these bounds excludes every competitor with at most five strictly convex corners. The resulting gap is greater than $1/1250$.

\section{A supporting line for thick sectors}
\label{I:sec:five-six-thick}
Let $n=5$ or $6$, $b=\pi/n$, $k=\tan b$, $t=k^2$,
$x=\tan a/k$, $y=x^2$ and $u=z^2$, where $a$ is the sector
half-angle and $z$ is the signed relative difference of its endpoint
radii. For the regular pair $(P,kQ)$ write $E,M_0$ for the reference
energy and mass in the centroid coordinates, and put
\[\fitmath{B=\frac{I(P_Y^2)}{tE},\qquad C=\frac{tI(Q_X^2)}E,\qquad
 v=\frac{tI(Q^2)}{M_0},\qquad
 \tau=\frac{I[(P_X-Q_Y)^2+(P_Y/k-kQ_X)^2]}E .}\]
Here the regular pair can be the polynomial or its projection onto
the first eigenspace. Define
\[\fitmath{f(x)=\frac{1+t}{k}[\arctan(kx)-b],\qquad
 h=1+2(C-B)-2v,\qquad T(x)=1+hf(x).}\]
For either pair the unmodified transport
$(p,q,r,\nu)=(P,Q,Q,P)$ satisfies $H_1\le T$ throughout the following
rectangles, with equality precisely at $x=1,u=0$:
\[\fitmath{\begin{array}{c|ccc|ccc|c}
 n&l&r_0&z_*&[B_-,B_+]&[C_-,C_+]&[v_-,v_+]&\tau_*\\\hline
5&.65&3&.41&[0,.00552]&[0,.00378]&[.1552,.1754]&.181\\
6&.77&2.32&.302&[.00167,.00787]&[.00089,.00606]
                       &[.09905,.11322]&.18266
\end{array}}\]
Here $l\le x\le r_0$, $0\le u\le z_*^2$, $0\le\tau\le\tau_*$.
The bounds above show that both the polynomial field and its projection onto the first eigenspace satisfy the tabulated inequalities.

\paragraph{Reduction to zero and maximal skew.}
Put $D=1-u$, $L=1+uty$, $J=y+u/t$, $M=(1-v)L+vJ$ and
\[\fitmath{N=(L-J)(BL-CJ)+\tau u(1+ty)^2/t .}\]
The transported field is obtained by applying the same
linear map $S=\bigl(\begin{smallmatrix}1&z/k\\ kxz&x\end{smallmatrix}\bigr)$
to the vector and its argument, up to a common factor. Its Gram matrix has diagonal $L,J$, off-diagonal
entry $z(1+ty)/k$ and determinant $yD^2$. The two gradient columns of
the transported field are obtained from $S\nabla U S^{-1}$ in physical
reference coordinates. After squaring, the odd integrands vanish by symmetry. This gives
the normalized energy $1+N/(yD^2)$ and mass $M$. The fan functional is
\[\fitmath{H_1=\frac{xD}{4M}\left(2+\frac{N}{yD^2}\right)^2.}\]
Increase $\tau$ to $\tau_*$ in this expression. This can only increase
$H_1$: the unsquared factor is bounded below by
$2-(B+C)LJ/(yD^2)>0$, since
$LJ/(yD^2)=1+u[(ty)^{-1}+2+ty]/D^2<3$ on the two rectangles.
The last inequality follows by taking the endpoints in $ty$, since
$1/T+T$ is convex for $T>0$. At $u=1$ the increased $N$ is positive,
because $N(1)/(1+ty)^2=(1-1/t)(B-C/t)+\tau_*/t>0$.
Consequently $N/D$ is affine plus the positive multiple $N(1)/D$,
and $2D+N/(yD)$ is positive and convex. Also $DM$ is positive and concave,
since $M$ is affine and increasing in $u$. The function
$(p,s)\mapsto p^2/s$ is convex, increasing in $p>0$ and decreasing in $s>0$.
Consequently the resulting upper bound for $H_1$ is convex in $u$. It suffices to prove
the supporting-line inequality at $0$ and $z_*^2$.

\paragraph{Zero skew.}
For $d=x-1$ put $K=(d-f)/d^2$ with continuous interpretation at $d=0$.
Expanding the preceding formula at $u=0$ gives the exact identity
\[\fitmath{\frac{M(T-H_1)}{d^2}=
 (v-C)(x+2)-2v(v+B-C)(x+1)-B/x-hMK
     -\frac{(x+1)^2(Cx^2-B)^2}{4x^3}.}\]
The last absolute value is below $.00025$, by bounding with the
larger of $C_+^2x(x+1)^2/4$ and $B_+^2(x+1)^2/(4x^3)$.
The first three terms increase with $v$ and decrease in $B,C$:
the $v$ derivative is $x+2-(4v+2B-2C)(x+1)>0$.
Integrating the angle derivative gives
$K=t\int_0^1 w R(1+wd)\,dw$, $R(s)=(1+s)/(1+ts^2)$.
The difference between its tangent at one and $R(s)$ is
\[\fitmath{\frac{t(s-1)^2[3-t+(1-3t)s]}{(1+t)^2(1+ts^2)}\ge0
                  \quad(0\le s\le3).}\]
Hence $K\le k_0+k_1(x-1)$ with
$k_0=t/(1+t)$, $k_1=t(1-3t)/[3(1+t)^2]$. Write
$h_+=1+2(C_+-B_-)-2v_-$. A lower bound for the zero-skew divided
difference above is
\[\fitmath{G(x)=(v_--C_+)(x+2)-2v_-(v_-+B_+-C_+)(x+1)-B_+/x
 -h_+(1-v_-+v_+x^2)(k_0+k_1(x-1))-.00025.}\]
For the cubic product before the last constant, its chord minus its
value equals $v_+(x-l)(r_0-x)[k_0+k_1(x+l+r_0-1)]$.
The bracket is positive by $k_0>0$ and $k_0+6k_1>0$; also $-B_+/x$
is concave. Thus $G$ lies above its chord. The two pairs of endpoint
values exceed $(.0809,.1250)$ and $(.0155,.0256)$ respectively.
This proves strictness for every $x\ne1$. At $x=1$ direct
substitution gives equality on this edge.

\paragraph{Maximal skew.}
Set $u=z_*^2$, $D=1-u$, $a_0=1-u/t$, $b_0=1-ut$,
$m=1-v+vu/t$, $s=v+(1-v)ut$. With $W=2+N/(yD^2)$ write
$W=A/y+B_0+C_0y$. Its coefficients are
\[\fitmath{\begin{aligned}
 A&=[a_0(B-Cu/t)+\tau_*u/t]/D^2,\\
 B_0&=2+[2\tau_*u-a_0(C-But)-b_0(B-Cu/t)]/D^2,\\
 C_0&=[b_0(C-But)+\tau_*ut]/D^2.
\end{aligned}}\]
Each coefficient is positive throughout the rectangles in the table. The first increases in $B$
and decreases in $C$, the last decreases in $B$ and increases in $C$,
and $B_0$ decreases in both since $2u/t,2ut<1$. These monotonicities give the extreme values of each coefficient. Write $R_0=r_0^2$. Differentiating
$H_1=DxW^2/[4(m+sy)]$ gives
\[\fitmath{\frac{4(m+sy)^2}D H_1'\le \alpha-\beta y-\lambda/y ,}\]
where
\[\fitmath{\begin{aligned}
 \alpha&=m(B_0^2+2AC_0),\qquad \lambda=2mAB_0,\\
 \beta&=s(B_0^2+2AC_0)-6mB_0C_0
    -R_0(5mC_0^2+2sB_0C_0+3sC_0^2R_0).
\end{aligned}}\]
Indeed $x y^j/(m+sy)$ has derivative numerator
$y^j[(2j+1)m+(2j-1)sy]$ over $(m+sy)^2$. For $j=-2,-1$
retain only $-2mAB_0/y$, and bound $y^2,y^3$ in the positive
contributions for $j=1,2$ by $R_0y,R_0^2y$.

Let subscripts minus and plus denote the lower and upper endpoints
on the rectangle. Bound $\alpha$ from above by its upper positive
products and $\beta$ from below by using lower endpoints in its positive
products and upper endpoints in its subtracted products. Write these
bounds $\bar\alpha,\underline\beta$. Put $\lambda_-=2m_- A_- (B_0)_-$,
$h_-=1+2(C_--B_+)-2v_+$, $g=4h_-(1+t)/D$. Multiplication of the
derivative difference $T'-H_1'$ by the positive number
$4(1+ty)(m+sy)^2/D$ gives the lower bound
$e_0+e_1y+e_2y^2+\lambda_-/y$, where
\[\fitmath{(e_0,e_1,e_2)=
 (gm_-^2-\bar\alpha+t\lambda_-,
  2gm_-s_--t\bar\alpha+\underline\beta,
  gs_-^2+t\underline\beta).}\]
The tangent $1/y\ge 2/3-y/9$ reduces positivity to the following
lower bounds, obtained by the substitutions just given using
$t=5-2\sqrt5,1/3$:
\[\fitmath{\begin{array}{c|ccc}
 n&e_2&e_0+2\lambda_-/3&e_1+l^2 e_2-\lambda_-/9\\\hline
 5&.4260&.0905&.3232\\ 6&.1760&.1239&.1428
\end{array}}\]
In particular $T-H_1$ is strictly increasing at maximal skew.

At $x=l$, $L>J$ and $f(l)<0$. If $F_1=W(L-J)/(2lDM)$, then
$\partial_B H_1=F_1L$, $\partial_C H_1=-F_1J$,
$\partial_v H_1=H_1(L-J)/M$. Bound $W$ by
$A_+/l^2+(B_0)_++(C_0)_+l^2$ and $M$ below by $(1-v_+)L+v_+J$.
The larger of $F_1 L,H_1(L-J)/M$, divided by $-2f(l)$, is then
less than $.848,.490$ respectively. Hence $T-H_1$ increases in $B,v$
and decreases in $C$; its minimum occurs at $B_-,C_+,v_-$.
The endpoint comparison has a direct algebraic form. Put
$D=1-z_*^2$, $y=l^2$, $s=(l-1)/(1+tl)$ and $f_-=(1+t)s(1-ts^2/3+t^2s^4/5)$.
The alternating arctangent series gives $f(l)\ge f_-$ since $ks<0$
and $|ks|<1$. Substituting $B_-,C_+,v_-,\tau_*$ at the indicated
endpoint gives
\[
\begin{array}{c|ccc}
 n&1+[1+2(C_+-B_-)-2v_-]f_->&
  lD[2+N/(yD^2)]^2/(4M)<& T-H_1>\\\hline
5&.725674&.716284&.00939\\
6&.802899&.802715&.000184
\end{array}
\]
At six sides every entry uses $t=1/3$ alone. At five sides insert
$t=5-2\sqrt5$ and square to compare in $\mathbb Q(\sqrt5)$.
Convexity in $u$ and the strictness on this
upper edge complete the proof, including $x=1$ with $u>0$.

\paragraph{The tangent line at the boundary of the narrow-sector range.}
Each choice of reference field uses its own $h$, constant across all
sectors of the polygon, and $\sum_iT(x_i)=n$ since $\sum_i a_i=nb$.
For the pentagonal polynomial its moments give
$1.02091< (1+t)h<1.02092$ and $(1+t)h/k<1.4053$. For
$x\le.65$ the line $T$ is consequently bounded below by
$1+1.4053[\arctan(.72654 x)-.628319]$.
For the short hexagonal polynomial the moment bounds above give
$h<.784161$. Thus on $x\le.77$ its own line dominates
$1+(.7842)(1+t)[\arctan(kx)-\pi/6]/k$.
These are the comparison lines used in the narrow-sector estimates.

\section{The pentagonal narrow-sector comparison}

Normalize the regular pentagon to circumradius one. Put $b=\pi/5$,
$k=\tan b$ and $t=k^2=5-2\sqrt5$. Let $2a$ be a sector angle at
the area centroid of a competitive convex pentagon of the same area,
and let $r_-,r_+$ be its two endpoint radii. Set
\[\fitmath{x=\tan(a)/k,\qquad z=(r_+-r_-)/(r_++r_-),\qquad u=z^2 .}\]
The geometric bounds used here are $x<2.943$, $|z|<.405$, and
\[\fitmath{|z|\le Z(x):=\frac{x(1.555-.5278x)}{1+1.553x}
                       \qquad(0<x\le .65).}\]
Indeed the joint covariance bound gives
$\tan a\le S$ and
$|z|\le(\tan a)(S-\tan a)/(1+S\tan a)$ with
$S=\sqrt{[70\pi/(9F)]^2-1}<2.13813$, where
$F=|P|\mu_1(P)\ge F(R_5)>10.3518$.
The right side increases with $S$, and
$1.553<2.13813 k<1.555$, $t>.5278$.

\paragraph{The reference field and its moments on a ray.}
Write $I$ for integration on the triangle $0<X<1,\ |Y|<X$, of area one.
For $\zeta=X+ikY$ use the regular polynomial
\[\fitmath{\begin{aligned}
 P_0+ikQ_0&=\zeta p_1(|\zeta|^2)-\bar\zeta^4p_4(|\zeta|^2)/20
                  +\zeta^6p_6(|\zeta|^2)/150,\\
 p_m(v)&=1-\frac{5v}{7(m+1)}+\frac{25v^2}{98(m+1)(m+2)} .
\end{aligned}}\]
Put $h(X)=P_0(X,X)=Q_0(X,X)$, $c=I(P_0)$,
\[\fitmath{M_0=I(P_0^2+tQ_0^2),\quad
 E_0=I(P_{0,X}^2+P_{0,Y}^2/t+tQ_{0,X}^2+Q_{0,Y}^2),
 \quad R=(1+t)E_0/M_0.}\]
Write $\langle f\rangle=2\int_0^1 Xf(X)\,dX$,
$H_2=\langle h^2\rangle$, $N=\langle(h'-h/X)^2\rangle$,
$p_1^*=h(1)$, $L_1=3\langle Xh\rangle/2$. For the thick-sector
supporting line put $B=I(P_{0,Y}^2)/(tE_0)$,
$C=tI(Q_{0,X}^2)/E_0$, $v=tI(Q_0^2)/M_0$,
$\gamma=(1+t)[1+2(C-B)-2v]/k$, and $T=1+\gamma(a-b)$.
The entries in the following table equal $(a_0+b_0t)/d_0$ exactly.
\[\fitmath{\begin{array}{c|r|r|r}
 &a_0&b_0&d_0\\\hline
 M_0&118927743533514713&19303938041129181&381289693676141250\\
 E_0&1570111015455908&1553164106566050&2475907101793125\\
 M_0v&5724564727334735&9298887641313633&190644846838070625\\
 E_0B&-220270406603&978127719453&110040315635250\\
 E_0C&-2035913003801&4817117252031&330120946905750\\
 c&7883585&-45903&15280650\\
 \langle h\rangle&76695&679&154350\\
 p_1^*&26703&9145&51450\\
 H_2&4723579468291&764962338005&19305318532500\\
 N&-74518831791&157262310676&56780348625\\
 L_1&47170311&1434725&88288200
\end{array}}\]
To obtain a row, expand the three harmonics by the binomial theorem.
Use $t^2=10t-5$, $1/t=2-t/5$,
$I(X^iY^j)=(1+(-1)^j)/[(j+1)(i+j+2)]$, and
$\langle X^j\rangle=2/(j+2)$. For the five ray rows one can instead set
$d_b=1/\cos b=\sqrt5-1$ and use the univariate polynomial
\[\fitmath{h(X)=Xp_1(d_b^2 X^2)
   +d_b^3 X^4p_4(d_b^2X^2)/20-d_b^5 X^6p_6(d_b^2X^2)/150.}\]
For the scalar bounds below $2.23606797749<\sqrt5<2.23606797751$
suffices, as follows by squaring. In particular
$4.3552444<R<4.3552445$, $1.02091<k\gamma<1.02092$,
$.51<c<.52$ and $\langle h\rangle>.49$. Fix $\eta=1.00032$.

\paragraph{A trial field with prescribed values on the rays.}
For $0<x\le.65$ choose the four profiles
\[\fitmath{p=h,\qquad q=r=Yh/X,\qquad \nu=h+\tfrac32X(X^2-Y^2).}\]
They are polynomials since $h(0)=0$. On the upper ray all four
agree with $h$; parity gives the matching trace on the lower ray.
Thus they agree there with the unmodified thick-sector trial, which
has $p=\nu=P_0$, $q=r=Q_0$. For any profile $f$ put
$I_c(f)=I(f^2)-(I(f)-c)^2$. With $s=x^2$, $D=1-u$ and $d=u/s$,
the centered mass factor is
\[\fitmath{A=\{I_c(p)+uI(r^2)+ts[I(q^2)+uI_c(\nu)]\}/M_0 .}\]
For the gradient put
$\mathcal J(f,g,w)=I((f-ug)^2+uw^2)$ and
\[\fitmath{\begin{aligned}
 J_1&=\mathcal J(p_X,r_Y,r_X-p_Y),&
 J_2&=\mathcal J(p_Y,r_X,r_Y-p_X),\\
 J_3&=\mathcal J(q_X,\nu_Y,\nu_X-q_Y),&
 J_4&=\mathcal J(q_Y,\nu_X,\nu_Y-q_X).
\end{aligned}}\]
This sector contributes
\[\fitmath{H_\eta=\frac{xW^2}{AD^3},\qquad
 W=\frac12\left[D^2+\frac{\eta}{E_0}
                       \left(J_1+J_4+\frac{J_2}{ts}+tsJ_3\right)\right].}\]
For fixed $x$, $A$ is positive, increasing and affine in $u$, since
$I(\nu)=\langle h\rangle+2/5>c/2$ and $I(p)>c/2$.
The energy term $J_1+J_4+J_2/(ts)+tsJ_3$ is quadratic
in $u$ with nonnegative value at $u=1$. Dividing by $D$ gives an
affine function plus a nonnegative multiple of $1/D$.
Thus $2W/D$ is nonnegative and convex in $u$, while $AD$ is positive and
concave. The function $(v,w)\mapsto v^2/w$ is convex, increasing
in $v\ge0$ and decreasing in $w>0$. This shows that $H_\eta$ is
convex in $u$, so it suffices to bound it at $u=0,Z(x)^2$.
The enlarged endpoint is admissible for these formulas:
$Z$ is increasing up to $.65$ because
$1.555-2(.5278)x-(1.553)(.5278)x^2>0$ there, and $Z(.65)<.393$.

At zero skew,
$A_0=[H_2-(\langle h\rangle-c)^2+tsH_2/3]/M_0$.
The coefficient of $u$ in $A$ is
$[H_2/3+tsI_c(\nu)]/M_0\ge A_0/3$ because
$I_c(\nu)\ge c(2I(\nu)-c)>H_2/9$. Thus
\[\fitmath{A\ge m(1+u/3),\qquad m=m_0+.138s,\quad m_0=.7836 .}\]
For the energy write $v'=Y/X$ temporarily and $b'=3X^2/2$.
If $g=h'-h/X$, the derivatives of $r$ are $v'g,h/X$, and those
of $\nu$ are $h'+b'(3-(v')^2),-2b'v'$.
Integration by parts gives $\langle h'h/X\rangle=(p_1^*)^2$ and
$\langle b'h'\rangle=3p_1^*-3L_1$. Averaging the even powers of
$v'$ on $[-1,1]$ and putting $S_1=N+2(p_1^*)^2$ gives
\[\fitmath{\begin{aligned}
 J_1+J_4={}&S_1+u(2N/3-4(p_1^*)^2+4p_1^*+1-32L_1/3)\\
           &+u^2(S_1+16p_1^*-16L_1+27/5),\\
 J_2={}&(u+u^2/3)N,\qquad
 J_3=N/3+u(N+20p_1^*-80L_1/3+27/5)+u^2 .
\end{aligned}}\]
For $W$, the computed upper coefficients of $1,s,d,u,ud,us,u^2,su^2$ are
\[\fitmath{.96671,\ .01364,\ .14685,\ -2.93870,\ .04895,\ .91084,\
 4.34485,\ .27352 .}\]
In particular, with $w=.967+.014s$ and $p_*=.147$,
\[\fitmath{W\le w+p_*d+u[-2.938+.0491d+.912s+(4.346+.275s)u].}\]
The identity $D^3(1+u/3)-(1-4u/3)^2=u^2(2/9-u^2/3)$ implies
$\sqrt{AD^3}\ge\sqrt m(1-4u/3)$.

\paragraph{The angular term.}
Write $\ell=.65,\ K=.72654,\ c_*=1-1.4053(.628319)-.005$ and
$g_*=1.4053K$. Since $K<k$, $b<.628319$, $\gamma<1.4053$ and
$a-b<0$, one has $T-.005\ge c_*+1.4053\arctan(Kx)$.
For $v'$ uniform on $[0,1]$ put $Z'=K^2x^2(v')^2$,
$\mu=\mathbb E Z'$. The exact tangent remainder of $(1+Z')^{-1}$
at $\mu$ yields
\[\fitmath{(1+\mu)\frac{\arctan(Kx)}{Kx}
   =1+\frac1{1+\mu}\mathbb E\frac{(Z'-\mu)^2}{1+Z'} .}\]
Now $m/m_0\ge1+\mu$ and $\operatorname{Var}((v')^2)=4/45$.
Since $4m_0g_* K^4/[45(1+K^2\ell^2/3)(1+K^2\ell^2)]>3/200$,
\[\fitmath{m(T-.005)\ge c_*m+m_0g_*x+3x^5/200 .}\]
The latter expression minus $xw^2$ has derivative at most
$m_0g_*-.967^2+2(.138)c_*\ell+5(3/200)\ell^4<-.1015$.
Its value at $\ell$ exceeds $.0008266$, proving the strict angular
comparison at $u=0$.

\paragraph{Maximal skew.}
At $u=Z(x)^2$, put $a_0=1.555,\ b_0=1.553,\ c_0=.5278$.
Then $d=[(a_0-c_0 x)/(1+b_0x)]^2$ is decreasing and convex:
its positive square root is decreasing and convex. Also $u\le.28x$,
since after multiplication by $(1+b_0x)^2/x$ this assertion becomes
\[\fitmath{.28-1.548345x+2.31676452x^2-.27857284x^3\ge0 .}\]
On $[0,.6]$ the quadratic minorant using $x^3\le.6x^2$ has minimum
above $.00118$. On $[.6,.65]$ the quadratic minorant using
$x^3\le.65x^2$ is increasing and positive at $.6$.

Let $X_0=\sqrt{p_*}$ and $B(x)=-2.938+4w/3+.0491d+.912x^2$.
We claim
\[\fitmath{\frac{W}{1-4u/3}\le w+(p_*-x^2)_+d .}\]
Subtract the proposed upper bound after multiplying by $1-4u/3$,
and divide by $u>0$. The resulting difference is bounded above,
on the two sides of $X_0$, respectively, by
\[\fitmath{\begin{aligned}
  &B(x)+1+4p_*d/3+(4.346+.275x^2-4/3)(.28x)
                   &&(0<x\le X_0),\\
  &B(x)+p_*/x^2+(4.346+.275x^2)(.28x)
                   &&(X_0\le x\le\ell).
 \end{aligned}}\]
Both functions are convex: besides constants they are nonnegative
multiples of $d,x,x^2,x^3$ and, on the second interval, $1/x^2$.
Their endpoint values are below $-.0560,-.0079$ and
$-.00566,-.07755$, respectively, using $.383405<X_0<.383407$.
This proves the claim. For $x\ge X_0$ it reduces the required
comparison to the zero-skew estimate already proved.

For $0<x\le X_0$ set $f=w+(p_*-x^2)d$, $v_0=x/(1+b_0x)$,
$D_0=c_0+a_0b_0$, $A_*=.967+p_*a_0^2$, and $B_*=2p_*a_0D_0$.
Direct expansion gives
\[\fitmath{f=A_*-B_*v_0-E(x)v_0^2,\qquad
 E(x)=(a_0-c_0x)^2-p_*D_0^2-.014(1+b_0x)^2 .}\]
Here $E$ decreases, $E(X_0)>.5210$; write $v_X=X_0/(1+b_0X_0)$
and $\xi=E(X_0)/(B_*+E(X_0)v_X)$. Then
$B_*\xi/(1-\xi v_0)\le E(X_0)$ and $b_0-\xi<1.199<6/5$.
Consequently
\[\fitmath{0<f\le A_*-\frac{B_*v_0}{1-\xi v_0}
   =A_*-\frac{B_*x}{1+(b_0-\xi)x}
   \le A_*-\frac{B_*x}{1+6x/5}.}\]
The angular lower bound is strictly greater than $m_0(c_*+g_*x)$.
After multiplying by $(1+6x/5)^2$, it suffices that
\[\fitmath{m_0(c_*+g_*x)(1+6x/5)^2-x[A_*+(6A_*/5-B_*)x]^2>0 .}\]
In ascending order the four coefficients exceed
$.08778,-.73814,1.40749,1.09370$. A positive minorant for $x\ge0$ is
\[\fitmath{\tfrac7{80}-\tfrac34x+\tfrac75x^2+x^3
 =(x-\tfrac15)^2(x+\tfrac25)+\tfrac75(x-\tfrac9{40})^2+\tfrac1{1600}.}\]
We have proved $H_\eta<T-.005$ on every pentagonal narrow sector,
including all intermediate skews by convexity.

\section{The hexagonal narrow-sector comparison}
For the five-frequency trial of degree nineteen, set
$K=1/\sqrt3$, $t=1/3$, $g=.7842(1+t)$ and $S=.77^2$. Then
on $0<x\le.77$, $|z|\le.52x$, the sector contribution $H_1$
satisfies
\[\fitmath{H_1<T_*-.00179,\qquad
 T_*=1+(g/K)[\arctan(Kx)-\pi/6].}\]
In particular the same strict bound holds relative to the polynomial field's tangent line, since $1+2(C-B)-2v<.784161<.7842$ and $x<1$.

We now compute the polynomial coefficients used in this comparison.
In the coefficient calculation write $Y=Xb$ (here $b$ is a local coordinate). With $\epsilon_m=1$ for $m=1,7,13$ and $-1$ for $m=5,11$,
write $P_0+iKQ_0$ as the sum of the terms
\[\fitmath{a_m X^{m+2l}(1+i\epsilon_m Kb)^m
  (1+tb^2)^l\frac{m!(-12129/16000)^l}{l!(m+l)!},
 \quad 0\le l\le3,\qquad
 (a_1,a_5,a_7,a_{11},a_{13})
 =(1,-1/34,1/170,1/1130,-399/940160).}\]
For $s=x^2$, $u=ks$ (in this section $k$ denotes squared skew divided by $s$), take the four profiles
$p=h+s(P_0-h)$, $r=bh$, $q=r+s(Q_0-r)$,
$\nu=h+3X^3(1-b^2)/2$, $h=P_0(X,X)$.
The four profiles have the same traces on $Y=X$.
Write $D=1-ks$, $m_0=I(P_0)$, $I_c(f)=I(f^2)-(I(f)-m_0)^2$,
where $I$ is area integration in $(X,Y)$ over $|Y|<X<1$. Put
\[\fitmath{\begin{aligned}
 \mathcal A&=[I_c(p)+uI(r^2)+ts(I(q^2)+uI_c(\nu))]/M_0,\qquad P=\mathcal A D^3,\\
 W&=\tfrac12[D^2+(J_1+J_4+J_2/(ts)+tsJ_3)/E_0],\qquad Q=W^2/P .
\end{aligned}}\]
Here $M_0=I(P_0^2+tQ_0^2)$, $E_0=I(P_{0X}^2+P_{0Y}^2/t+tQ_{0X}^2+Q_{0Y}^2)$,
and $J(f,g,w)=I((f-ug)^2+uw^2)$ is applied to
$(p_X,r_Y,r_X-p_Y)$, $(p_Y,r_X,r_Y-p_X)$,
$(q_X,\nu_Y,\nu_X-q_Y)$, $(q_Y,\nu_X,\nu_Y-q_X)$ in order.
The second of these has a factor $s$. The patched fan identity gives
$H_1=xQ$ on each of the two edges $k=0,169/625$.
For completeness, fix $x>0$ before differentiating in $u=z^2$.
All four profiles are then fixed. The coefficient rule below gives
$.5134<m_0<.5135$ and $I(P_{0Y}^2)<.00151$. Integration by parts on
the vertical fibers and $I(Y^2)=1/6$ give
$|I(h)-m_0|=|I(YP_{0Y})|<.016$. Thus $I(p),I(\nu)>m_0/2$
since $0\le s\le S<1$ and the correction in $\nu$ is nonnegative.
The identity $I_c(f)=\operatorname{Var}(f)+m_0(2I(f)-m_0)$
gives $\mathcal A(u)>0$, affine with strictly positive slope.

Write $D=1-u$, $B=D\mathcal A$ and
$N=(J_1+J_4+J_2/(ts)+tsJ_3)/E_0$. Each $J_i$ is quadratic in $u$
and $J_i(1)\ge0$, directly from its defining squares.
Then $H_1=(x/4)C^2/B$, with $C=D+N/D>0$.
Polynomial division gives $C''=2N(1)/D^3\ge0$ and
$B''=-2\mathcal A'<0$. Consequently
\[\fitmath{\frac4x\frac{d^2H_1}{du^2}
 =\frac2B\left(C'-\frac{CB'}B\right)^2
   +\frac{2CC''}B-\frac{C^2B''}{B^2}>0.}\]
This proves the reduction to $u=0$ and $u=(169/625)x^2$ for
each fixed aperture. At the other ray $Y=-X$, parity gives
$p=\nu=h$ and $q=r=-h$, as required for the global trace.

All the following coefficients have rational expressions. For the
products in the definitions above use
\[\fitmath{I(X^i b^jX^{i'}b^{j'})=
 \frac{1+(-1)^{j+j'}}{(i+i'+2)(j+j'+1)}.}\]
The two Cartesian derivatives of $X^i b^j$ are
$(i-j)X^{i-1}b^j$ and $jX^{i-1}b^{j-1}$. The coefficients in ascending powers of $s$ are enclosed below.
Each nonzero entry is a lower bound. Add $.0001$ to obtain the upper bound.
A printed zero is exact.
\[
\begin{array}{c|rrrrr}
 k=0,\ \mathcal A&.8359&.1504&.0127&.0007&0\\
 W&.9679&.0331&-.0055&.0042&0\\
 k=169/625,\ \mathcal A&.8359&.2258&.2632&.0007&0\\
 W&1.0364&-.7783&.4789&.0328&.0004
\end{array}
\]
For example the constant mass term is
$[I(h^2)-(I(h)-I(P_0))^2]/M_0$.
Writing $p_1=P_0-h$, $q_1=Q_0-r$, the last nonzero mass term
is $t I(q_1^2)/M_0$ at either skew. All nine skew energy
coefficients are obtained by the same four squares and moment rule.

At zero skew $P=\mathcal A$ and $P(1)=W(1)=1$ exactly, by 
$P_0,Q_0$ being the regular trial and by the absence of a spectral
factor in $W$. Put $U=W-P$, $V=2W-P$ and
$D_U=U/(s-1)=U_1+(s+1)U_2+(s^2+s+1)U_3$ in coefficient notation.
From the two rows at $k=0$,
\[\fitmath{U_3>.0033,\quad U_2+3U_3<-.007,\quad
 D_U(0)<-.13,\qquad V_3>.007 .}\]
Thus $D_U$ is negative and decreasing on $[0,1]$. Since
$Q=V+U^2/P$ and $P$ increases, the second divided difference
of $Q$ at one is at most $M=V_2+(2+S)V_3+D_U(S)^2/P(0)$.
Let $F=T_*/x$, $c_*=1-g\pi/(6K)$. Then
$F'(1)=[g/(1+t)-1]/2$ and
$F''(1)=[3-g(3+5t)/(1+t)^2]/4$. The difference in slopes contributes at most
$E=[F'(1)-V'(1)]_+/(1-S)$ in divided-difference units.
The intervals in the table give $M<.0205$, $E<.0030$, while $F''(1)/2=.0319125$.
Since $.05174<c_*<.05175$, the remaining coefficient exceeds $3c_*/20$.
In $F=c_*/\sqrt s+g\int_0^1(1+tsv^2)^{-1}dv$, the integral
lies above its second Taylor polynomial at one for $s\le1$.
Retaining the exact remainder of the first term gives
\[
 T_*-H_1\ge c_*[(1-x)^3(8+9x+3x^2)/8+3x(1-x^2)^2/20]>.0023 .
\]
Indeed the derivative of the first term in brackets is
$-15(1-x^2)^2/8$, while the derivative of $x(1-x^2)^2$ is at most
$(1-x^2)^2$. It suffices to insert $x=.77$.

At $k=169/625$, the coefficient bounds in the table give $W''>.957$, $W'(S)<-.175$
and $\mathcal A'(S)<.54<.67<3k\mathcal A(0)$.
Thus $W$ decreases, as does $P$ by $P'=D^2(D\mathcal A'-3k\mathcal A)$.
Instead of combining the smallest numerator with the largest denominator
over the whole interval, compare them on $[0,S/2]$ and $[S/2,S]$.
On either $[a,b]$ use $W(b)/P(a)\le \rho=W/P\le W(a)/P(b)$.
This gives $1.01<\rho<1.44$ on the first half and $1.04<\rho<1.35$
on the second. Since the mass coefficients are positive,
\[
 P''\le \mathcal A''(S)-6k(1-kS)^2\mathcal A'(0)+6k^2\mathcal A(S)<.738 .
\]
Differentiation gives
\[
 Q''=2(W'-\rho P')^2/P+2\rho W''-\rho^2P''
 \ge1.914\rho-.738\rho^2>1.17,
\]
as is checked at $1$ and $1.45$ by concavity.
Moreover $Q(0)<1.286$, $Q(S)<.896$. Use the chord of
$Q-(1.17/2)s^2$ and $\arctan(Kx)\ge Kx-(Kx)^3/3$ to obtain
\[
 T_*-H_1> .05174-.2404x+.888x^3-.585x^5 .
\]
Here $1.286-g=.2404$ and $(1.286-.896)/S+.585S-g/9>.888$.
For $x\le .32$ the cubic and quintic together exceed or equal
$(.888-.585(.32)^2)x^3$. The minimum of $ux^3-vx$ on $x\ge0$ is
$-2v\sqrt{v/(3u)}/3$; the resulting bound exceeds $.00188$.
For $x\ge .32$ the derivative of the displayed quintic is a concave
quadratic in $x^2$, positive at both $.32^2$ and $S$. Thus the same
lower bound applies to the full interval at maximal skew.
The estimate at zero skew and convexity now give the bound for every intermediate skew.

\section{Summing the sector estimates for pentagons and hexagons}\label{paper-small-sum}
If all sectors are thick, use the first eigenfield throughout. Its Rayleigh quotient on $R_n$ is $\mu_1(R_n)$ and its tangent has the same slope on the whole fan. The thick-sector inequality gives $\sum_i H_{1,i}\le n$, strictly unless $x_i=1$ and $z_i=0$ for all $i$. The variational inequality in Section~\ref{I:sec:fan} therefore gives a strict inequality unless these equalities hold.

If some sector is narrow, use the polynomial field on the regular polygon as the reference for all sectors. Transport it unchanged to each thick sector. On each narrow sector use the modified polynomial with the same boundary trace on the two rays. Every thick sector has $H_1\le T$. A narrow sector has $H_1\le H_{1.00032}<T-.005$ when $n=5$ and $H_1<T-.0017$ when $n=6$. Thus $\sum_i H_{1,i}<n-b_n$, where $b_n$ is the gain from just one narrow sector. Let $L_n$ be the lower bound for the eigenvalue of the regular polygon of circumradius one, and let $R_{\max}$ bound the Rayleigh quotient of the polynomial used on this fan. The estimates above give
\[\fitmath{\begin{array}{c|cccc|c}
 n&L_n&R_{\max}&b_n& n(1-L_n/R_{\max})\text{ is below}
 &b_n-n(1-L_n/R_{\max})\text{ exceeds}\\\hline
 5&4.35415&4.355245&.005&.001258&.003742\\
 6&4.04315&4.043369&.0017&.000325&.001375
\end{array}}\]
Hence $\mu_1(P)\le R\sum_iH_{1,i}/n<L_n<\mu_1(R_n)$ in this case.

When equality holds, every angle at the centroid is $2\pi/n$ and every pair of consecutive radii is equal. The vertices form a regular $n$-gon. Conversely, scaling and isometries preserve $F$, so every similar copy attains the bound. This proves the cases $n=5,6$. Section~\ref{I:sec:lowertri} and Theorem~\ref{I:thm:all-real-quads} establish $n=3,4$.

\clearpage\part{The centroid comparison for seven through nine sides}\label{part:seven}
\section{Area normalization and the fan}
Throughout this part $n=7,8,9$. Scale the polygons so that $|P|=|R_n|$ and the regular polygon has circumradius one. Call $P$ competitive if
$\mu_1(P)\ge\mu_1(R_n)$. Theorem~\ref{IV:thm:corner-gap} excludes fewer than $n$
strictly convex vertices. Therefore a competitive $P$ is a strictly
convex $n$-gon. Let $a_i$ denote half the angle subtended by an edge at the
area centroid, $b=\pi/n$, $k=\tan b$, $x_i=\tan a_i/k$ and
$z_i=(r_{i+1}-r_i)/(r_{i+1}+r_i)$ for consecutive vertex radii. Then
\[\fitmath{\sum_i\arctan(kx_i)=nb,\qquad u_i=z_i^2 .}\]
We compare $H_{1,i}=j_i(1+e_i)^2/(4m_i)$ to a tangent
$T_i=1+(\gamma/k)(a_i-b)$ with the same $\gamma$ on all sectors
of a given polygon. Here $j_i=x_i(1-u_i)$, $e_i$ is the sector energy
per unit area divided by the reference energy, and $m_i$ is the corresponding
centered-mass lower estimate before rescaling the sector. In particular
$\sum T_i=n$.

For either reference field with Rayleigh quotient $R$, apply Section~\ref{I:sec:fan} to obtain $\mu_1(P)/R\le n^{-1}\sum_iH_{1,i}$. The tangent for the chosen reference is the same on every sector.
\section{Bounds on the centroid sectors and a polynomial trial field}

For $n=7,8,9$ respectively, the bounds on the centroid sectors proved in Section~\ref{seven-window-proof} are
\[
 x<r_n=(1.78,1.64,1.59),\qquad |z|<(.16,.105,.075),\qquad
 |z|<\rho_n x,\quad \rho_n=(.230,.129,.085).
\]
\subsection{Moving the center of the fan to the centroid}\label{seven-transport}
The change of center in Section~\ref{seven-window-proof} uses the following three differential bounds.
Fix the base endpoints $a,b$ and move the center $y$ of the fan along the segment
from the conformal center $O$ to the area centroid $G$.
Write $d$ for its distance to the base line,
$\alpha$ for the sector half-angle and $\delta$ for half the
difference of the base angles. In directions measured from the edge
normal the endpoint rays have arguments $\delta\pm\alpha$, up to
a common sign. In particular $|\delta|+\alpha<\pi/2$.
From the derivatives of the two arguments and the sine rule,
\[\fitmath{|\nabla\alpha|=|a-b|/(2|a-y||b-y|)
       =\sin(2\alpha)/(2d),\qquad
 |\nabla\operatorname{atanh}z|=\sin(2\alpha)/(2d),}\]
where $z=(|b-y|-|a-y|)/(|b-y|+|a-y|)$ and
$\operatorname{atanh}z=\frac12\log(|b-y|/|a-y|)$.
The magnitude of the sum of the two argument gradients gives
\[\fitmath{4d^2|\nabla\delta|^2=
     1+2\cos(2\alpha)\cos(2\delta)+\cos^2(2\alpha)
     \le 4\cos^2\delta .}\]
The difference on the right factors as
$(1-\cos2\alpha)(1+\cos2\alpha+2\cos2\delta)\ge0$,
since $\cos2\delta\ge-\cos2\alpha$. Consequently
$|\nabla\log\tan\alpha|\le1/d$ and
$|\nabla\operatorname{arsinh}(\tan\delta)|\le1/d$.
Put $D=\operatorname{dist}(O,\partial P)$, $h=|G-O|/D$ and $p=1-h$.
At displacement $s$ along the segment we have $d\ge D-s$.
Integrating the three bounds yields
\[\fitmath{\tan\alpha_G\le\tan\alpha_O/p,\qquad
 |z_G|\le\tanh\!\left(\operatorname{atanh}|z_O|-\tfrac12\log p\right),
 \qquad
 |\tan\delta_G|\le\sinh\!\left(\operatorname{arsinh}|\tan\delta_O|-\log p\right).}\]
Since $|z_G|/x_G=\tan(\pi/n)|\tan\delta_G|$, these give all three
bounds needed when the center moves from the conformal center to the centroid.
\subsection{Sector polynomials and exact moments}\label{seven-moments}\paragraph{Seven sides.}
Use $b=\pi/7$, $c=\cos b$, $k=\tan b$ and $t=k^2$, in the fan coordinates
of the centroid lemma. Put $z=c(X+ikY)$ and $Z=(X^2-Y^2)^2$. Take
\[
 V=z+\{ -465 z|z|^2+15 z|z|^4-52\overline z^6
                  +Z(-46+79X-27ikY)\}/1024 .
\]
The first three corrections have rotation character one.
The last correction and its first derivatives vanish on both radial edges,
so rotation and reflection extend $V$ to an admissible field.

\paragraph{Eight and nine sides.}
For the octagon take $c=\cos(\pi/8)$, $k=\tan(\pi/8)$, $t=k^2$.
The following fifteen ordered pairs form the degree-five basis used
at eight and nine sides. Start with $(X^j,X^{j-1}Y)$, $1\le j\le5$.
Then for $d=0,1,2,3$ in increasing order and $j=0,\ldots,d$ in
increasing order put $B_{dj}=(X^2-Y^2)X^{d-j}Y^j$ and append
$(B_{dj},0)$ for even $j$ and $(0,B_{dj})$ for odd $j$. Every pair
has the required parity and equal component traces on $Y=X$.
Thus a linear combination, continued by reflection and rotation, is
in $H^1$ on the regular polygon. Continuity of the trial field and the weak
eigenfunction equation suffice everywhere below.
In the degree-five sector basis let $V=c(P_0+ikQ_0)=P+ikQ$ have coefficient
numerators $(1024,3,-493,15,78,2,49,255,54,-644,34,-74,548,-42,22)$
over $1024$.

Set $c=\cos(\pi/9)$, $k=\tan(\pi/9)$, $t=k^2$ and
$l=.67$, $r=1.59$. Use $V=c(P_0+ikQ_0)=P+ikQ$ in the same
degree-five sector basis as the octagon. Its coefficient numerators are
\[
 (4096,-8,-1843,-77,350,25,42,1170,377,-3184,97,-339,2570,-120,77)
\]
over $4096$. The traces agree across the sector boundaries by construction.

\paragraph{Exact integration of the reference field.}
Let $\psi=X(X^2-Y^2)$, $c=\cos(\pi/n)$, $t=c^{-2}-1$.
Write the reference field as $V=\sigma(p+i\sqrt t\,q)$ with
$\sigma=1$ at seven and $\sigma=c$ at eight and nine. Thus $p,q$ at
eight and nine are the rational quintics specified by their fifteen
integer coefficients. At seven sides, $p,q$ are the real and scaled
imaginary components of the displayed sextic.
For monomials define a bilinear form by
\[\fitmath{[X^iY^j,X^kY^s]=
 \frac{1+(-1)^{j+s}}{(j+s+1)(i+j+k+s+2)} .}\]
This is exactly the triangle integral of their product, by a single
integration in each coordinate. Use the following ordered list of moments:
\[\fitmath{\begin{aligned}
 m&=[p,p],\quad z=[q,q],\quad u=[p_Y,p_Y],\quad v_x=[q_X,q_X],
 \quad w=[p_X,p_X]+[q_Y,q_Y],\\
 d&=[p_X-q_Y,p_X-q_Y],\quad h_s=[p_Y,q_X],\\
 l_x&=[p_X,\psi_X],\quad l_y=[p_Y,\psi_Y],\quad
 l_Q=[q_Y,\psi_X]+[q_X,\psi_Y],\quad f=[p,\psi],\quad p_0=[p,1].
\end{aligned}}\]
In the heptagonal column $(a,b,d)/N$ means $(a+bc+dc^2)/N$.
Thus the heptagonal column lies in $\mathbb Q(c)$, with basis
$1,c,c^2$; both columns for the quintics, after the common factor
$\sigma$ is divided out, lie in $\mathbb Q$. Termwise integration gives
\[\fitmath{\begin{array}{c|l}
&7\\\hline
m&(8836164924,-11597944585,249444819030)/871995801600\\
z&(402031002,-854929883,40830974578)/406931374080\\
u&(4701146047,-2949079185,-1572886440)/871995801600\\
v_x&(-79362673,213401295,1122518280)/96888422400\\
w&(4386558889,-2855853141,80963692320)/87199580160\\
d&(-1376134913,5784249525,25487095560)/217998950400\\
h_s&(354341427,-537371345,1188315040)/290665267200\\
l_x&(108661,338438,-274056)/516096\\
l_y&(16257,38542,-50568)/516096\\
l_Q&(-2489,1347410,11400)/1290240\\
f&(237103,4184190,-715780)/23654400\\
p_0&(5986,176655,-18990)/322560\\
\end{array}}\]
\[\fitmath{\begin{array}{c|ll}
&8&9\\\hline
m&61433686441/217998950400&981999977239/3487983206400\\
z&335753333/3460300800&13423369217/139519328256\\
u&958649/1238630400&1418071/2642411520\\
v_x&113540743/6606028800&29276095/1409286144\\
w&3754074491/3963617280&60049913633/63417876480\\
d&912007829/6606028800&14850513287/105696460800\\
h_s&4257929/1238630400&125648359/39636172800\\
l_x&10957/26880&8259/20480\\
l_y&73/4032&6511/430080\\
l_Q&4265/4032&689363/645120\\
f&44243/276480&6185363/38707200\\
p_0&66251/129024&331109/645120\\
\end{array}}\]
For example, if $p=\sum p_{ij}X^iY^j$ and
$q=\sum q_{ks}X^kY^s$, the mixed derivative entry in every column is
\[\fitmath{h_s=\sum_{i,j,k,s:\;jk>0}
 \frac{jk[1+(-1)^{j+s-1}]p_{ij}q_{ks}}
 {(j+s)(i+j+k+s)}.}\]
Here $i,j$ index the coefficients of $p$ and $k,s$ those of $q$. The other products use the same form. At seven
$8c^3=4c^2+4c-1$, and
$\Re\overline z^6=\sum_{j\ {\rm even}}\binom6j(-1)^{j/2}
 (cX)^{6-j}(1-c^2)^{j/2}Y^j$ for $z=c(X+i\sqrt{t}Y)$.
The scaled imaginary part uses odd $j$ and factors
$-c(-1)^{(j-1)/2}(1-c^2)^{(j-1)/2}$ in place of the even powers.
These reductions give the first column in the cubic field. For the angle
variable $t$ the defining equations and isolating intervals are
\[
\begin{array}{c|l|c}
n&\text{polynomial vanishing at }t&10^5t\ \text{lies in}\\\hline
7&t^3-21t^2+35t-7&(23191,23192)\\
8&t^2-6t+1&(17157,17158)\\
9&t^3-33t^2+27t-3&(13247,13248).
\end{array}
\]
The first identity follows from $\sin(7\pi/7)=0$ by dividing the
multiple-angle formula by $\sin(\pi/7)\cos^6(\pi/7)$. The middle
root is $3-2\sqrt2$. At nine sides, $\tan(3\pi/9)=\sqrt3$ gives
$t(3-t)^2=3(1-3t)^2$. Each polynomial is strictly monotone on
$[0,1/4]$ and has opposite signs at the two indicated rational
endpoints. The required angles have $t<1/4$, for example by the cosine
bounds below. At seven and nine sides $\mathbb Q(t)=\mathbb Q(c)$,
since $t=c^{-2}-1$ and the displayed cubics for $t$ are irreducible
by the rational-root test. At eight sides $c$ lies in a quartic field.
In particular the quadratic energy and mass at all three side counts
lie in $\mathbb Q(t)$; moments carrying a single factor of $1/\sigma$
use $\mathbb Q(c)$.

Set $M'=m+tz$, $E'=w+u/t+tv_x$. Then $M=\sigma^2 M'$,
$E=\sigma^2E'$ and the Rayleigh quotient is $(1+t)E'/M'$.
For the corrected trial field, the coefficients in the angular estimate are
\[\fitmath{\begin{gathered}
 B=u/(tE'),\quad C=tv_x/E',\quad v=tz/M',\quad
 b=.03l_y/(\sigma t E'),\quad e=.03 f/(\sigma M'),\\
 D=B-2b+.0004/(tE),\qquad \kappa=.03l_x/(\sigma E')+b-e,\qquad
 d_g=(d-2h_s)/E',\qquad \|g_1-g_4\|^2=d/E'.
 \end{gathered}}\]
The three remaining linear moments are $L_X=l_x/(\sigma E')$,
$L_Y=l_y/(\sigma E')$ and $L_Q=l_Q/(\sigma E')$.
In particular the two gradient quantities in the shear estimate follow
directly from the positive energy and the two mixed moments $d,h_s$.

To check the sign using rational arithmetic, $10^{12}c$ lies between $N$ and
$N+1$ for \[\fitmath{N=900968867902,\quad 923879532511,\quad 939692620785,}\] respectively.
Indeed substitution gives opposite signs in $8c^3-4c^2-4c+1$,
$8c^4-8c^2+1$ and $8c^3-6c-1$, respectively. Each increases on
$[.9,1]$ and has just the root $\cos(\pi/n)$ there.
The preceding exact columns then give
\[\fitmath{\begin{array}{c|rrr}
 &7&8&9\\\hline
 \sqrt B>&.06971&.06874&.06517\\
 \sqrt C<&.05604&.05559&.05372\\
 \sqrt D<&.02629&.02020&.01350\\
 \sqrt{d/E'}<&.37960&.38030&.38383\\
 d_g<&.13674&.13743&.14068\\
 |\kappa|<&.000029&.000046&.000011\\
 p_0/\sqrt{M'}>&.93018&.93990&.94612
\end{array}}\]

The endpoint bounds for the remaining linear moments are
$7/16<L_X<31/64$, $L_Y<1/40$ and $L_Q>19/16$. Substituting the exact moments gives the following bounds.
\[
\begin{array}{c|cccc}
 n&E&M&\text{transverse fraction}&e\\ \hline
 7&(.7807,.7808)&(.2490,.2491)&(.074978,.075026)&(.017442,.017450)\\
 8&(.814,.815)&(.254,.256)&(.055777,.055781)&(.017410,.017411)\\
 9&(.8421,.8422)&(.2598,.2599)&(.043307,.043314)&(.017335,.017337)
\end{array}
\]

\subsection{Eigenvalue and eigenfield bounds for the regular polygons}\label{seven-spectral}
Normalize the regular circumradius to one. Let $\bar R_n$ be the upper bound for the Rayleigh quotient obtained from the exact moments in $(1+t)E'/M'$, $U_n$ the constant $U(n)$ of Lemma~\ref{IV:lem:uniform-regular-remainder}, and $p_0=3.3899577$. With $s_n=\operatorname{sinc}(2\pi/n)$ this lemma gives
\[
 \mu_1(R_n)>\frac{p_0}{s_n[1+(4.809+U_n/n^2)/(n^3(p_0-1))]}>L_n.
\]
For the last inequality use $s_n<.87102643,.90031632,.92072544$ in order, by the alternating sinc series through degree eight and $3.14159265358<\pi<3.14159265360$.
Put $b=\pi/n$. For $n\ge7$ one may take $G=28/(1+22/(7n))$. To see this, map a disk sector by $(r,\theta)\mapsto r(\cos b,s\theta)$ with $s=\sin(b)/b$. The Jacobian in orthonormal polar frames is constant $s\cos b$. The squared operator norm is at most $1+b$, since the two diagonal entries of the column Gram matrix are at most $1$ and $s^2$ and its off-diagonal has modulus at most $bs^2$. These sector maps agree on their common edges and commute with rotation by $2\pi/n$. By the min--max principle, the second eigenvalue in rotation character one on the polygon is at least the corresponding disk eigenvalue multiplied by $1/(1+b)$. On the disk the second radial eigenvalue in angular order one exceeds 28, and all other angular orders congruent to one modulo $n$ have absolute value $\ge n-1$ and eigenvalue $>28$, by Lemma~\ref{IV:lem:disk-angular-floors}. Now use $b<22/(7n)$.

Normalize both fields to unit energy and choose the first eigenfield $u$ with positive inner product with the unit reference $v$. If $z$ is the complementary energy fraction, then $1/R\le(1-z)/\mu_1(R_n)+z/G_n$, where $G_n=28/(1+22/(7n))$. As $\mu_1(R_n)\ge L_n$ and $R\le \bar R_n<G_n$,
\[
 z\le p=\frac{1-L_n/\bar R_n}{1-L_n/G_n},\qquad
 \|u-v\|_E^2\le \frac{2p}{1+\sqrt{1-p}}<\delta_n^2 .
\]
After normalization to unit mass, the squared sine of the angle between the fields is at most $z R/G_n$ by spectral expansion. The distance between these unit fields is at most $\epsilon_n=1.0001\delta_n\sqrt{\bar R_n/G_n}$, since the sine is less than $.02$. We use the following bounds both on central sectors and when summing over the fan.
\[
\begin{array}{c|rrrr}
n& U_n& L_n&\bar R_n&\delta_n\\\hline
7&92.5&3.86037&112/29&.02380\\
8&85.25&3.74649&3.74722&.01580\\
9&80.75&3.66960&3.67006&.01274
\end{array}
\]
In the remaining estimates the two fields have equal energy: their masses satisfy $M_*=(R/\mu_1(R_n))M_V$. Projections onto individual gradient components contract in energy. For the mass fraction, the distance bound for fields of unit mass gives $\sqrt{v_*}\le\sqrt{v_V}+\epsilon_n$. For a linear moment $e=I(P\Phi)/M$ the mass normalization factor is $\sqrt{\mu_1(R_n)/R}$, so $|e_*-e_V|\le\epsilon_n\sqrt{I(\Phi^2)/M_V}+2\delta_n^2|e_V|$. The centered variance scales as $U_*/U_V=\mu_1(R_n)/R$. Also $I(P_*)/\sqrt{M_*}\ge I(P_V)/\sqrt{M_V}-\epsilon_n$ by Cauchy--Schwarz. These estimates fix all normalizations in the shear test.
\section{Angular comparison for seven, eight, and nine sides}
\label{seven-central}
With $x=\tan(a)/\tan(\pi/n)$ write $s=(T-H_1)/(x-1)^2$ at zero skew, continued at $x=1$. We first bound the transverse mass of the normalized first field.
\paragraph{Transverse mass.} Write $a=2/n$, $\sigma=\operatorname{sinc}(\pi a)$ and $A_j=(a)_j/j!$. Section~\ref{seven-fourier} gives the radial inequality and the two adjacent-order estimates. In the adjacent sine estimate retain the factor $1-n/(n^2-1)$. Indeed $u_{-1}-u_1\ge(u_{-1}+u_1)/n$ by the positive Green-kernel ratio, $U_1/2\le nu_0/(n^2-1)$ for $U_1=u_{-1}+u_1$, and $t_0-t_2\ge(u_0-U_1/2)(W_0-W_2)$ for the Fourier coefficients of $W\operatorname{Re}(e^{-i\theta}Y)$. At the boundary, the beta integral gives $M_2(1)=a(a+1)/[(1-a)(2-a)]$ by two successive beta recurrences, and $M_1(1)=a/(1-a)$. Thus the multiplier of the adjacent pair is at least
\[\fitmath{J=1-\left(\frac a{1-a}\right)^2+\frac12\left(1-\frac n{n^2-1}\right)
       \left(1-\frac{a(a+1)}{(1-a)(2-a)}\right).}\]
Here the factor $1/2=1/(na)$ is the ratio of the sine contribution to the factor $\sigma a^2/(1-a^2)$ in the cosine contribution. The radial order required in this argument follows from $m=|1+nj|\ge n-1\ge6$ for nonzero $j$. $rH_m'/H_m\ge 2(1-r^2)\ge1-r^3\ge rg'/g$.

Use $Z_+=[1-(4.809+13.6/n^2)/n^3]^{-1}$, $D=Z_+-1/3.39$ and
\[\fitmath{H_+=\frac a{1-a}\left[1+\left(\frac a{1-a}\right)^2\left(1+\frac1{1-2a}\right)\right],\quad
 K=1-\frac1{3.39Z_+}+\frac{3.29}{(n-1)^2}
 +\frac{2H_+}{n}\left[\frac{6.58}{n-1}+\frac{14.436}{(n-1)^2}\right].}\]
Indeed $Z\le Z_+$ by the regular gamma-product bound. The adjacent coefficient sum satisfies $\sum_{p\ge0}A_pA_{p+1}=[a/(1-a)]\sum_{p\ge0}A_p^2\le H_+$, since $A_p\le ae^a p^{a-1}$ for $p\ge1$, $e^a\le1/(1-a)$, and the decreasing power tail is bounded by its integral. The mass estimate in Section~\ref{seven-fourier} now gives $M/u_0(1)^2\le K$. For regular polygons of area $\pi$, the eigenvalue lower bound exceeds $3.35$ in every row, by multiplying the lower bound at circumradius one by $\operatorname{sinc}(2\pi/n)$.

Set $k_0=n+3$, $R_m=2(k_0+m+3)/[m(m+3)(k_0+3)(k_0+m)]$. In the equations for the coefficients of angular orders $n\pm1$ we keep only the term $W_1u_0\ge au_0(1)r^{n+1}/Z$ of the positive right-hand side and solve the radial equation explicitly, as in Section~\ref{seven-fourier}. The contribution of these two orders is then at least
\[\fitmath{A_* = \frac{\sigma a^3 J}{1-a^2}\frac{3.35(R_{n-1}+R_{n+1})}{Z_+^2K}.}\]
For $h(z)=[1-(1.2+2.4/(z+4))/((z+2)(z+3))]/(z+1)$, the radial integral estimate in Section~\ref{seven-fourier} bounds $\int_0^1 t^zg(\sqrt t)^2dt$ below by $h(z)$. Retain $j=1,\ldots,b$ and $p=0,\ldots,b-1$ in the nonnegative radial series, with $b=3,3,5$ in order. Then
\[\fitmath{v\ge (1-\sigma)/2+\sum_{j=1}^b F_j+A_*,\qquad
 F_j=\frac{\sigma a^2}{D(j^2-a^2)}\sum_{p=0}^{b-1} A_p A_{p+j}h((2p+j)/a).}\]
The following are integer lower bounds after multiplication by $10^6$. Each column is a sum of at most five rational terms after replacing $\sigma$ by its lower sine bound in the positive terms and by its upper bound in $(1-\sigma)/2$. At $h_0=3.14160a$, use the alternating lower bound through degree six for sinc, and at $l_0=3.14159a$ the upper bound through degree eight.
\[\fitmath{\begin{array}{c|r|rrrrr|r|c}
 n&(1-\sigma)/2&F_1&F_2&F_3&F_4&F_5&A_*& v>\\\hline
 7&64486&7172&660&161&-&-&417&.07285\\
 8&49841&4371&388&92&-&-&257&.05491\\
 9&39637&2835&247&58&20&9&154&.04295
\end{array}}\]
A dash denotes a term omitted by choice of $b$. For example $F_1$ uses $A_{p+1}=A_p(p+a)/(p+1)$ with $A_0=1$.

\paragraph{Angular lower bound on central sectors.} From Section~\ref{seven-moments} take
$\underline{\sqrt{B_V}}=(.065,.068,.065)$,
$\overline{\sqrt{C_V}}=(.057,.056,.0538)$ and
$\overline{\sqrt{D_V}}=(.02653,.02021,.013526)$.
Write $\delta=\delta_n$ from Section~\ref{seven-spectral}. Set
\[\fitmath{B=(\underline{\sqrt{B_V}}-\delta)^2,\quad C=(\overline{\sqrt{C_V}}+\delta)^2,\quad
 D_0=(\overline{\sqrt{D_V}}+\delta)^2.}\]

Use $v$ equal to the transverse lower bound in the last column, $e=.0178,.0177,.0175$ and $U=.00020$. The inequalities are in the required direction by the following derivative test. The formulas for $m,A_J,p$ follow by expansion in Section~\ref{seven-angular}. Write $l=.67$.
Write $z=x-1$, $J=t/(1+t)+zt(1-3t)/[3(1+t)^2]$ and
$\zeta=xp/(2m)$. Eliminating $b$ by $2b=B+H-D$, let
\[
\begin{aligned}
 m&=1+(x^2-1)v+2ze+z^2U,\qquad K=1+2(C-B-v),\\
 A_J&=(x+2)(v-C)-KJ+B+H-(1+1/x)D-x(G-U),\\
 p&=(x+1)(2v-C)+2e+(1+1/x)B+(1-1/x)H+(x^{-2}-1)D-z(G-2U).
\end{aligned}
\]
For the first eigenfield, the reference moments and projection bounds give
$v<.0813,.0595,.0459$, using $\sqrt{v_*}\le\sqrt{v_V}+\epsilon$;
$ .0171<e<.0178,.0177,.0175$ in order, with $e>.0172$
at eight and nine; and $.00020<U<.00023$, since
$U_V=.03^2(14/225)/M_V$ and $U_*/U_V=\mu/R$.
The same intervals contain the values for the polynomial reference field. The eigenfield bounds on $D$
follow from $\sqrt{D_*}\le\sqrt{D_V}+\delta$.
For both fields the derivative of the quotient $s$ with respect to $v$ is positive:
after discarding the nonnegative $2j$, it is at least
\[
 2-r[(r+1)(.54)/.93-1]
 -r(1-l^2)(.54)^2/[4(.93)^2]>.8 .
\]
Here $m>.93$, $0<p<.54$ hold throughout all successive replacements.
For example, in the upper bound for $p$ use $x+1\le2.78$,
$v<.0813$, $C\ge0$, $B<.009$, $0\le D\le.0026$,
$.0171<e\le.0178$, $.00020\le U\le.00023$,
$G<.0028$, $H<.0036$, and $|z|\le.78$.
For the lower bound use $(1+l)(2(.04295)-.0066)$ for its first term.

Set $S_J=A_J-xp^2/(4m)$. Its derivatives in $B,C,D$ respectively are
$1+2J-(1+1/x)\zeta$, $-x-2-2J+(x+1)\zeta$ and
$-(1+1/x)[1-(1-1/x)\zeta]$.
Its derivatives in $e,U$ are
$-2\zeta[1-(1-1/x)\zeta]$ and $x[1-(1-1/x)\zeta]^2$.
Thus $S_J$ increases with $B,U$ and decreases with $C,D,e$:
in particular $(1+1/x)\zeta\le(r+1)(.54)/(2(.93))<1$.
Let $t=\tan^2(\pi/n)$, $G_0=.00216/E$, $H_0=.0004/(tE)$, $l=.67$, $r=1.78,1.64,1.59$ and $w=1/80,1/100,3/400$ respectively. The arctangent remainder from Section~\ref{seven-angular} supplies $w(x-1)^2$: indeed $[1-2(B_++v)]t^2[3-t+(1-3t)l]/[4(1+t)^2(1+tr^2)]>w$ with $B_+=.009,.0075,.0062$. Define
\begin{gather*}
 \eta=2v-C+e+B,\quad K_0=1+2(C-B-v),\quad P=2v-C-G_0+2U,\quad d=(1+\eta)D_0,\\
 S=3(v-C)-K_0t/(1+t)+B+H_0-2D_0-G_0+U-\eta^2,\\
 L=v-C-K_0t(1-3t)/(3(1+t)^2)+D_0-G_0+U
       -\eta(P-B+H_0-2D_0)-\eta^2(1-2v-2e).
\end{gather*}
Completing the square at $x=1$ gives $s\ge S+(x-1)L+(x-1)^2 Q$, where
\[\fitmath{Q=w-\eta P+\eta^2[2(e-U)+(x+2)(v+U)]-d/x-x R^2/(4m),}\]
\[\fitmath{m=1+(x^2-1)v+2(x-1)e+(x-1)^2U,\quad
 R=P-(B-H_0)/x-D_0(x+1)/x^2-2\eta[(x+1)v+2e+(x-1)U].}\]
The identity follows from $p-2\eta m=(x-1)R$ in the formula for $s$.

\paragraph{Endpoint bounds on the central interval.}
Since $m\ge1+(l^2-1)v+2(l-1)e$, take
$m_*=.94,.95,.96$, respectively. The term involving $R$ comes from the remaining square after matching the value and slope at $x=1$.
We bound it on the whole interval. Indeed
\[
 \begin{aligned}
 R&\ge P-(B-H_0)_+/l-D_0(1+l)/l^2
              -2\eta[(1+r)v+2e+(r-1)U]>.046,\\
 R&\le P+(H_0-B)_+/l-D_0(1+r)/r^2
              -2\eta[(1+l)v+2e+(l-1)U]<R_* ,
 \end{aligned}
\]
where $R_*=.086,.070,.059$, in order. Here $u_+=\max(0,u)$.
Set
\[
 C_*(x)=w-\eta P+\eta^2[2(e-U)+(x+2)(v+U)]-d/x
                           -xR_*^2/(4m_*).
\]
Then $Q\ge C_*(x)$. Since $C_*''=-2d/x^3<0$, the lesser of
$C_*(l),C_*(r)$ bounds the quadratic coefficient on the whole
aperture interval.

For the arithmetic write $a=.0004$, $g=.00216$, $y=t/(1+t)$,
$P_0=2v-C+2U$, and put
\[
 S_0=3(v-C)+B-2D_0+U-\eta^2,\qquad
 L_0=v-C+D_0+U-\eta(P_0-B-2D_0)-\eta^2(1-2v-2e).
\]
Thus $P=P_0-g/E$ and
\[
 S=S_0-K_0 y+(a/t-g)/E,\qquad
 L=L_0-K_0y(1-4y)/3-[(1-\eta)g+\eta a/t]/E .
\]
Thus $S>s_0$ is equivalent to
\[
 p_{s_0}(t,E)=tE[(1+t)(S_0-s_0)-K_0t]+(1+t)(a-gt)>0.
\]
It is affine in $E$ and concave in $t$, since $S_0-s_0<K_0$.
It suffices to test the four corners of the rational $t,E$ intervals
in Section~\ref{seven-moments}. The third column is an integer lower
bound for $10^9p_{s_0}$ at all four corners.
\[
\begin{array}{c|r|r|c|rr}
 n&s_0&10^9p_{s_0}>&L\ \text{lies in}&C_*(l)>&C_*(r)>\\\hline
 7&.00740&2960&(.01128,.01132)&-.0092&-.0066\\
 8&.00307&1320&(.00460,.00464)&-.0031&-.0022\\
 9&.00046&620 &(.00069,.00072)&-.0008&-.0006
\end{array}
\]
For example, for the last row the floors of $10^9p_{s_0}$ are
$1581$, $1566$, $643$ and $628$, at the corners
$(t_-,E_-)$, $(t_-,E_+)$, $(t_+,E_-)$ and $(t_+,E_+)$, respectively.
The last three columns use $y_\pm=t_\pm/(1+t_\pm)$:
bound $y(1-4y)$ between $y_-(1-4y_+)$ and $y_+(1-4y_-)$,
and use $P\le P_0-g/E_+$ in $C_*$. All these replacements involve
positive denominators.

Let $c_*=.0092,.0031,.0008$ be the negatives of the penultimate
column. Then $s\ge S+(x-1)L-c_*(x-1)^2$.
For $x\le1$ this exceeds
$s_0-.33L_+-.33^2 c_*>.00263,.00118,.00012$, respectively. For
$x\ge1$ one has $L_--(r-1)c_*>0$, so $s\ge S$.
This proves the strict comparison at zero skew for the first eigenfield whenever the central aperture differs from one. For nonzero skew we use the three convex upper bounds in Section~\ref{seven-shear}. The energy upper bounds $.02,.14,.482$ lie below the corresponding mass lower bounds $.23,.23,.494$ for nonzero shear.

\paragraph{The polynomial field on central sectors.}
The same estimates apply to every central sector when the polynomial trial is used throughout the fan. Its transverse fraction is at least $.07497,.05577,.04330$ respectively by direct integration. Each value exceeds the corresponding lower bound for $v$ in the transverse mass table above. The bounds on $B,C,D_0,e,U$ and the three shear majorants include the polynomial field by taking its projection displacement to be zero. The polynomial tangent additionally has $\kappa=I(P_X\Phi_X)/E+b-e$. At fixed original moments put $\zeta=xp/(2m)$ with $p$ evaluated at $\kappa=0$. The resulting change in $s$ is
\[\fitmath{-2(1+j-\zeta)\kappa-x\kappa^2/m.}\]
The same rectangle has $0<\zeta<.52$, $j<.212$ and $m>.93$, so the loss is at most $2.5|\kappa|$. The exact moments in Section~\ref{seven-moments} yield $|\kappa|<.000037,.000046,.000012$. Consequently the angular lower bounds to the left of one remain above $.00253,.00106,.00009$ and the bounds to the right of one are larger. At $x=1$ the tangent and the zero-skew estimate coincide. Nonzero skew gives strict inequality. This establishes $T_i-H_{1,i}\ge0$ for each complementary central sector, with the same polynomial tangent on all sectors.

\section{A lower bound on one narrow sector suffices for the whole fan}\label{sec:seven-narrow}
Write $k=\tan(\pi/n)$, $t=k^2$, $h_*=67/100$, $q_*=449/1000$,
$\tau=x^2\le h_*^2<q_*$ and $u=r\tau$. In this section $r=0$
or $r=\rho_n^2$, with $\rho_n=.230,.129,.085$, respectively.
These values correspond to the endpoints of the allowed skew interval. We prove
\[\fitmath{T(x)-H_1(x,u)>b_n,\qquad
 (b_7,b_8,b_9)=(.00368,.00219,.00138).}\]
The additional bound $|z|<.16,.105,.075$ is not needed on this branch.

\paragraph{A quadratic lower bound for the mass.}
Let $P+ikQ$ be the regular polynomial, $h=P(X,X)$,
$Z=X^2-Y^2$, $j=Yh/X$. Use the four profiles
\[\fitmath{p=h+\tau tZ/4,\quad q=j+\tau tXYZ,\quad r_0=j,\quad
 \nu=p+3XZ/2.}\]
On both rays the correction terms vanish. All these profiles have
the same upper-ray trace $h$, with the required sign on the lower ray.
Write $I$ for integration on the reference triangle, $c_0=I(P)$,
$I_c(f)=I(f^2)-(I(f)-c_0)^2$ and $y=t\tau$. With
$E=I(P_X^2+P_Y^2/t+tQ_X^2+Q_Y^2)$ and $M=I(P^2+tQ^2)$, put
\[\fitmath{m=\{I_c(p)+uI(r_0^2)+y[I(q^2)+uI_c(\nu)]\}/M.}\]
For $J(f,g,w)=I[(f-ug)^2+uw^2]$, define
\[\fitmath{\begin{aligned}
 J_1&=J(p_X,r_{0,Y},r_{0,X}-p_Y),&
 J_2&=J(p_Y,r_{0,X},r_{0,Y}-p_X),\\
 J_3&=J(q_X,\nu_Y,\nu_X-q_Y),&
 J_4&=J(q_Y,\nu_X,\nu_Y-q_X).
\end{aligned}}\]
Then $D=(1-u)^2+(J_1+J_2/y+yJ_3+J_4)/E$ is a polynomial in $\tau$
(the numerator of $J_2/y$ is divisible by $\tau$). The fan formula is
$H_1=xD^2/[4m(1-u)^3]$. Let $m_0=m(0)$ and $D_0=D(0)$,
evaluating after $u=r\tau$.

We compare the coefficients as follows. $A$ bounds
$D_0^2/(4m_0)$ from above. The entries in columns four and five are $10^6b_1$ and $10^6b_2$, where $b_1,b_2$ bound the first two nonconstant coefficients of $m/m_0$ from below.
Higher nonzero mass coefficients are positive. The vector in the
last column gives upper bounds for the nonconstant coefficients of $D/D_0$
in millionths, in increasing degree. An omitted trailing coefficient is zero.
\[\fitmath{\begin{array}{cr|rrr|l}
 n&r&A&10^6b_1&10^6b_2&10^6[\tau^{1,2,\ldots}](D/D_0)\text{ upper bounds}\\\hline
7&0&1.083596&170263&12662&43171,1519,3988\\
7&\rho_7^2&1.126985&187908&56924&-130226,45439,11027,303,1\\
8&0&1.061451&123975&6742&31356,668,1541\\
8&\rho_8^2&1.079349&129525&16791&-23585,8580,2596,34,1\\
9&0&1.047179&94511&3933&23867,324,684\\
9&\rho_9^2&1.057081&96920&7235&16,2661,939,7,1
\end{array}}\]
Here and below the decimal endpoints are rational. One way to reproduce
the table uses only univariate trace polynomials. In the three rows
$h/\cos(\pi/n)$ is respectively
\[\fitmath{X+(-465X^3+15X^5+52X^6)/1024,\quad
 X+(3X^2-493X^3+15X^4+78X^5)/1024,\quad
 X+(-8X^2-1843X^3-77X^4+350X^5)/4096 .}\]
These traces also verify positivity of the mass for every admissible skew.
On $0<X\le1$, retaining only the negative coefficients gives
$h(X)>.46X$ in all three cases, since $\cos(\pi/n)>.9$.
Thus $I(h)>.46(2/3)>.30$. By the moment bounds in
Section~\ref{seven-moments}, $0<c_0\le\sqrt M<.51$. As
$Z=X^2-Y^2\ge0$ and $\nu\ge p\ge h$, we have
$I(\nu)\ge I(p)>c_0/2$. It follows from
$I_c(f)=\operatorname{Var}(f)+2c_0 I(f)-c_0^2$ that $m>0$.
Every column follows from the same integration formula
\[\fitmath{[\tau^a]I(fg)=
 \sum_{i,j,k,l,\ b+d=a}f_{ijb}g_{kld}
 \frac{1+(-1)^{j+l}}{(j+l+1)(i+j+k+l+2)}
 \quad (f=\textstyle\sum f_{ijb}X^iY^j\tau^b).}\]
Apply it in the displayed $m,D$, with
$[\tau^a]I_c(f)=[\tau^a](I(f^2)-(I(f)-c_0)^2)$.
For instance at $r=0$, with $p_1=tZ/4$, $q_1=tXYZ$, the constant
and first two coefficients of $Mm$ are $I_c(h)$,
$2I(hp_1)-2(Ih-c_0)I(p_1)+tI(j^2)$, and
$I(p_1^2)-I(p_1)^2+2tI(jq_1)$. For $r>0$ add
$r\tau I(j^2)+tr\tau^2 I_c(p+3XZ/2)$ to this same mass polynomial.
In order down the table the cubic mass coefficients, divided by $m_0$,
exceed $.00087,.00263,.00033,.00062,.00014,.00022$. The respective
quartic coefficients for $r>0$ exceed $.000012,.0000014,.00000028$.
All other higher mass coefficients vanish. The moment sum also gives
$c_0,E,M$ from the sector polynomials in Section~\ref{seven-moments}. To bound the rational
functions of $c=\cos(\pi/n)$ one may insert the isolating endpoints
given in the regular-moment calculation, using $t=c^{-2}-1$.
No value of the trial at a spatial node enters this identity.

In the six rows put respectively
\[\fitmath{10^6 d=(44660,-107570,31970,-19200,24155,1405),\qquad
 S=(.079,.244,.059,.1175,.0455,.072).}\]
If $a_1,\ldots,a_s$ are the upper coefficient bounds of $D/D_0$,
then $d>a_1+\sum_{j\ge2}a_j q_*^{j-1}$. Hence $D/D_0\le1+d\tau$.
Bernoulli's inequality gives
\[\fitmath{(1-r\tau)^3m/m_0\ge 1+B\tau+C\tau^2,\quad
 B=b_1-3r,\qquad C=b_2(1-3rq_*)-3rb_1.}\]
In every row $C>0$. To obtain $H_1\le Ax(1-S\tau)$ it suffices that
\[\fitmath{B-S-2d+v(C-SB-d^2)-SCv^2\ge0 \quad(0\le v\le q_*).}\]
This polynomial is concave. The smaller of its two endpoint values,
in units of $10^{-6}$, exceeds $480,340,230,360,230,195$ in row order.
This proves the bound on the whole interval at both endpoint values of the skew.

\paragraph{A lower bound uniform in the skew.}
The polynomial tangent has the form
$T=1-\gamma\int_x^1(1+ts^2)^{-1}ds$, where $0<\gamma<g$ with
$g=1.04284,1.03716,1.03131$, respectively. These upper bounds follow from $\gamma=(1+t)[1+2(C-B+\kappa-v)]$ using the exact moment columns in Section~\ref{seven-moments}.
Use
\[\fitmath{(l,t_+)=(.23190,.23192),\ (35/204,169/985),\ (51/385,31/234)}\]
for the respective lower and upper endpoints for $t$. Define
$U(v)=\sum_{j=0}^4(-1)^j t_j^j(1-v^{2j+1})/(2j+1)$ with
$t_j=t_+$ for even $j$ and $t_j=l$ for odd $j$.
The finite geometric expansion bounds the integral from $v$ to one by
$U(v)$. At zero skew the comparison function
\[\fitmath{f(x)=1-g\int_x^1(1+ts^2)^{-1}ds-Ax(1-Sx^2)}\]
is convex, since $f''(x)=2x[3AS-gt/(1+tx^2)^2]\ge0$.
Indeed $3AS-gt_+>.0149,.0099,.0063$. Its derivative at $h_*$
is at most $g/(1+lh_*^2)-A+3ASh_*^2<0$ in each row.
Its minimum is therefore at least $1-gU(h_*)-Ah_*(1-Sh_*^2)$. The latter is greater than
$.003682,.002199,.001386$.

At the other endpoint $r=\rho_n^2$, use
$\int_0^x (1+ts^2)^{-1}ds\ge x-t_+x^3/3$.
Then $f(x)\ge Z-U_1x+Vx^3$, where
$Z=1-gU(0)$, $U_1=A-g>0$ and $V=AS-gt_+/3>0$.
For each $b=.0067,.0038,.0018$ respectively,
$Z>b$ and $27V(Z-b)^2>4U_1^3$. These inequalities say that the
minimum of the cubic on the nonnegative axis exceeds $b$.

Concavity extends the estimate to intermediate skews. At
fixed $x$, let $e$ denote the normalized uncentered energy and set
$(1-u)^2e=e_0+n_1u+n_2u^2$, $m=m'_0+m'_1u$.
With $v=u/(1-u)$, the condition $H_1\le T-b_n$ is equivalent to
\[\fitmath{2\sqrt{(T-b_n)/x}\sqrt{m'_0+(m'_0+m'_1)v}
    -(1+e_0)-(2e_0+n_1)v-(e_0+n_1+n_2)v^2\ge0 .}\]
The coefficient $e_0+n_1+n_2$ is nonnegative by the sum-of-squares
energy at $u=1$ after its denominator is cleared. Both endpoint masses are positive, so their affine interpolation is positive throughout. The endpoint estimates also show that the target $T-b_n$ is positive and that the displayed left side is strictly positive at both ends. The left side
is concave, giving strict positivity throughout the interval.

\section{Summation and the extremizer}\label{paper-seven-sum}
For the first eigenfield on a central sector, the zero-skew inequality
is strict for $x\ne1$. The signed shear inequality of Section~\ref{seven-shear}
is strict for $z\ne0$. Together with the bounds on the centroid sectors, these inequalities give
$H_{1,i}\le T_i$, with equality only at $(x_i,z_i)=(1,0)$.
The same central-sector inequalities hold for the polynomial tangent,
using its own slope correction $\kappa$.
If every aperture satisfies $x_i\ge.67$, use the first eigenfield of the regular polygon on the entire fan. Then $R=\mu_1(R_n)$.
Unless all sectors have $x_i=1,z_i=0$, the fan inequality gives
$\mu_1(P)<\mu_1(R_n)$. In the exceptional configuration the angles
between successive rays are $2\pi/n$ and all vertex radii are equal,
which determines a regular polygon.

If at least one aperture has $x_i<.67$, select the polynomial reference on
the entire polygon. Its central and narrow profiles have the same trace on their common rays.
Each narrow sector contributes $T_i-H_{1,i}>b_n$; all remaining
sectors have $T_i-H_{1,i}\ge0$. The bounds $L_n,\bar R_n$ from Section~\ref{seven-spectral} control the Rayleigh quotient error over the entire fan, as the table shows.
\[
\begin{array}{c|r|c|c}
n&b_n& n(1-L_n/\bar R_n)&b_n-n(1-L_n/\bar R_n)\\\hline
7&.00368&4927/1600000&961/1600000\\
8&.00219&292/187361&11832059/18736100000\\
9&.00138&207/183503&2311707/9175150000
\end{array}
\]
It follows that $\sum_iH_{1,i}<n-b_n<nL_n/R$.
The fan inequality for this field gives $\mu_1(P)<L_n<\mu_1(R_n)$,
contrary to competitiveness. Similar copies of $R_n$ give equality by
scaling. This proves the comparison at seven, eight and nine sides.

\clearpage\part{Geometric, shear and Fourier estimates for the centroid trials}\label{part:sector}
\section{Bounds on centroid sectors for seven, eight, and nine sides}\label{seven-window-proof}
Normalize the area of the strictly convex polygon $P$ and its regular
comparison polygon to $\pi$. A polygon is competitive if its first positive
Neumann eigenvalue is at least that of the regular polygon of the same area.
Write $a$ for the sector half-angle at the area centroid, $b=\pi/n$,
$x=\tan(a)/\tan b$, and $z=(r_+-r_-)/(r_++r_-)$ for the two
centroid-to-vertex radii $r_+,r_-$. We derive the aperture and skew bounds simultaneously for $n=7,8,9$. The first
eigenfield of $R_n$ has rotation character one, by the
proof of Lemma~\ref{IV:lem:regular-character}. All geometric estimates in
this section concern convex sectors.

For a convex fan of area $\pi$ put $T=r(P,O)^{-2}-1$, $q=np$,
$y=(n-q)/(n-1)$ and $D_j(q)=q^j+(n-1)y^j-n$. Put
\[
 \begin{gathered}
 a_2=144/n^3,\quad a_4=4.8n-216/n^3,\quad a_6=12.4/n+144/n^3,\\
 A(q)=36/n^2+\sum_{j=2,4,6}a_jD_j(q).
 \end{gathered}
\]
For $\eta_i=b_i/p_i-1$ the one-sector inequality is
\begin{equation}\label{V:eq:finite-angle-budget}
 A(q)+n^5\sum_i p_i[\eta_i-\log(1+\eta_i)+W(p_i,\delta_i/\pi)]
 \le n^5\log(1+T)-4.808n^2-12.4 .
\end{equation}
Apply \eqref{IV:eq:single-aperture-budget} and use $v^8\ge4v^6-6v^4+4v^2-1$ for each normalized aperture $v$, since the difference is $(v^2-1)^4$. This yields the displayed $a_j$ and $A$.
The bound on skew in terms of aperture follows from Corollary~\ref{IV:cor:aperture-skew-product}. The sine rule gives $|z|=|\tan(\pi p)\tan\delta|$, and triangle area gives $r_ir_{i+1}\sin(2\pi p_i)=2\pi b_i$ in area-$\pi$ units.
\subsection{Harmonic tests with radial profiles vanishing on the boundary}\label{V:sec:geo-bubble}
We use the normalization of Proposition~\ref{IV:prop:harmonic-budget-65}. Thus $Q=f'/f'(0)=\sum_{i\ge0}\gamma_i z^i$,
$\gamma_0=1$, $w_i=|\gamma_i|^2/(i+1)$, $X=\sum_{i\ge1}|\gamma_i|^2/(i+2)$,
$\lambda=(j'_{1,1})^2$, $d_i=1/\lambda-H-(1-M_{i+1})$, and the Rayleigh
multiplier is at least $s$. Write
\[
 A_j=\frac{239j^2+2529j+7090}{200(j+2)(j+3)(j+4)(j+5)}.
\]
Lemma~\ref{IV:lem:uniform-harmonic-correlations} gives $1-M_{j+1}\le A_j<1.2/(j+5/3)^2$.
We shall use an arbitrary multiplier $s>12/5$ and coefficient $0\le\beta\le1$ with $5s/12>\beta$.
For each shift $j\ge3$ and each of its signed frequencies $1\pm j$
take the radial profile $v_m(r)=r^m(1-r^2)$, where $m=|1\pm j|$.
Its Dirichlet energy for angular mass $dA/\pi$ is $4/(m+2)$
by integration. Centering subtracts a constant and the signed
frequencies are pairwise distinct, so their Dirichlet energies decouple.
The leading mass pairing $2\int_0^1 g(r)r^{j+m+1}(1-r^2)\,dr$ is bounded
below using $g(r)\ge r[1+(1-r^2)/2]$, giving
\[
 V_{j-1}=\frac{j+4}{(j+1)(j+2)(j+3)},\qquad
 V_{j+1}=\frac{j+5}{(j+2)(j+3)(j+4)}.
\]
Set $G_j=(j+1)\sum_{m=j\pm1}(m+2)V_m^2/4$. Then
\[
 \begin{gathered} G_j=\frac{(j+4)^2+(j+1)(j+3)(j+5)^2/(j+4)^2}
            {4(j+2)^2(j+3)^2},\\
 \frac{G_j}{A_j}\ge
 \frac{100(j+5)^2(j+1)(j+4)}{(j+2)^2(239j^2+2529j+7090)}
 \ge\frac5{12}.\end{gathered}
\]
For the first inequality apply the arithmetic--geometric mean bound
to the two numerator terms, then use
$\sqrt{(j+1)(j+3)}\ge(j+1)(j+3)/(j+2)$. The second reduces to
$j^4+115j^3+798j^2+1124j-4360>0$, already positive at $j=3$
and increasing.

For $k>j$, $g\le1$ bounds the mass pairing from index $k-j$ by
$4/[(2k+m-j+2)(2k+m-j+4)]$.
For $0\le\beta\le1$ put $t_j=\beta A_j/(sG_j)<1$ and
$c=s\beta/(5s/12-\beta)$. In the simultaneous Rayleigh inequality
split each squared mass pairing into its leading term $p$ and tail $z$
using $|p+z|^2\ge t_j|p|^2-t_j|z|^2/(1-t_j)$.
The coefficient of $w_j$ increases by at least $\beta A_j$, and
$s t_j/(1-t_j)\le c$. Applying Cauchy--Schwarz to the tail with weight sum $X$ bounds the coefficient of each
$w_k$ by $(k-j+2)(k+1)$ times the squared mass pairing.
It follows that the diagonal coefficient is at least $1/\lambda-H-L_k$, where
\[
 L_k\le (1-\beta)A_k+cX S_k,\qquad
 S_k=\frac{4(k+1)}{(2k+3)^2}\sum_{j=3}^{k-1}(k-j+2)
 \left[\frac{j+1}{(2k+1)^2}+\frac{j+3}{(2k+5)^2}\right].
\]
For a finite upper cutoff $J$ these estimates hold at $k\le J$.
For $k>J$ the omitted leading contribution is at most $A_k$.
Its total error is at most $\sup_{k>J}A_k\sum_\ell w_\ell\to0$.
Thus the same summed inequality holds for all indices.

The two sums of numerators in $S_k$ are $(k-3)(k^2+12k+8)/6$
and $(k-3)(k^2+18k+20)/6$. Replacing both by their arithmetic mean
increases $S_k$, because the larger is paired with the smaller reciprocal.
With $x=1/(2k+3)\le1/9$ this gives
\[
 S_k\le\frac{(1-x)^2(1-9x)(1+25x)(1+4x^2)}{12(1-4x^2)^2}
 \le\frac{1+16x-225x^2}{12},\qquad A_k\le4.8x^2 .
\]
Indeed $(1-4x^2)^2-(1-x)^2(1+4x^2)=x(2-13x+8x^2+12x^3)\ge0$.

\subsection{Geometric bounds for seven, eight, and nine sides}\label{V:sec:geo-seven-nine}
Normalize the area to $\pi$ and take a strictly convex $n$-gon.
We use a lower bound $\mu_*\le\mu_1(P)$ for which
$H=(\lambda-1)/\mu_*-(1-1/\lambda)\le H_+$. The eigenvalue bound for $R_n$
in Lemma~\ref{IV:lem:uniform-regular-remainder} supplies this bound for
a competitive polygon, and Proposition~\ref{prop:quantitative-window-entry}
supplies it for a small positive deficit. Set $H_+=(4.809+\mathcal R_n/n^2)/(3.3899n^3)$,
where $\mathcal R_n=92.5,85.25,80.75$ at $n=7,8,9$.
Put $U_1=H_+(6.812+37H_+)$, $U_2=H_+(9.13+37.5H_+)$.
Lemma~\ref{IV:lem:two-harmonic-secants} and Proposition~\ref{IV:prop:harmonic-budget-65} give $T\le P(D)$ whenever the coefficient of $w_i$ is at least $D$ for every $i\ge3$, where
\[
 P(D)=H_+/D+
 [1-(.1583-H_+)/D]U_1^2/2+
 [1-(.2169-H_+)/D]U_2^2/3 .
\]
To see the direction of each substitution, write $P(D)=G+E/D$, where
\[
 G=U_1^2/2+U_2^2/3,\qquad
 E=H_+-(.1583-H_+)U_1^2/2-(.2169-H_+)U_2^2/3.
\]
On $D>.23$ the two bracketed coefficients in $P$ are positive, so replacing
$H$ and the two harmonics by their upper bounds increases the estimate.
Here $E>.93H_+$ by substitution, so $P$ decreases with $D$.
Initially $D=.2448-H_+$ gives $T<.025$; at the largest $H_+$ (seven sides)
the lower bound $\lambda/[(1+\lambda H_+/(\lambda-1))(1+.025)]>3.19$
therefore permits $s=3.19$ throughout.

Put $a=1/3.39-H_+$, $X=H_+/a+7H_+^2$ and $\beta=1-10/n^2$. The coefficient estimate following Lemma~\ref{IV:lem:two-harmonic-secants} bounds the weighted tail by $X$.
For $x=1/(2i+3)$ the radial test functions of Section~\ref{V:sec:geo-bubble} bound
the loss by $J(1+16x)-Kx^2$, where
\[
 J=\frac{s\beta X}{12(5s/12-\beta)},\qquad K=225J-4.8(1-\beta).
\]
Completing the square gives the diagonal $D=a-J-64J^2/K$.
To prove that $C$ is an upper bound in the sector inequality, set
$\ell=(4.808n^2+12.4+C)/n^5$. For a proposed lower bound $D_*$ test
\[
 I=K(a-J-D_*)-64J^2,\qquad
 N_0=n^5\{D_*(\ell+\ell^2/2-G)-E\}.
\]
If $K,I,N_0>0$, then $D>D_*$ and $T\le P(D)<\ell+\ell^2/2<e^\ell-1$.
Thus the right side of \eqref{V:eq:finite-angle-budget} is less than $C$.
All entries in this test are rational. In the table $T_+$ also bounds
$\ell+\ell^2/2$, so it may be used in every subsequent estimate.
\[
\begin{array}{c|rrr|rrr|r}
 n&\mathcal R_n&C&D_*&K>&10^5 I>&N_0>&T_+\\\hline
 7&92.5&109&.2762&.81&3&.09&.0215\\
 8&85.25&96.2&.2818&.52&4&.06&.0128\\
 9&80.75&88.6&.2852&.34&1.8&.07&.00835
\end{array}
\]
In particular $A(q)+n^5\sum p_iW_i<C$ for every sector $q$, since
$\eta-\log(1+\eta)\ge0$.

It remains to locate the centroid. With $r=.9$, $v=U_2+U_1^2/2$,
$\Delta=(1-2U_1-4v)/3$, put
\[
\begin{aligned}
 N&=r-U_1r^2/2-U_2r^3/3
  -\sqrt{T_+(1.662-r^2-r^4/2-r^6/3)},\\
 L&=\sqrt{T_+(5(.694)/\Delta-3/2)}
      +[64\Delta(1-U_1-v)]^{-1}+2^{-15}.
\end{aligned}
\]
Lemma~\ref{IV:lem:convex-boundary-two-heads} and Lemma~\ref{IV:lem:physical} imply
$h=|O-G|/\operatorname{dist}(O,\partial P)\le T_+(5/16+L)/N$.
Indeed dropping all divisors $\sqrt{1+T}$ on the resulting upper
bound suffices; use $e^w\le 1/(1-w)$ for $0\le w<1$.
Here $\log 2<.694$ and $-\log(.19)<1.662$ follow, for example,
from the first six and eight positive terms of the exponential series
at $.694$ and $1.662$, respectively.
The pairs $(U_1,U_2)$ are bounded by
$(.0405,.0539),(.0246,.0328),(.0163,.0217)$.
Thus $L<.611,.44,.35$, $N>.784,.815,.834$, and
\[
                 h<.0254,.012,.007
\]
in order. These are the only displacement bounds needed below.

\paragraph{A bound on half the absolute difference between the base angles.}
Use the bounds $C=109,96.2,88.6$ in the three cases.
At $q=.1$ the margins $A(q)-C$ exceed $22,43,61$, respectively;
at $q=u=1.55,1.50,1.49$ they exceed $8,4,15$. These are rational substitutions. For example,
with $v=q-1$, cancellation in $D_j$ can be avoided by using
\[
 D_j(1+v)=\sum_{b=2}^{j}\binom jb \big[1+(-1)^b/(n-1)^{b-1}\big]v^b .
\]
Since $A$ is strictly convex with minimum at one, $.1<q<u$.

Write $s=\max_i|\delta_i|$ for the largest absolute half-difference
of the base angles at the conformal center. We may suppose $s>0$.
Set $L(s)=1-4s^2/\pi^2$ and $M=C\pi^4/(32n^5)$.
In every sector the skew-square estimate gives
$|z_{O,i}|^2\le M/L(s)$. The positive angular term in the sector
realizing $s$ gives $C\ge8n^3 q^2(1-L)/L$, hence
\[
 L(s)\ge\frac{8n^3q^2}{C+8n^3q^2}>.2>3M ,
\]
using only $q>.1$, $n\ge7$, $C\le109$, $\pi<22/7$.
In particular the squared increments of the logarithm of the radius all have coefficient at least
\[
 k(s)=8n L(s)(1-M/L(s))^2/\pi^4 .
\]
Let $q$ now refer to a sector realizing $s$. Its $n-1$ complementary
logarithmic radius increments sum to the negative of its own. By
Cauchy--Schwarz their contribution is at least
$4n^4 k(s)(\operatorname{atanh}z_O)^2/(n-1)$.
Since $|z_O|=\tan(\pi q/n)\tan s\ge(\pi q/n)\tan s$,
retaining the first aperture correction in the skew estimate gives
\[
 \begin{split}
 C&\ge A(q)+B(s)q^2+(256ns^2/\pi^2)q^4,\\
 B(s)&=\frac{n^3s^2}{\pi^2}\left(\frac{32}{L(s)}+\frac{41}{25}\right)
       +\frac{4\pi^2 n^2}{n-1}k(s)\tan^2s .
 \end{split}
\]
Indeed the first term in the skew bound of Lemma~\ref{IV:lem:skew-fan} retains the factor
$(1-4p^2)^{-2}\ge1+8p^2$ on $32/L(s)$; use $p=q/n$ and $L(s)\le1$ in the added term.
The last term of $B$ comes from the other radius increments. The function $B$
increases up to the largest admissible $s$. Its last term is the only one
requiring a check: $L\tan s$ increases since
$\pi^2-4s^2\ge4s\sin(2s)$, and $L^{-1/2}-ML^{-3/2}$
increases as $L$ decreases on $L>3M$. Their product is proportional
to $\sqrt{k(s)}\tan s$.

To exclude $s\ge d$, choose integers $b<B(d)$, $c<256nd^2/\pi^2$
and put $F(q)=A(q)+cq^4$. Convexity gives
\[
 C\ge F(q)+bq^2\ge F(w)-wF'(w)-F'(w)^2/(4b).
\]
The following rational contact points suffice. The last two columns only
require $D'_j(w)=j[w^{j-1}-((n-w)/(n-1))^{j-1}]$.
\[
\begin{array}{c|rr|r|rr|rr}
n&d&w&k(d)>&b&c&F(w)-wF'(w)>&|F'(w)|<\\\hline
7&.42&.38&.51&254&32&143&185\\
8&.29&.53&.62&171&17&138&167\\
9&.22&.62&.71&136&11&135&158
\end{array}
\]
For example $143-185^2/(4\cdot254)>109$; the other two rows give
$138-167^2/(4\cdot171)>96.2$ and $135-158^2/(4\cdot136)>88.6$.
Hence $s<d$ in every sector. For the $k(d)$ column insert
$1-4d^2/(3.14159)^2$ for its lower $L$ and $C(3.14160)^4/(32n^5)$
for its upper $M$ in $8nL(1-M/L)^2/\pi^4$. With this column in $B(d)$, using $3.14159<\pi<3.14160$ gives respectively $B(d)>254.1,171.6,136.7$.

To finish the transfer, use $k_*=(.51,.62,.71)$ from this table and the
displacement $h$ above. The sum of the squares of an increment and its
$n-1$ complementary increments is at least $n/(n-1)$ times the
square of the first. Therefore
$|z_O|\le\tfrac12\sqrt{(n-1)C/(n^5k_*)}<Z$, with $Z=.14,.092,.066$.
The following rational bounds suffice for the three remaining estimates.
On $0\le v\le.70$ the alternating sine and cosine series give
\[
 T_-(v)=\frac{v-v^3/6+v^5/120-v^7/5040}{1-v^2/2+v^4/24}
 \le\tan v\le
 T_+(v)=\frac{v-v^3/6+v^5/120}{1-v^2/2+v^4/24-v^6/720}.
\]
Use $T_-(d)$ also in the lower estimate of $B(d)$. Set $p=1-h$, $t=T_+(d)$,
$\pi_-=3.14159$, $\pi_+=3.14160$. Section~\ref{seven-transport} bounds
the three quantities on the left by the corresponding expressions on the right:
\[
 x\le\frac{T_+(\pi_+u/n)}{pT_-(\pi_-/n)},\qquad
 |z|\le\frac{Z+h/(2-h)}{1+Zh/(2-h)},\qquad
 \frac{|z|}{x}\le T_+(\pi_+/n)\,
 \frac{(1+p^2)t+(1-p^2)(1+t^2/2)}{2p}.
\]
Using $\sqrt{1+t^2}\le1+t^2/2$ gives the three respective bounds
$(1.779,1.633,1.585)$, $(.153,.099,.070)$ and $(.229,.1289,.0841)$,
with each triple ordered by $n$. These prove the bounds on the centroid sectors in
Part~\ref{part:seven}.

\section{Corrected energy and mass on a central sector}\label{seven-angular}

For normalized energy $e$ and mass $m>0$ at zero skew put $d=x-1$,
$H=x(1+e)^2/(4m)$ and $p=(2m-1-e)/d$. For any tangent comparison $T$ one has
\begin{equation}\label{V:eq:general-gap}
 \frac{T-H}{d^2}=A-\frac{xp^2}{4m},\qquad A=\frac{T-x(1+e-m)}{d^2}.
\end{equation}
Indeed $(1+e)^2/(4m)=1+e-m+(2m-1-e)^2/(4m)$. At $x=1$ use the continuous limits whenever the numerators have the indicated zeros.
For the shear estimates, let $a=y^{-1}+y+2$, $y,t>0$, and define
\[
 B_2=2(1+t)+U(ty-1)+\tfrac32(t/y-1),\qquad
 B_3=2(1+1/t)+U(1/(ty)-1)+\tfrac32(y/t-1).
\]
The common quadratic form on four gradient components is
\begin{equation}\label{V:eq:common-shear}
 g^\top Sg=a(g_1-g_4)^2+B_2g_2^2+B_3g_3^2-2ag_2g_3 .
\end{equation}
The first block has largest eigenvalue $2a$; bounding the mixed term in the other block by $a(g_2^2+g_3^2)$ bounds its largest eigenvalue by $a+\max(B_2,B_3)$. The linear term depends on the particular correction field.
Write $b_* = \pi/n$, $k=\tan b_*$, $t=k^2$ and $\sqrt{U_s}=.16,.105,.075$ for $n=7,8,9$, respectively. The geometric bounds give $u<U_s$ and $x<r=r_n$. Use the field $V=P+ikQ$ in Section~\ref{seven-moments} and write
$M=I(|V|^2)$, $E=I(|V_X|^2+|V_Y|^2/t)$. To compare with the first eigenfield
match the energy to $E$, and in all mass ratios use the mass
of the field in question. Write $d=x-1$, $\psi=X(X^2-Y^2)$,
$\Phi=.03\psi$. At zero skew take $P+d\Phi,Q$ for the two profiles.
Define
\begin{gather*}
 B=I(P_Y^2)/(tE),\quad C=tI(Q_X^2)/E,\quad v=tI(Q^2)/M,\quad
 b=I(P_Y\Phi_Y)/(tE),\quad e=I(P\Phi)/M,\\
 G=I(\Phi_X^2)/E,\quad H=I(\Phi_Y^2)/(tE),\quad
 U=[I(\Phi^2)-I(\Phi)^2]/M,\\
 D=I((P_Y-\Phi_Y)^2)/(tE)=B-2b+H,\quad
 \kappa=I(P_X\Phi_X)/E+b-e .
\end{gather*}
For the first eigenfield $\kappa=0$ by testing its weak equation against
the real correction extended by rotation on all sectors. Its traces on
the radial edges vanish. The tangent and its divided remainder are
\[
 T=1+\tfrac{1+t}{k}K[\arctan(kx)-b_*],\qquad
 j=t\int_0^1\frac{p(1+Y)}{1+tY^2}\,dp,\quad Y=1+pd.
\]
Here
\begin{align*}
 m&=1+(x^2-1)v+2de+d^2U,\qquad K=1+2(C-B-v+\kappa),\\
 A&=(x+2)(v-C)-Kj+2b-D/x-2\kappa-x(G-U),\\
 p&=(x+1)(2v-C)+2(e+b)+(B-H)/x+D/x^2-d(G-2U)-2\kappa .
\end{align*}
The formula for $s$ follows from \eqref{V:eq:general-gap} with normalized zero-skew energy
\[
 W=1+(x^2-1)C+(x^{-2}-1)B+2d(\kappa+e-b+b/x^2)+d^2(G+H/x^2).
\]
Insertion of $T=1+K(d-d^2j)$ gives the displayed $A,p$.

Put $l=.67$, $J=j_0+j_1d$, where $j_0=t/(1+t)$ and
$j_1=t(1-3t)/[3(1+t)^2]$. The one-variable tangent remainder gives
\[
 j\le J-q d^2,\qquad
 q=\frac{t^2[3-t+(1-3t)l]}{4(1+t)^2(1+tr^2)} .
\]
Indeed the tangent at $1$ minus $(1+Y)/(1+tY^2)$ equals
$t(Y-1)^2[3-t+(1-3t)Y]/[(1+t)^2(1+tY^2)]$.
Use the angular expression $s_J$ obtained by replacing $j$ with $J$.
Thus $s\ge s_J+Kq d^2$.
The derivative of $s$ with respect to $v$ is
\[
 x+2+2j-x(x+1)p/m+x(x^2-1)p^2/(4m^2).
\]

For the first eigenfield, the three squared gradient quantities above have the following bounds, respectively:
$(\sqrt{B_V}\pm\delta)^2$, $(\sqrt{C_V}\pm\delta)^2$ and
$(\sqrt{D_V}\pm\delta)^2$, with the difference replaced by
its positive part in each lower bound. The same projection gives
$|b-b_V|\le\delta\sqrt H$, and the normalized mass projection gives
$|e-e_V|\le\epsilon\sqrt{I(\Phi^2)/M_V}+2\delta^2|e_V|$
and $(1-\delta^2/.7)U_V\le U\le U_V$. In fact $U/U_V=\mu/R$
and the mass normalization factor for $e$ is $\sqrt{\mu/R}$.
These estimates and the moments of each fixed field will be used below.

\subsection{A strict decrease under shear}\label{seven-shear}
At fixed $x$ write $L(u)=(1-u)^2 e=e_0+n_1u+n_2u^2$, $m=m_0+m_1u$. Direct division in $H_1=x(1-u)(1+e)^2/(4m)$ gives
\[
\sqrt{H_1(x,u)/H_1(x,0)}=
\frac{1+u[n_1-2+(n_2+1)u]/(1+e_0)}{(1-u)^{3/2}\sqrt{1+(m_1/m_0)u}} .
\]
On the central interval use $p=P+\alpha\psi$, $q=r_0=Q$,
$\nu=P+\beta\psi$, with $\alpha=.03d$, $\beta=\alpha+b_0$ and
$b_0=(n-4)/3$ for $n=7,8,9$. Here $r_0$ is the odd shear trial,
distinct from the aperture endpoint $r$. Put $y=tx^2$, $g=(P_X,P_Y/k,kQ_X,Q_Y)/\sqrt E$,
$G_0=12/(5E)$, $H_0=4/(9E)$, $v_0=\alpha/y+\beta y$, and form the symmetric matrix
$S$ by
\[S=S(t,y,U_s)\quad\text{as in \eqref{V:eq:common-shear}},\qquad a_2=B_2,\quad a_3=B_3,\quad a=y^{-1}+y+2.\]
Set
\begin{align*}
 \sqrt E\,\ell=\big(& (v_0+U_s\beta+3\alpha/2)\psi_X, (\alpha+\beta+U_s\beta y+3\alpha/(2y))k\psi_Y,\\
 &-(v_0+\alpha+\beta)\psi_Y/k, -(v_0+\alpha+\beta)\psi_X\big).
\end{align*}
Write $(1-u)^2 e(u)=e_0+n_1u+n_2u^2$, $m=m_0+m_1u$ in the centroid
fan formula and $\mathcal F=n_1+U_s(n_2+1)+3e_0/2-1/2$. Expansion gives
\begin{align*}
 \mathcal F={}&-1+2U_s+I(g^\top Sg+2\ell\cdot g)\\
 &+G_0(\alpha^2/y+\beta^2y+U_s\beta^2+3\alpha^2/2)
 +H_0(\alpha^2+\beta^2+U_s\beta^2y+3\alpha^2/(2y)).
\end{align*}

The mass increases since $I(\nu)\ge I(P)>0$ and $I_{I(P)}(\nu)\ge I(P)^2$.
The square-root ratio is less than one for nonzero skew whenever
\[
 \mathcal F<(1+e_0)(1-U_s)^{3/2}
 \frac{m_1/m_0}{1+\sqrt{1+U_s m_1/m_0}} .
\]
Indeed the numerator of that ratio is at most
$1-3u/2+u\mathcal F/(1+e_0)$; its denominator exceeds
$1-3u/2+u(1-U_s)^{3/2}(m_1/m_0)/(1+\sqrt{1+U_sm_1/m_0})$
for $0<u\le U_s$, by rationalizing the last square root.
The fraction on the right increases with $m_1/m_0$, so we can use a lower bound for this mass ratio.

\paragraph{Three coefficients in the shear estimate.}
The following argument applies to the reference field and the first eigenfield normalized to the same energy. Write $t=\tan^2(\pi/n)$, $y=tx^2$, $u_s=U_s$, $w=(n-4)/3$,
$\alpha=.03(x-1)$ and $\beta=w+\alpha$. In the shear expansion above
define $d_g=I[(g_1-g_4)^2-2g_2g_3]$ and
\[
 L_X=I(P_X\psi_X)/E,\quad L_Y=I(P_Y\psi_Y)/E,\quad
 L_Q=I(Q_Y\psi_X+Q_X\psi_Y)/E .
\]
The identities for the reference field give
\[
 \|g_{V,1}-g_{V,4}\|<.39,\quad d_{g_V}<1/7,\quad
 7/16<L_{X,V}<31/64,\quad L_{Y,V}<1/40,\quad L_{Q,V}>19/16 .
\]
Let $p=.070,.069,.0652$ and $q=.057,.056,.0538$ be the
reference upper bounds on $\sqrt B,\sqrt C$, respectively.
Write $G=12/(5E_-)$, $H=4/(9E_-)$ with
$E_-=.7807,.814,.8421$ and use the energy norm bounds
$\delta=\delta_n$ of Section~\ref{seven-spectral}. For either field put
\[
\begin{aligned}
 D_0&=1/7+2\delta\sqrt{2(.39)^2+p^2+q^2}+2\delta^2,
 & b&=(p+\delta)^2,& c&=(q+\delta)^2,\\
 L_-&=7/16-\delta\sqrt G,& l&=31/64+\delta\sqrt G,\\
 m&=1/40+\delta\sqrt{tH},& j&=19/16-\delta\sqrt{G+H/t}.
\end{aligned}
\]
Then $d_g\le D_0$, $B\le b$, $C\le c$, $L_-\le L_X\le l$,
$L_Y\le m$, $L_Q\ge j$. Indeed the operator of $d_g$ maps
$g$ to $(g_1-g_4,-g_3,-g_2,g_4-g_1)$ and has norm two.
The three linear errors follow from Cauchy--Schwarz and
$I(\psi_X^2)=12/5$, $I(\psi_Y^2)=4/9$.

In the linear part of the shear expansion the factors of
$2L_X,2L_Y,-2L_Q$ are respectively
\[
 A_1=\alpha/y+\beta y+u_s\beta+3\alpha/2,\quad
 A_2=\alpha+\beta+u_s\beta y+3\alpha/(2y),\quad
 A_4=\alpha/y+\beta y+\alpha+\beta .
\]
Here $y>.059$, $\beta\ge .99$, $\alpha\ge-.0099$, whence
$A_2,A_4>0$ and $A_1\ge-.01/y$. The coefficients of $B,C$
in the quadratic part are also positive: each is the sum of
a constant $\ge 1/2-u_s$ and positive $t,y$ terms. Thus
$A_1L_X\le A_1l+.01(l-L_-)/y$ and all other transfers have a fixed
sign. Collecting powers of $y$ gives $\mathcal F\le A/y+B_0+C_0y$, where
\[
\begin{aligned}
 A={}&D_0+\tfrac32 tb+u_sc/t+(G+\tfrac32 H)\alpha^2
       +2\alpha(l+\tfrac32m-j)+.02(l-L_-),\\
 B_0={}&-1+2u_s+2D_0+(\tfrac12+2t-u_s)b+(\tfrac12+2/t-u_s)c\\
 &+G(u_s\beta^2+\tfrac32\alpha^2)+H(\alpha^2+\beta^2)
       +2[(u_s\beta+\tfrac32\alpha)l+(\alpha+\beta)(m-j)],\\
 C_0={}&D_0+u_stb+\tfrac32c/t+(G+u_sH)\beta^2
                   +2\beta(l+u_sm-j).
\end{aligned}
\]
Each coefficient is a convex quadratic in $\alpha$, so one uses
$\alpha=-.0099,.03(r-1)$ in each. In order $n=7,8,9$ the
following simple majorants result:
\[
 \mathcal F\le .190/y-2.04+2.22y,\qquad
 .180/y-2.50+3.90y,\qquad .174/y-2.81+6.14y ,
\]
respectively. These bounds use only
$.23190<t_7<.23192$, $35/204<t_8<169/985$,
$51/385<t_9<31/234$ in the displayed coefficients.

Each final majorant is convex in $y$, so its largest value on
$[.67^2t,r^2t]$ is attained at an endpoint. The three upper bounds
are $.02,.14,.482$, respectively.

The projected mass bounds in Section~\ref{seven-spectral} and the reference moments in Section~\ref{seven-moments} give the following bounds for the first eigenfield. In particular $v_*\le(\sqrt{v_V}+\epsilon_n)^2<.0813,.0595,.0459$, and for the nine-sided field $M_*<.261$, $I(P_*)/\sqrt{M_*}>.93$, and $e_*>.0172$ using $I(\Phi^2)=.00012$. The lower bounds for $v$ used below are the Fourier estimates in Section~\ref{seven-central}. We now bound the mass increase from below in the shear inequality. Write $c_P=I(P)>0$.
At $n=7,8$ one has $m_1\ge v/t>.309$, since
$I_{c_P}(P+\beta\psi)=\operatorname{Var}(P+\beta\psi)+c_P^2+2c_P\beta I(\psi)>0$.
Also $m_0<1.205$ from $x<1.78$, $v<.0813$, $e<.0178$, $U<.00023$.
In both cases $1+e_0\ge1+(\sqrt{1-B-C}-.0234\sqrt G)^2>1.89$,
using $b<.009,c<.0066$. The right side of the shear criterion is
therefore at least
$1.89(1-.16^2)^{3/2}(.256)/(1+\sqrt{1+.16^2(.256)})>.23$.
For $n=9$ we use the same moments to obtain a sharper bound. After normalization to the same energy, both fields satisfy
$.2598<M<.261$, $c_P/\sqrt M>.93$, $v\ge.04295$,
$.0172<e<.0175$ and $.00020<U<.00023$. In particular
\[
 \operatorname{Cov}(P,\psi)=Me/.03-c_P(4/15)>0,
\]
as follows already from $M>.249$, $e>.0171$, $M<.27$,
$c_P\le\sqrt M$. Since $\operatorname{Var}\psi=14/225$,
retaining two terms of the squared variance yields
\[
 m_1\ge v/t+tx^2\{c_P^2/M+2c_P\beta(4/15)/M+
                              \beta^2(14/225)/M\}
       > .321+.412x^2 .
\]
Here $\beta>1.656$ and the expression in braces is $>3.12$.
Also $m_0<.924+.0462x^2+.035x$ by $v<.0459$.
At $.67$ the first polynomial exceeds $.52$ times the latter, and their
difference increases from there. Thus $m_1/m_0>.52$.
Furthermore
$1+e_0\ge 1+(\sqrt{1-B-C}-.0177\sqrt G)^2>1.92$,
using $B<.0062,C<.0046$.
The right side in the criterion is therefore at least
$1.92(1-.075^2)^{3/2}(.52)/(1+\sqrt{1+.075^2(.52)})>.494$.
This proves strict decrease at every nonzero central shear for each of the three side counts.

\section{Ordering of the Fourier coefficients and the transverse mass fraction}\label{seven-fourier}
Fix $7\le n\le9$, $a=2/n$, $b=\pi/n$ and $s=\sin(2b)/(2b)$. Let $v$ be the transverse mass fraction of the first eigenfield in a regular side sector. Take the polygon to have area $\pi$.

On its conformal disk the area weight and its Fourier coefficients in $\phi=n\theta$ are
\[
 W=Z^{-1}|1-z^n|^{-2a},\quad
 W_j=Z^{-1}\sum_{p\ge0} A_p A_{p+|j|}r^{n(2p+|j|)},\quad
 Z=\sum_{p\ge0}\frac{A_p^2}{np+1},\quad A_p=(a)_p/p!.
\]
These formulas follow by the binomial series for the regular polygon's Schwarz--Christoffel derivative, followed by area normalization. Write $\langle h\rangle=2\int_0^1h(r)r\,dr$.
Let $g(r)=J_1(\sqrt\lambda r)/J_1(\sqrt\lambda)$ be the disk dipole, with $3.3899<\lambda<3.39$, $g'(1)=0$ and $g(1)=1$.

\paragraph{Radial and angular order.}
On the disk the first field is $Y=e^{i\theta}\sum u_j(r)e^{ij\phi}$. It may be chosen with
\[
 u_0\ge u_{-1}\ge u_1\ge u_{-2}\ge u_2\ge\cdots\ge0.
\]
Indeed the inverse radial operator in order $m=|1+nj|$ has positive kernel
$G_m(r,t)=[(rt)^m+(r_</r_>)^m]/(4m)$ against $2t\,dt$, where $r_<,r_>$ denote the smaller and larger of $r,t$. Multiplication by $W$ convolves with a symmetric decreasing sequence; convolution preserves this interlacing order, by decomposing its arguments into positive sums of indicators of integer intervals centered at $0$ or $-1/2$. The kernels then decrease with $m$. Iterate from $g e^{i\theta}$. This converges, after normalization, in energy to the first eigenfield by compactness of the inverse and simplicity of the first eigenvalue in rotation character one. For nonzero overlap, taking the modulus of every radial Fourier coefficient decreases the energy and increases the weighted mass of a maximizer, since all $W_j$ are positive. The Green equation therefore gives a positive order-one coefficient of the first field.

Each $u_j/g$ is nondecreasing. The kernel $G_m$ is proportional to $r_<^m(r_>^m+r_>^{-m})$, and $t^m/(t^m+t^{-m})$ increases. Thus its ordered two-by-two determinants are nonnegative. Expanding such a determinant in a double integral shows that $G_mF/G_m h$ increases when $F/h$ increases, for positive $h$. If each input coefficient divided by $g$ increases, so do the coefficients of $W Y$ divided by $g$, since the series defining the $W_j$ have nonnegative coefficients. At $m=1$ use $G_1g=g/\lambda$. For $m\ge n-1\ge6$ the source itself increases; compare with $G_m1$, proportional to $H_m=(mr^2-2r^m)/(m-2)$. In fact
\[
 rH_m'/H_m=2m(1-r^{m-2})/(m-2r^{m-2})\ge1-r^3 .
\]
The Bessel product~\cite[Section 10.21(iii)]{DLMF} gives $rg'/g=1-\sum_{j\ge1}b_jr^{2j}$ with all $b_j>0$, $\sum b_j=1$ and $\sum j b_j=(\lambda-1)/2$, using the ODE and $g'(1)=0$. Jensen's inequality yields $rg'/g\le1-r^{\lambda-1}\le1-r^3$. Thus $H_m/g$ is nondecreasing. For $m\ge6$ the left side is at least $2(1-r^2)\ge1-r^3$, since $r^{m-2}\le r^2$. This proves the assertion by iteration. Energy convergence gives uniform convergence of each radial coefficient on compact subintervals of $(0,1)$ by one-dimensional Sobolev embedding. Monotone convergence of the nonnegative sums then justifies the radial integrations.

Write $Y=e^{i\theta}(p+iq)$, $C(\theta)=\sin^2(\theta-b)$ on $[0,2b]$, periodically extended. The transverse squared component minus $C|Y|^2$ is $q^2\cos(2\theta-2b)+pq\sin(2\theta-2b)$. The first term is nonnegative. In the second, $-q$ has nonnegative sine coefficients $u_{-j}-u_j$, and the sine coefficients of $-\sin(a(\phi-\pi))$ are $2jas/(j^2-a^2)>0$. The Fourier coefficients $t_j$ of $Wp$ are symmetric decreasing. Multiplying the latter sine series by $Wp$ preserves nonnegative sine coefficients, by $t_{j-k}-t_{j+k}\ge0$ for $j,k\ge1$. Thus the mean of the second term against $W$ is nonnegative.

At fixed radius, put $M_j=W_j/W_0$. These are moments of a probability measure on $[0,1]$: $A_{p+j}/A_p$ is the $j$th moment of the beta law with parameters $p+a,1-a$, and the weight formula mixes these laws after multiplying their variable by $r^n$. Consequently $M_{|p-q|}\ge M_{p+q}\ge M_p M_q$. The second inequality is the elementary covariance inequality for increasing functions of the same variable. Hence under the angular density $W/W_0$ all pairs $\cos p\phi,\cos q\phi$ have nonnegative covariance and $\operatorname{Var}(\cos\phi)\ge(1-M_1^2)/2$.
Moreover $M_j$ increases with radius: $A_{p+j}/A_p$ increases in $p$, and the probabilities proportional to $A_p^2r^{2np}$ shift monotonically towards larger $p$. In particular $M_1\le a/(1-a)$, its value at one. For this value Euler's beta integral in the two binomial sums gives $\sum A_p^2=\Gamma(1-2a)/\Gamma(1-a)^2$ and $\sum A_pA_{p+1}=[a/(1-a)]\sum A_p^2$.

The constant Fourier coefficient of $C$ is $(1-s)/2$; the coefficient at $j\ne0$ is $sa^2/[2(j^2-a^2)]$. Those of $|Y|^2/g^2$ are nonnegative and nondecreasing in radius by the radial order just proved. Thus the angular expectations of $C$ and $|Y|^2/g^2$ are nondecreasing in radius. Apply the one-variable covariance inequality with radial density proportional to $2rg^2 W_0$, and retain also the angular covariance coming from the two first cosine coefficients. The latter coefficient of $|Y|^2$, in exponential-series normalization, is at least $u_0(u_1+u_{-1})$. Writing $M=\pi^{-1}\int W|Y|^2$, this gives
\[
 v\ge \frac{\int g^2 W C}{\int g^2 W}
 +\frac{sa^2}{1-a^2}\left(1-\frac{a^2}{(1-a)^2}\right)
       \frac{\langle W_0u_0(u_1+u_{-1})\rangle}{M}.
\]

\paragraph{The first adjacent orders.}
Let $c=u_0(1)$. The kernel comparison gives $u_j\le u_0/|1+nj|$ for $j\ne0$, and $u_0\ge cr$ because $G_1(r,t)\ge rG_1(1,t)$. Also $u_0\le cg$ by radial order. Write $H=\sum_pA_pA_{p+1}$. The weight formula gives
$\langle W_j\rangle\le2H/(nZ|j|)$ for $j\ne0$.
Expanding the mass in Fourier pairs, using $|1+nj|\ge(n-1)|j|$, gives
\[
 M/c^2\le K=1-1/(\lambda Z)+\frac{2\zeta(2)}{(n-1)^2}
 +\frac{2H}{nZ}\left(\frac{4\zeta(2)}{n-1}+\frac{12\zeta(3)}{(n-1)^2}\right).
\]
Indeed the zeroth diagonal pair uses $\langle W_0g^2\rangle\le1-1/(\lambda Z)$, since $\langle g^2\rangle=1-1/\lambda$. The off-diagonal pairs with nonzero indices use $\sum_{ij(i-j)\ne0}1/|ij(i-j)|=12\zeta(3)$, obtained by splitting signs and summing $\sum_{j,k>0}1/[jk(j+k)]=2\zeta(3)$.
In orders $d=n\pm1$ keep only the term $W_1u_0\ge ac r^{n+1}/Z$ of the positive right-hand side. The solution of the radial equation in order $d$ with right-hand side $r^{n+1}$ and vanishing derivative at $r=1$ is $[(k/d)r^d-r^k]/(k^2-d^2)$, $k=n+3$. Pairing with $r$ and using $W_0\ge1/Z$ gives
\[
 \frac{\langle W_0u_0(u_1+u_{-1})\rangle}{M}
 \ge\frac{\mu a}{Z^2 K}(R_{n-1}+R_{n+1}),\qquad
 R_d=\frac{2(k+d+3)}{d(d+3)(k+3)(k+d)} .
\]
Here $\mu$ is the regular eigenvalue at area $\pi$; Lemma~\ref{IV:lem:uniform-regular-remainder} implies $\mu>3.35$ throughout this range.
The sine term in the transverse component also contributes a positive
amount. At fixed radius put $U_1=u_{-1}+u_1$. The interlacing order
holds also for the source $WY$. The displayed Green kernels satisfy
$G_{n-1}\ge[(n+1)/(n-1)]G_{n+1}$ pointwise. Hence
$u_{-1}-u_1\ge U_1/n$. The cosine series of $p$ has exponential
coefficients $u_0,U_1/2,\ldots$, symmetric and decreasing.
Subtracting $(u_0-U_1/2)$ from its constant coefficient preserves
these properties. Convolution with $W$ therefore gives
$t_0-t_2\ge (u_0-U_1/2)(W_0-W_2)$.
In the angular mean of $Wpq\sin(2\theta-2b)$ keep the $j=k=1$
term of the nonnegative sine pairing above. Its value is at least
$as(u_{-1}-u_1)(t_0-t_2)/(1-a^2)$.
Applying the Green bound $u_j\le u_0/|1+nj|$ at $j=-1,1$ and
adding gives $U_1/2\le nu_0/(n^2-1)$. Together with the boundary value from the beta integral
$M_2(1)=a(a+1)/[(1-a)(2-a)]$, this gives the multiplier $J$ used in
Section~\ref{seven-central}.
For the radial integral set $y=1-t$. The derivative of $1-g(\sqrt{1-y})^2$ is the product of $g(\sqrt{1-y})/\sqrt{1-y}$ and $g'(\sqrt{1-y})$. Their Taylor coefficients are nonnegative by the Bessel products with first zero at or beyond the boundary radius. This difference has quadratic coefficient $(\lambda-1)/4<.6$, no constant or linear term, and all coefficients sum to one. Thus
\[
 1-g(\sqrt{1-y})^2\le .6y^2+.4y^3 .
\]
Integrating this inequality against $t^z$ gives the rational lower bound $h(z)$ in Section~\ref{seven-central}.

\clearpage\part{The projected field and the comparison from ten sides}\label{part:large}
\section{Regular eigenspaces and signed turning angles}
Normalize by similarity so that $R_n$ has circumradius one and $|P|=|R_n|$. We call $P$ competitive when $F(P)\ge F(R_n)$.

\subsection{Symmetries of the first eigenspace of a regular polygon}
\label{IV:sec:single-small-turn-base}
\begin{lemma}\label{IV:lem:regular-character}
For a centered regular polygon $K=R_N$, $N\ge3$, the gradient map
$u\mapsto\nabla u(0)$ is an equivariant isomorphism from the first positive Neumann eigenspace onto $\mathbb R^2$.
\end{lemma}
\begin{proof}
Let $\phi>0$ be the Dirichlet ground state and $w=e^{ik\cdot x}$, $k$ parallel
to a side and $|k|^2=\lambda_D$. By integration by parts $w-c\phi$ has zero
defect at $\lambda_D$. Choose $c=\int_Kw/\int_K\phi$ so that the function has mean zero. Equality
would impose the Neumann condition. Hopf's boundary point lemma on the parallel side would then give $c=0$, but the plane wave has nonzero normal derivative on a nonparallel side. Thus $\mu_1<\lambda_D$.

A real first mode has two nodal domains by Courant's theorem. Each meets an open edge, since otherwise the corresponding positive or negative part has Dirichlet quotient $\mu_1$. If the mode were invariant under rotation by $2\pi/N$, a positive and a negative boundary point would have interlacing orbits on the boundary. Two positive points separating two negative ones in cyclic order cannot be joined through the positive domain while the negative points are joined through the negative domain: the two arcs must intersect. If $\nabla u(0)=0$, each nonzero
reflection-odd part $u-u\circ S$ is a mixed ground state on the half-polygon
by odd extension. Hopf's lemma at the center contradicts its zero gradient. Thus $u$
is fixed by every reflection, hence by rotation, so $u=0$. The injective
equivariant gradient has nonzero rotation-invariant image and is therefore onto.
Consequently the complex-valued first eigenfunctions $u$ with $u(\omega z)=\omega u(z)$, $\omega=e^{2\pi i/N}$, form a complex line, and so do those with $u(\omega z)=\bar\omega u(z)$.
On the polygon of circumradius one, $R_N$, the projection of $U=g(r)e^{i\theta}$ onto the first eigenspace is nonzero
because its imaginary part pairs positively with the positive odd first mode
on the upper half-polygon.
\end{proof}

\begin{lemma}[Disk estimates in a fixed angular order]\label{IV:lem:disk-angular-floors}
Let $\nu_m$ be the lowest Neumann eigenvalue on the unit disk in
angular order $m\ge1$. If $b$ is the positive root of
$2b^2=(m+4)b+m$, then
\[
             \nu_m\ge m^2+2m+m^2/b .
\]
In particular $\nu_2>9.2$, $\nu_4>24+32/9$ and $\nu_m\ge40$ for $m\ge5$.
The second Neumann eigenvalue in disk order one exceeds $28$.
On the circumradius-one regular $n$-gon set
$K=(\sqrt{4+\tan^2(\pi/n)}+\tan(\pi/n))^2/4$.
Among functions of rotation character one, that is, with $u(\omega z)=\omega u(z)$ for $\omega=e^{2\pi i/n}$, its second eigenvalue is at least $28/K$
for $n\ge6$, and exceeds $(24+32/9)/K$ for $n=5$.
\end{lemma}
\begin{proof}
Set $y=x^2>0$ and $a_j=xJ_j(x)/J_{j-1}(x)$. Fix an integer $k>\sqrt y$.
With $A=k-\sqrt{k^2-y}$ the recurrence
$a_j=y/(2j-a_{j+1})$ gives $0<a_j<A$ for $j\ge k$:
this holds at large $j$ by the Bessel series, where $J_j>0$ and
$a_j\to0$, and $y/(2j-A)\le A$ gives the induction step from one index to the preceding one.
Every denominator in the recurrence is positive.

Take $k=m+2$. We have
$a_{m+1}<y/[m+\sqrt{(m+2)^2-y}]\le m$ up to
$y=m^2+\delta$, where $\delta$ is the positive root of
$\delta^2+m^2\delta=4m^2(m+1)$. This is precisely
$\delta=2m+m^2/b$ by the defining quadratic for $b$.
Since $xJ'_m=J_m(m-a_{m+1})>0$ throughout, the first stationary
point lies beyond this interval. At $m=2,4,5$ one has respectively
$b<10/3$, $b=2+\sqrt6<9/2$ and $b=5$. The radial Rayleigh
quotient increases with $m$ because its angular potential is $m^2/r^2$.

For the second mode in order one take $y\le28$ and $k=7$.
Then $a_7<5/2$, hence $a_6<3$, $a_5<4$, $a_4<7$ and $J_3>0$.
Two more recurrences give
$x^3J'_1(x)=J_3(x)[24-7y+(y-4)a_4]$.
For $4\le y\le28$ the bracket is at most $-4$. On $0<x\le2$
the alternating series gives $J_1(x)/(x/2)\ge1-y/8>0$.
Since $(xJ'_1)'=(1/x-x)J_1$, the function $xJ'_1$ increases up
to one and then decreases through zero before two. Thus $J'_1$
has exactly one zero through $\sqrt{28}$.

The equivariant radial disk-to-polygon map has factor
$\rho=\cos(\pi/n)/\cos\theta$ on a bisector sector. Its polar derivative
is $\rho$ times a unit shear of entry $\tan\theta$. Thus its Jacobian is
$\rho^2\le1$ and its singular-value ratio is at most $K$. By min--max the
eigenvalues of the polygon in each rotation character are at least the disk values divided
by $K$. In rotation character one the angular orders on the disk are $1+jn$, and those other than one have
absolute value at least $n-1$.
\end{proof}

The second eigenvalue in rotation character one also exceeds $55/3$ for $n\ge5$. Put $b=\pi/n$, $s=\operatorname{sinc}b$.
The map $(r,\theta)\mapsto r(\cos b,s\theta)$ has constant Jacobian; its largest
squared stretch occurs at $|\theta|=b$, where the Gram matrix has diagonal $1,s^2$
and off-diagonal $bs^2$. For $n\ge6$ this stretch is at most $1+b<32/21$.
For the pentagon, $s<15/16$, $b<22/35$ give $(3/2-s^2)/2>b^2s^4$, so the squared stretch is below $3/2$.
Divide the disk eigenvalue lower bounds $28,248/9$, respectively, by these bounds.

Normalize $g(1)=1$ and write $\lambda=(j'_{1,1})^2$.
In $2J'_1(\sqrt x)=\sum_{j\ge0}(-x/4)^j(2j+1)/(j!(j+1)!)$ the alternating
sums through $j=9,8$ have signs $+,-$ at $x=3.389957716,3.389957717$, respectively.
The zero count above locates $\lambda$ between these endpoints.
Since $g$ increases to one and $(r^3(g/r)')'=-\lambda r^3(g/r)$, $g/r$
decreases to one. In particular $r\le g\le1$ and
$0\le g'\le g'(0)\le(1-\lambda/8)^{-1}<2$, using the first two terms of
$2J_1(\sqrt\lambda)/\sqrt\lambda$. Differentiating
$r^2g'^2+(\lambda r^2-1)g^2$ gives $2\lambda rg^2$, hence
$2\int_0^1rg^2dr=1-1/\lambda$.
\subsection{Centering the disk and counting convex vertices}
For a Riemann map $f_0:\mathbb D\to P$ and continuous
$U(re^{i\theta})=g(r)e^{i\theta}$, $g(0)=0,g(1)=1$, choose coordinates with
$\int U|f'|^2=0$, $f'(0)>0$. Indeed the normalized mean of
$U((z-a)/(1-\bar a z))$ against $|f_0'|^2dA$ extends continuously
to $-a$ on the unit circle by dominated convergence. Degree one gives an
interior zero; precompose $f_0$ by the inverse automorphism and rotate.
We call the resulting map $f$ \emph{balanced} and $O=f(0)$ a \emph{conformal center} of $P$.

\begin{lemma}[Convex vertices and the circle mean]\label{IV:lem:signed-circle}
For a conformal map from the disk onto a simple polygon, put $Q=f'/f'(0)=1+\sum_{k\ge1}\gamma_k z^k$
and $H(s)=\sum_{k\ge1}|\gamma_k|^2e^{-ks}$. If at most $m$ corners ($m\ge3$) of
$P$ are strictly convex, then for every $s>0$, $b\ge0$,
\[
            H''/4+bH'+b^2(1+H)\ge(1+b)^2e^{-ms}.
\]
Put $U=(e^{ms}-1)^{-1}$, $\Phi=(1-e^{-ms})^{4/m}$ and
$\ell=2/[m(m-2)]$. For $0\le W(s)\le s^2/2$,
\begin{equation}\label{IV:eq:integrated-turn-functional}
 4\int_0^\infty WU/\Phi\le\int_0^\infty WH''+
             \ell\int_0^\infty e^{-ms}H .
\end{equation}
\end{lemma}
\begin{proof}
Let $2\pi\nu_j$ be the nonzero signed exterior angles and $\xi_j$ their
prevertices. The Schwarz--Christoffel formula \cite{DT} gives
$S=zQ'/(2Q)=\sum_j\nu_j/(1-z/\xi_j)-1$, with $\sum_j\nu_j=1$.
Let $p\le m$ be the number of positive weights and use normalized
circle measure. The signed-measure estimate from the author's earlier
work \cite[proof of Lemma~3.3, (3.8), and proof of Theorem~1.3]{Wang}
gives, for $0<r<1$,
\[
 \int_{\mathbb T}\log|S(r\zeta)-b|\,d\sigma(\zeta)
 \ge \log(1+b)+p\log r.
\]
This estimate uses only the signed weights and their prevertices, so it
applies to the balanced conformal map. The circle
mean of $\log|Q|$ vanishes since $Q$ is zero-free and $Q(0)=1$.
The arithmetic--geometric mean inequality therefore gives
$\|Q(r\cdot)(S(r\cdot)-b)\|_2^2\ge(1+b)^2r^{2m}$.
At the conformal center the first Fourier coefficient is retained, as
in Dong and Wang \cite[Lemma~3]{DW}. With $r=e^{-s/2}$, Parseval
applied to $Q(S-b)=zQ'/2-bQ$ gives the left side explicitly as
\[
 b^2+\sum_{k\ge1}(k/2-b)^2|\gamma_k|^2 e^{-ks}
 =H''/4+bH'+b^2(1+H).
\]

At $b=U$ the inequality becomes $H''\ge4U(1-H')-4U^2H$, hence
$(\Phi(1-H'))'\le4\Phi U^2H$. Integrate from $s$ to infinity, where
$\Phi(1-H')\to1$, and substitute in the first inequality.
After integration against $W$ the total error is
\[
 \int_0^\infty U(v)^2H(v)
 \left[4W(v)+16\Phi(v)\int_0^v W(s)U(s)/\Phi(s)\,ds\right]dv .
\]
Here $U(s)\le1/(ms)$ and $\Phi(v)/\Phi(s)\le(v/s)^{4/m}$ for $s\le v$.
The bracket is at most $[2+4/(m-2)]v^2$; now
$v^2U(v)^2\le e^{-mv}/m^2$ gives \eqref{IV:eq:integrated-turn-functional}.

Set $W_0(s)=(s-1+e^{-s}){\bf1}_{s\ge0}$. Its Laplace transform is
$1/[k^2(k+1)]$, so $\int W_0 H''=\int e^{-s}H$ term by term.
\end{proof}

\section{Polygons with fewer than the prescribed number of strictly convex vertices}
\label{IV:sec:conformal}

Write $F(P)=|P|\mu_1(P)$, $\lambda=(j'_{1,1})^2$ and $B=4\zeta(3)$.
Lemma~\ref{IV:lem:signed-circle} supplies the circle-mean bound at the conformal center. Dong and Yang \cite[Section~9, Corollary~9.2]{DY} obtained a spectral gap for a missing convex vertex. For six sides we integrate the circle mean against a rational weight. For seven or more sides we combine its integrated form with the following eigenvalue estimate for regular polygons in Theorem~\ref{IV:thm:corner-gap}.

\begin{lemma}[A lower bound for $J_m$]\label{IV:lem:positive-turn-kernel}
Put $h=1/m$, $c_j=\prod_{i<j}(1+4h/i)$, $B=4\zeta(3)$,
$r_k=(M_{k+1}-M_1)/(1-M_1)$ for $M_j=2j\int_0^1g^2r^{2j-1}dr$, and
$m^3J_m=4\sum_{j\ge1}c_jr_{mj}/[j^3(1+h/j)]$.
For $m\ge6$, $m^3J_m>B-10h^2$; for $m\ge9$ one may add $50h^3$. For $m=6,7,8$ the respective strict lower bounds are
$4.832,4.804,4.792$. Thus $m^3J_m>4.60$ throughout and
$J_m>B/n^3+3B/n^4$ for $n=m+1\ge8$.
\end{lemma}
\begin{proof}
Write $A_j=\sum_{i<j}1/i$, $e_d(j)=[z^d]\prod_{i<j}(1+z/i)$.
Two identities sum the linear and higher terms:
\[
 4\sum_j A_j/j^3=\zeta(4),\qquad
 \sum_j e_d(j)/j^3=\sum_j A_j/j^{d+2}\quad(d\ge1).
\]
For the first, twice the sum on the left without the $4$ equals
$\sum_{a,b\ge1}1/[ab(a+b)^2]=\zeta(4)/2$ by
$2/[ab(a+b)^2]=1/(ab)^2-(a^{-2}+b^{-2})/(a+b)^2$.
For the second set $L=-\log t$, $A=-\log(1-t)$. Successive integration gives
$\sum_j e_d(j)t^j/j=A^{d+1}/(d+1)!$. Integrate against $L\,dt/t$, integrate
by parts, then reflect $t\mapsto1-t$ to obtain $\int_0^1A^2 L^d\,dt/(2d!t)$.

Let $S$ denote $m^3J_m$ with $r_{mj}=1$. Removing the constant and
linear terms from $(1+4hA_j)/(1+h/j)$ leaves
$h^2(j^{-2}-4A_j/j)/(1+h/j)$, negative except at $j=1$.
All remaining terms have $d\ge2$, $j\ge3$; use $1+h/j\le1+h/3$
and the second identity to sum their geometric series. Thus
\[
 \frac{S-B}{4h^2}>
 A(h):=\frac1{1+h}+.036+
 \left(\frac{16}{1+h/3}-4\right).095+
 \frac{64h}{1+h/3}\sum_{j=2}^5\frac{A_j}{j^5(1-4h/j)}.
\]
Here $\zeta(5)-1>.036$ and $\sum A_j/j^4>.095$ already through $j=9$.
By Lemma~\ref{IV:lem:uniform-harmonic-correlations},
$1-r_k\le E_k:=3.39(4.78+.78/(2k+5))/[(2k+3)(2k+4)]<3.39(4.78)/(4k^2)$.
Separating the first two terms of this error and using $c_j/j\le c_2/2$ for $j\ge3$ yields
\[
 m^3J_m>B+4h^2 A(h)-\frac{4E_m}{1+h}
 -\frac{c_2E_{2m}}{2(1+h/2)}-.1621c_2h^2.
\]
Indeed $\sum_{j\ge3}j^{-4}<.020$ (separate $j=3,4$ and integrate
from $9/2$).

For $m=6,7,8$ a shorter calculation uses the first five summands of
$J_m$ itself. Put $\alpha=4/m$ and
$u_j=4c_j(1-E_{mj})/(1+h/j)$. Each of the five summands is
at least $u_j/j^3$. Concavity gives
$1+\alpha/i\ge(1+1/i)^\alpha$, so that
$c_j\ge c_6(j/6)^\alpha$ for $j\ge6$.
The factor $(1-E_{mj})/(1+h/j)$ increases with $j$.
For the convex function $f(t)=t^{-3+\alpha}$ the trapezoid
inequality on successive unit intervals gives
$\sum_{j\ge6}f(j)\ge f(6)/2+\int_6^\infty f(t)\,dt$.
Thus the infinite tail is bounded below by the rational expression
$u_6[1/(2\cdot6^3)+1/(6^2(2-\alpha))]$.
The table gives lower bounds, in units of $10^{-4}$, using
$c_1=1$, $c_{j+1}=c_j(1+\alpha/j)$.
\[
\begin{array}{c|rrrrr|r|r}
m&j=1&2&3&4&5&j\ge6&m^3J_m>\\\hline
6&31948&7526&3086&1619&977&3226&4.8382\\
7&33130&7212&2834&1444&852&2598&4.8070\\
8&34027&6967&2650&1321&766&2206&4.7937
\end{array}
\]
Only the six integer columns are scaled.
For $h\le1/5$, $4A(h)>7.65$ and the three subtractions divided by
$h^2$ are below $16.4,.925,.294$. For $h\le1/9$ the first is at most
$3.39(4.78+.39h)/[(1+h)(1+1.5h)(1+2h)]$. Its factors are
positive decreasing convex, pairing the numerator with $1+h$,
so its endpoint chord lies below $16.4-50h$.
From $m=7$ this exceeds $4.60$, while $(1-1/n)^{-3}\ge1+3/n+6/n^2$
and $10(8/7)^5<6B$ give the last assertion.
\end{proof}

\begin{lemma}[An upper bound for $d_R$ for seven or more sides]\label{IV:lem:uniform-regular-remainder}\label{IV:prop:temple-regular}\label{IV:cor:block-regular}

For a regular $n$-gon of area $\pi$, put $B=4\zeta(3)$,
$e=\lambda-1$ and $d_R=e(\lambda/\mu_1(R_n)-1)$. Here $4.808<B<4.809$. For the lower bound sum through twenty and add $4\int_{21}^\infty t^{-3}\,dt$; for the upper bound sum through five and use midpoint convexity to bound the remaining tail by $4\int_{11/2}^\infty t^{-3}\,dt$.
For every $n\ge7$,
\begin{equation}\label{IV:eq:regular-continuous-band}
 d_R\le B/n^3+U(n)/n^5,\qquad
 U(n)=\begin{cases}
 92.5,&n=7,\\ 85.25,&n=8,\\ 80.75,&n=9,\\
 59.93+197/n,&10\le n\le14,\\ 74,&n\ge15 .
 \end{cases}
\end{equation}
In particular $U(n)\le79.63<84$ for $n\ge10$.

\end{lemma}
\begin{proof}
Set $x=1/n$, $a=2x$. On the conformal disk the area weight is
$Z^{-1}|1-z^n|^{-2a}$, where $a_j=(a)_j/j!$ and
$Z=1+\sum_{j\ge1}a_j^2/(nj+1)$. With the unscaled map normalized by $f'(0)=1$, the beta integral gives its circumradius $f(1)=\Gamma(1+x)\Gamma(1-2x)/\Gamma(1-x)$. The polygon area formula then gives $Z=\operatorname{sinc}(2\pi x)f(1)^2$. Euler's reflection formula yields
$Z=\Gamma(1+x)^2\Gamma(1-2x)/[\Gamma(1-x)^2\Gamma(1+2x)]$.
The expansion of $\log\Gamma(1+z)$ therefore gives
\[
 L=\log Z=2\sum_{k\ge3,\ k\ {\rm odd}}(2^k-2)\zeta(k)x^k/k
 \le B x^3+b(x)x^5,\qquad b(x)=12(1.037)+\frac{2(1.01)2^7x^2}{7(1-4x^2)}.
\]
Indeed the fifth coefficient is $12\zeta(5)$, and the remaining odd
terms are bounded by the geometric tail starting at $k=7$.
By Lemma~\ref{IV:lem:regular-character}, $1/\mu_1(R_n)$ is the largest
ratio of mass to energy among functions of rotation character one. Use energy and mass
divided by $\pi$. Decompose such a function into a multiple of the normalized disk eigenfunction
$u=g(r)e^{i\theta}/\sqrt e$, a component in its radial complement in the same angular order, and all
higher orders. Write $q$ for the diagonal on $u$ and $d_1,d_2$ for upper
bounds on the two complementary diagonals. If their mutual coupling
has norm $c$ and their squared couplings with $u$ are bounded by
$Px^6,Vx^5$, bound the cross term by Young's inequality with coefficients $4c$ on the squared norm of the radial complement and $c/4$ on that of the higher orders. Completing the squares gives
\[
 1/\mu_1(R_n)\le q+Px^6/(q-d_1-4c)+Vx^5/(q-d_2-c/4).
\]
The lower bounds below show that both denominators are positive. Radial covariance gives $q\ge1/\lambda$.
More precisely, with $D=\sum_{j\ge1}a_j^2\int_0^1(1-g^2)2r^{2nj+1}dr$,
$\lambda e(q-1/\lambda)=1-(1+\lambda D)/Z$.
Set $z=1-r^2$, $K(X)=g(\sqrt X)/\sqrt X$ and $S=K+2XK'$.
By Lemma~\ref{IV:lem:wall-zero-products} both $K(1-z)$ and $S(1-z)/z$
have nonnegative coefficients. Since $[1-g(\sqrt{1-z})^2]'=K(1-z)S(1-z)$,
the loss has nonnegative coefficients summing to one, starting at degree two.
The disk equation gives $c_2=(\lambda-1)/4$, $c_3=(2\lambda-3)/12$.
Bound its remaining powers between zero and $(1-c_2-c_3)z^4$.
Beta integration against $j(1-z)^{j-1}$, with $M_j=2j\int_0^1g^2r^{2j-1}dr$
and $(x)_k=x\cdots(x+k-1)$, gives
\begin{equation}\label{IV:eq:regular-beta-band}
 0\le 1-M_j-\frac{2c_2}{(j+1)_2}-\frac{6c_3}{(j+1)_3}
 \le \frac{24(1-c_2-c_3)}{(j+1)_4}.
\end{equation}
With $y=1-r^2$ the lower polynomial itself gives
\begin{equation}\label{IV:eq:positive-regular-mass-loss}
 1-g^2\ge ey^2/4+(2\lambda-3)y^3/12 .
\end{equation}
The first Fourier coefficient $a_1=2x$ and the beta integrals for
$y^2,y^3$ now give
\begin{equation}\label{IV:eq:regular-first-power-loss}
 n^5\lambda D\ge K(x):=\frac{2\lambda e}{(1+x)(1+2x)(1+3x)}
 \left(1+\frac{(2\lambda-3)x}{e(1+4x)}\right).
\end{equation}
For $j\ge2$, $a_j\ge(2x/j)(1+2x)$ and
$\prod_{i=1}^3(1+ix/j)\le(1+x)^3\le(1+2x)^2$. Thus the quadratic term
in \eqref{IV:eq:positive-regular-mass-loss} contributes at least $2\lambda e/j^5$
to $n^5\lambda D$; summing $2\le j\le5$ gives $n^5\lambda D\ge K(x)+.59$.
For a supporting line write $\lambda_*=3.3899$, $e_*=\lambda_*-1$
and let $K_*$ be $K$ at $\lambda_*$. Both factors increase with $\lambda$,
so $K\ge K_*$. Apart from a positive factor, $K_*$ equals
$(1+\kappa x)F(x)$ with $\kappa=6-1/e_*$ and
$F=\prod_{i=1}^4(1+ix)^{-1}$. Put $h_i=i/(1+ix)$ and $s=\sum h_i$.
Since $F''/F=s^2+\sum h_i^2\ge5s^2/4$ and
$(1+\kappa x)s\ge10>8\kappa/5$, the product is convex.
At $t=1/20$, $K_*(t)-tK_*'(t)+.59>16.42$ and $K_*'(t)>-56.55$.
Put $h=4.809+b(x)x^2$, $L_+=hx^3$, $L_-=4.808x^3$, $Z_-=1+L_-$.
The mass identity and $1-e^{-L}\le L-(1-L/3)L^2/2$ give
\[
 A:=x^{-5}[\lambda e(q-1/\lambda)-B x^3]
 \le b-(1-L_+)(K_*+.59)-(1-L_+/3)4.808^2x/2.
\]
Here $L_-\le L\le L_+$. On $[0,1/10]$ the supporting line gives
$A\le A_*:=b-16.42+45x+16.82L_+$.

In angular order $m\ne0$, radial integration gives
$\|v_m\|_\infty^2\le E_m/(2|m|)$ for normalized energy $E_m$.
In order one the Neumann evaluation kernel is
$(rt+\min(r,t)/\max(r,t))/4$ (its radial flux jumps by $-1/2$).
Projecting off $u$ gives evaluation norm squared $(1+r^2)/4-g^2/e\le1/4$
as $g\ge r$ and $e<4$. The excess radial weight
$\sum_{j\ge1}a_j^2r^{2nj}2r\,dr$ has mass $Z-1$. In each power
$\mathbb E y^2/\mathbb E y=2/(nj+3)$, and
$g^2\ge(1-y)(1+y/2)^2\ge1-y^2$.
Hence its mean squared evaluation norm is at most $1/2-1/e$ for $n\ge7$.
Cauchy--Schwarz, $g\le1$, and $1-1/Z\le L_+$ give
$P=h^2(1/2-1/2.39)/e_*$.

Put $d_{jp}=a_pa_{p+j}/a$, $T_j=\sum_{p\ge0}d_{jp}/(p+j+1)$ and $U_j=(j+1)T_j$.
Gauss summation gives $T_1=B_0/(2-a)$, $T_{j+1}/T_j=(j+a)/(j+2-a)$,
where $B_0=\Gamma(2-2a)/\Gamma(2-a)^2\le[1-.645a^2-.405a^3/(1-a)]^{-1}$.
The bound holds coefficientwise: the coefficients $(2^k-2)(\zeta(k)-1)/k$ of $\log B_0$ are below
$.645$ for $k=2$ and $.405$ for $k\ge3$, using $\zeta(k)-1\le2^{3-k}(\zeta(3)-1)$.

The squared dual norm of $2\int r^{l-1}v\,dr$ in angular order $m>0$ is
$(2/m+1/l)/(l+m)^2$, by the Neumann representative
$[(l/m)r^m-r^l]/(l^2-m^2)$ and continuity. Thus for the power $(2p+j)n$
and the two orders $jn\pm1$ the norm is $n^{-3/2}L_{jp}$, with
\[
 L_{jp}^2=\sum_{m=j\pm x}\frac{2/m+1/Y}{(Y+m)^2},\quad Y=2p+j+2x,
 \qquad H=L_{10}+d_{11}L_{11}+4L_{12}(T_1-1/2-d_{11}/3).
\]
Here $d_{11}=a(1+a)/2$. Each $(p+j+1)L_{jp}$ decreases with $p$ since
$(p+j+1)/(Y+m)$ decreases. Also $L_{jp}(x)\le L_{jp}(0)$ for $p\le2$:
their squared ratio is at most $[(1-x^2/j^2)(1+x/(2p+2j))^2]^{-1}\le1$.
In particular $L_{j0}(0)^2=3/(2j^3)$. For $j\ge2$ we have
$U_j\le U_2(2/j)^{1-2a}$ since
$\log[(1+a/j)/(1-a/(j+2))]\le a(1/j+1/(j+2-a))
\le2a/(j+1/2)\le2a\log(1+1/j)$.
Use the triangle inequality for fixed $j$ and orthogonality across $j$.
Bounding $\sum_{j\ge3}j^{-5+4a}$ by its midpoint integral from $5/2$,
and $(5/4)^{4a}\le(1-.9a)^{-1}$, gives
\[
 S=H^2+(U_2L_{20})^2+\frac{3U_2^2(4/5)^4}{8(4-4a)(1-.9a)},
 \qquad V=\frac{4S}{e_*Z_-^2},\quad c=x^{5/2}\sqrt S/Z_- .
\]
The squared coupling to the higher orders is at most $Vx^5$. Evaluation with norm $1/2$
gives the bound $c$ for their coupling to the radial complement, since the radial Green kernel is positive.

The integrated Fourier weight at displacement $j$ is at most $2aw_j/(nj)$, where
$w_j=jT_j+(a_j/a)j(1-a)/[(j+a)(j+1)]$: separate $p=0$ and use
$2p+j+a\ge p+j+1$ otherwise. Both terms decrease with $j$ by their ratios.
For $0\le a\le2/7$, $w_1\le1$ and $w_2\le3/4$ by testing the convex
majorants $(1+a)^{-1}+T_1-1/2$ and $2T_2+1/6-a/12$ at the endpoints.
With Schur weights $|ni+1|^{-1/2}$, the convolution
$\sum_{k\ne0,i}(|k|\,|i-k|)^{-1}=2(H_{|i|-1}^{\rm harm}+H_{|i|}^{\rm harm})/|i|\le5/2$.
For $|i|=1$ this sum is $2$, of which $1/2$ has $|i-k|=1$.
Thus the corresponding weighted sum is at most $3.25x^3/(1-x)$. For $|i|\ge2$ the
disk eigenvalue is at least $(2n-1)^2$. The additional terms satisfy the same bound because
$(2-x)^{-2}+1.75x/(1-x)<.59<(1+x/2)^{-1}$.
Lemma~\ref{IV:lem:disk-angular-floors} gives the remaining eigenvalue lower bound $n^2+n/2$.
Including the contribution of the excess radial weight gives the denominators
\[
 D_1=1/3.39-1/28-L_+/2-4c,\quad
 D_2=1/3.39-\frac{x^2}{1+x/2}-\frac{xL_+}{2(1-x)}-\frac{3.25x^3}{1-x}-c/4 .
\]
The Schur bound is $n^5(d_R-B x^3)\le A+3.39(2.39)(xP/D_1+V/D_2)$.
At $n=7,8,9$ retain the first radial moment in the evaluation bound. In the order-one
radial complement the squared evaluation norm is at most
$(1+r^2)/4-g(r)^2/e$, where $e=\lambda-1$.
Averaging over the excess radial measure and putting $y=1-r^2$
therefore gains, compared with $1/2-1/e$, at least
\[
 \left(\frac14-\frac{2}{e_*(n+3)}\right)
 \frac{a^2(1-L_+)}{(n+1)(n+2)L_+}=:\eta_n .
\]
Indeed $1-g^2\le y^2$, and in every term of the measure
$\mathbb E y^2/\mathbb E y=2/(nj+3)$. The contribution of its first
term to the unnormalized $y$-moment is $a^2/[(n+1)(n+2)]$.
The total mass $Z-1\le L_+/(1-L_+)$.
Thus one can replace $P$ by $P'=P-h^2\eta_n/e_*$
without altering any high-mode bound.

For each fixed $x=1/n$ put $H_2=U_2L_{20}$, $H_3=S-H^2-H_2^2$,
$D'_2=Z_-^2D_2$, and evaluate the single expression
\[
 W_{\rm fin}(x)=A+3.39(2.39)\left(\frac{xP'}{D_1}
       +\frac{4(H^2+H_2^2+H_3)}{2.3899D'_2}\right).
\]
Here $A$ may be replaced by the upper bound following
\eqref{IV:eq:regular-first-power-loss}. Each row displays the entries to insert;
all upper endpoints round up and both denominators round down.
\[
\begin{array}{c|rrr|rr|rr|r}
x&H<&H_2<&H_3<&A<&P'<&D_1>&D'_2>&W_{\rm fin}<\\\hline
1/7&1.231&.353&.057&2.35&.746&.212&.268&92.29\\
1/8&1.214&.329&.045&1.65&.730&.226&.275&85.16\\
1/9&1.205&.313&.037&1.09&.722&.235&.280&80.72
\end{array}
\]
To recover the first three columns use
$T_1\le[(2-a)(1-.645a^2-.405a^3/(1-a))]^{-1}$,
$U_2=3T_1(1+a)/(3-a)$ and the four square roots
$L_{10},L_{11},L_{12},L_{20}$ defined above. For instance the last row
uses $1.09+3.39(2.39)[.722/(9\cdot.235)+
4(1.205^2+.313^2+.037)/(2.3899\cdot.280)]<80.72$.
The rounded values give the stated $U(7),U(8),U(9)$.

For $0\le x\le1/10$ take $A_*$ from the supporting line for $K_*$ above. Use the following convex majorant for the remaining side counts.
In $H$ use $L_{11}(0),L_{12}(0)$.
For the second component of $\sqrt S$ use
\[
 H_2=\tfrac12L_{20}(x)+(U_2-\tfrac12)L_{20}(0).
\]
Denote the resulting upper bound on $S$ by $S_c$. Its square root is
convex on $[0,1/10]$. Indeed $L_{j0}$ is log-convex, since positive
sums and products of reciprocals of positive affine functions are
log-convex. The remaining summands of the first two components have
nonnegative Taylor coefficients; this holds for $T_1-1/2-d_{11}/3$
by the terms $p\ge2$ in its defining series and the coefficientwise
majorant for $B_0$. The third component is a product of functions with
nonnegative coefficients. The Euclidean norm of these three nonnegative
convex functions is convex.

In $D_1,D_2$ take $c=1.4x^{5/2}$. Indeed $\sqrt{S_c}<1.255$ at either
endpoint and $Z_-\ge1$. The resulting $D_2$ decreases to a value above
$.280$. Thus $Z_-^2D_2\ge D_2+2(.265)L_-=:D'_2$. Both $D_1,D'_2$
are positive, decreasing and concave. For the latter note that the combined
cubic subtractions have coefficient $3.25/(1-x)-2(.265)(4.808)$
with nonnegative power series, and that $x^2/(1+x/2)$ is convex.
The Schur estimate therefore gives
\[
 n^5(d_R-4\zeta(3)x^3)
 \le W(x)= A_*+3.39(2.39)
 \left(\frac{xP}{D_1}+\frac{4S_c}{e_*D'_2}\right).
\]
This majorant is convex. In the first quotient the numerator and
reciprocal denominator are nonnegative, increasing and convex.
For the second use convexity of the quadratic-over-linear function,
monotone increasing in $\sqrt{S_c}\ge0$ and decreasing in $D'_2$.

At $0$, $1/13$ and $1/10$ the three entries
$(H,H_2,S_c-H^2-H_2^2)$ are bounded above, respectively, by
\[
 (1.22475,.21651,.00961),\quad
 (1.19622,.28000,.02366),\quad
 (1.20322,.30587,.03174).
\]
These follow from $T_1=[(2-a)(1-.645a^2-.405a^3/(1-a))]^{-1}$,
$U_2=3T_1(1+a)/(3-a)$ and the displayed two-term formula for
$L_{jp}^2$, using $a=2x$. Thus $\sqrt{S_c}$ is less than
$1.248,1.239,1.255$, respectively. At zero use
$A_*<-3.97$ and $D'_2>.2949$. At $1/13$ use
$A_*<-.25$, $P<.816$, $D_1>.248$, $D'_2>.288$;
at $1/10$ use $A_*<1$, $P<.833$, $D_1>.239$, $D'_2>.283$.
Substitution gives
\[
 W(0)<68,\qquad W(1/13)<74.2,\qquad W(1/10)<79.4 .
\]
On $[1/13,1/10]$ the chord between these two upper bounds lies below
$59.93+197x$ by evaluation at the two ends. On $[0,1/13]$, convexity
gives $W(x)\le68+80.6x$. At $n=14$ this is below $59.93+197/n$;
at $n\ge15$ it is below $74$. Consequently one may take
\[
 U_n=59.93+197/n\quad(10\le n\le14),\qquad
 U_n=74\quad(n\ge15).
\]

\end{proof}

\begin{lemma}[The upper bound for $d_R$ at five and six sides]
\label{small-schur-regular}
In area $\pi$ normalization,
$d_R\le 4\zeta(3)/n^3+U/n^5$ with $U=94.3$ at five sides and $U=98$
at six sides.
\end{lemma}
\begin{proof}
At five sides Lemma~\ref{lem:degree-six-pentagon-floor} gives
$F(R_5)>4.35415\,(5/2)\sin(2\pi/5)>10.3526$,
using $\sin(2\pi/5)=\sqrt{10+2\sqrt5}/4$. Hence
\[
 5^5(d_R-B/5^3)<
 3125\left[(3.389958-1)\left(\frac{3.1416(3.389958)}{10.3526}-1\right)
              -\frac{4.808}{125}\right]<94.3 .
\]
For $n=6$ use the preceding lemma through its formulas for the dual norm $L_{jp}$
and the covariance $P'$ at $x=1/6$. There $\eta_6>0$ since
$1/4-2/(9e_*)>0$.
Treat the two angular orders $1-n$ and $1+n$
separately: the order $1+n$ has the larger disk eigenvalue.

Use the following upper bound for $T_1$
\[
 \begin{gathered}
 G=.645a^2+.405a^3+.33a^4/(1-a),\\
 T^+_1=\frac{1+G+G^2/[2(1-G/3)]}{2-a},
 \qquad T^+_{j+1}=T^+_j(j+a)/(j+2-a).
 \end{gathered}
\]
Indeed for $k\ge4$ the coefficient of $a^k$ in $\log B_0$ is at most
$4(\zeta(4)-1)<.33$ by the same termwise comparison as before
($\zeta(4)=\pi^4/90$). Sum the exponential using
$k!\ge2\cdot3^{k-2}$ for $k\ge2$.

For $1\le j\le3$ enclose the norm in each paired angular order by
\[
 B_j=d_{j0}L_{j0}+d_{j1}L_{j1}+(j+3)L_{j2}
  \left(T^+_j-\frac{d_{j0}}{j+1}-\frac{d_{j1}}{j+2}\right).
\]
For $B_1^\pm$ use the same expression at $j=1$, keeping only the
summand $m=1\pm x$ under each square root $L_{1p}$.
These bounds use the decreasing sequence $(p+j+1)L_{jp}$ separately
in each sign. For $j\ge4$ put
\[
 T_{\rm hi}=\frac{3(5T^+_4)^2\,16}{2(7/2)^4(4-4a)(1+a/2)},\qquad
 S_0=B_1^2+B_2^2+B_3^2+T_{\rm hi},\quad c=x^{5/2}\sqrt{S_0}/Z_- .
\]
In fact $U_j\le U_4(4/j)^{1-2a}$: the logarithmic ratio argument above
requires only $j\ge2$, $a\le2/5$, since
$1/j+1/(j+2-a)\le2/(j+1/2)$. Multiply its square by
$L_{j0}(0)^2=3/(2j^3)$, and integrate from $7/2$ as before;
$(7/8)^{4a}\le(1+a/2)^{-1}$.
Here $L_{j0}(x)\le L_{j0}(0)$ follows from
$(1-x^2/j^2)(1+x/(2j))^2\ge1$ for $x/j\le1/5$.

Write $v_1=(1+a)^{-1}+T^+_1-1/2$ and $v_2=2T^+_2+1/6-a/12$.
These bound $w_1,w_2$ by their definitions; $w_j$ decreases.
In the Schur bound for the high orders $m=|ni+1|$ use the weights
$m^{-1/2}$ at $|i|=1$ and $m^{-1/2}/2$ at $|i|\ge2$. The respective
off-diagonal row sums are at most
\[
 N=\frac{2x^3}{Z_-}\left[\frac{v_2}{2(1-x)}
       +\frac{v_1/2+v_2}{2(1-x/2)}\right],\qquad
 F=\frac{2x^3v_1}{Z_-(1-x)}\left(\frac52+\frac43\right).
\]
Indeed each unweighted summand is bounded by
$2x^3 w_{|i-k|}/(Z_-|k+x|\,|i-k|)$ after applying the weights
$|nk+1|^{-1/2}$. At $|i|=1$ the contribution from the other of the two orders $1\pm n$ is bounded by
$v_2/[2(1-x)]$. On $|k|\ge2$ the convolution weights
$1/(|k|\,|i-k|)$ sum to $3/2$, of which $1/2$ has $|i-k|=1$.
At $|i|\ge2$ the two orders $1\pm n$ have combined convolution weight
$1/(|i|-1)+1/(|i|+1)\le4/3$. All other weights remain unchanged.

For the diagonal of the disk field, the positive excess radial measure
has mean $1-g^2$ at most $1-M_{n+1}$ by monotonicity of the moments.
Its mass is $Z-1$. The mass identity in the preceding lemma therefore gives
\[
 q\ge q_0=1/3.39+
 \frac{L_-}{(1+L_-)3.39(2.39)}
 \left[1-\frac{3.39(4.78+.78/(2n+5))}{(2n+3)(2n+4)}\right].
\]
In order $m$, Lemma~\ref{IV:lem:disk-angular-floors} gives the disk eigenvalue lower bound $f_m=m^2+2m+2m^2/(m+6)$ by inserting $(m+6)/2$ as an upper
bound for its positive root. In this order the excess weight contributes at
most $L_+/(2m)$ by radial evaluation. Thus use the four denominators
\begin{align*}
 D&=q_0-1/(28Z_-)-L_+/2-4c,\\
 E_\pm&=Z_-^2[q_0-1/(Z_- f_{n\pm1})-L_+/(2(n\pm1))-N-c/4],\\
 E_t&=Z_-^2[q_0-1/(Z_-[(2n-1)^2+2(2n-1)])-L_+/(2(2n-1))-F-c/4].
\end{align*}
Completing the same squares as in the preceding lemma in the three
high-order groups gives
\[
 n^5(d_R-Bx^3)\le A+3.39(2.39)\left[\frac{xP'}D+
 \frac4{e_*}\left(\frac{(B_1^-)^2}{E_-}+\frac{(B_1^+)^2}{E_+}
                  +\frac{B_2^2+B_3^2+T_{\rm hi}}{E_t}\right)\right].
\]
With $d_{jp}=a_pa_{j+p}/a$ and $a_{k+1}=a_k(k+a)/(k+1)$ the numerator
bounds are
\[
 (B_1^-,B_1^+,B_2,B_3,T_{\rm hi})<(.995,.770,.364,.175,.0288),
 \qquad (A,P')<(3.263,.776).
\]
The denominators satisfy $(D,E_-,E_+,E_t)>(.191,.2676,.2797,.2591)$
componentwise. The one final substitution is
\[
 3.263+3.39(2.39)\left[\frac{.776}{6(.191)}
 +\frac{4}{2.3899}\left(\frac{.995^2}{.2676}+\frac{.770^2}{.2797}
 +\frac{.364^2+.175^2+.0288}{.2591}\right)\right]<97.71<98 .
\]
\end{proof}

\begin{theorem}\label{IV:thm:corner-gap}
If a simple polygon $P$ has at most $n-1$ strictly convex corners, then
\[
 F(R_6)-F(P)>1/1250\quad(n=6),\qquad
 F(R_n)-F(P)>5/n^4\quad(n\ge7).
\]
For $n\ge8$ the gap exceeds $(63-350/n)/n^4$; at $n=7$ it exceeds $7.2/n^4$.

\end{theorem}

\begin{proof}
Normalize area to $\pi$. Let $B=4\zeta(3)$, $4.808<B<4.809$,
$3.3899<\lambda<3.389958$, and write $d_R=(\lambda-1)(\pi\lambda/F(R_n)-1)$.
The upper bound \eqref{IV:eq:regular-continuous-band} for $d_R$ gives
$d_R\le B/n^3+84/n^5$ for $n\ge10$, in particular
$\lambda_0\operatorname{sinc}(2\pi/n)\ge\lambda-34/n^3$ in
circumradius-one normalization. The latter follows by inverting
$1+d_R/(\lambda-1)$ and using $(1+y)^{-1}\ge1-y$.

\emph{One missing convex vertex.} Choose a Riemann map $f$ for which
the two disk Neumann coordinate functions have zero pullback mean
(as shown before Lemma~\ref{IV:lem:signed-circle}). Put
$f'/f'(0)=1+\sum_{k\ge1}\gamma_k z^k$,
$M_j=2j\int_0^1g^2r^{2j-1}dr$, $M=M_1$, $D_0=1-M=1/\lambda$,
\[
 T=\sum_{k\ge1}\frac{|\gamma_k|^2}{k+1},\qquad
 G=\sum_{k\ge1}(M_{k+1}-M)\frac{|\gamma_k|^2}{k+1},\qquad
 d_P=(\lambda-1)(\lambda\pi/F(P)-1).
\]
The centered Rayleigh estimate and monotonicity of the moments give
$G/D_0\le d_P(1+T)$ and
$(5-\lambda)T/3\le G/D_0$; here
$M_2=2(\lambda+1)/(3\lambda)$ by integrating Bessel's equation.

Let $m=n-1$, $w_k=|\gamma_k|^2/(k+1)$, $t_k=(k+1)/(k+m)$,
$r_k=(M_{k+1}-M)/D_0$.

At $n=6$ take $m=5$. For $x=e^{-s}$ the shifted circle mean
reads
\[
 4[(1+b)^2x^5-b^2]\le
 \sum_{k\ge1}(k+1)(k-2b)^2 x^k w_k .
\]
The weight $(1-x)^2$ offsets the quadratic growth in the order $k$ after
integration, since $\int_0^1x^{k-1}(1-x)^2\,dx=2/[k(k+1)(k+2)]$.
Choose the shift $b=a x^3/(1-x)$, $a=9/40$, and integrate with
weight $(1-x)^2/x$. For $p>0$ write
$M_{p,j}=j!/[p(p+1)\cdots(p+j)]$. Integrating each nonnegative term
on the right gives $J\le\sum C_kw_k$, where
\begin{align*}
 J&=4[M_{5,2}+2aM_{8,1}+a^2(M_{11,0}-M_{6,0})]=353/7392,\\
 C_k&=(k+1)[k^2M_{k,2}-4akM_{k+3,1}+4a^2M_{k+6,0}].
\end{align*}
We show $1+C_k/J<(57/2)r_k$ simultaneously for all $k$.
At $k=1$ use $r_1=(5-\lambda)/3>8/15$; the difference with $8/15$
in place of $r_1$ is $1611/1765$. For $k\ge2$, Lemma~\ref{IV:lem:uniform-harmonic-correlations}
gives $r_k\ge s_k=1-17/[(2k+3)(2k+4)]$ since
$3.39[4.78+.78/(2k+5)]<17$. Multiply $(57/2)s_k-1-C_k/J$ by
$35300(k+2)(k+3)(k+4)(k+6)(2k+3)/10^6$. The resulting polynomial
is bounded below by
\[
 k^3(.015 k^2-.311 k+2.73)+44k^2+16.9k-208.
\]
Indeed its coefficients in descending degree are
$.015884$, $-.310794$, $2.732485$, $44.088555$, $16.981506$ and $-207.116136$.
The quadratic in parentheses has negative discriminant, and the
remaining quadratic increases from a positive value at $k=2$.
Thus $1+T<(57/2)\sum r_kw_k\le(57/2)d_P(1+T)$, giving
$d_P>2/57$.

For the regular hexagon,
Lemma~\ref{small-schur-regular} gives
$d_R<4.809/6^3+98/6^5$. The function
$F=\pi\lambda/(1+d/(\lambda-1))$ decreases in $d\ge0$
and increases in $\lambda>1$. For explicit endpoints put
$\lambda_-=3.389957716$ and $\lambda_+=3.389957717$; the disk
series brackets $\lambda$ between them. Then
\[
\begin{aligned}
 F(R_6)&>
  \frac{3.1415926\,\lambda_-}
  {1+(4.809/6^3+98/6^5)/(\lambda_--1)}>10.4967,\\
 F(P)&<
  \frac{3.1415927\,\lambda_+}
  {1+(2/57)/(\lambda_+-1)}<10.4959 .
\end{aligned}
\]
Their difference exceeds $.0008=1/1250$.

For $n\ge7$, in \eqref{IV:eq:integrated-turn-functional} take
\[
 W(s)=D_0^{-1}\int e^{-v}[W_0(s)-W_0(s-v)]\,d\nu(v),\qquad
 d\nu(v)=-d[g(e^{-v/2})^2].
\]
Integration by parts gives $\int e^{-v}d\nu=D_0$ and
$k^2(k+1)\mathcal L W(k)=r_k$, so $0\le W\le W_0\le s^2/2$.
Expanding $U/\Phi$ and using the centered Rayleigh test gives
\begin{equation}\label{IV:eq:missing-turn}
 J_m\le\sum_k(r_k+\ell t_k)w_k,\qquad
 \sum_k r_kw_k\le d_P(1+\textstyle\sum_k w_k).
\end{equation}
For $m\ge6$ one has $1-r_k\le.70(1-t_k)$: at $k=1$ the ratio is
$(\lambda-2)(m+1)/[3(m-1)]$. For $k\ge2$, Lemma~\ref{IV:lem:uniform-harmonic-correlations} instead gives
\[
 \frac{1-r_k}{1-t_k}\le\frac{4.92(3.39)(k+m)}
 {(2k+3)(2k+4)(m-1)}\le\frac{4.92(3.39)(7)}{7\cdot8\cdot4}<.53 .
\]

Choose any $0<j\le J_m$ with $j\le\ell/3$. Then
$(\ell+j)(1-r_k)\le\ell(1-t_k)$, or
$r_k+\ell t_k\le(1+\ell+j)r_k-j$. Summing gives $j(1+T)\le(1+\ell+j)d_P(1+T)$ and hence
\begin{equation}\label{IV:eq:turn-direct-lower}
                        d_P\ge j/(1+\ell+j).
\end{equation}
Lemma~\ref{IV:lem:positive-turn-kernel} gives an admissible $j$ upon
dividing its lower bounds by $m^3$: each numerator is below $5$, so
$j/\ell\le5(m-2)/(2m^2)<1/3$. 

At $n=7,8,9$ use the same lower bounds for $J_m$ and the constants
$U(n)=92.5,85.25,80.75$ from Lemma~\ref{IV:lem:uniform-regular-remainder}.
Inverting their endpoints gives $n^4(F(R_n)-F(P))>(7.2,19.4,27.1)$.

For the remaining counts we show $d_P>B(1+3/n)/n^3$.
For $m=n-1\ge9$ put $h=1/m$ and take $j=h^3[B-(10-50h)h^2]$.
By \eqref{IV:eq:turn-direct-lower} the desired comparison follows from
\[
 B(6+4h)-B(1+4h)\bigl[2/(1-2h)+Bh\bigr]>10-10h .
\]
Indeed $j\le Bh^3$, $(1+h)^4\ge1+4h+(6+4h)h^2$ and
$(10-50h)(1+h)^4\le10-10h$ by concavity. The displayed difference
is concave, positive at $0,1/9$. Hence for $n\ge10$, with $x=1/n$,
$d_P>d_R+\delta$, $\delta=x^4(3B-79.63x)$ by the upper bound \eqref{IV:eq:regular-continuous-band} for $d_R$.
Since $d_R+\delta\le Bx^3(1+3x)<.008$, inversion multiplies $\delta$
by at least $\pi\lambda(\lambda-1)/(\lambda-1+.008)^2>4.37$.
This gives the stated gap.

\end{proof}

\begin{remark}[A stronger hexagon bound from a polynomial trial field]
\label{rem:hexagon-missing-trial}
The degree-nineteen field in Section~\ref{I:sec:five-six-regular-data}
gives a second lower bound for the regular hexagon of circumradius one, $\mu_1(R_6)>4.04315$.
Since $|R_6|=3\sqrt3/2$, it follows that $F(R_6)>10.504$. Together
with the competitor bound $F(P)<10.4959$ established in
Theorem~\ref{IV:thm:corner-gap}, this also gives
$F(R_6)-F(P)>1/200$ for a polygon with at most five strictly convex corners.
\end{remark}

\section{The projected disk field on the triangles of a polygon}
\label{IV:sec:background}

For the construction and identities in this section assume $n\ge3$. A single sector has $0<hq<\pi$ and $R>|w|/2$; the centered fan of a convex polygon satisfies these conditions. The endpoint radii of the regular sector are $(1,1)$. Normalize the competitor to the regular
area $A_0=n\sin h/2$, $h=2\pi/n$. We use bisector coordinates on the reference and target triangles.

\paragraph{The trial field and its normalization.}
Let $j$ be the first positive zero of $J_1'$, $\lambda=j^2$,
$g(r)=J_1(jr)/J_1(j)$ and $U(x)=g(|x|)x/|x|$, with $U(0)=0$.
The same formula defines the field for $r>1$.
Let $\lambda_0=\mu_1(R_n)$ and let $Y$ be the componentwise $L^2(R_n)$
orthogonal projection of $U$ onto the first positive eigenspace. Set $D=Y-U$.
On one regular fan triangle $C_*$ put
$E_D=\int_{C_*}|\nabla D|^2$, $M_D=\int_{C_*}|D|^2$ and $S_D=\sqrt{E_D}$.
Use coordinates relative to the chosen fan center. A target triangle has opening $hq$ and endpoint
radii $R-w/2,R+w/2$, where $R=1+e$. We call $q$ the normalized aperture of the sector, $e$ the mean displacement and $w$ the skew displacement of its endpoint radii. In bisector coordinates put
$c=\cos(h/2)$, $s_0=\sin(h/2)$, $C=\cos(hq/2)$, $S=\sin(hq/2)$,
$A=C/c$, $B=S/s_0$. The affine map between the sectors, its determinant and the trial field are
\[
 T=\begin{pmatrix}RA&w C/(2s_0)\\w S/(2c)&RB\end{pmatrix},
 \quad d_R=(R^2-w^2/4)AB,\qquad V(Tx)=U(Tx)+T D(x).
\]
Thus $D$ is the correction of least $L^2$ norm for which each component of $U+D$ lies in the
first eigenspace of $R_n$. By reflection and rotation, $D$ is radial on the vertex rays
with common trace at equal fractional radii. The sector maps $T$ carry each ray
to its target ray, so the sector fields agree on their common rays and define a field in $H^1$.
Set
\[
 G=\int_{TC_*}(|\nabla V|^2-\lambda_0|V|^2+c_0),\qquad
 c_0=\lambda_0 A_0^{-1}\int_{R_n}|Y|^2 .
\]
At $(1,0,0)$, $V=Y$. To sum the sector contributions and express their variation with the endpoint radii, we use Lemma~\ref{IV:lem:projected-fan-identities}.
\begin{lemma}[Mass and boundary traces of the regular-sector correction]\label{IV:lem:bg-spectral}
Fix $n\ge10$, $h=2\pi/n$, $0<c\le\cos(h/2)$, $a\ge\tan(h/2)$,
$K=(\sqrt{4+a^2}+a)^2/4$, $G=28/K$. Suppose $\lambda_0\le\Lambda<28$,
$|\lambda_0-\lambda|\le d h^2$, $0<\sqrt{s}\le S/h^2\le z$,
$S^2=E_D$. Choose $b_1\ge3.59/(c\sqrt{120})$,
$b_2\ge1.59d/(12c\sqrt2)$. Set
\[
\begin{aligned}
 P&=[zh^{3/2}/G+b_1h^2/(c\sqrt{2G})]/.592,\qquad B=2.39\sqrt h/(12c\sqrt s),\\
 V&=\sqrt{2Ka^3/(3h)}+\Lambda h^2/(4c\sqrt{28G})+b_2h^{5/2}/\sqrt{56s},\\
 W_E&=(V+.405B)/(1-\Lambda/28),\qquad
 W_M=(V/\sqrt{28}+.05613B)/(1-\Lambda/28).
\end{aligned}
\]
For $W_E^2<K$, upper bounds for $M_D/(hS^2)$ and the total and transverse
boundary trace norms in $L^2(d\theta/h)$ divided by $S$ are respectively
\[
 \begin{aligned}
 m&=[P^2+(1-28/(n-1)^2)W_M^2+K/(n-1)^2]/h,\\
 T^2&=[P/(.592\sqrt h)+.405W_E/\sqrt h]^2+(K-W_E^2)/(2\pi-h),\\
 T_2&=P\sqrt h/(.592\sqrt{12})+\sqrt{K[.405^2 h/12+1/(2\pi-h)]}.
 \end{aligned}
\]
\end{lemma}
\begin{proof}
Pull back by $X=r\rho e^{i\theta}$, $\rho=\cos(h/2)/\cos\theta$.
The Jacobian determinant lies in $[c^2,1]$ and the energy distortion is
at most $K$. Projection gives $\int U\cdot D=-M_D$,
$\int\nabla U:\nabla D=-S^2$ in the sector. By
Lemma~\ref{IV:lem:disk-angular-floors} the next eigenvalue in rotation
character one is at least $G$, so $M_D\le S^2/G$.
Write $D\circ X=p+w+v$ for its components along the disk eigenfunction, in the radial complement of angular order one, and in the higher orders. These components are orthogonal in both disk energy and disk mass. The first has $\|p\|/S\le P$: use
$\|g e^{i\theta}\|>.592\sqrt h$, $\|D\circ X\|\le S/(c\sqrt G)$
and the projection identity. Indeed
$1-\rho\le h^2[1/4-(\theta/h)^2]/(2c)$ has angular mean
$h^2/(12c)$ and root mean square $h^2/(c\sqrt{120})$; since $|g'|<1.59$,
the root mean square of $g(r)-\rho^2g(r\rho)$ is at most $b_1h^2$.

On the radial complement the function representing the boundary trace in the energy inner product is $r-2g/(\lambda-1)$.
Its squared energy and mass in the radial measures are respectively
$1-2/(\lambda-1)<.405^2$ and $1/4-2/[\lambda(\lambda-1)]<.05613^2$:
integration against the disk equation gives the projection of
$r$ onto $g$ in the energy inner product. After multiplication by $e^{-i\theta}$ the boundary flux is
$\rho g'(\rho)-ig(\rho)\tan\theta$; its imaginary mean vanishes and
the Bessel equation and $g'(1)=0$ give
\[
 r g'(r)=\int_r^1(\lambda t-1/t)g(t)\,dt\le(\lambda-1)(1-r)
 \qquad(0<r\le1),
\]
since $0\le g\le1$. Taking $r=\rho$ gives the stated bound $B$.
In the orthonormal polar frame the pullback energy matrix differs from
the identity by
$\bigl(\begin{smallmatrix}\tau^2&-\tau\\-\tau&0\end{smallmatrix}\bigr)$,
where $\tau=\tan\theta$. For a test $\phi(r)e^{i\theta}$, the two
gradient components have the same angular factor. In the squared norm
of the matrix difference applied to this gradient, the mixed term
integrates to zero on $[-h/2,h/2]$ because its coefficient is $-2\tau^3$.
The two diagonal coefficients are $\tau^2+\tau^4$ and $\tau^2$.
Their angular averages are at most $2a^3/(3h)$, since
$\int_{-h/2}^{h/2}(\tau^2+\tau^4)\,d\theta=2\tan^3(h/2)/3$.
This gives the claimed bound in terms of the energy of the test.
The mass distortion uses $1-\rho^2\le h^2/4$; for the volume source,
radial orthogonality removes $\rho^2g(r)$, leaving mean variation
at most $1.59rh^2/(12c)$. These give the three terms of $V$.
The disk eigenvalues on this radial complement exceed $28$, so the shifted inverse on this
complement gives $\|w\|_E/S\le W_E$, $\|w\|/S\le W_M$.

For the higher orders $|k|\ge n-1$, the mass is at most $1/(n-1)^2$ times the energy
and the squared boundary trace norm is at most $1/(2\pi-h)$ times the energy. Their energy
is at most $KS^2-\|w\|_E^2$, where $\|w\|_E^2\ge28\|w\|^2$.
The total squared trace increases in $\|w\|_E/S$ because
$.405^2/h>1/(2\pi-h)$. By reflection symmetry the boundary traces of
the angular order-one terms $p,w$ are real multiples of $e^{i\theta}$.
Their transverse components therefore carry the factor
$\|\sin\theta\|_{L^2(d\theta/h)}\le h/\sqrt{12}$. For the higher
orders $v$ use the full boundary trace bound above. Orthogonality gives
$\|w\|_E^2+\|v\|_E^2\le KS^2$, so Cauchy--Schwarz yields
\[
 .405\sqrt{h/12}\,\frac{\|w\|_E}{S}
 +\frac{\|v\|_E}{S\sqrt{2\pi-h}}
 \le \sqrt{K\bigl(.405^2h/12+1/(2\pi-h)\bigr)}.
\]
Adding the contribution $P\sqrt h/(.592\sqrt{12})$ from $p$ gives $T_2$.
\end{proof}

\paragraph{Sector energy at equal endpoint radii.}
The following changes of variables and integrations by parts use only $0<hq<\pi$ and $R>|w|/2$. For the equal-radius calculation take $(R,w)=(1,0)$.
At equal endpoint radii let $d=AB$, $\omega=dB^2$, $C_q$ the cross
energy minus $\lambda_0$ times cross mass of $U,TD$, and
$\Phi'(r)=r[g'^2+g^2/r^2-\lambda_0g^2+c_0]$. The nonlinear disk
term is $hq\int[\Phi(\cos(hq/2)/\cos(hqs))-\Phi(1)]ds$.
Projection gives $C_1=-E_D+\lambda_0M_D$ and
$A_0-\int|Y|^2=\pi/\lambda-\int_{\mathbb D\setminus R_n}(1-g^2)+nM_D$.
We call the energy minus $\lambda_0$ times the mass of $TD$ alone the \emph{pure term}, in contrast with the cross term $C_q$.
Absorbing $2\omega C_1$ in the pure term changes its gradient weights to
$d(1,A^2/B^2,B^2/A^2,1)-2\omega$ and its $\lambda_0D_i^2$ weights to
$(2\omega-dA^2,\omega)$. By Green's formula, the remaining $C_q-\omega C_1$
has boundary term with vector coefficient $F_q-\omega F_1$ and volume factor
$(\lambda-\lambda_0)d[U(Tx)\cdot TD-B^2U(x)\cdot D]$:
the ray flux is transverse to $TD$.
On the chord set $v=s_0t$, $p=\sin^2(hq/2)$, $p_1=s_0^2$, $z=1-t^2$,
$K(X)=g(\sqrt X)/\sqrt X$, $k(x)=g'(\sqrt{1-x})/x$, $j(x)=2K'(1-x)$.
With $\Delta k=k(pz)-k(p_1z)$, and likewise $\Delta j$, the two components of the boundary term are $dp(z\Delta k-t^2\Delta j)$ and $cv\omega\Delta j$.
Here $k(x)=-\int_0^1(3K'+2XK'')(1-vx)\,dv$, using $g'(1)=0$.
For $H_A=dA^2,H_B=dB^2,K_q=K(|Tx|^2)$ the integrated volume divided
by $\lambda-\lambda_0$ is
\[
 (H_B-H_A)\left(M_D+\int K_1x_2D_2\right)
 +\int(K_q-K_1)(H_Ax_1D_1+H_Bx_2D_2)
\]
by $\int U\cdot D=-M_D$.

For the radial first variation differentiate $G_e(p,0,0)/p$ and
average in $p^2$ between $1$ and $q^2$, using $G_e(1,0,0)=0$.
Its disk integrand is $h\rho_p\Phi'(\rho_p)$; the radial powers
in the pure energy and mass are $R^2,R^4$. With $D_0=AB/p$,
$W=4K+2XK'$, the differentiated cross volume divided by
$2(\lambda-\lambda_0)p$ is $\int D_0W(X)(A^2x_1D_1+B^2x_2D_2)$.
For the boundary term this derivative divided by $p$ has components
$D_0(2g'+rg''-ut^2H_1),D_0cB^2vH_1$, where $H_1=-\lambda g/r$ by
the disk equation and $u=\sin^2(hp/2)$.

\paragraph{Energy and mass of the correction under the affine map.}
These identities hold for every $q>0$ with $hq<\pi$ and every positive pair of endpoint radii.
Put $k=w/(hR)$, $\delta=h^2/4$, $T=R(A_1+kB_1)$ with
$A_1=\operatorname{diag}(A,B)$ and $B_1$ off-diagonal with entries
$c_1=hC/(2s_0),d_1=hS/(2c)$. By reflection the correction energy and mass are
\begin{equation}\label{IV:eq:any-n-affine}
 R^2\frac{P_0+P_2k^2+P_4k^4}{1-\delta k^2},\qquad
 R^4AB(1-\delta k^2)\int(|A_1D|^2+k^2|B_1D|^2).
\end{equation}
Here $\sum P_{2j}k^{2j}=(AB)^{-1}\int|(A_1+kB_1)\nabla D\,
\operatorname{adj}(A_1+kB_1)|^2$. For $\tau=B/A$ and
$k_0=c_1^2/\tau+\tau d_1^2$ the diagonal-gradient block of $P_2$
has eigenvalues $-2c_1d_1,2k_0+2c_1d_1$; the off-diagonal block has
diagonal $(2d_1^2/\tau,2\tau c_1^2)$ and off-diagonal $-k_0-2c_1d_1$.

\begin{proposition}[A chord identity for the regular-sector correction]
\label{IV:prop:multiharmonic-correction}\label{IV:cor:sharp-correction-energy}\label{IV:lem:scalar-chord-energy}
For $n\ge3$, put $b=\sin(h/2)$, $S=\operatorname{sinc}h$, $L(p)=1-g(\sqrt{1-p})^2$, and $s_D=(E_D-\lambda_0M_D)/h^4$. Let $C=\int_{\mathrm{cap}}L(1-|x|^2)\,dx$ on the circular segment between the regular chord and the unit circle. Then
\begin{equation}\label{IV:eq:any-n-chord}
 h^4s_D=hS\int_0^1 L'(b^2(1-t^2))\,dt
 -(\lambda_0-\lambda)\{h(S-1/\lambda)/2+C\}.
\end{equation} The bounds from ten sides are given in Section~\ref{mid-response}. As $h\to0$, $E_D/h^4=\zeta(3)/(4\pi^3)+O(h)$.
\end{proposition}
\begin{proof}
Put $Y=h^2$, $b=\sin(h/2)$, $c=\cos(h/2)$, $S=\operatorname{sinc}h$,
$e=\lambda-1$ and $L(p)=1-g(\sqrt{1-p})^2=\sum_{j\ge2}\ell_jp^j$.
Orthogonality to $Y_0=U+D$ gives $A_{\lambda_0}(D)=A_{\lambda_0}(U)$.
On the chord $\partial_\nu |U|^2/2=cL'(1-|x|^2)$; on the rays it vanishes.
Put $J=\int_0^1L'(b^2(1-t^2))\,dt$. Subtracting the cap from the circular sector in Green's identity yields \eqref{IV:eq:any-n-chord}.

For the asymptotic formula the upper and lower bounds for $d_R$ give
$d_R=4\zeta(3)n^{-3}+O(n^{-5})$. Section~\ref{mid-response} rewrites the chord term as $G_c/h^3=O(h)$, hence $s_D=\zeta(3)/(4\pi^3)+O(h)$. Its uniform mass estimate gives $M_D/E_D=O(h)$, proving the assertion.
\end{proof}

On the reference triangle $|Tx|^2=A-pZ$ at zero offsets, where
$A=(x_1/c)^2$, $Z=A-(x_2/s_0)^2\in[0,A]$ and $p=\sin^2(hq/2)$.
Integrating on $0\le x_1\le c$, $|x_2|\le s_0x_1/c$ gives
$\|x_1\|^2=c^3s_0/2\le h/4$, $\|x_2\|^2=cs_0^3/6\le h^3/48$.
Since $D_2=0$ on the bisector, the Dirichlet--Neumann Poincar\'e bound
on each vertical half-fiber gives $\|D_2\|\le2s_0\sqrt{E_D}/\pi$.

\begin{lemma}[Flux across the chord and its radial derivatives]\label{IV:thm:radial}\label{IV:lem:radial-forty}
With $R=1+e$, $u=\sin^2(hq/2)$, $a=1-u+ut^2$, $b=u+(1-u)t^2$,
$\kappa=\operatorname{sinc}(hq)/\operatorname{sinc}h$ and
$k=h\cot(h/2)$, put
\[
 \begin{aligned} P&=R^2-w^2/4,&X&=aR^2+tRw+bw^2/4,\\
  V&=(1-u)R^2+tRw/2+uw^2/4,&W&=utR^2+Rw/2+(1-u)tw^2/4.\end{aligned}
\]
The coefficients $(F_1,hF_2)/q$ on the reference chord $|t|\le1$ are
$\kappa P(K+2VK',2kWK')$, with the derivatives of $K$ taken at $X$.
For $w=0$ their second $e$-derivatives are
$\kappa(A+2VB,2kWB)$, where $A=2K+10XK'+4X^2K''$ and $B=A'$.
In their pairing with $D$, reflection selects the even and odd
parts in $t$, respectively, and makes both integrals even in $w$.
\end{lemma}
\begin{proof}
Green's formula gives $T^t DU(Tx)(\operatorname{cof}T)e_1$ and
$DU(y)=KI+2K'yy^t$. Substitute the sector matrix and differentiate.
For $d=\det T$, $v=\dot Tx$, $b=TD$, the volume integrand has second variation
\[
 \ddot d\,U\cdot b+2\dot d(DU[v]\cdot b+U\cdot\dot b)
             +d(D^2U[v,v]\cdot b+2DU[v]\cdot\dot b).
\]
Here $T$ is affine in the endpoint radii. For $j=\sqrt\lambda$,
$U(y)=J_1(j)^{-1}\langle\xi\sin(j\xi\cdot y)\rangle_{|\xi|=1}$
gives $\|D^aU\|\le j^a/(2J_1(j))$ for $a\ge1$ by Cauchy--Schwarz
on two factors. The alternating Bessel series gives $J_1(j)>.581$,
so $\|DU\|<1.59$, $\|D^2U\|<2.92$. Also $|U|\le1$: beyond
$r=1$, $p=\lambda-r^{-2}$ gives $(g^2+g'^2/p)'=-(2/(rp)+p'/p^2)g'^2$.
\end{proof}

\begin{lemma}[Identities for the projected field on a convex polygon]\label{IV:lem:projected-fan-identities}
Let $n\ge3$, let $P$ be a convex $n$-gon of area $A_0$ and take any interior fan center. Use $h,q,e,w$ and $G$ as defined above, and write $\overline V=|P|^{-1}\int_P V$. Put
\[
 D(q)=G(1,0,0)+(q-1)G_q(1,0,0)-G(q,0,0),\qquad
 C(q,e,w)=G(q,e,w)-G(q,0,0).
\]
For the glued field one has the exact identity
\[
 \int_P\bigl(|\nabla V|^2-\lambda_0|V-\overline V|^2\bigr)
 =\sum_i\bigl(C(q_i,e_i,w_i)-D(q_i)\bigr)
       +\lambda_0 A_0|\overline V|^2.
\]
At the regular sector, $G_e(1,0,0)=0$. By reflection, $G_w(q,e,0)=0$ throughout the allowed range of sector parameters. The radial increment and the coefficient of its linear term satisfy
\begin{align*}
 C(q,e,w)={}& eG_e(q,0,0)
 +e^2\int_0^1(1-t)G_{ee}(q,te,0)\,dt\\
 &+w^2\int_0^1(1-t)G_{ww}(q,e,tw)\,dt,\\
 G_e(q,0,0)={}& (q-1)\int_0^1G_{qe}(1+t(q-1),0,0)\,dt .
\end{align*}
\end{lemma}
\begin{proof}
All sectors of the regular fan have equal mass and energy by equivariance. Hence $nG(1,0,0)=c_0A_0$ by the regular weak equation. Since $\sum_i q_i=n$, subtracting the affine part of $G$ at $q=1$ cancels exactly this area term. Subtracting the mean of the field gives the last term in the identity.
To vary the regular radii by a common factor $R$, pull the field back to the fixed polygon as $U(Rx)+R D(x)$. In two dimensions the gradient energy has no power of $R$ outside its integral, and the mass and area each carry the factor $R^2$. Differentiating at $R=1$, the weak equation for $Y$ cancels the derivative of the field. The remaining terms are $2c_0 A_0-2\lambda_0\int_{R_n}|Y|^2=0$. The sectors are congruent, so $G_e(1,0,0)=0$. Reflection in the bisector sends $w$ to $-w$ and preserves energy and mass, giving the other vanishing derivative. Integrate first in $e$ at $w=0$ and then in $w$ at fixed $e$, using Taylor's formula with integral remainder. The segment lies in $R>|w|/2$, and differentiation under the fixed reference-triangle integrals follows from invertibility of $T$, smoothness of $U$, and $D\in H^1$.
\end{proof}

\paragraph{Changes in energy with the opening angle and endpoint radii.}
Give $P$ the area of the circumradius-one $R_n$, put $h=2\pi/n$ and
write $\lambda_0=\mu_1(R_n)$. By Theorem~\ref{IV:thm:corner-gap}
a competitor is strictly convex with $n$ sides. Use the glued projected field
$V(Q+Tx)=U(Tx)+TD(x)$ of Section~\ref{IV:sec:background}, with
$D=Y-U$ and center $Q$ chosen so that $\int_P V=0$. Let $G(q,e,w)$ be
the sector functional defined there. The opening angle is $hq$ and the endpoint
radii are $1+e\pm w/2$. Set $d=h^2$, $v=q^2-1$, $x=e/d$, $y=w/d$,
$J=G_e(q,0,0)/(qh^3v)$. In the radial Taylor formula of
Lemma~\ref{IV:lem:projected-fan-identities}, denote the averages of
$-G_{ee}/(qh)$ and $-G_{ww}/(qh)$ against $2(1-t)\,dt$ by $K_e,K_w$.
Then
\begin{equation}\label{mid-radial-identity}
 [G(q,e,w)-G(q,0,0)]/(qh^5)=Jvx-K_ex^2/2-K_wy^2/2 .
\end{equation}
The preceding formula for the first radial derivative gives the continuous extension of $J$ at $q=1$. Put $\Delta(q)=G(1,0,0)+(q-1)G_q(1,0,0)-G(q,0,0)$ and $C(q,e,w)=G(q,e,w)-G(q,0,0)$.

\section{Bounds for the regular-sector correction for ten or more sides}\label{mid-response}

Here $X=h^2$, $b=\sin(h/2)$, $c=\cos(h/2)$, $S=\operatorname{sinc}h$, $J_0=bc$, and $e_0=\lambda-1$. Let $E_D,M_D$ be the energy and mass of the correction $D$ on one regular sector. Put $d_*=(\lambda_0-\lambda)/X$, $s_D=(E_D-\lambda_0M_D)/h^4$, $z=E_D/h^4$ and $\rho=\lambda_0 M_D/E_D$. The upper bound of Lemma~\ref{IV:lem:uniform-regular-remainder} and the lower bound of Lemma~\ref{IV:lem:conformal-sideband-energy} give the following bounds for the normalized spectral difference $d_R=e_0(\pi\lambda/F(R_n)-1)$:
\[
 {4.808+41.7/n^2\over n^3}<d_R<
 {4.809+U(1/n)/n^2\over n^3},\qquad
 \lambda_0={\lambda\over S(1+d_R/e_0)}.
\]
Here $U(x)=\max(74,59.93+197x)$. Put $x=1/n$ and
$D_-=(4.808+41.7x^2)x^3$. On $0<x\le1/10$ the regular polygon has
\begin{equation}\label{regular-common-ceiling}
 \lambda_0\le \overline\lambda(x)=\frac{\lambda}{S(1+D_-/e_0)} .
\end{equation}
This follows from the lower bound for the normalized spectral difference in Lemma~\ref{IV:lem:conformal-sideband-energy}. The upper bound increases with $x$. Indeed, the sine product gives $1/v-\cot v\ge v/3$, so the logarithmic derivative of $1/S$ is at least $4\pi^2x/3$, while $D_-'/ (e_0+D_-)<.70x$. For fifteen or more sides this gives the uniform upper bound $\overline\lambda(1/15)<3.489$.
For the regular chord expression set $G_* =e_0[S(4-\cos h)-3]/6+(1-S)^2/(2S)$ and $A=2(2\lambda-3)/15$. Let $L_0(p)=1-(1-p)K(1-p)^2$, $C=\int_{\rm cap}(1-|U|^2)$, and put
\[
 J=\int_0^1L_0'\bigl(b^2(1-t^2)\bigr)\,dt,\qquad
 G_c=SJ-\tfrac12(S^{-1}-1)(\lambda S-1).
\]
Substituting $\lambda_0=\lambda/[S(1+d_R/e_0)]$ in
\eqref{IV:eq:any-n-chord} and collecting the terms in $d_R$ gives
\[
s_D={G_c\over h^3}
  +{\lambda-1/S\over 2(e_0+d_R)}{d_R\over h^3}
  -{(\lambda_0-\lambda) C\over h^4},\qquad 0\le C\le {3hX^3\over768}.
\]
On the cap, $r\ge c=\cos(h/2)$. Integrating
$tg'(t)\le(\lambda-1)(1-t)$ and using $g(1)=1$ gives
$1-g(r)^2\le(\lambda-1)(1-r)^2/r\le3(1-r)^2$ for $n\ge10$.
Now $1-c\le X/8$ and the cap area is at most $hX/12$.

For a lower bound on $s_D$, retain the first four coefficients of
$L_0(p)=1-(1-p)K(1-p)^2$. The Bessel product formulas in the proof of Lemma~\ref{IV:lem:wall-zero-products} give $L_0'(p)=s_1pB(p)T(p)$ with $s_1>0$ and both factors having
nonnegative coefficients. Write $K(1-p)=\sum k_jp^j$.
The Bessel equation gives $k_0=1$, $k_1=1/2$ and
$k_{j+2}=k_{j+1}-\lambda k_j/[4(j+1)(j+2)]$. In particular
\[\fitmath{[p^2]L_0=(\lambda-1)/4,\qquad
 [p^3]L_0=(2\lambda-3)/12,\qquad
 [p^4]L_0=(6-\lambda)(\lambda-2)/48>.0755.}\]

The chord average of $L_0'(b^2(1-t^2))$ and the ratio bound on its positive coefficients give
\[
 S[A b^4+4(.0755)(16/35)b^6]\le G_c-G_*\le S A b^4/(1-.4b^2).
\]
These bounds give $\lambda_0<\Lambda$, $d_-<d_*<d_+$ and $s_-<s_D<s_+$ at ten and eleven sides:
\[
\begin{array}{c|r|rr|rr|rrrrr}
 n&\Lambda&d_-&d_+&s_-&s_+&m_*<&\rho<&z<&T_1<&T_2<\\\hline
10&3.61588&.5706&.5722&.01284&.01408&.0540&.1229&.0161&.654&.524\\
11&3.57582&.56846&.56953&.01267&.01362&.0500&.1026&.01520&.642&.510\\
\end{array}
\]
For endpoint evaluation at the disk limit let $X=h^2$ and set
\[
 \phi=(1-S)/X,\quad b_0=b^2/X,\quad R_0=[S(4-\cos h)-3]/X^2,
\]
\[
 Q=d_R/x^3\in[4.808+41.7x^2,\ 4.809+U(x)x^2].
\]
Their three removable limits are $1/6,1/4,-1/10$.
To evaluate the spectral difference near zero, write it in the form
$d_*=\lambda_0[\phi-SQx/(4\pi^2 e_0)]$. The two bounds on $G_c/h^4$ have common term
$e_0R_0/6+\phi^2/(2S)$; add
$S[Ab_0^2+4(.0755)(16/35)Xb_0^3]$ below and
$SAb_0^2/(1-.4b^2)$ above. Thus the expression for $s_D$ is regular at $h=0$, using $d_R/h^3=Q/(2\pi)^3$.
For the last five columns use Lemma~\ref{IV:lem:bg-spectral} with $K=(\sqrt{4+t^2}+t)^2/4$, $t=\tan(h/2)$ and $G=28/K$.
Since the second eigenvalue in rotation character one is at least $G$, we have $\sqrt{s_D}\le\sqrt z\le z_0=\sqrt{s_D/(1-\lambda_0/G)}$.
Use these bounds in $P,B,V,W_E,W_M$ there, with $b_1=3.59/(c\sqrt{120})$, $b_2=1.59d_* /(12c\sqrt2)$.
For $p=P/\sqrt h$, $w=W_E/\sqrt h$, $y=W_M/\sqrt h$ its three conclusions read
\[
 m=p^2+[1-28x^2/(1-x)^2]y^2+Kh/(2\pi-h)^2,\qquad
 T_1^2=(p/.592+.405w)^2+(K-hw^2)/(2\pi-h),
\]
\[
 T_2=ph/(.592\sqrt{12})+\sqrt{K[.405^2h/12+1/(2\pi-h)]}.
\]
These expressions extend continuously to zero: in $V/\sqrt h$ the first term is $\sqrt{2Kh(t/h)^3/3}$, where $t/h=\sqrt{b_0}/c$. The second and third terms use $h^{3/2}$ and $h^2$. Use $\rho\le\lambda_0hm$ and $z\le s_D/(1-\lambda_0hm)$. The trace bound increases with $w$ since $.405^2>h/(2\pi-h)$, so the same upper endpoint can be inserted in its two occurrences.

For ten and eleven sides substitute $x=1/10,1/11$, respectively. For $S,\phi,b_0,R_0$ the coefficients of $X^i$ are respectively
\[
 \frac{(-1)^i}{(2i+1)!},\qquad \frac{(-1)^i}{(2i+3)!},\qquad
 \frac{(-1)^i}{2(2i+2)!},\qquad \frac{(-1)^i(4-4^{i+2})}{(2i+5)!}.
\]
They give alternating terms of decreasing magnitude for $X\le.4$, so truncating at indices seven and eight bounds both endpoints. In these two substitutions use the upper endpoint for $w$ only after checking $hw^2<K$. The next proposition gives a separate endpoint argument on the entire remaining interval.

\begin{proposition}[Uniform bounds for the mass, energy and boundary traces of the correction]\label{continuous-response}
Put $H=\pi/6$, $c_H=\cos(H/2)$, $t_H=\tan(H/2)$ and set
\[
 \begin{aligned}
 \alpha_D&=\frac{\zeta(3)}{4\pi^3},\quad
 K_H=(\sqrt{4+t_H^2}+t_H)^2/4,\quad G_H=28/K_H,\\
 p_H&=\frac1{.592}\left(\frac{\sqrt{.015}\,H}{G_H}
       +\frac{3.59H^{3/2}}{c_H^2\sqrt{240G_H}}\right),\\
 \mathcal T_H&=\frac{p_H H}{.592\sqrt{12}}
       +\sqrt{K_H\!\left(\frac{.405^2H}{12}+\frac1{2\pi-H}\right)} .
 \end{aligned}
\]
For $0<h=2\pi/n\le H$ the quantities defined above for the regular sector satisfy
\begin{gather*}
 .56<d_*<.575,\qquad \alpha_D+h/200<s_D<.0135,\\
 m_*<.048,\quad \rho<.09,\quad z<.015,\quad T_1<.65,\quad T_2\le\mathcal T_H .
\end{gather*}
The positive projection on the cubic direction $j$ of
Lemma~\ref{IV:lem:common-diagonal-cubic} satisfies
$\langle\nabla D/h^2,j\rangle^2>.0072h$.
\end{proposition}
\begin{proof}
Write $X=h^2\le.275$, $x=h/(2\pi)$, $Q=d_R/x^3$, $A=2(2\lambda-3)/15$.
The lower estimate in Lemma~\ref{IV:lem:conformal-sideband-energy}
and the upper estimate in Lemma~\ref{IV:lem:uniform-regular-remainder} give
$4\zeta(3)+41.7x^2<Q<4.809+76.35x^2$ and
\eqref{regular-common-ceiling} gives $\lambda_0<3.546$ at $x=1/12$.
From $\lambda_0 S(1+d_R/e_0)=\lambda$,
\[
 d_*={\lambda\over X}(S^{-1}-1)-{\lambda_0Qx\over4\pi^2 e_0}.
\]
Use $1+X/6+7X^2/360\le1/S\le1+X/6+.021X^2$.
The lower bound is the positive cosecant series and the upper bound follows on multiplying by $1-X/6+X^2/120-X^3/5040$.
The upper bound $\lambda_0 Q/(4\pi^2e_0)<.202$ gives
$d_*>.56498+2.60x^2-.202x>.56$ by square completion.
Then $\lambda_0>\lambda$, so $\lambda_0 Q/(4\pi^2e_0)>.17$ and
$d_*<.565+2.82x^2-.17x<.575$ by convexity.

To bound the chord contribution, put
$b_-=1/4-X/48$, $b_+=b_-+X^2/1440$, $S_-=1-X/6$ and $S_+=S_-+X^2/120$.
Alternating sine sums give $b_-\le b^2/X\le b_+$.
The following two polynomials bound $G_c/h^4$ from below and above:
\begin{align*}
 P_-={}&{e_0\over6}(-1/10+X/84-X^2/1440)
   +\tfrac12(1/6-X/120)^2(1+X/6)\\
 &+S_-[A b_-^2+4(.0755)(16/35)Xb_-^3],\\
 P_+={}&{e_0\over6}(-1/10+X/84)
   +\tfrac12(1/6-X/120+X^2/5040)^2(1+.175X)\\
 &+S_+ A b_+^2(1+.103X).
\end{align*}
Indeed $1/S\le1+.175X$ and $(1-.4b^2)^{-1}\le1+.103X$ by $b^2\le X/4$.
Both polynomials have constant coefficient $1/180$.  In powers starting at
degree one, the coefficients of $P_-$ exceed
$(-.00268,-.00028,0,-.000009,0)$.
The respective upper bounds for $P_+$ are
\[
 (-.00147,.000275,.000045,0,.00000046,0,.000000001).
\]
Substitution of $X\le.275$ therefore gives
$1/180-.0028X<P_-\le P_+\le1/180$ for $X>0$.
The remaining positive term of $s_D$ is
$C_0=(\lambda-1/S)Q/[2(e_0+d_R)(2\pi)^3]$.
Since $d_R/X<.012$ and $(.175+.012)/e_0<.08$,
\[
 { (1-.08X)(4\zeta(3)+41.7X/(2\pi)^2)\over 2(2\pi)^3}
 \le C_0\le
 {(1-X/(6e_0))(4.809+76.35X/(2\pi)^2)\over2(2\pi)^3}.
\]
The lower bound is at least $\alpha_D+.00130X$. Indeed its
coefficient above the constant, divided by $X$, is bounded below by
\[
 \frac{(1-.08\cdot .275)41.7/(2\pi)^2-.08(4.809)}
      {2(2\pi)^3}>.00130 .
\]
The upper bound is below $.009695+.00323X$. The cap term is negative and has magnitude
at most $3(.575)h^5/768$. Consequently
\[
 \frac{s_D-\alpha_D}{h}\ge
 \frac1{180}+.00130h-.0028h^2-\frac{3(.575)}{768}h^4>\frac1{200}.
\]
The polynomial is concave for $h\ge0$ and exceeds $.00529$ at
both $0$ and $.524$. The upper bound is
$.009695+.00323H^2+H/180<.0135$.

Put $u=\sqrt{h/H}$. In the notation of the lemma on mass and boundary traces, let
$p=P/\sqrt h$, $v=V/\sqrt h$, $w=W_E/\sqrt h$, $y=W_M/\sqrt h$.
Use the increasing functions $\tan(h/2)/h$, $1/c$, $K$ and $1/G$
at $H$. Their needed values are $c>.9659$, $K(H)<1.307$,
$G(H)>21.43$. In the formulas below the actual endpoint $H=\pi/6$
may be used without rounding those three values. Thus $p<.0397u^2$ using
$z_0=\sqrt{.0135/(1-3.546/G(H))}$, and the three terms of $v$ are bounded
by $.24720u,.01420u^3,.02076u^4$ using $s_D>\alpha_D> .00969$.
Convexity of $s\mapsto1/\sqrt{s}$ and the affine lower bound on $s_D$
give the secant bound
\[
 B/\sqrt h\le {2.39\over12\cos(H/2)}
 \left[{1-u^2\over\sqrt{\alpha_D}}+{u^2\over\sqrt{\alpha_D+H/200}}\right]
 <2.095-.220u^2 .
\]
Consequently $v\le.2823u$ and
\[
 y\le Y=.135+.0612u-.014u^2,\qquad
 w\le W=.973+.324u-.101u^2.
\]
Both increase on $[0,1]$, and $hW^2<1\le K$.
For the mass coefficient, retaining the negative term from the higher orders gives
\[
 m_*\le .0397^2u^4+.02064u^2+
 \left[1-{28u^4\over(12-u^2)^2}\right]Y^2.
\]
This upper bound increases. The only negative contribution to its derivative has magnitude at most
$(1344/11^3)(.1822)^2u^3<.034u^3$, while the middle term has derivative
$.04128u$. Its endpoint is below $.048$, giving
$\rho\le3.546H(.048)<.09$ and $z\le .0135/(1-.09)<.015$.

In the trace estimate, convexity of $K$ in $h$ gives
$K\le1+.307u^2$. Indeed $K=1+t^2/2+t\sqrt{1+t^2/4}$ is increasing
and convex for $t\ge0$, as is $t=\tan(h/2)$. The trace formula increases
in $w$ since $.405^2>1/11$. Hence
\[
 T_1^2\le (.0676u^2+.405W)^2+(1+.307u^2)/5.759-u^2W^2/12 .
\]
The seven coefficients of the polynomial on the right in the degree-six
Bernstein basis on $[0,1]$ are bounded above by
$.330,.347,.365,.381,.395,.405,.413$.
For a polynomial $\sum a_i u^i$ the $j$th such coefficient is
$\sum_{i=0}^j a_i\binom ji/\binom6i$, and the basis is nonnegative
with sum one. Thus $T_1<.65$. In the transverse trace substitute
the energy bound $z<.015$ just obtained. Each positive factor in the
resulting bound on $p$ and $T_2$ increases with $h$; substitution at
$H$ gives $p\le p_H$ and $T_2\le\mathcal T_H$. The two comparisons
$.49599<\mathcal T_H<.49601$ follow directly from the displayed
formula: use $t_H=2-\sqrt3$, $c_H^2=(2+\sqrt3)/4$,
$3.1415926<\pi<3.1415927$, and $1.7320508<\sqrt3<1.7320509$,
rounding square roots outwards by squaring. Thus the distance to $1/2$ exceeds $.00399$. In the mixed
estimate we insert the displayed expression for $\mathcal T_H$.

Finally put $a=1/6+1/24+.076/20$, $b_H=\sin(H/2)$ and $S_H=\sinc H$.
The cubic mass satisfies $B/h\ge (S/2)[a-(a-.1)b^2]$.
For $\|U\|^2$, use the smaller bound in
Lemma~\ref{IV:lem:common-diagonal-cubic}, obtained by subtracting the integral over the circular cap. For the squared residual
$(\|\Psi\|^2-2B/3+\|U\|^2/9)/h$ it gives the upper bound
\[
 S\left[\tfrac1{18}+\tfrac12(\tfrac1{12}-2a/3)(1-b^2)
   +\tfrac12(\tfrac{29}{420}-\tfrac1{15})b^2\right]
 -{1\over18(3.39)}+{3H^6\over9(768)}<.0122.
\]
The bracket is positive and increases in $b^2$, so $S\le1$ and
$b^2\le b_H^2$ suffice. Divide the signed cubic estimate by $\sqrt h$
and use $J/h=S/2$. Dropping the positive fifth-order flux term gives
\begin{align*}
 {\langle\nabla D/h^2,j\rangle\over\sqrt h}\ge{}&
 \sqrt{5S_H/16}\left[.56\{a-(a-.1)b_H^2\}
              +\tfrac4{15}(b_H/H)^2\right]\\
 &-{3.546\over\sqrt{4S_H/5}}
 \left[.048(.015)H^2/3+\sqrt{.048(.015)(.0122)}\right]>.085 .
\end{align*}
Here $\sinc(h/2)$ decreases. Squaring proves the last assertion.
\end{proof}

\paragraph{The moments $I$, $\nu$ and $z_d$.}
Put $I=-h^{-3}\int D_1\,dv$, $\nu=h^{-4}\int vD_2\,dv$ on the regular
chord and $z_d=h^{-4}(\|\partial_1 D_1\|^2+\|\partial_2 D_2\|^2)$.
Apply \eqref{op-coordinate} and \eqref{op-torque} with the preceding bounds for the regular sector.
For the diagonal use \eqref{IV:eq:common-diagonal-test} at $\alpha=1/3$ and
the polynomial $B$ in Lemma~\ref{IV:lem:common-diagonal-cubic}. At ten and eleven sides retain the exact triangle moments and the smaller disk mass bound $H_U$ there. The last column bounds the square of the
positive projection on $j=\nabla\Psi/\|\nabla\Psi\|$ from below, and hence also $z_d$.
\[
\begin{array}{c|rrrr}
n&I>&I<&\nu>&z_d>\\\hline
10&.0212&.0434&.0033&.00470\\
11&.0230&.0425&.0042&.00445\\
\end{array}
\]

The same endpoint bounds give $ .0226<I<.0433$ for all $n\ge12$.
In \eqref{op-coordinate} the common error is less than $v=.005341$ by
$.05613(3.546)\sqrt{.048(.015)}<v$. Keeping the two occurrences of $X$
together, its lower and upper polynomials are
\begin{align*}
 I_-={}&(2.3899/12)(1-X/4)\\
 &\quad-.575[1/3.3899-(1-X/4)(X/12-X^2/240)]-v,\\
 I_+={}&(2.39/12)(1-3X/16+.020X^2)+.20X/112\\
       &-.56[1/3.39-X/12]+v+.84(3.546)(.048)(.015)X .
\end{align*}
The derivative of the first is negative and that of the second positive
on $[0,.275]$. At $.275$ they give the stated bounds.
In inequality \eqref{op-torque} for $\nu$, the subtracted coefficient
is at most
$(41/144)(.0433)+.65\sqrt{.015/.9659}/40+
H\mathcal T_H\sqrt{.015}/(24\sqrt{12})<.015$.
Since $\cos(h/2)(1+\sin^2(h/2)/6)\le1$,
this yields $\nu>\alpha_D-.010h$ on the same interval.

\paragraph{Bounds for the cap density.} Write $s=\sinc h$, let $r_-,r_+$ be the two bounds for $d_R$ above and set $\lambda_-=3.3899$, $\lambda_+=3.39$. The cap mass identity $-L=1/[s^2(1+d_R/(\lambda-1))]-1+
2\lambda_0(M_D-C)/(hs)$, with $0\le C\le3hX^3/768$, gives
\begin{equation}\label{common-cap-height}
 \ell_-={1\over s^2(1+r_+/(\lambda_--1))}-1-{6\Lambda X^3\over768s}
 < -L <
 {1\over s^2(1+r_-/(\lambda_+-1))}-1+{2\Lambda mZX^2\over s}=\ell_+ .
\end{equation}
At ten and eleven sides this gives, respectively,
$.3498X<-L<X/3+.062X^2$ and $.3472X<-L<X/3+.0595X^2$.
For every $n\ge12$ Proposition~\ref{IV:prop:unconditional-cap-twelve} gives
$.32X<-L<X/3+.077X^2$.

\section{Geometric bounds for a competitive polygon}
Fix a balanced conformal map $f$ and write $O=f(0)$ for its conformal center.
We center the fan at an interior point $Q$ where the projected trial field has mean zero. We will bound its distance from the conformal center.

\begin{lemma}[Conformal radius on a polygonal fan]\label{IV:lem:skew-fan}\label{IV:lem:gamma-full-skew}
Let $P$ be a convex polygon of area $\pi$, $O\in P$, and $r(P,O)$ the
conformal radius. Let $p_i,b_i$ be the fractions of the full angle and area in the $i$th triangle, and $\delta_i$ half the difference of its base angles. Then,
\[
 -\log r(P,O)^2\ge \sum_i p_i\log(p_i/b_i)+4\zeta(3)\sum_i p_i^4
   +\frac{28\zeta(3)}{\pi^2}\sum_i p_i^2\delta_i^2 .
\]
For $0<p<1/2$ and
$|t|<1/2-p$ put $c=1/2-p$, $d=1/2+p$, $x=t^2$, and write
\[
 W(p,t)=\sum_{j=0}^{\infty}\log
 \frac{1-x/(j+d)^2}{1-x/(j+c)^2}.
\]
Write $Q_{\rm iso}(p)$ for the isosceles factor computed below.
Then
\[
 -\log r(P,O)^2\ \ge\
 \sum_i p_i\left[\log(p_i/b_i)-\log Q_{\rm iso}(p_i)
             +W(p_i,\delta_i/\pi)\right].
\]
The excess inside each triangle over its isosceles value equals $W$. In particular it diverges as $|t|\uparrow c$ at
fixed $p$. It has the elementary lower bounds
\begin{align*}
 W(p,t)&\ge \log\frac{1-x/d^2}{1-x/c^2}+\frac{41}{25}px,\\
 W(p,t)&\ge
 px\left[\frac4{2(1/4-p^2)^2-(1/2+2p^2)x}+\frac{41}{25}\right].
\end{align*}
It also admits the uniform series estimate
\[
 W(p,t)\ge4p\sum_{k\ge1}
 (2^{2k+1}-1)\zeta(2k+1)t^{2k}.
\]
\end{lemma}
\begin{proof}
Write the triangle angles as $\pi\alpha=2\pi p,\pi\beta,\pi\gamma$, so
$\delta=\pi(\beta-\gamma)/2$. The Schwarz--Christoffel map
$F(z)=C\int_0^z t^{\alpha-1}(1-t)^{\beta-1}dt$ from the upper half-plane,
with $C>0$ and prevertices $0,1,\infty$, gives apex sides $CB(\alpha,\beta),CB(\alpha,\gamma)$.
Composing with $z=4\xi/(1+\xi)^2$, $\xi=y^{1/\alpha}$, maps the unit sector
of angle $\pi\alpha$ to the triangle, its circular arc to the base, and has
apex expansion $\rho y(1+o(1))$, $\rho=4^\alpha C/\alpha$.

Put $L=\log(1/\varepsilon)$. The condenser from $B_\varepsilon(O)$ to
$\partial P$ has capacity $2\pi/(L+\log r(P,O)+o(1))$ by conformal invariance.
Cutting along the rays with free traces lowers this capacity to the sum of the
sector capacities $2\pi p_i/(L+\log\rho_i+o(1))$: the inner arcs have radii
$(\varepsilon/\rho_i)(1+o(1))$, and the logarithmic minimizer on concentric
arcs has zero normal derivative on the rays. Comparing the $L^{-2}$ terms gives
$\log r(P,O)\le\sum_i p_i\log\rho_i$.
If the triangle has area $\pi b$, its two side lengths give
\[
 Q:=\frac{\rho^2p}{b}
 =\frac{\pi16^\alpha}
        {\alpha\sin(\pi\alpha)B(\alpha,\beta)B(\alpha,\gamma)} .
\]
For $h(u)=\log\Gamma(u)-\log\Gamma(1-u)$ the skew-dependent part of
$-\log Q$ is $h(c+t)+h(c-t)$, since $\alpha+\beta=1-\gamma$.
Euler's gamma product gives the central second difference of $\log\Gamma(u)$
as $-\sum_{j\ge0}\log(1-t^2/(u+j)^2)$; subtraction at $c,d$ gives $W$.
Write $z=j+1/2$. For $j\ge1$ the first powers in the logarithms give at least
$4pxz/(z^2-p^2)^2\ge4px/z^3$, and
$\sum_{j=1}^4(j+1/2)^{-3}+\int_5^\infty(u+1/2)^{-3}du>41/100$.
Retaining $j=0$ thus gives the first elementary bound; the second follows
from $\log R\ge2(R-1)/(R+1)$ for $R\ge1$ by differentiation.
For the full series, integration of the convex function $v^{-2k-1}$
on $[z-p,z+p]$ gives $(z-p)^{-2k}-(z+p)^{-2k}\ge4kp/z^{2k+1}$.
Summing over $j$ in each positive logarithmic coefficient gives the odd zeta sums.
At $\alpha=2p$ and $\beta=\gamma=1/2-p$, the gamma duplication
formula simplifies $Q_{\rm iso}$ to
\[
 Q_{\rm iso}=
 \frac{2\pi p}{\sin(2\pi p)}
 \left(\frac{\Gamma(1-p)}
 {\Gamma(1+p)\Gamma(1-2p)}\right)^2 .
\]
Expanding the gamma functions on $0<p<1/2$ gives
\[
 -\log Q_{\rm iso}=2\sum_{k\ge3,\ k\ {\rm odd}}
          (2^k-2)\zeta(k)p^k/k
 \ge4\zeta(3)p^3+12\zeta(5)p^5 .
\]
Using $-2\log\rho=\log(p/b)-\log Q$ in the sector radius comparison gives the full inequality. The $k=1$ term in the series bound for $W$ is $28p\zeta(3)\delta^2/\pi^2$, giving the stated quadratic estimate.
\end{proof}
\begin{lemma}[Moments of the radial disk profile]\label{IV:lem:uniform-harmonic-correlations}
For the first disk profile and $M_j=2j\int_0^1g^2r^{2j-1}\,dr$,
$|g''|<2.42$ on $[0,1]$, $M_2=2(\lambda+1)/(3\lambda)$ and
\[
 1-M_{i+1}\le\frac{4.78+.78/(2i+5)}{(2i+3)(2i+4)},\qquad i\ge1 .
\]
\end{lemma}
\begin{proof}
The disk equation gives $(r^3(g/r)')'=-\lambda r^2g$ and
$g'''=(3/r^2-\lambda)(g'-g/r)$. Hence $g''$ decreases to $r_0=\sqrt{3/\lambda}>15/16$
and then increases to $1-\lambda$. On $[r_0,1]$ one has
$g'''\le(\lambda-3/r^2)/r\le\lambda-3$, so $|g''|\le\lambda-1+(1-r_0)(\lambda-3)<2.42$.
In \eqref{IV:eq:regular-beta-band} increasing $\lambda$ replaces fourth powers
by lower powers. At $\lambda=3.39$ its upper bound for $1-M_{i+1}$ is
$A_i=(239i^2+2529i+7090)/[200(i+2)(i+3)(i+4)(i+5)]$;
subtracting from the stated bound gives
$3(713i^2+7327i+13910)/[200(i+2)_4(2i+3)(2i+5)]>0$.
Also $M_2=2(\lambda+1)/(3\lambda)$. Indeed, multiplying the radial
Bessel equation $g''+g'/r+(\lambda-r^{-2})g=0$ by $r^3g$ and $r^4g'$
and integrating by parts gives respectively
$\int r^3(g')^2=-1+\int rg^2+\lambda\int r^3g^2$
and $\int r^3(g')^2=(\lambda-1)/2+\int rg^2-2\lambda\int r^3g^2$,
with all integrals over $(0,1)$.
\end{proof}

\begin{proposition}[Fourier coefficient bounds for a balanced conformal map]\label{IV:prop:harmonic-budget-65}
Let $P$ be a simply connected Lipschitz Jordan domain of area $\pi$,
and $0<\mu_*\le\mu_1(P)$. For the disk profile $g$ and
$\lambda=(j'_{1,1})^2$ write $E=\lambda-1$,
$M_j=2j\int_0^1 g^2r^{2j-1}dr$, $M=M_1=E/\lambda$,
$b=E/\mu_*=M+H$ and $d_i=M_{i+1}-b$ for $i\ge1$.
Assume $d_1>0$. For a balanced map $f$, let $\gamma_i$ be the coefficient of $z^i$ in $f'/f'(0)$ and put $w_i=|\gamma_i|^2/(i+1)$. Then $\sum_{i\ge1}d_iw_i\le H$. Testing also with the centered harmonics of frequencies $1\pm k$ gives the simultaneous bounds on their mass pairings below.
\end{proposition}
\begin{proof}
Expand $Q=f'/f'(0)=\sum_{i\ge0}\gamma_i z^i$, where $\gamma_0=1$,
let $w_i=|\gamma_i|^2/(i+1)$, $X_i^2=d_iw_i$, and put
$W=\sum_{i\ge1}X_i^2$. Area gives $T=\sum_{i\ge1}w_i$.
Choose the balanced conformal map by the argument before Lemma~\ref{IV:lem:signed-circle}. Then $U=g(r)e^{i\theta}$ has zero
mean for the pulled-back area. Its Rayleigh bound gives
$W\le H$, hence $T\le H/d_1$. Set $s=\mu_* /(1+T)$. Centering $U$ and applying the Rayleigh bound to every real projection
of $U$ give, by angular integration,
\[
 \sum_{i\ge0}B_i\gamma_i\overline{\gamma_{i+1}}=0,\qquad
 \left|\sum_{i\ge0}H_i\gamma_i\overline{\gamma_{i+2}}\right|\le H-W,
\]
where $B_i=2\int_0^1gr^{2i+2}\,dr$,
$H_i=2\int_0^1g^2r^{2i+3}\,dr$. 
The moments increase with $j$ because $g$ increases.
Here $f'(0)>0$ and $T=f'(0)^{-2}-1$; area measure on the disk is $dA/\pi$.
Fix $0<s_0\le s$.
For $3\le k\le K$ include both centered harmonics of frequencies
$1\pm k$, writing $m=|1\pm k|$, and
$v_i^{(k,m)}=2\int_0^1g r^{2i+k+m+1}\,dr$. Their signed frequencies are
distinct. Their Dirichlet energy matrix is diagonal, with entries $2m$.
Index each signed harmonic by $\ell=1\pm k$ and put $m=|\ell|$. Let $C_\ell$ denote its mass pairing with $U$ using the weight $|Q|^2$. Subtracting its weighted mean from the harmonic
preserves this coupling. The Rayleigh inequality for $U$ plus a linear
combination of the centered harmonics implies
$E-sM_U+\sum_\ell 2|\ell|\,|a_\ell|^2-2s\operatorname{Re}\sum_\ell\bar a_\ell C_\ell\ge0$,
where $M_U$ is the weighted mass of $U$: the omitted quadratic
auxiliary mass is nonnegative. Minimize in the independent coefficients $a_\ell$,
use $(E/s)-M_U=H-W$, and compute each $C_\ell$ by angular integration, to obtain
\[
 H-W\ge s_0\sum_{k=3}^{K}\sum_{\ell=1\pm k}\frac1{2|\ell|}
       \left|\sum_{i\ge0}v_i^{(k,|\ell|)}\gamma_{i+k}\overline{\gamma_i}\right|^2 .
\]
Put $X=\sum_{i\ge1}|\gamma_i|^2/(i+2)$.
For either sign write $v_i=v_i^{(k,m)}$. Since $v_i\le1/(i+k+1/2)$, the squared tail ($i\ge1$) of either
mass pairing is at most $X\sum_{j>k}(j-k+2)w_j/j$, by Cauchy--Schwarz, whose second factor has coefficient
$(j-k+2)(j+1)/(j+1/2)^2\le(j-k+2)/j$.
For each $k$, choose $0<t_k<1$ and separate the leading term with
$|p+z|^2\ge t_k|p|^2-t_k|z|^2/(1-t_k)$.
Since $g\ge r$, the leading terms contribute at least the following coefficient of $w_k$:
\[
 \frac{t_ks_0}{2}\big[(k^2-1)^{-1}+(k+2)^{-2}\big].
\]
For each $w_j$, $j>k$, the two remainders together reduce its coefficient by at most
$s_0 X[k/(k^2-1)][t_k/(1-t_k)](j-k+2)/j$. The remainders involve no lower indices.
\end{proof}
\subsection{The conformal center and the area centroid}
\label{IV:sec:center-transfer}
\begin{lemma}[Distance from the conformal center to the centroid]\label{IV:lem:physical}
Let $f$ be the balanced conformal map of the disk onto a polygon of area $\pi$
with centroid $G$. Put $O=f(0)$, $a=f'(0)>0$, $v=f-O-az$ and $T=a^{-2}-1$.
Then
\begin{equation}\label{paper-centroid-shift}
 |G-O|\le\frac{T}{1+T}\left(\frac5{16}+\|v\|_\infty\right).
\end{equation}
Hence
$|G-O|\le(21/16+R_0)T$ if $|f-O|\le R_0$.
\end{lemma}
\begin{proof}
Write $f-O=\sum_{k\ge1}a_k z^k$, $d\sigma=dA/\pi$ and
$B_j=2\int_0^1g(r)r^{2j+2}\,dr$, so $B_0=2/\lambda$.
In the centroid integral $\int(az+v)|a+v'|^2d\sigma$, subtract from
$a\int z|f'|^2d\sigma$ the zero sum
$(a/2B_0)\sum_{k\ge1}k(k+1)B_{k-1}a_k\bar a_{k+1}$ from the zero-mean condition for the disk field.
Orthogonality then gives
\[
 G-O=a\sum_{k\ge2}d_k a_k\bar a_{k+1}+\int v|v'|^2d\sigma,\quad
 d_k=1+k[1-(k+1)B_{k-1}/(2B_0)] .
\]
The averages $(k+1)B_{k-1}$ of $g/r$ decrease from $2B_0$, since
$g/r$ decreases; $g/r\ge1+(1-r^2)/2$ gives $1\le d_k\le
s_k=1+k[1-(\lambda/4)(1+1/(2k+4))]$.
Put $x_k=\sqrt k\,|a_k|$, $b_k=s_k/\sqrt{k(k+1)}$.
Both $s_k/k=1-\lambda/4+1/k-\lambda/[8(k+2)]$ and
$s_k/(k+1)=1-\lambda/4+\lambda(k+4)/[8(k+1)(k+2)]$ decrease;
hence so does $b_k$. Substitution of $\lambda>61/18$ gives
$b_2^2<1/5$, $b_3^2<1/8$, $b_4^2<9/100$ and $b_5^2<75/1024$.
In any finite chain, completing the square for the mixed term $bxy$ with diagonal coefficient $d>0$ reduces the preceding diagonal coefficient by
$b^2/(4d)$. Starting with $5/16$, the coefficients for $k\ge5$
remain at least $15/64$ since $5/16-(75/1024)/(4\cdot15/64)=15/64$.
Backward substitution at $k=4,3,2$ leaves at least $3/14,1/6,1/80$.
Pass to the limit and use $a\le1$ and $\sum_{k\ge2}x_k^2=T/(1+T)$.
Finally $\|v\|_\infty\le R_0+1$. Scaling to circumradius-one
reference area multiplies the displacement by $\sqrt{\operatorname{sinc}h}\le1$.
\end{proof}

For a convex fan of area $\pi$ put $T=r(P,O)^{-2}-1$, $q=np$,
$y=(n-q)/(n-1)$ and $D_j(q)=q^j+(n-1)y^j-n$. Set
$b(q)=4.8nD_4(q)+(12.4/n)D_6(q)$ and
$A_n(q)=b(q)+(36/n^3)[D_8(q)+n]$ and $\eta_i=b_i/p_i-1$. For each choice of sector, its aperture and the two sums below satisfy
\begin{equation}\label{IV:eq:single-aperture-budget}
 A_n(q)+n^5\!\left\{\sum_i p_i[\eta_i-\log(1+\eta_i)]
              +\sum_i p_iW(p_i,\delta_i/\pi)\right\}
 \le n^5\log(1+T)-4.808n^2-12.4 .
\end{equation}
Indeed $\sum p_i\eta_i=0$, so the first sum equals $\sum_ip_i\log(p_i/b_i)$, the relative entropy of the angle fractions with respect to the area fractions. Keep the quartic, sextic and octic powers
in Lemma~\ref{IV:lem:skew-fan} and apply Jensen on the other $n-1$
apertures; their coefficients exceed $4.808,12.4,36$, respectively,
and $D_j\ge0$. Each $D_j$ is strictly convex with minimum at one.
With $B=4\zeta(3)$ the same argument retains the exact quartic
baseline $Bn^2$ on the right. The skew term is at least
$3.4\sum_i p_i^2\delta_i^2$ before multiplication by $n^5$.

The sine rule and triangle area give the relative endpoint difference
$|z|=|\tan(\pi p)\tan\delta|$ and mean radius
$R_O^2=(1+\eta)\operatorname{sinc}h/[(1-z^2)\operatorname{sinc}(hq)]$
in reference-area units. For two interior centers $O,G$, put $a=a_O$, $a'=a_G$, $z=z_O$, $z'=z_G$.
Let $d_O,d_G$ be their distances to the same side line. Then
\begin{equation}\label{IV:eq:chord-center-lipschitz}
 \sinh^2\!\left(\tfrac12\log\frac{\tan a'}{\tan a}\right)
 +\frac{\sinh^2(\operatorname{atanh}z'-\operatorname{atanh}z)}
       {\sin(2a')\sin(2a)}
 =\frac{|O-G|^2}{4d_Od_G}.
\end{equation}
Indeed $Z=(X-v_+)/(X-v_-)$ maps the supporting half-plane isometrically
to a half-plane, with $\log|Z|=\pm2\operatorname{atanh}z$ and
$\arg Z=\pm2a$ up to an additive constant. Substitute in the
half-plane distance formula, using
$\sin^2(a'-a)/(\sin2a'\sin2a)=\sinh^2(\frac12\log(\tan a'/\tan a))$.
To control the bisector tilt, normalize the side to $[-1,1]$ and set
$f(X)=\sqrt{(X-1)(X+1)}$, taking the average of the two upper half-plane
arguments. The signed distance of $f(X)$ from the imaginary axis is
$\pm\operatorname{arsinh}(\tan\delta_X)$. Schwarz--Pick bounds its change
by the same hyperbolic distance $D=2\operatorname{arsinh}(|O-G|/(2\sqrt{d_Od_G}))$.
Thus the changes in $\log\tan a$, $2\operatorname{atanh}z$ and
$\operatorname{arsinh}(\tan\delta)$ are at most $D\le-\log(1-h)$ for
$h=|O-G|/\operatorname{dist}(O,\partial P)<1$.

\subsection{Fourier coefficients of the balanced conformal map}
\label{IV:sec:pure-allocation}
\begin{lemma}[The first two harmonic coefficients]\label{IV:lem:two-harmonic-secants}
In the balanced conformal normalization of
Proposition~\ref{IV:prop:harmonic-budget-65}, put $w_i=|\gamma_i|^2/(i+1)$
and $d_i=1/\lambda-H-(1-M_{i+1})$. If $0\le H\le1/100$, then
\[
 |\gamma_1|\le H(6.812+37H),\qquad
 |\gamma_2|\le H(9.13+37.5H).
\]
For $0<D\le1/\lambda$, the terms with $i=1,2$ contribute
at most $(19.03+380H)H^2$ to $\sum_i(1-d_i/D)_+w_i$.
\end{lemma}
\begin{proof}
For $0\le H\le.01$ put $L=3.94+16H$. In fact
\[
 |\gamma_1|\le (\lambda/2)LH,\qquad
 |\gamma_2|\le [3\lambda/(\lambda+1)]LH .
\]
By \eqref{IV:eq:regular-beta-band}, for $i=1,2,3$ and $i\ge4$, respectively, this gives
\[
 Ld_i\ \ge\ 3/5,\quad 5/6,\quad 19/20,\quad 41/40 .
\]
Here one may check just $M_2=2(\lambda+1)/(3\lambda)$ and $j=3,4,5$,
as $M_j$ increases. The beta estimate gives
$d_i+H>.1583,.2169,.2448,.2602$, respectively, using
$3.3899<\lambda<3.39$. Multiply each minus $H$ by $L$ and take the
two endpoints by concavity.

Let $u_i=|\gamma_i|/\sqrt{i+1}$. The two tails in
Proposition~\ref{IV:prop:harmonic-budget-65} are bounded by
$\sum_{i\ge1} b_i u_i u_{i+k}$ for $k=1,2$ with
$b_i=B_i\sqrt{(i+1)(i+2)}$ or $H_i\sqrt{(i+1)(i+3)}$.
All $b_i\le1$ by $g\le1$ and the arithmetic--geometric mean inequality.
Moreover $b_1^2\le8/9$: at step one Cauchy--Schwarz gives
$6B_1^2\le M_2$, and at step two $H_1\le1/3$.
Apply the same backward square completion as in Lemma~\ref{IV:lem:physical} on each progression of step $k$. The constant tail $41/40$
leaves at least $5/8$ since $5/8+1/(4\cdot5/8)=41/40$.
The third and second sites thus leave at least
$19/20-2/5=11/20$ and $5/6-5/11=25/66>3/8$.
At the first site $i=1$ the subtraction is at most
$(8/9)/(4\cdot3/8)=16/27<3/5$. In the progression starting at $i=2$, the coefficient of $u_2^2$ before completing squares is at least $5/6>1/(4\cdot3/8)$. Truncation and passage
to the limit now bound both tails by $L\sum d_i u_i^2=LW$.
The leading coefficients are $B_0=2/\lambda$ and
$H_0=M_2/2=(1+1/\lambda)/3$. The additional step-two error $H-W$ is
absorbed since $L>1$ and $W\le H$. This proves the stronger bounds
above, and hence the two in the statement.
Finally use $D\le1/\lambda<.295$ and put $D_1=.1583-H$, $D_2=.2169-H$.
With $A=6.812+37H$, $B=9.13+37.5H$,
the first two terms divided by $H^2$ are at most
$(1-D_1/.295)A^2/2+(1-D_2/.295)B^2/3$.
This polynomial has nonnegative coefficients and its endpoint values
are below $19.03,22.83$. Its secant proves the claim.
\end{proof}

For $H\le1/100$ set $a=1/\lambda-H$. The weighted coefficient
sum $X$ in Proposition~\ref{IV:prop:harmonic-budget-65} satisfies
$X\le H/a+7H^2$.
Indeed $d_i\ge a(i+1)/(i+2)$ for $i\ge3$ by
Lemma~\ref{IV:lem:uniform-harmonic-correlations}, since $4.92/[2(2i+3)]<a$.
For the first two indices the preceding lemma bounds the discrepancy
divided by $H^2$ by
$(2/3-.1483/.295)7.182^2/2+(3/4-.1971/.295)9.505^2/3<7$, using $d_1>.1483$, $d_2>.1971$. Sum and use $\sum d_iw_i\le H$.

\paragraph{A second-moment identity for the radial function.}
Let $|P|=\pi a$, $a\le1$, and write $\rho_X$ for its radial graph
about an interior point $X$ of the convex polygon. Use normalized
angular measure and set $D_X=\|1-\rho_X^2\|_2$. Polar integration
and the parallel-axis identity give, with centroid $G$,
\[
 D_X^2=1-2a+\frac2\pi\int_P |x-X|^2\,dx=D_G^2+2a|X-G|^2.
\]
Let $\theta_i$ be local aperture, $w_i$ the endpoint radius difference,
$p_i=\theta_i/(2\pi)$ and $b_i=p_i(1+y_i)$ the area fraction at $X$.
With $C_i=3\sin\theta_i/(2+\cos\theta_i)$ one has
\begin{equation}\label{IV:eq:polar-fan-rms}
 D_X^2=(1-a)^2+a^2\sum p_i y_i^2+
 a^2\sum p_i(1+y_i)^2(\theta_i/C_i-1)
 +\tfrac a3\sum p_i(1+y_i)w_i^2 .
\end{equation}
Indeed a triangle of area $A_i$ and vertices $0,v,v'$ has polar
moment $A_i(|v|^2+v\cdot v'+|v'|^2)/6=A_i w_i^2/6+A_i^2/C_i$.
Sum over the fan and use $\sum p_i y_i=0$. The angular factor is
at most $\theta_i^4/(180\operatorname{sinc}\theta_i)$: the third
derivative of $t(2+\cos t)-3\sin t$ is $t\sin t\le t^2$ and the
first three initial values vanish.

Suppose $.74<\rho_G<1.175$, $D_G<.0905$ and $|\beta|a<.035$.
There is a fan center $Q$ balancing $U(x-Q)+\beta(x-Q)$ with
\begin{equation}\label{IV:eq:moment-center-ball}
                         |Q-G|\le 2 D_G^2/5 .
\end{equation}
To see this, on $.73<r<1$ one has $r\le g\le1$, $0\le g'\le g/r$.
For $1<r<1.18$ one has $0<g\le1$ and $|g'|\le3.1(r-1)$ by
integrating $-(rg')'=(\lambda r-r^{-1})g$ from one.
Thus $B(r)=\int_1^r(g(t)-t)t\,dt$ satisfies
$|B(r)|\le(1-r)^2/2$ inside and
$|B(1+s)|\le s^2(1+1.7s+.775s^2)/2$ outside.
In particular $|B(\rho_G)|\le (1-\rho_G^2)^2/6$.
The centroid condition then gives
$|J(G)|\le2\pi\langle |B(\rho_G)|\rangle\le\pi D_G^2/3$,
where $J(Q)=\int_P U(x-Q)\,dx$.
On the indicated ball $.73<\rho_Q<1.18$ and the normalized area
of $(P-Q)\triangle\mathbb D$ is at most
$D_Q\le (D_G^2+2a(2D_G^2/5)^2)^{1/2}<.091$.
There $\|DU\|\le7/5$. Since $\int_{\mathbb D} DU=\pi I$,
the negative derivative of $J(Q)+\beta |P|(G-Q)$ is at least
$\pi[1-(7/5)(.091)-.035]I>(5\pi/6)I$.
The vector field points inward on the boundary, and degree gives
\eqref{IV:eq:moment-center-ball}. If $D_G=0$, the centered disk already
has $J(G)=0$.

\subsection{Centering the projected trial field}\label{mid-geometry}

First take $n=10,11$, $h=2\pi/n$, $t=h/2$, $a=\operatorname{sinc}h$.
For a strictly convex polygon the only spectral input here is
$(1-1/\lambda)(\pi\lambda/F(P)-1)\le H$, with $H$ specified below.
Thus the lower frequency used in the harmonic tests may be taken
to be $\lambda/(1+\lambda H/(\lambda-1))$.
Write $k$ for the normalized aperture at the conformal center $O$ and $q$ for the corresponding aperture at a center $Q$ where the projected trial field has mean zero.
Use $A_n$ and $\eta$ from \eqref{IV:eq:single-aperture-budget}. The first four columns after $n$ give an upper bound $C$ for the left side of the conformal inequality, the aperture range $l<k<u$, and a bound $|\eta|<Y$ on the relative area error. The remaining bounds hold for the same polygon with the fan centered at $Q$, with $|w|/h^2<W$ everywhere and $|w|/(qh^2)<W_-$ for $q\le1$.
\begin{equation}\label{paper-geometric-window}
\begin{array}{c|rccc|cc|cc|cc}
n&C&l&u&Y&q>&q<&e>&e<&W&W_-\\\hline
10&83.1&.408&1.430&.15&.426&1.413&-.063&.064&.227&.330\\
11&80&.44&1.411&.12&.463&1.391&-.049&.049&.209&.281\\
\end{array}
\end{equation}
Near the largest aperture the following upper bounds on $e$ will also be used. Each pair denotes a lower threshold on $q$ and its corresponding upper bound on $e$:
\[
\begin{array}{c|llll}
10&(1.35,.053)&(1.37,.049)&(1.40,.043)&(1.41,.036)\\
11&(1.35,.038)&(1.37,.033)&(1.39,.024)\\
\end{array}
\]
For the decagon, the same ellipse bounds $|w|/h^2$ by the second row in each half of the following display:
\[
\begin{array}{c|ccc@{\qquad}c|cccc}
 q\le&.48&.52&.60 & q\ge&1.25&1.30&1.35&1.40\\
 |w|/h^2<&.087&.113&.150 & |w|/h^2<&.192&.171&.137&.073
\end{array}
\]

\paragraph{An upper bound from the conformal map.}
Set $H=(4.809+(59.93+197/n)/n^2)/(3.3899n^3)$, $s=3.35$, $\gamma=1-50/n^3$, $a_0=1/3.39-H$, and
\[
 P(D)=H/D+
 [1-(.1583-H)/D]H^2(6.812+37H)^2/2+
 [1-(.2169-H)/D]H^2(9.13+37.5H)^2/3 .
\]
Write $w_i=|\gamma_i|^2/(i+1)$. If the coefficient of $w_i$ is at least $D>.23$ for each $i\ge3$, the bounds on $w_1,w_2$ in Lemma~\ref{IV:lem:two-harmonic-secants} give $T\le P(D)$. Section~\ref{V:sec:geo-seven-nine} explains why we may use the indicated upper bound on $H$. Before adding the harmonic tests, the coefficient is $M_{i+1}-M-H\ge M_4-M-H\ge .2448-H$. It follows that the coefficient $\mu_* /(1+T)$ multiplying the mass in the Rayleigh inequality is at least
$3.3899/[(1+3.39H/2.3899)(1+P(.2448-H))]>s$. In Section~\ref{V:sec:geo-bubble} take
\[
 X=H/a_0+7H^2,\quad J=\frac{s\gamma X}{12(5s/12-\gamma)},\quad
 L=225J-4.8(1-\gamma)>0,\quad D=a_0-J-64J^2/L .
\]
The tail loss there is at most $J(1+16x)-Lx^2\le J+64J^2/L$. We may therefore use $T_*=P(D)$ in both rows. The right side of the conformal inequality is bounded by
$n^5(T_*-T_*^2/2+T_*^3/3)-4.808n^2-12.4<C$.

\paragraph{A quadratic bound for the endpoint radii.}
Write $g(k)=-\log\operatorname{sinc}(hk)$, $b(k)=g(k)-g(1)$, $A=4.8n(3+2l+l^2)$.
The strict convexity and endpoint tests $A_n(l),A_n(u)>C$ give $l<k<u$.
The quartic term is at least $A(k-1)^2$. Subtract this lower bound from the total upper bound. The remainder divided by $k$ is at most $2C/(1+\sqrt{1-C/A})$, and divided by $k^2$ is at most $C/(1-C/A)$ by square completion.
The first bound is below $n^4[Y-\log(1+Y)]$, giving $|\eta|<Y$.
The skew term of \eqref{IV:eq:single-aperture-budget} gives
$k|\delta|\le c=\sqrt{C/(3.4n^3)}$ and $|\delta|\le d=c/\sqrt{1-C/A}$, where $\delta$ is the angle from the side normal to its bisector. With $\tau(x)=\tan(x)/x$,
\[
 |z|=|\tan(tk)\tan\delta|\le Z=t\tau(tu)\tau(d)c.
\]
Here $tu<.45$, $d<.36$, $Z<.055$. Corollary~\ref{IV:cor:aperture-skew-product} therefore gives the lower bound $K(\log r_i-\log r_{i+1})^2$ for the skew term on each edge, with $K=8n^5/\pi^4$.

For a path of $n-1$ edges denote its endpoint coefficients from Lemma~\ref{IV:lem:middle-conformal-entry} by $M_\pm(B,K)$; explicitly they are $\sqrt{BK}(1\mp r^{n-1})/(1\pm r^{n-1})$ with $r=(\sqrt K-\sqrt B)/(\sqrt K+\sqrt B)$. For $v=\log(1+\eta)$ the summand $\eta-\log(1+\eta)=e^v-1-v$ of the relative entropy satisfies $e^v-1-v\ge[1+\log(1-Y)/3]v^2/2$, by the cubic minorant on $v<0$ and the quadratic minorant on $v\ge0$. Set
\[
 E=\tfrac12 n^4[1+\log(1-Y)/3],\qquad
 \alpha^{-1}=(El)^{-1}+g'(u)^2/A .
\]
Cauchy--Schwarz and the preceding path lemma give an ellipse for the half-sum and half-difference of the adjacent log radii with squared semiaxes
$p^2=C/[4(\alpha+M_+(\alpha,K))]$ and $s_0^2=C/[4(K+M_-(\alpha,K))]$.
Thus at $O$ we have $|e_O|\le E_O=e^p-1$ and $|w_O|\le W_O=2e^{p^2/2}\sinh s_0$. Indeed $s_0^2<p$, and $e^{p x}\cosh(s_0\sqrt{1-x^2})\le e^p$ on $[-1,1]$.

For the skew bound above use $\sinh(s_0\sqrt{1-x^2})\le\sqrt{1-x^2}\sinh s_0$ and $e^{px}\sqrt{1-x^2}\le e^{p^2/2}$, by $\log(1-x^2)\le-x^2$.

Fix the aperture $k$ of one sector to sharpen the bound on its endpoint radii. Put $y=(n-k)/(n-1)$, $y_0=(n-u)/(n-1)$,
\[
 B=4.8n(l^2+2ly_0+3y_0^2),\quad m^{-1}=(El)^{-1}+g'(u)^2/B,\quad
 M=M_+(m,K),\quad N=M_-(m,K).
\]
Subtract the lower bound $A_n(k)$ from Jensen's inequality. On the other edges the quartic remainder about $y$ is at least $B\sum_{j\ne i}(k_j-y)^2$ by expansion of the fourth power. The relative entropy term is at least $El\sum_{j\ne i}\log^2(1+\eta_j)$, while $|g(k_j)-g(y)|\le g'(u)|k_j-y|$. By Cauchy--Schwarz, after shifting the log radii by $b(y)/2$, the coefficient of the squared sum of adjacent log radii along this path is at least $m$. The coefficient of their squared difference is $K$. Write $P_0,S_0$ for the half-sum and half-difference of the selected endpoint log radii. Adding the relative entropy and skew terms for the selected edge gives
\begin{equation}\label{paper-joint-ellipse}
 A_n(k)+Ek(2P_0-b(k))^2+M(2P_0-b(y))^2+4(K+N)S_0^2\le C .
\end{equation}
Alternatively retain $3.4n^3 k^2\delta^2$ for the selected skew. Since $2|S_0|=2|\operatorname{atanh}z|\ge2tk|\delta|$, the complementary difference term yields
$|\delta|\le\delta_*=[C/((3.4n^3+4t^2N)(1-C/A))]^{1/2}$.

\paragraph{Choosing a center where the trial field has mean zero.}
Take $c_1^*=H(6.812+37H)/2$ and $c_2^*=H(9.13+37.5H)+2(c_1^*)^2$ in Lemma~\ref{IV:lem:convex-boundary-two-heads}. Let $v_*$ denote its right side with these absolute coefficient bounds and $T_*$. Lemma~\ref{IV:lem:physical} gives $|G-O|\le g_*=T_*(5/16+v_*)$ at the area centroid $G$. In Corollary~\ref{IV:cor:rms-fan-center} insert $C,Y,u$ to get
\[
 N_0^2=(1-a)^2/h^4+\frac{a^2(1+Y)^2}{180\operatorname{sinc}(hu)}
       +\frac{2a^2(1+Y)C}{(2\pi)^4n}.
\]
Its two coefficient ratios are bounded by $2(1+E_O)^2[\pi\tau(tu)\tau(d)/\sqrt{3.4}]^2u/(3an)$ and $(1+Y)(2\pi)^4/[864n\operatorname{sinc}(hu)]$, each below one. Thus $D_G^2\le h^4N_0^2$. The support at $O$ is at least
$R=(1-E_O)(1-Z^2)\cos(tu)\cos\delta_*$. The outer radius is at most $1+E_O+W_O/2$. After allowing $g_*$ these imply $.74<\rho_G<1.175$, and $h^2N_0<.0905$. The same endpoints give $H<.01$. Here $Z_s=\pi\tau(tu)\tau(d)/\sqrt{3.4}<2$ bounds the sine-rule factor used in the moment estimate. The quantity $Z$ above bounds the skew $|z|$ itself. The integral mean of $T_iD$ is $\beta$ times the triangle centroid by reflection and affine transport. Its magnitude satisfies $|\beta|a\le3\sqrt{aM_D/(2h)}/\cos t<.035$ using $M_D\le m_*h^5 z$ from Section~\ref{mid-response}.
Formula~\eqref{IV:eq:moment-center-ball} gives a center $Q$ where the projected trial field has mean zero, with
$|Q-O|\le d_Q=g_*+2h^4N_0^2/5$. Along the segment joining the centers, the distance to each supporting line is at least $R-d_Q$. Put $f=\exp(d_Q/(R-d_Q))$. By \eqref{IV:eq:chord-center-lipschitz}, $k/f\le q\le kf$.
The changes in mean radius and skew are bounded by $d_Q$ and $d_Qhkf$, respectively, since the difference of the two radial unit vectors has length $2\sin(tq)$. The bisector rotates by at most $\log f$. Hence for $q\le1$, $|w|/(qh^2)\le(1+\sup e)\tau(t)\tan(\delta_*+\log f)/h$.

\paragraph{Bounds on the aperture and radii.}
To bound the aperture, put
$H(k)=EkM/(Ek+M)$ and $\theta(k)=Ek/(Ek+M)$. Completing the square in \eqref{paper-joint-ellipse} leaves
$C-A_n(k)-H(k)[b(k)-b(y)]^2$ and log center $[\theta(k)b(k)+(1-\theta(k))b(y)]/2$.
To exclude $q\le q_-$ test $k_0=q_-f$ unless $k_0<l$, when the result is immediate. Both $A_n(k)$ and $[b(k)-b(y)]^2$ decrease on $k<1$, so $A_n(k_0)+H(l)[b(k_0)-b(y(k_0))]^2>C$ suffices. To exclude $q\ge q_+$ test $k_0=q_+/f$, using $H(k_0)$ and the increasing right branch, or directly $k_0>u$.

\paragraph{Supporting tangents for ten and eleven sides.}
At ten and eleven sides a proposed endpoint radius can be tested directly
against the joint ellipse. Write $y=(n-k)/(n-1)$, $B_k=Ek+M$,
$\theta=Ek/B_k$, $\Delta=b(k)-b(y)$ and
\[
 \mu=(\theta b(k)+(1-\theta)b(y))/2,\qquad
 R=C-A_n(k)-M\theta\Delta^2 .
\]
On the admissible set $R\ge0$, the log-radius semiaxes are
$p^2=R/(4B_k)$, $s^2=R/[4(K+N)]$. They satisfy $s^2\le p$ since
$C(Eu+M)<4(K+N)^2$. Thus the mean radius at $O$ is between
$e^{\mu-p}$ and $e^{\mu+p}$ by the preceding cosh inequality.

To prove an upper (respectively lower) bound $e_*$ for the mean-radius displacement at the new center, put
$r=\log(1+e_*-d_Q)$ (respectively $\log(1+e_*+d_Q)$). Define
\[
 L(k)=A_n(k)+Ek[2r-b(k)]^2+M[2r-b(y)]^2 .
\]
Whenever $L(k)>C$, the log threshold $r$ is outside the ellipse:
the identity is $L(k)-C=4B_k(r-\mu)^2-R$.
For the upper thresholds used below $2r>\max\{\theta(u)b(u),b(y(l))\}$;
for the lower thresholds $2r<\min\{b(l),b(y(u))\}$. These give
the required orientation since $b(k)$ and $b(y)$ have opposite signs.

On $v\le k\le u$ take the tangent to $L$ at $z$. With
$D=2r-b(z)$, $D'=2r-b(y(z))$ and $g''_u=h^2/\sin^2(hu)-1/u^2$,
a sufficient condition that it supports $L$ is
\[
 4.8n\,{v^2+2vz+3z^2\over v+2z}
 -E D_+ g''_u-{M(D')_+g''_u\over(n-1)^2(v+2z)}\ge0 .
\]
Indeed linearize each square at its value at $z$, retaining the
nonnegative square of its increment. The divided difference $[k,z,z]A_n$
dominates $4.8n(k^2+2kz+3z^2)$ by convexity of the other powers.
The sine product gives $0<g'(k)\le k g''_u$, $0<g''(k)\le g''_u$.
Consequently $[k,z,z](k b(k))\le g''_u(k+2z)/2$ and
$[k,z,z]b(y)\le g''_u/[2(n-1)^2]$, both nonnegative.
Finally $(k^2+2kz+3z^2)/(k+2z)=k+3z^2/(k+2z)$ increases.
The displayed condition has left side greater than $19$ in every
substitution below. The lower bound for $L$ is therefore just the
smaller of the two endpoint values of its tangent. For differentiation,
\[
 L'(z)=A_n'(z)+E[D^2-2zD b'(z)]+2M D'b'(y(z))/(n-1).
\]
For the global lower and upper bounds take $v=l$ and the following $z$.
For each upper bound that applies above an aperture threshold $q_0$, take $v=z=q_0/f$.
The table lists the bounds in order of increasing threshold.
\[
\begin{array}{c|cc|cc|l}
 n&z_{\rm lower}&z_{\rm upper}&\inf L_{\rm lower}>&\inf L_{\rm upper}>&
 \text{lower bounds for }\inf L\text{ by threshold}\\ \hline
 10&.812&1.13 &92&96&88,86,87,84\\
 11&.83&1.117&89&90&85,81,80.6
\end{array}
\]
Every entry exceeds $C=83.1,80$ in its row.

\paragraph{Skew from the same ellipse.}
At fixed $k$ its skew, after the center shift, is at most
$U(k)=2e^{\mu+p^2/2}\sinh s+d_Qhkf$.
For the global estimate put $C_0=C-36/n^2$, $a=4.8(6n+5)$ and
\[
 w_0=2\exp\{C_0/[8(El+M)]\}\sinh\sqrt{C_0/[4(K+N)]},\qquad
 \gamma={\exp(\theta(u)b'(u)(u-1)/2)-1\over u-1}.
\]
For $k\le1$ one has $\mu\le0$ since
$(1-\theta(l))b'(u)/(n-1)<\theta(l)h^2(1+l)/6$ and
$|b(k)|\ge h^2(1-k^2)/6$. Thus $U\le w_0+d_Qhf$.
For $k\ge1$, write $t=k-1$. Then $e^\mu\le1+\gamma t$ by convexity
of the exponential and $R\le C_0-at^2$ by the quartic.
Use that $\sinh s/s$ increases and
$\sqrt{1-at^2/C_0}\le1-at^2/(2C_0)$. Maximizing the resulting
quadratic gives
\[
 U\le w_0+d_Q h f+{C_0(w_0\gamma+d_Qhf)^2\over2w_0a}.
\]
After division by $h^2$ this is less than $.226,.2084$ at ten and
eleven sides, respectively.

For the seven aperture ranges in the decagon skew estimate, $U$ increases on the admissible portion
of $[l,.61]$ and decreases on that of $[1.24,u]$. To check the signs
put $B=Ek+M$ and $\Theta=\mu+R/(8B)$.
On the whole window $|\Delta|<.08$, $b'(u)<.20$, $E/M<1.43$,
$E<4730$ and $B>5200$, so $|\mu'|<.17$ by differentiation.
On the left, $R'\ge-A_n'(.61)-E(.08)^2>150$, using
$0<(M\theta)'\le E$, $\Delta<0$, $\Delta'>0$. Hence
$\Theta'>-.173$ and the logarithmic derivative of
$e^\Theta\sinh s$ is at least $-.173+150/(2C)>0$.
On the right $R'\le-A_n'(1.24)<-205$ and $\Theta'\le.17$.
There the negative derivative of $2e^\Theta\sinh s$ is at least
$.96[205/2-.17C]/\sqrt{(K+N)C}>d_Qhf$, since $e^\mu>.96$ and
$s\coth s\ge1$. These derivative tests are used where $R>0$,
with the endpoints following by continuity. Thus for each upper threshold on a narrow aperture
evaluate $U$ at $q_0 f$, and for each lower threshold on a wide aperture at $q_0/f$.
In the order of the skew table their upper bounds after division by
$h^2$ are $.0860,.1121,.1494,.1916,.1698,.1364,.0724$.

\begin{proposition}[Sector angles and endpoint radii for twelve or more sides]\label{IV:prop:entry-twentyfour}\label{IV:prop:entry-22-23}\label{IV:cor:taper-lower-geometry}
Let $N=12$, $B=4\zeta(3)$ and $U=\max(74,59.93+197/N)$.
For every $n\ge N$, normalize $R_n$ to have circumradius one and $P$ to
have the same area $\pi\operatorname{sinc}(2\pi/n)$. Each strictly convex $n$-gon $P$ satisfying
\[
 (\lambda-1)(\pi\lambda/F(P)-1)\le (B+U/n^2)/n^3
\]
admits a fan center where the projected trial field has mean zero and the following bounds hold.
The last column applies below unit aperture.
\[
\begin{array}{c|cc|ccc}
 N&q>&q<& |e|/h^2<& |w|/h^2<& |w|/(qh^2)<\\ \hline
 12&.475&1.403&.160&.202&.259\\
\end{array}
\]
\end{proposition}
\begin{proof}
Scale out $x=1/n$ in the preceding ellipse and center formulas.
Put $B=4\zeta(3)$, $U=\max(74,59.93+197/N)$, $Q=B+Ux^2$,
$H=x^3\eta$, $\eta=Q/\lambda$ and $s=3.35$.
In the conformal estimate divide $J,L$ by $x^3$, and write
\[
 j={s\gamma\over12(5s/12-\gamma)}
 \left({\eta\over1/\lambda-H}+7\eta^2x^3\right),\quad
 \gamma=1-50x^3,\quad V=\eta+j+64j^2/(225j-240).
\]
For the estimates below the argument of an increasing factor is
$x_N=1/N$. Use $1-50x_N^3\le\gamma\le1$ and
$4.808<B<4.809$. For an upper bound on $V$, insert the upper endpoint of $j$ in the numerator and its lower endpoint in $225j-240$. In particular $225j-240>450$ for $N=12$,
using $\eta\ge4.808/3.39$. Here $H<.000912$. Taking $D=.2448-H$ in $P(D)$ above gives $T<.00379$, since
\[
 P(D)<.000912/.2438+.000912^2\big[(6.812+37\cdot.000912)^2/2+(9.13+37.5\cdot.000912)^2/3\big]<.00379.
\]
Hence $\mu_* /(1+T)>3.3899/[(1+3.39\cdot.000912/2.3899)(1+.00379)]>3.37$, which permits $s=3.35$ in the additional harmonic tests. With these tests the coefficient of $w_i$ for every $i\ge3$ is at least $1/\lambda-x^3V>.290$. The bounds on the Fourier coefficients now give
\[
 T/x^3\le T_3={Q\over1-\lambda x^3V}+(19.03+380H)\eta^2x^3,
\quad
 C_* =U-12.4+x\left[{Q\lambda V\over1-\lambda x^3V}
                         +(19.03+380H)\eta^2\right].
\]
Here $C_*$ bounds the conformal excess by $\log(1+T)\le T$ and
cancellation of $Bx^3$ before division by $x^5$.
The following substitutions give the initial bounds on the aperture at the conformal center, $l<k<u$;
$Y=Y_0x^2$ bounds $|\eta_i|$.
\[
\begin{array}{c|rr|r|c|c|c|ccc}
N&T_3<&C_*<&C&l&u&Y_0&E_O/x^2<& W_O/x^2<&d_Q/x^2<\\\hline
12&5.45&79.28&79.3&.44&1.415&14&6.24&7.59&.436\\
\end{array}
\]
Expansion of the quartic term gives the lower bound $4.8[N(k^2+2k+3)+5](k-1)^2$ at $l,u$. Write $A=4.8N(3+2l+l^2)$. The condition from the relative entropy term is
$2C/(1+\sqrt{1-C/A})<(Y_0^2/2)(1-2Y/3)$.
We compute the last three columns as follows.
In the formulas for the ellipse divide $\alpha$ and $K$ by $n^4$.
For $g(k)=-\log\sinc(hk)$ the function $g'(u)/x^2$ increases in $x$
by the sine product. Consequently the least scaled coefficients are
bounded below by
\[
 a=\left[2/\{l(1+\log(1-Y_N)/3)\}
           +(g_N'(u)/x_N^2)^2/A\right]^{-1},
 \qquad b=8N/\pi^4 ,
\]
where a subscript $N$ means evaluation at $x_N$.
Use $M_+(a,b;N-1)$ and $\sqrt{ab}$ for the two scaled path coefficients.
Indeed path energy increases in each coefficient, and
$(1-r^L)/(1+r^L)$ increases with $L$ for $0\le r<1$.
With $p^2=C/[4(a+M_+)]$, $s_0^2=C/[4(b+\sqrt{ab})]$,
the bounds for $E_O/x^2,W_O/x^2$ are $p\exp(px_N^2)$ and
$2\exp(p^2x_N^4/2)\sinh(s_0x_N^2)/x_N^2$.
The scaled center shift is bounded by
$x_NT_3(5/16+v_*)+(2/5)(2\pi)^4x_N^2N_0^2$; in the formula above
for $N_0^2$ use $(1-\sinc h)/h^2\le1/6$ and $\sinc h\le1$ in
numerators. The two coefficient ratios for the moment estimate are
below $.30$ and $.19$, respectively, for $N=12$;
the support, outer radius and polar moment tests give
$\rho_G>.87$, $\rho_G<1.08$, $D_G^2<.0034$.
Use $m_*<.048$, $z<.015$ from Proposition~\ref{continuous-response}
in the balancing test.

For the last skew bound keep the complementary path in the selected
edge estimate. Put $y_0=(N-u)/(N-1)$ and define $a'$ by replacing
$A$ in $a$ with $4.8N(l^2+2ly_0+3y_0^2)$. Thus
$M_-/n^4\ge\sqrt{a'b}$, and the estimate after
\eqref{paper-joint-ellipse} gives
\[
 (\delta_*/x)^2\le
 {Cx\over(3.4+4\pi^2 x\sqrt{a'b})(1-C/A)}.
\]
This increases up to $x_N$. Insert it in the tilt estimate for the center $Q$. With $d_1$ the resulting upper bound for
$\delta_*/x$ and $d_2=(d_Q/x^2)/(R-d_Q)$,
the bound on $|w|/(qh^2)$ is
$(1+\sup e)\tau(\pi x_N)\tau(d_1x_N+d_2x_N^2)(d_1+d_2x_N)/(2\pi)$.
The same substitution gives $f<1.0035$ and already gives
both stated skew bounds. To bound the selected aperture and mean radius,
write $E_0=(1+\log(1-Y_N)/3)/2>.4829$ and
$M_0=M_+(a',b;N-1)>.392$.
Since $g'(k)\ge h^2 k/3$, the term involving the squared difference between the
selected and complementary log-radius centers in \eqref{paper-joint-ellipse}
is at least $(2\pi)^4(1-k^2)^2 H_0(v)/36$,
where $H_0(v)=E_0vM_0/(E_0v+M_0)$, $v=l$ below one and $v=k$ above.
Let $A_4=4.8(k^2+2k+3)(k-1)^2$ and $A_6=12.4(k^6-1-6(k-1))$.
Expansion at one and convexity give
$A_n(k)\ge nA_4+A_6/n+4.8c(k-1)^2$ with $c=6$ below one and $c=5$
above. At $k=.475(1.0035),1.403/1.0035$, $A_6<N^2A_4$, so the minimum
of the bound as a function of $n$ occurs at $N$. Adding the squared-difference term gives respectively
more than $79.6,88$, both exceeding $C$. Each term increases as the aperture moves farther from one, which excludes the two outer intervals.

For the mean displacement put $k\in[v,w]$, with consecutive endpoints
$.44,.72,.95,1.415$. Since $g'(k)/(2kx^2)\le g'_N(u)/(2ux_N^2)<6.84$
by the sine product, and $b(k),b(y)$ have opposite signs in the center
of the selected ellipse, its absolute center divided by $x^2$ is at most
$\Gamma|k-1|$, with
\[
 \Gamma=3.42\max\left\{{w(w+1)\over w+2(.392)},2.06/11\right\}.
\]
Indeed $|y-1|=|k-1|/(n-1)$ and $1+y<2.06$. Set
$A_v=4.8[N(v^2+2v+3)+5]$, $D_v=.4829v+.392$.
With $s=|k-1|$ the total log bound divided by $x^2$ is
\[
 T_v(s)=\Gamma s+\sqrt{(C-A_vs^2)/(4D_v)}.
\]
The cosh estimate for the mean applies to the enlarged ellipse as well:
its coefficients are $D_v,b$, its right side is at most $C$,
and $C D_vx_N^4<4b^2$.
Maximize $T_v$ at $s=\sqrt{C/A_v}\,\Gamma/\sqrt{\Gamma^2+A_v/(4D_v)}$,
clipped to the range of $|k-1|$. The three resulting bounds are
$5.74,5.51,5.45$. Thus
$|e|/h^2\le [(\exp(5.74x_N^2)-1)/x_N^2+.436]/(2\pi)^2<.160$.
This completes all endpoint tests.
\end{proof}

\section{Aperture, skew, and logarithmic endpoint radii}
\begin{corollary}[Aperture and skew at the conformal center]\label{IV:cor:aperture-skew-product}
Let $z=(r_+-r_-)/(r_++r_-)$ be the signed relative difference of
the endpoint radii in a sector of half-aperture $\pi p$ and base-angle
half-difference $\delta$. The skew excess of Lemma~\ref{IV:lem:gamma-full-skew}
satisfies
\[
 pW(p,\delta/\pi)\ \ge\ \frac{32}{\pi^4}
              (1-4\delta^2/\pi^2)z^2 .
\]
On the range $a=\pi p\le.45$, $d=|\delta|\le.36$, $|z|\le.055$, the
logarithmic endpoint radii satisfy the stronger estimate
\[
 pW(p,\delta/\pi)\ \ge\ \frac8{\pi^4}
       \left(\log\frac{r_+}{r_-}\right)^2 .
\]
\end{corollary}
\begin{proof}
Put $v=\delta^2/\pi^2$. The rational denominator in that lemma yields
\[
 pW\ge\frac{32p^2v}{(1-4p^2)^2(1-4v)} ,
\]
because $(1+4p^2)/(1-4p^2)^2\ge1$.
For $0<s<\pi/2$ put $c=4/\pi^2$ and
$F(s)=s\cos s-(1-cs^2)\sin s$. Both endpoint values vanish, while
$F'(s)/(s\sin s)=2c-1+cs\cot s$ strictly decreases, because
$\sin s\cos s<s$. Hence $F\ge0$, or
$\tan(s)/s\le(1-cs^2)^{-1}$.
Apply the bound at $s=\pi p$ and $s=|\delta|$, then use the sine
rule $|z|=\tan(\pi p)\tan|\delta|$.

For the logarithmic bound keep also the $41px/25$ term of
Lemma~\ref{IV:lem:gamma-full-skew}. Dividing by $32 z^2/\pi^4$ gives the lower bound
\[
 (d\cot d)^2\left\{
 \frac{(a\cot a)^2}{(1-4p^2)^2(1-4d^2/\pi^2)}
       +\frac{41}{800}(a\cot a)^2\right\}.
\]
At $z=0$ the result holds directly. For $0\le s\le.45$,
$s\cot s\ge1-.35s^2$ by dividing $1-s^2/2$ by $1-s^2/6+s^4/120$.
As $4/\pi^2>.4$, the ratio $(a\cot a)^2/(1-4p^2)^2$ is at least
$1+.10a^2$; the quadratic term here exceeds the decrease $(41/800)(.70a^2)$ in the second
term. The display is therefore at least
$(1-.70d^2)(1+.40d^2+41/800)>1.003$, by its endpoint values.
But $(\operatorname{atanh}z/z)^2\le[1+z^2/(3(1-z^2))]^2<1.003$,
and $\log(r_+/r_-)=2\operatorname{atanh}z$.
\end{proof}

\section{Estimates for the disk term near the unit circle}
\begin{proposition}\label{IV:prop:unconditional-cap-twelve}
For $n\ge12$, $h=2\pi/n$, $s=\operatorname{sinc}h$ and $a=\mu_1(R_n)$,
the cap density $L=1-a+c_0$ satisfies $\lambda<a\le\lambda/s$ and
$.32h^2<-L<h^2/3+.077h^4$.
\end{proposition}
\begin{proof}
Put $\delta=1-as/\lambda$. The upper and lower bounds for $d_R$ give
$0\le\delta<5.64/(2.3899n^3)<1-s$ by $1-s\ge h^2/6-h^4/120$.
The projected mass identity gives
\[
 -L=(1-\delta)/s^2-1+
 \frac a{\pi s}\left[nM_D-\int_{\mathbb D\setminus R_n}(1-g^2)\right].
\]
Proposition~\ref{continuous-response} gives $M_D<.048(.015)h^5$.
Since $a<3.546$ and $s\ge\sinc(\pi/6)$, its positive contribution
$anM_D/(\pi s)=2aM_D/(hs)$ is less than $.006h^4$ throughout.
For all $n\ge12$ the positive series for $s^{-2}$ bounds
$(s^{-2}-1-h^2/3)/h^4$ between zero and its value $<.070$ at $h=\pi/6$.
On the caps $1-g^2\le3(h^2/8)^2$ by $rg'(r)\le(\lambda-1)(1-r)$.
Their area is $\pi(1-s)$, so the sum of the two terms being subtracted is less than
$.0055h^2+.030h^6$. This proves both bounds.
\end{proof}

\begin{lemma}[Expansion at the disk boundary from the Bessel zeros]\label{IV:lem:wall-zero-products}\label{IV:lem:lower-signed-cap}
Let $x=j'_{1,1}$, $\lambda=x^2$, $K(X)=J_1(x\sqrt X)/(J_1(x)\sqrt X)$
and $S(X)=K(X)+2XK'(X)$. For $a>0$ put
\[
 W=S^2+(1-aX)K^2+c_0,\qquad
 W(1-z)=L+\sum_{k\ge1}d_kz^k,\qquad R(X)=XW(X).
\]
Then $d_1=1$, all $d_k$ are positive, and
\[
 d_3<\tfrac35d_2,\qquad d_{k+1}<\tfrac4{15}d_k\quad(k\ge3).
\]
In particular, for $0\le z\le1/2$,
\begin{equation}
 2(d_2-1)-6z d_2\le R''(1-z)\le 2(d_2-1).
 \label{IV:eq:common-wall-curvature}
\end{equation}
\end{lemma}
\begin{proof}
Put $B=K(1-z)$, $S(1-z)=s_1zT(z)$, $T(0)=1$ and
$M=1-(1-z)B^2$. Then $M'=s_1zBT$ and $d_k=[z^k](s_1^2z^2T^2+B^2+aM)$
for $k\ge1$. Both $B,T$ are normalized products of $1+\alpha z$,
$\alpha=x^2/(j^2-x^2)>0$, over the remaining positive zeros $j$
of $J_1,J'_1$, respectively: the zeros are real
\cite[Section~10.21(i),(iii)]{DLMF}, and the series in the squared argument
have order $1/2$. For such a product $H=\sum h_jz^j$ the ratios
$(j+1)h_{j+1}/h_j$ decrease by $(Q')^2\ge QQ''$ at zero for each
derivative $Q$ of a finite product (Rolle), then passage to the limit.
This applies to $B^2,T^2,BT$; in $M$ the successive ratios from degree
two are $k((k-1)[z^{k-1}]BT/[z^{k-2}]BT)/[(k+1)(k-1)]$ and hence decrease.

Write $B=1+z/2+b_2z^2+\cdots$, $T=1+rz+rtz^2+\cdots$. The equation
$4XK''+8K'+\lambda K=0$ gives $s_1=(\lambda-1)/2$,
$b_2=(4-\lambda)/8<1/12$, $2/7<r=(\lambda-2)/(2\lambda-2)<3/10$,
$t=(8-\lambda)(\lambda-3)/(12(\lambda-2))<7/60$. At degrees
two and three the successive coefficient ratios of $B^2$ are at most
$1/3,1/4$; those of $z^2T^2$ are $2r,r/2+t$; those of $M$ are
$2(1/2+r)/3,3(b_2+r/2+rt)/(4(1/2+r))$. Each pair lies below
$(3/5,4/15)$. Positive addition proves the ratio bounds and $d_1=1$.
The same recurrence gives
\[
 \begin{aligned}
 4d_2&=a\lambda-a+\lambda^2-3\lambda+6,\\
 12d_3&=2a\lambda-3a+3\lambda^2-14\lambda+24,\\
 -48d_4&=a\lambda^2-8a\lambda+12a+\lambda^3-16\lambda^2+75\lambda-120 .
 \end{aligned}
\]
Finally $[z^j]R(1-z)=d_j-d_{j-1}<0$ for $j\ge3$. For $z\le1/2$
the curvature loss is at most $6z\sum_{j\ge3}(d_{j-1}-d_j)=6zd_2$,
as $j(j-1)z^{j-2}\le6z$.
\end{proof}

\begin{lemma}[Cap and chord curvature]\label{IV:lem:exact-real-cap}
\label{IV:lem:all-side-chord-remainder}\label{IV:lem:weighted-disk-curvature}\label{IV:lem:exact-skew-chord}
Use $W,L,d_j$ from Lemma~\ref{IV:lem:wall-zero-products}, with $a>0$ and
$\Phi'(r)=rW(r^2)$. Set $c_k=2(k+1)d_{k+1}-(2k+1)d_k$. Then $c_k<0$
for $k\ge2$, $2|c_{k+1}|\le |c_k|$ for $k\ge3$, and, if
$-2\le x=1-r^2\le S<1$, $S\ge0$, then
\[
 -\Phi''(r)\ge 2-L+c_1x+c_2x^2-|c_3|(x_+)^3/(1-S/2).
\]
For aperture $0<t<\pi$ write $S_t=\sin^2(t/2)$, $I=\operatorname{sinc}t$,
$\rho_0(s)=\cos(t/2)/\cos(ts)$, $|s|\le1/2$, and $I_j=\int_0^1(1-v^2)^j\,dv$.
The signed cap integral $F(t)=(t/2)\int_{-1/2}^{1/2}\int_1^{\rho_0^2}W(X)\,dX\,ds$
satisfies
\[
 -F''(t)=\tfrac12\sin t\left[L+\sum_{k\ge1}(k+1)I_kd_k S_t^k\right].
\]
Thus $F''$ changes sign at most once and any retained subset of the
positive terms gives a lower bound for $-F''$.
For $X=1-S_t(1-v^2)$ the same substitution yields
\[
 \int_{-1/2}^{1/2}\rho_0^2 H(1-\rho_0^2)\,ds= I\int_0^1 H(1-X)\,dv .
\]
If $D(e,w)=\int_{-1/2}^{1/2}\Phi(\rho)\,ds$ is the chord contribution
with endpoint radii $R\pm w/2$, $R=1+e$, then
\[
 -D_{ww}(e,0)/I=\int_0^1\left[-\tfrac14v^2\Phi''(R\sqrt X)
                  +\tfrac12(1-v^2)W(R^2X)\right]dv .
\]
\end{lemma}
\begin{proof}
The $c_k$ are the coefficients of $-\Phi''$ after its constant term.
Lemma~\ref{IV:lem:wall-zero-products} gives $d_3/d_2<3/5$ and, for $k\ge3$,
$d_{k+1}/d_k<4/15<(2k+1)/(6k+8)$, proving both sign and ratio bounds.
The tail starting at $k=3$ is nonnegative by alternation for $x\le0$;
for $x\ge0$ sum the geometric bound. The first two $d_k$ are given
by the recurrence at the disk boundary above.

Take $v=\tan(ts)/\tan(t/2)$, so $\rho_0^2=X$ and
$\rho_0^2ds=I\,dv/2$ on $[-1,1]$. For
$A_j(t)=(2j)^{-1}\int_{-t/2}^{t/2}[1-\cos^2(t/2)\sec^2\theta]^j\,d\theta$,
the endpoints vanish and the same substitution gives $A'_j=I_{j-1}S_t^j$.
Differentiate $F=-LA_1-\sum_{k\ge1}d_kA_{k+1}$ twice; the entire series
converges locally with its derivatives. Its displayed bracket increases.
Finally $\rho=R\rho_0(1-z^2)/(1-vz)$, $z=w/(2R)$, gives
$\rho_w=v\rho_0/2$, $\rho_{ww}=-(1-v^2)\rho_0/(2R)$ at zero skew.
Differentiation and reflection give the claimed chord curvature.
\end{proof}

\begin{corollary}[Disk curvature along a chord]\label{IV:cor:disk-annulus-product}
Suppose $\lambda\le a\le3.63$, $.782\le r\le1.10$ and
$\Phi'(r)=rW(r^2)$. Put $p=1-r^2$ and $C=4d_2-3$. Then
\[
 \begin{array}{ll}
 C-15.28p-31.41p^2-18p^3/(1-p/2)^2\le \Phi'''(r)/(2r)\le C
       & (r\le1),\\
 C\le \Phi'''(r)/(2r)\le C+15.28(r^2-1)&(r\ge1).
 \end{array}
\]
Here $12.42<C<13$, and $C<12.57$ when $a\le3.447$.
The lower bounds are positive and increase with $r$ after multiplication
by $2r$, as do the upper bounds.
\end{corollary}
\begin{proof}
For $H(p)=-\Phi''(r)$, $\Phi'''(r)=2rH'(p)$; the preceding lemma
gives the negative coefficients and their ratio bound. The expansion at the disk boundary in Lemma~\ref{IV:lem:lower-signed-cap}, with $3.3899<\lambda<3.39$, gives the stated bounds for $C=c_1$, as well as
$7.36<|c_2|<7.64$, $|c_3|<10.47$ and $|c_4|<9d_4<4.5$.
For $0\le p<.389$, the tail of $H'$ after its quadratic term
has magnitude at most $4|c_4|p^3/(1-p/2)^2$.
The resulting lower bound is positive even at $p=.389$.
For $p=-z\le0$ the tail after the linear term alternates with
decreasing magnitudes. Hence
$C+(2|c_2|-3|c_3|z)z\le H'(-z)\le C+2|c_2|z$.
The parenthesis is positive for $z\le.21$.
\end{proof}

\begin{proposition}[Integral identities on the sector and its chord]\label{IV:prop:lower-radial-profile}\label{IV:cor:lower-radial-twentyfour}
With $c=\cos(h/2)$ and $v$ the Cartesian coordinate along the chord, put
$F_*=c^2g'(r)/r^2+v^2g(r)/r^3$ and $J_*=2K'(r^2)$. Then
\[
 \int D_1\,dv=(\lambda_0-\lambda)\int U\cdot x+
 \lambda_0\int D\cdot x-\int cg'(r)\,dv,\qquad
 E_D-\lambda M_D=-\int(F_*D_1+cvJ_*D_2)\,dv .
\]
\end{proposition}
\begin{proof}
For every equivariant vector test $\Psi$ on the regular triangle with radial
traces on the vertex rays, the weak equation for $Y$ and integration by parts for $U$ give
\begin{equation}\label{op-weak-correction}
 \int\nabla D:\nabla\Psi=(\lambda_0-\lambda)\int U\cdot\Psi+
 \lambda_0\int D\cdot\Psi-\int_{\rm chord}\partial_\nu U\cdot\Psi .
\end{equation}
Indeed the reflected tests glue, the ray flux of $U$ is perpendicular to these
traces and all sectors give equal integrals. Take $\Psi=x,D$ and use
$\int U\cdot D=-M_D$. On the disk sector $\int U\cdot x=h/\lambda$ by the
Bessel recurrence and $J'_1(\sqrt\lambda)=0$.
In the radial variation of the chord integral write $p$ for the intermediate aperture,
$y=h^2/4$, $f(x)=\cos\sqrt x\operatorname{sinc}^3\sqrt x$,
$D_0=AB/p$, $N_p=D_0(2g'+rg''-u_0t^2H_1)$ and
$B_p=D_0B^2H_1$, where $H_1=-\lambda g/r$.
Then $D_0B^2=p^2f(yp^2)/f(y)$, and the preceding identity extracts
the principal term:
\[
 -2\lambda(E_D-\lambda M_D)+\int(N'_p/p-2\lambda F_*)D_1\,dv
       +c\int v(B'_p/p-2\lambda J_*)D_2\,dv .
\]
\end{proof}

\section{Conformal sector bounds and eigenvalue estimates for regular polygons}
\begin{lemma}[Endpoint bounds from a path]\label{IV:lem:middle-conformal-entry}
Let $p_i=q_i/n$, let $b_i$ be the area fractions, and put
$a_i=\log(b_i/p_i)$, $v_i=\log r_i$, $f(q)=-\log\operatorname{sinc}(hq)$.
Suppose $l<q_i<u$, $a_i>-.17$, and the skew terms in the conformal estimate with upper bound $C$ dominate
$\beta\sum(v_{i+1}-v_i)^2$. Including the relative entropy and quartic terms then gives
\[
 \sum_i[\alpha(v_i+v_{i+1})^2+\beta(v_i-v_{i+1})^2]\le C,\qquad
 \alpha^{-1}=2e^{.17}/(ln^4)+f'(u)^2/[4.8n(3+2l+l^2)] .
\]
Indeed $v_i+v_{i+1}=a_i+f(q_i)-f(1)$ by area, $f'$ increases,
and the relative entropy and quartic terms are bounded below by $ln^4e^{-.17}a_i^2/2$ and
$4.8n(3+2l+l^2)(q_i-1)^2$. For a path of $L$ edges
with coefficients $A,B$, let $c=\sqrt A+\sqrt B$,
$\rho=(\sqrt B-\sqrt A)/c$. The contribution of each edge equals
$c^2(v-\rho u)^2+2\sqrt{AB}(u^2-v^2)$. With endpoints fixed the
increments pair with $(\rho^{L-k})_{k=1}^L$ to give $v_L-\rho^Lv_0$;
Cauchy--Schwarz is attained by proportional increments. Thus the minimum
is $M_+(v_0+v_L)^2+M_-(v_0-v_L)^2$,
$M_\pm=\sqrt{AB}(1\mp\rho^L)/(1\pm\rho^L)$.
With $(A,B,L)=(\alpha,\beta,n-1)$, add the selected edge to get an ellipse
of squared semiaxes $p^2=C/[4(\alpha+M_+)]$, $s^2=C/[4(\beta+M_-)]$ for
$(v_i+v_{i+1})/2,(v_i-v_{i+1})/2$. Hence
$|e_O|\le e^p\cosh s-1$, $|w_O|\le2e^{p^2/2}\sinh s$ by
$e^{pv}\sqrt{1-v^2}\le e^{p^2/2}$ and convexity of $\sinh$.
\end{lemma}

\begin{lemma}[A lower bound for $d_R$]\label{IV:lem:conformal-sideband-energy}
For every $n\ge3$, $x=1/n$, the normalized spectral difference of the regular polygon satisfies
$d_R>4\zeta(3)x^3+41.7x^5$.
\end{lemma}
\begin{proof}
Use the normalized disk notation $Z,e=\lambda-1$, $a_i=(2x)_i/i!$,
$u=g(r)e^{i\theta}/\sqrt e$ and $q$ of Lemma~\ref{IV:lem:uniform-regular-remainder}.
Suppose first that $0<x=1/n\le1/5$.
Take the Neumann inverse of $r^{n+1}(e^{i(1-n)\theta}+e^{i(1+n)\theta})$
as a second trial direction. Its two radial factors for $m=|1\pm n|$ are
$v_m=t_m r^m[1+m(1-r^{N-m})/(N-m)]$ with $N=n+3$, $t_m=1/[m(m+N)]$,
as follows by differentiation. Its normalized Dirichlet energy, also its
pairing with its source, is $E=\sum_m t_m^2(2m+m^2/N)$. Normalize the
direction by $\sqrt E$. All terms in its mass pairing with $u$ are positive.
The first cross term of the conformal density has radial coefficient
$a_0a_1r^n/Z$. Keeping this term and using $g\ge r$ gives exactly
the source pairing above, hence a lower bound $b_0=2x\sqrt{E/e}/Z$.
Writing $F_m(p)=2[1+m/(p+N-m)]/p$ records the integral of the bracket
in $v_m$ against $2r^{p-1}dr$.

In its integral the density relative to $r^{p-1}dr$ decreases, where $p=n+m+3$;
hence the mean of $y=1-r^2$ is at least $2/(p+2)$. The stronger
$g/r\ge1+y/2$ therefore improves this term by a factor $1+x/(2+6x)$.
For the next term the Fourier coefficient relative to the first is $a_2=x(1+2x)$.
In $F_m(p+2n)/F_m(p)$ the factor $p/(p+2n)$ is at least $(1+x)/(2+x)$
and $(1+m/(4n+6))/(1+m/(2n+6))\ge5/6$ by $m\le n+1$.
Thus the coupling exceeds $b_*=b_0(1+7x/8)$: indeed
$(2+6x)^{-1}+5(1+2x)(1+x)/[6(2+x)]\ge7/8$, by
$120x^3+157x^2-15x+2>0$.

Write $Y=1+3x$. The energy formula gives
$2(2+1/Y)/(1+Y)^2\le E/x^3<2$ by convexity in $m/n$ for the lower bound.
The positive log series gives $Bx^3+12.44x^5<\log Z<5.38x^3$,
so $Z\le1+5.65x^3\le1+.226x$. Consequently
$(1+7x/8)/Z\ge1+.62x$ and
\[
 \frac{b_*^2 e}{4x^5}\ge 2(2+1/Y)\left(\frac{1+.62x}{2+3x}\right)^2
 \ge 1.305-1.620x .
\]
Each positive factor in the middle is decreasing and convex (the ratio is
a constant plus a positive reciprocal). Their product lies above its tangent
at $x=1/5$, whose intercept exceeds $1.305$ and slope exceeds $-1.620$.
The scalar disk entry satisfies $1/\lambda\le q\le[1+(1-1/Z)/e]/\lambda$;
$E<2x^3$ now gives $q+b_*^2/q<.306$. The two-dimensional trial thus improves the lower bound by
at least $b_*^2/.306$ over $q$.

For the disk loss, $1-g^2\le(3y^2+y^3)/4$. Beta integration and $a_j\le a_2$
for $j\ge2$ give
\[
 D/x^5\le6\left[\frac{1+x/(1+4x)}{(1+x)(1+2x)(1+3x)}
 +\frac{.203(1+2x)^2}4\left(1+\frac{x}{2+4x}\right)\right]
 <6.64-16.60x .
\]
Here $(1+x/(1+4x))/(1+x)\le1$; use the endpoint chord of the convex
$1/[(1+2x)(1+3x)]$ and bound the last term at $1/5$.
The mass identity and $1-1/Z\ge\log Z-(\log Z)^2/2$ now give
\[
 (d_R-Bx^3)/x^5\ge12.44-5.38^2x/2-3.39(6.64-16.60x)
       +4(3.3899)(1.305-1.620x)/.306>41.7
\]
by its two endpoints. The explicit triangular and square spectra give
$d_R>.95,.18$ respectively, which also exceed $Bx^3+41.7x^5$.
\end{proof}

\begin{proposition}[Two moments along the chord of a regular sector]
\label{IV:prop:refined-middle-entry}\label{IV:prop:regular-lower-calibration}\label{IV:lem:projected-wall-pair}
For $n\ge10$ let $I=-h^{-3}\int D_1\,dv$ and $\nu=h^{-4}\int vD_2\,dv$ on the regular chord. These moments satisfy \eqref{op-coordinate} and \eqref{op-torque}. Combining these inequalities with the bounds in Section~\ref{mid-response} gives the two moment columns in that section.
\end{proposition}
\begin{proof}
Write $X=h^2$, $z=E_D/h^4$, $b=\sin(h/2)$, $c=\cos(h/2)$;
$d_\pm$ bound $(\lambda_0-\lambda)/h^2$. Taking $\Psi=x$ in
\eqref{op-weak-correction} gives
$I=A-dX_U-\lambda_0h^{-3}\int D\cdot x$, where
$A=h^{-3}\int c g'(r)\,dv$, $r^2=c^2+v^2$ and
$X_U=h^{-1}\int U\cdot x$. The removed circular cap has area
$h(1-\operatorname{sinc}h)/2$ and $c^2\le rg(r)\le1$ there.
Thus $1/\lambda-(1-\operatorname{sinc}h)/2\le X_U\le
1/\lambda-c^2(1-\operatorname{sinc}h)/2$.
Projection gives $\int D\cdot U=-M_D$ and
Lemma~\ref{IV:lem:bg-spectral} gives
$\|x-2U/(\lambda-1)\|_2<.05613\sqrt h$ on the disk sector.
Thus the last term of $I$ lies in
$[-v,v+.84\lambda_0mzX]$, $v=.05613\lambda_0\sqrt{mz}$.

Here $c>.951$ and $0\le g'\le g/r$, so
$0\le g'''=(3/r^2-\lambda)(g'-g/r)\le(\lambda-3/r^2)/r<.39$.
Two integrations from one give
$0\le g'(r)-(\lambda-1)(1-r)\le.20(1-r)^2$.
Use $1-r=(b^2-v^2)/(1+r)$, whose numerator integrates to $4b^3/3$
and whose squared numerator integrates to $16b^5/15$. This yields
\[
 \frac{2.3899}{12}(1-X/4)\le A\le
 \frac{2.39}{12}(1-3X/16+.020X^2)+.20X/112.
\]
Indeed $c\operatorname{sinc}^3(h/2)\ge(1-X/8)^2$, while
$2c/(1+c)\le1-X/16+X^2/768$ and
$\operatorname{sinc}^3(h/2)\le1-X/8+.007X^2$ by the first
three Taylor terms. The square term uses $c/(1+c)^2\le1/4$.
For $d_-<d<d_+$ the resulting polynomial bounds are
\begin{equation}\label{op-coordinate}\begin{aligned}
 I&\ge\tfrac{2.3899}{12}(1-X/4)
 -d_+[1/3.3899-(1-X/4)(X/12-X^2/240)]-v,\\
 I&\le\tfrac{2.39}{12}(1-3X/16+.020X^2)+.20X/112
 -d_-(1/3.39-X/12)+v+.84\lambda_0mzX .
\end{aligned}\end{equation}
Indeed $X/6-X^2/120\le1-\operatorname{sinc}h\le X/6$.
For the second moment, use the chord identity
$E_D-\lambda M_D=-\int(FD_1+2cvK'(p)D_2)\,dv$ with $t=v/b$,
$p=1-b^2(1-t^2)$ and $F=K(p)+2c^2K'(p)$. On $[c^2,1]$ one has
$.152<K''<.16$: $K''(1)=(4-\lambda)/4$ and $|K'''|<.04$.
For the latter, the positive Taylor coefficients $a_j$ of $K(1-u)$
satisfy $2(j+1)a_{j+1}\le a_j$ by the product in
Lemma~\ref{IV:lem:wall-zero-products}, and
$6a_3=(24-7\lambda)/8$, so $|K'''(1-u)|\le6a_3 e^{u/2}<.04$.
Since $F(1)=b^2$ and $-K'(1)=1/2$, integration gives
\[
 F/h^2=(b/h)^2[1+B(t)(1-t^2)],\qquad
 1/6<B(t)<1/4 .
\]
Indeed $B$ averages $-K'-2c^2K''$, which is between $.5-2(.16)$
and $.5+.16b^2-2c^2(.152)$. Put $B_0=5/24$ and
$f=(b/h)^2(1+2B_0/3)\le41/144$. The rms of $F/h^2-f$ in
normalized Lebesgue measure on the chord is at most
$[B_0(2/\sqrt{45})+(1/24)\sqrt{8/15}]/4<1/40$.
Also $-2K'$ differs from $1+b^2/6$ by at most $b^2/6$.
Cauchy--Schwarz, $2b\le h$, $dv\le(h/c)ds$ and
$\int v^2\,dv\le h^3/12$ give
\begin{equation}\label{op-torque}
 c(1+b^2/6)\nu\ge s_D-
 h\left[fI+T_1\sqrt{z/c}/40+hT_2\sqrt z/(24\sqrt{12})\right].
\end{equation}
The bounds for both chord moments in Section~\ref{mid-response} follow by the substitutions stated there.
\end{proof}

\begin{proposition}[An estimate for the disk cap using three moments]\label{IV:prop:expanded-cap-eighteen-nineteen}
Set $H_j(q)=\sum_{i=0}^{j-2}(j-1-i)q^i$, $P(q)=(q-1)(3q^2+9q+8)/360$,
and let $\mathscr Tf$ denote the angular Taylor quotient at one. Use the
regular-polygon data $a=\lambda_0=\mu_1(R_n)$ and
$c_0=\lambda_0|R_n|^{-1}\int_{R_n}|Y|^2$ in $W$ and $L=1-a+c_0$,
with $R_n$ of circumradius one and $Y$ the projected field of Section~\ref{IV:sec:background}.
Suppose $n\ge10$, $1\le Q\le1.6$, $hQ\le1$, and $bh^2\le-L\le h^2/3+eh^4$
with $b\ge0$. Write $z=h^2H_9(Q)/H_7(Q)$. For $0\le q\le Q$ the disk
term of Lemma~\ref{IV:lem:exact-real-cap} satisfies
\[
 \mathscr T G_U(q)\le-h^5P(q)-
 h^7[kH_7(q)+bH_5(q)/240-eH_3(q)/12].
\]
With $d_1=1$ and $d_2,d_3$ from Lemma~\ref{IV:lem:wall-zero-products},
one may take
\begin{align*}
 k={}&\frac{L-10+12d_2}{10080}
 +\frac{-L+42-168d_2+144d_3}{725760}z\\
 &+\frac{L-170+1764d_2-4320d_3}{79833600}z^2
 -\frac{d_3}{1587600}z^4 .
\end{align*}
The elementary lower bounds
$.00370-.000460z-.000114z^2$ for $n=10,11$ and
$.003584-.000467z-.000114z^2$ for $n\ge12$ may also be used.
In particular $k>.003$ for $Q\le1.42$. For $n\ge12$ one may take
$b=.32$, $e=.077$; the bracket then has positive coefficients for $Q\le1.6$.
\end{proposition}
\begin{proof}
Write $I_j=\prod_{i=1}^j2i/(2i+1)$. In $-F''(t)$, $0\le t\le1$, retain
$L\sin(t)/2$ and the three positive terms
$(j+1)d_jI_j\sin(t)\sin^{2j}(t/2)/2$, $1\le j\le3$.
For independent uniforms $V_i$ on $[-1,1]$, the characteristic function gives
\[
 \sin(t)\sin^{2j}(t/2)=4^{-j}t^{2j+1}\mathbb E\cos(tY_j),
 \qquad Y_j=V_0+\tfrac12\sum_{i=1}^{2j} V_i .
\]
Use the degree-six cosine lower polynomial and the degree-nine sine
upper polynomial in the term with $L<0$. Equivalently expand
$\sin(t)[(1-\cos(t))/2]^j$ through degree $2j+7$ in each positive term.
Twice integrating the resulting differential inequality in the nonnegative
Taylor kernel gives coefficients $p_m$ of the lower polynomial for $-F$.
The first two are $p_3=L/12$ and $p_5=(2-L)/240$.
The coefficients $p_7,p_9,p_{11},p_{15}$ are the respective coefficients
of $1,z,z^2,z^4$ in the displayed formula for $k$.
The remaining coefficient is
$p_{13}=(-32d_2+171d_3)/23587200>0$.
These identities follow by dividing the degree-$(m-2)$ coefficient by $m(m-1)$.
Here $p_9,p_{11},p_{15}<0$, since $3.855<d_2<4$, $1.98<d_3<2.07$ and $-.145<L<0$.
For $j\ge7$ the identity
\[
 {H_{j+2}(q)\over H_j(q)}=q^2+{j(1+q)+1\over H_j(q)}
 \le{H_9(Q)\over H_7(Q)}
\]
follows since each coefficient $(j-1-i)/[j(1+q)+1]$ increases in $j$
and the coefficient ratios $(8-i)/(6-i)$ for $H_9/H_7$ increase in $i$.
Discard $p_{13}$ and use this bound in each negative term.
The first two terms give the remaining coefficients since $H_3/36-H_5/120=-P$.

To check the simpler constants observe that $H_7(Q)\ge6+5Q$, so
$z\le h^2(Q^2+7/5)<1.56$ since $hQ\le1$ and $h\le\pi/5$. For $n=10,11$ the recurrences at the disk boundary give
$3.966<d_2<3.992$, $2.044<d_3<2.058$; thus $p_7>.003714$,
$p_9>-.000461$, and $p_{11}+p_{15}z^2>-.000032$. For $n\ge12$ use
$3.855<d_2<3.963$, $1.98<d_3<2.05$, and $-L<.13$ from
Proposition~\ref{IV:prop:unconditional-cap-twelve}.
The corresponding bounds are $.003584,-.000467,-.000032$.
They imply both simpler lower estimates.
At $Q=1.42$ one has $z<.96$ at ten sides and $z<.67$ at twelve;
for $n\ge12$, $Q\le1.6$, the second polynomial exceeds $.00314$.
This proves the remaining assertions.
\end{proof}

\begin{lemma}[First variation and correction]\label{IV:lem:middle-radial-slope}\label{IV:lem:one-sided-pure-radius}
Let $y=h^2/4$, $c=\cos(h/2)$, $s=\operatorname{sinc}h$,
$Q=\max(1,q)$, $v_*=yQ^2\le.222$, $z=E_D/h^4$,
$z_y=h^{-4}\int|\partial_2D|^2$, $M_D\le mh E_D$ and
$f(x)=\cos\sqrt{x}\operatorname{sinc}^3\sqrt{x}$.
In $J=G_e/[qh^3(q^2-1)]$ at zero offsets the principal boundary term is
$-2\lambda(E_D-\lambda M_D)/h^3$. For the remaining boundary terms use
$N=\overline{N'_p/(ph^2)}-2\lambda F_*/h^2$,
$T=\overline{B'_p/p}-2\lambda J_*$, averaging in $p^2$ between $1$ and $q^2$.
If $-.32y\le J_*+1\le0$, then
\begin{equation}
 |T|\le2\lambda\max\{1/f(y)-1,1-[1-y(1+q^2)]/f(y)+.32y\}. \label{paper-wall-slope}
\end{equation}
The absolute volume and the signed pure correction have respective upper bounds
\begin{equation}
 \frac{\lambda_0-\lambda}{h^2}\sqrt{mz}
 \left[\frac{h^4(5(6.36)/6+4.06/2)}{2sc^2}
       +\frac{h^3(2(6.36)+4.06h^2Q^2/2)}{\sqrt{48}f(y)}\right],             \label{paper-wall-volume}
\end{equation}
\begin{equation}
 2h\left[1+\frac{1.2v_*^2-.8v_*^3}{(1-v_*)^2}\right](z-z_y)
            +5\lambda_0m h^4z/(3sc^2).                                   \label{paper-wall-remainder}
\end{equation}
\end{lemma}
\begin{proof}
The decomposition of the boundary integral follows from Proposition~\ref{IV:prop:lower-radial-profile}.
Put $x=yp^2$, $S=\sin^2\sqrt x$, $d=1-t^2$, $X=1-Sd$.
Since $B_p=-\lambda p^2 f(x)K(X)/f(y)$, $B'_p/p+2\lambda
=2\lambda[1-A/f(y)]$ with $A=(f+xf')K-xfS'dK'$.
The alternating series of $8r^3f(r^2)=2\sin(2r)-\sin(4r)$ gives
$1-2x\le f+xf'\le1-2x+1.2x^2$ and hence
$1-x\le f\le1-x+.4x^2$. By Lemma~\ref{IV:lem:normal-wall-mean},
$1\le K\le1+.54x$, $-.54\le K'\le0$. Thus
$1-2x\le A\le(1-2x+1.2x^2)(1+.54x)+.54x\le1$;
averaging proves \eqref{paper-wall-slope}.

Write $A_0=\cos(hp/2)/c$, $B_0=\sin(hp/2)/\sin(h/2)$, $D_0=A_0B_0/p$.
The two volume weights are $D_0(A_0^2,B_0^2)(4K+2XK')$ with
$X=A_0^2x_1^2+B_0^2x_2^2$, $|X'_p|\le ph^2/2$. By the sine series
and $D_0B_0^2=p^2f(x)/f(y)$ their first factors have derivatives bounded
by $5ph^2/(6sc^2),2p/f(y)$ and values bounded by $1/(sc^2),Q^2/f(y)$.
The Bessel series gives $|4K+2XK'|<6.36$, $|6K'+2XK''|<4.06$ on
$[0,1]$. Use $\int x_1^2\le h/4$, $\int x_2^2\le h^3/48$ for \eqref{paper-wall-volume}.

In the pure term let $E_{ij}=\int|\partial_j D_i|^2$, $M_i=\int|D_i|^2$.
Before radial differentiation its energy and mass weights divided by $q$ are
$R^2D_0(1,(A_0/B_0)^2,(B_0/A_0)^2,1)$,
$-\lambda_0 R^4D_0(A_0^2,B_0^2)$, respectively. In the $q^2$ difference
quotient all but $E_{21}$ and $M_1$ have favorable signs by monotonicity,
using $f+xf'>0$. For $j=f/\cos^2\sqrt{x}$ one has
$D_0(B_0/A_0)^2=q^2j(yq^2)/j(y)$, $j(y)\ge1$ by
$(1-y/6)^3\ge1-y/2+y^2/24$. The preceding $f,f+xf'$ bounds,
$\cos^2\sqrt x\ge1-x$ and $S'\le1$ give
$j+xj'\le1+(1.2x^2-.8x^3)/(1-x)^2$ (less than $1.015$ at $x\le.10$).
Use $E_{21}\le h^4(z-z_y)$ and the bound on $(D_0A_0^2)'$ for \eqref{paper-wall-remainder}.

A lower pure bound is $-2Phz-4\lambda_0M_D/(h^3 f(y))$, where
\[
 P=\min(1,q)^{-2}\left[1+1.5y/((1-v_*/6)(1-y/2)^3)\right].
\]
Indeed $D_0(A_0/B_0)^2=q^{-2}v(yq^2)/v(y)$ with
$v(x)=\cos^3\sqrt x/\operatorname{sinc}\sqrt x$,
$v(y)\ge(1-y/2)^3$, $v'\ge-1.5/(1-v_*/6)$, $v$ decreasing
(logarithmic differentiation). Its difference quotient is at least $-P$:
if $q\ge1$ expand at $v(y)$, and if $q<1$ multiply the negative
variation by $q^{-2}$. The magnitude for $D_0$ is at most $h^2/(6s)<P$;
the $D_0B_0^2$ quotient is at most $1/f(y)$.
\end{proof}

\begin{lemma}\label{IV:lem:convex-radius-disk}
For $n\ge10$, take the regular-polygon data $a=\mu_1(R_n)$ and $c_0$
of Proposition~\ref{IV:prop:expanded-cap-eighteen-nineteen}. Suppose $-.145<L=L_{\rm cap}<0$ and let
$0\le(hq)^2\le.89$. The disk contribution to $J$ is the divided
difference of the convex function $D$ below, with $D'(0)=(1-L)/6$.
With $\ell=-L$, $d_2<4$, $1.98<d_3<2.07$ and $d_4\le1/2$ one has
\begin{align}
 D'(x)\le P(x)={}&(1+\ell)/6+(d_2/15-1/12-\ell/60)x\nonumber\\
 &+(\ell/1680+1/80-d_2/30+3d_3/140)x^2+.0049x^3+.0059x^4 .
 \label{IV:eq:common-disk-slope}
\end{align}
\end{lemma}
\begin{proof}
Put $w=\operatorname{sinc}\sqrt z$, $S=(1-\cos\sqrt z)/2$ ($S'=w/4$).
The chord change of variable cancels the factor $X$ in $XW(X)$, giving
the secant from $h^2$ to $(hq)^2$ of $D(z)=w\sum_{j\ge0}d_jI_jS^j$,
$I_0=1$, $I_j/I_{j-1}=2j/(2j+1)$, $d_0=L$. For $z\le.89$ the sine
series gives $S<.21$, $w>.85$, $-1/6\le w'\le0$, $0\le w''\le1/60$.
For $j\ge1$, $(wS^j)'=S^{j-1}(w'S+jw^2/4)\ge0$ and
\[
 (wS^j)''=w''S^j+(3j/4)ww'S^{j-1}+j(j-1)w^3S^{j-2}/16 .
\]
This is nonnegative for $j\ge2$ since $(j-1)w^2+12w'S>0$;
for $j=1,2,3$ it is at least $-w/8,w(.85^2/8-S/4),.19wS$.
Using $d_2>3.85,d_3>1.98$ from Lemma~\ref{IV:lem:lower-signed-cap}
with $a=\lambda_0\ge\lambda$
gives $D''\ge-.145/60+w[-1/12+2.05(.85^2/8-S/4)+.90(.19)S]>0$.

For $j\le3$ write $wS^j=(z/4)^j\mathbb E\cos(\sqrt z\,Z_j)$, with
$Z_j$ the sum of one uniform on $[-1,1]$ and $2j$ independent uniforms
on $[-1/2,1/2]$. The successive even-moment ratios are at most $(j+1)^2$,
so derivative-series terms decrease from index $m=2$: their ratio is at
most $.89(j+1)^2(j+m+1)/[(2m+1)(2m+2)(j+m)]<1$.
Expanding these finite products of sinc series gives the upper coefficients
of $I_j(wS^j)'$, in increasing powers of $z$, as
\[
 \begin{array}{c|rrrrr}
 j=1&1/6&-1/12&1/80&-17/18144&31/725760\\
 2&0&1/15&-1/30&7/1080&0\\
 3&0&0&3/140&-1/84&19/6720
 \end{array}.
\]
For $j\ge4$ the decreasing sinc factors bound each derivative by
$jz^{j-1}/4^j$. With $r=.89/15$ the weighted tail is at most
$(I_4/512)[4/(1-r)+r/(1-r)^2]z^3<.003429z^3$
by $d_4\le1/2$ and the $4/15$ ratio. Add $-w'\le1/6-z/60+z^2/1680$
with its factor $-L$ to obtain \eqref{IV:eq:common-disk-slope}.
\end{proof}

\begin{lemma}[Mean of the normal coefficient along the chord]\label{IV:lem:normal-wall-mean}
Put $y=h^2/4$, $Y=yp^2$, $S=\sin^2\sqrt Y$, $d=1-t^2$,
$X=1-Sd$, $s=\operatorname{sinc}h$, $c=\cos(h/2)$ and
$u_1=\sin^2(h/2)$. Average in $p^2$ between $1$ and $q^2$.
For $Y\le.222$ the longitudinal coefficient before subtracting the reference-chord term satisfies
$\overline{N'_p/(ph^2)}\ge [4.9831-7.56\,y(1+q^2)/2]/(2s)$.
If $K''\le k$ on the reference chord, the chord mean of $F_*/h^2$
with respect to $dv$ is at most $f_0+2b_+/3$, where
$f_0=u_1/h^2$ and $b_+=f_0(.5+ku_1-2c^2(.1525))$.
\end{lemma}
\begin{proof}
Write $w(Y)=\operatorname{sinc}(2\sqrt Y)$, so $S'=w$.
The normal coefficient in Proposition~\ref{IV:prop:lower-radial-profile} is
$N_p=wH/s$, $H=(2-\lambda+\lambda S)K+2XK'$.
Set $A=(1-S)K'+K/2$. The Bessel equation
$2XK''+4K'+\lambda K/2=0$ gives exactly
\[
 H_S=\lambda(K+dA),\qquad H_d=\lambda S A,\qquad
 N'_p/(ph^2)=(w'H+w^2H_S)/(2s).
\]
At $S=0$, $H=1-\lambda$, $A=0$. On this chord
$.1525<K''<.172$, since $K''(1)=(4-\lambda)/4$ and
$K''(1-Sd)\le K''(1)e^{S/2}$ by the product for $K$.
Thus $-.54<K'\le-.5$, $1\le K\le1+.54S$ and
$0<A_S=-(1+d/2)K'-(1-S)dK''<.80$.
In particular $\lambda\le H_S\le\lambda(1+1.35S)$ and
$H\le1-\lambda+\lambda Y(1+.675Y)<0$.
Use $-2/3\le w'\le-2/3+.267Y$ and $w^2\ge1-4Y/3$:
\[
 w'H+w^2H_S\ge
 (2/3-.267Y)[\lambda-1-\lambda Y(1+.675Y)]
       +\lambda(1-4Y/3)\ge4.9831-7.56Y .
\]
Indeed the linear loss is at most $2\lambda+.267(\lambda-1)$
and the quadratic loss at most $.183\lambda Y^2$.
For the upper bound $w'H+w^2H_S\le4.984$, use
$S\ge(1-Y/3)Y$, $H\ge1-\lambda+\lambda S$ and
$w^2\le1-4Y/3+(4/9+4/15)Y^2\le1-1.17Y$.
Indeed $(1-1.17Y)(1+1.35Y)\le1+.18Y$.
On the reference chord $F_*=K(X)+2c^2K'(X)$ at $S=u_1$.
There $F_*(d=0)=u_1$, and
$\partial_d(F_*/h^2)=f_0(-K'-2c^2K'')\le b_+$.
The chord mean of $d$ is $2/3$.
\end{proof}

\begin{lemma}[Integrating the disk field on the affine triangle]\label{IV:lem:affine-disk-skew}
Fix an aperture $t\in(0,\pi)$ and endpoint radii $R\pm w/2$.
Write $D=\int_{-1/2}^{1/2}\Phi(\rho)\,ds$ for the disk term,
$S=\sin^2(t/2)$, $z=w/(2R)$ and $I=\operatorname{sinc}t$.
Suppose $\lambda\le a\le3.63$, $3/4\le R^2\le32/25$,
$S\le2/9$, and $|z|\le\min(1/14,\sqrt{S(1-S)})$.
At fixed $R,t$ one has $D_{wwww}<5$. Consequently
\[
 D_{ww}(w)\le D_{ww}(0)+5w^2/2,\qquad
 2[D(w)-D(0)]/w^2\le D_{ww}(0)+5w^2/12
\]
with the second expression interpreted by continuity at zero.
\end{lemma}
\begin{proof}
In affine coordinates on the chord, for $y$ uniform on $[-1,1]$,
$X=R^2(1-S+Sy^2+2yz+Bz^2)$, $B=S+(1-S)y^2$.
The triangle Jacobian gives $D-\Phi(0)=IR^2\mathbb E[(1-z^2)V(X)]/2$,
$V(X)=\int_0^1W(vX)dv$. By
$W-c_0=\frac12\Delta(g^2)-(a-\lambda)g^2$, Bessel convolution gives
$V=\mathrm{const}+\sum_{j\ge1}(-1)^jv_jX^j$, where
\[
 h_j=(\lambda/4)^j\binom{2j+2}{j+1}/[j!(j+2)!],\qquad
 v_j=K(0)^2[2(j+1)h_j+(a-\lambda)h_{j-1}/(j+1)] .
\]
The scaled successive ratios $(j+1)v_{j+1}/v_j$ are below $.64$ at
$j\ge3$: for the two summands they equal
$\lambda(2j+3)/[2(j+1)(j+3)]$ and
$\lambda(2j+1)(j+1)/[2j(j+2)^2]$, both decreasing.
Hence $(-1)^kV^{(k)}\ge0$ for $k\ge3$ and $X\le3/2$ by alternation.
Since $R^2(1-S)(1-|z|)^2\le X\le R^2(1+|z|)^2$, here $1/2<X<3/2$.
The first five coefficients give
$V''(3/2)\ge2v_2-9v_3+27v_4-67.5v_5>0$, $V''\le 2v_2<10$,
$V^{(4)}\le24v_4<5$, $-V'\le v_1-v_2+3v_3/4<11/2$;
use $1.58<K(0)<1.59$ by the alternating Bessel series.

Put $v=R^2B$ and $P=X'=2R^2(y+zB)$. In the fourth $z$-derivative
of $(1-z^2)V(X)$ the $V'''$ coefficient is
$P^2[12v(1-z^2)-8zP]\ge0$ because $B\ge2|zy|$.
Completing the square bounds the $V''$ coefficient by
$12(1+3z^2)v^2$; the $V'$ coefficient is $-24v$.
With $Z=1/14$, $\partial_w=(2R)^{-1}\partial_z$ now gives
\[
 D_{wwww}\le (5/2)R^6\mathbb E(|y|+Z)^4+
 (15/4)R^2(1+3Z^2)\mathbb E B^2+(33/8)\mathbb E B<5 ,
\]
since the respective uniform moments are below $.283,.286,.482$.
Evenness gives both stated Taylor bounds. In the geometric windows for
$10\le n\le19$, $h<.629$, $q\ge1/3$, $hq<.941$, $|w|\le.30h^2$
and $.878<R<1.122$ give $S<2/9$, $|z|<1/14$ and
$2|z|<5hq/6\le\sin(hq)$, as required.
\end{proof}

\begin{lemma}\label{IV:lem:exact-trial-centering}
For the projected trial on a fan based at $Q$,
\[
 \int_P V_Q=\int_P U(x-Q)\,dx+\alpha|P|(G-Q),\qquad
 \alpha=\frac{3\int_{C_*}D_1}{2\cos(h/2)|C_*|},
\]
where $G$ is the area centroid and $C_*$ the reference triangle.
\end{lemma}
\begin{proof}
Bisector reflection and the triangle centroid formula give
$\int_{C_*}D=\alpha\int_{C_*}x$. Under $x\mapsto Q+Tx$ both
sides transform by $|\det T|T$. Summing gives the formula.
\end{proof}

\begin{lemma}[Boundary control from the two leading conformal coefficients]
\label{IV:lem:convex-boundary-two-heads}
Let $f$ map the unit disk conformally onto a convex polygon of area $\pi$,
with $O=f(0)$ and $a=f'(0)>0$. Write $f'/a=1+\sum_{j\ge1}\gamma_jz^j$,
$T=a^{-2}-1$ and $v=f-O-az$. Put
\[
 c_1=\gamma_1/2,\qquad c_2=\gamma_2-\gamma_1^2/2,\qquad
 \delta=[1-4(|c_1|+|c_2|)]/3,\qquad C=e^{2|c_1|+|c_2|}.
\]
If $\delta>0$ then
\begin{equation}
 \|v\|_\infty/a\le
 \sqrt{[5(\log 2)/\delta-3/2]T}+C/(64\delta)+2^{-15}.                 \label{paper-wall-centering}
\end{equation}
For $T\le1/24$ the disk $B_{2/3}(O)$ lies in the polygon.

\end{lemma}
\begin{proof}
Convexity gives $\Re p>0$ for $p=1+zf''/f'$. Equivalently,
the Schwarz--Christoffel formula expresses $p$ as a convex combination of
$(1+\bar\xi z)/(1-\bar\xi z)$, with weights given by the exterior
angles divided by $2\pi$. The first terms of $p$ are
$1+2c_1z+2c_2z^2$. With $s=2(|c_1|+|c_2|)$ the function
$P=(p+s-2c_1z-2c_2z^2)/(1+s)$ has positive real part and
$P=1+O(z^3)$. Schwarz's lemma gives
$|(P-1)/(P+1)|\le |z|^3$. Integrating
$\partial_r\log|f'(re^{it})|=\Re(p-1)/r$ yields
$|f'(re^{it})|/a\le C(1-r^3)^{\delta-1}$.
On the circle of radius $r$, Cauchy--Schwarz and the area identity give
$|v|/a\le\sqrt{T[-\log(1-r^2)-r^2]}$. To reach the boundary, integrate
$|f'|/a+1$ on the remaining radial segment. Choose $1-r=2^{-5/\delta}$
so that $(1-r)^\delta=1/32$. As $0<\delta\le1/3$, $r> .99$ and
$(1+r+r^2)^{1-\delta}>2$. The integral is therefore at most
$2^{-15}+C/(64\delta)$. In the interior estimate retain
$\log(1+r)+r^2>3/2$. This proves \eqref{paper-wall-centering}.
Evaluating the same circle bound at $r=7/8$ in every support direction
gives an inscribed radius at least
$a[7/8-\sqrt{3T/4}]>2/3$ for $T\le1/24$; here
$-\log(1-r^2)-r^2=\log(64/15)-49/64<3/4$.
\end{proof}

\begin{lemma}[Mean-square control of the fan center]\label{IV:cor:rms-fan-center}
Normalize the area of the convex $n$-gon to $\pi a$, where $h=2\pi/n$
and $a=\sinc h$. At the conformal center $O$, let
$Y\ge\max_i(b_i/p_i-1)$, $u_0=\max_i q_i\le1.47$, and
$\max_i|\delta_i|\le d<\pi/2$. Suppose each sector's mean endpoint radius
differs from one by at most $Eh^2$. Write $\tau(v)=\tan(v)/v$,
$\tau(0)=1$, and $Z_s=\pi\tau(\pi u_0/n)\tau(d)/\sqrt{3.4}$.
Let $C$ bound the sum of the relative entropy, angle and skew terms in the sector estimate.
If the two ratios
\[
 \frac{2(1+Eh^2)^2 Z_s^2 u_0}{3an},\qquad
 \frac{(1+Y)(2\pi)^4}{864n\sinc(hu_0)}
\]
are at most one, then
\[
 \|1-\rho_O^2\|_2^2/h^4\le N^2=
 (1-a)^2/h^4+\frac{a^2(1+Y)^2}{180\sinc(hu_0)}
 +\frac{2a^2(1+Y)C}{(2\pi)^4 n}.
\]
\end{lemma}
\begin{proof}
The relative entropy, angle and skew contributions satisfy $C_e+C_q+C_s\le C$.
For $y_i=b_i/p_i-1\le Y$ the relative entropy gives
$\sum p_i y_i^2\le2(1+Y)C_e/n^5$. The sine rule gives
$|w_i|\le2(1+Eh^2)\tan(\pi q_i/n)\tan|\delta_i|$.
Since $p_i=q_i/n$ and $C_s\ge3.4n^3\sum_i q_i^2\delta_i^2$,
squaring and summing gives
$\sum_i p_i w_i^2\le4(1+Eh^2)^2Z_s^2u_0 C_s/n^6$.
Since $q^5-5q+4\le2(q^4-4q+3)$ on $[0,1.47]$,
$\sum p_iq_i^4\le1+2C_q/(4.8n^2)$. In \eqref{IV:eq:polar-fan-rms} the
coefficients of $C_s$ and $C_q$, divided by the coefficient of $C_e$, are bounded by the first and second ratios in the hypothesis, respectively. The claim follows. Equation~\eqref{IV:eq:moment-center-ball}
then gives the shift of at most $2N^2h^4/5$ from the area centroid under the hypotheses of that estimate.
\end{proof}

\section{Boundary integrals and chord inequalities}
\begin{lemma}[Two chord primitives]\label{IV:lem:two-chord-primitives}
On the triangle $x_1=c_0r$, $x_2=b_0rt$ with $c_0=\cos(h/2)$,
$b_0=\sin(h/2)$, take $w\in L^2(-1,1)$ with $\int_{-1}^1w(t)\,dt=0$ and set $W(t)=\int_{-1}^t w$. For $k>0$ one has
\[
 \begin{aligned}
 \left|\textstyle\int_{\rm chord} w(t)D_i\,dv\right|
 &\le h\,N_k(w)\|\nabla D_i\|,\\
 N_k(w)^2&=\frac{b_0c_0}{2kh^2}\int_{-1}^1
 [w^2+(b_0/c_0)^2(kW-tw)^2]\,dt .
 \end{aligned}
\]
In addition, for $K(X)=g(\sqrt X)/\sqrt X$ one has
\[
\begin{aligned}
 |\textstyle\int K_1x_2D_2|&\le h^{5/2}L_y\|\partial_2D_2\|,\\
 L_y^2&=\frac{c_0}{120}\left(\frac{2b_0}{h}\right)^5
 \left[\frac{A^2}{6}-\frac{AB}{4}+\frac{B^2}{10}\right],
 \qquad A=1.59,\quad B=.59(1+c_0^2)/2.
\end{aligned}\]
\end{lemma}
\begin{proof}
For $\psi=f(r)S(t)$ with $f(1)=1$, $f(0)=0$, write $\tau=b_0/c_0$,
$A=\|S\|^2$ and $B=\|tS'\|^2+\tau^{-2}\|S'\|^2$. The Cartesian derivatives are
$(f'S-tfS'/r)/c_0,fS'/(b_0r)$. Therefore
$\|\nabla\psi\|^2/\tau=A\int r(f')^2+B\int f^2/r-\Re\langle S,tS'\rangle$.
Completing $\int r(f'-kf/r)^2$ gives the minimum at $f=r^k$, $k^2=B/A$
when $A>0$. With $S=b_0W$ the rotated gradient has normal trace $w$ on the chord
and zero on the rays, so integration by parts and Cauchy--Schwarz give $N_k$.
On the fiber $|v|\le b_0r$ the primitive of $vK(x_1^2+v^2)$ vanishing at
both ends has magnitude $\frac12\int_{c_0^2r^2+v^2}^{r^2}K(X)\,dX$.
The Bessel series gives $K''>0,K(0)<1.59,K(1)=1$, hence $K(X)\le1.59-.59X$.
The primitive is at most $(b_0^2r^2-v^2)(1.59-Br^2)/2$ with $B$ as in the
statement. Integrating its square in $v$ gives at most
$4(b_0r)^5(1.59-Br^2)^2/15$. Integrate in $x_1=c_0r$ and apply
Cauchy--Schwarz to the fiberwise integration by parts.
\end{proof}

\begin{corollary}\label{IV:cor:optimal-polynomial-chord-lift}
For $n\ge10$ the two chord weights
\[
 w_N=.308(1-t^2)+.047(1-t^2)^2-.2304,\qquad
 w_T=.308t(t^2-3/5)
\]
have dual constants $\inf_{k>0}N_k(w_N)<.044$ and
$\inf_{k>0}N_k(w_T)<.0156$ in Lemma~\ref{IV:lem:two-chord-primitives}.
\end{corollary}
\begin{proof}
Put $A=\int w^2$, $B=\int t^2w^2$, $C=\int W^2$ and $x=\pi/n$.
Since $2\int twW=-C$, the preceding minimum gives
$4\inf_k N_k^2=\operatorname{sinc}^2 x[\sqrt{C(A+B\tan^2x)}+C\tan x/2]$.
Integrating the two polynomial weights using $\int_{-1}^1t^{2j}dt=2/(2j+1)$ gives
$B/A<.52$ and $\sqrt{AC}+C/6<.00762,.000948$, respectively.
For $x\le\pi/10$, the sine series gives $\operatorname{sinc}x\le1-.16x^2$, hence
$\operatorname{sinc}^{-4}x\ge1+.64x^2\ge1+(B/A)\tan^2x$ since $\cos x>.95$.
Thus $4\inf_k N_k^2\le\sqrt{AC}+C/6$ by $\tan x<1/3$.
\end{proof}

\begin{lemma}[Angular gain from regular-sector data]\label{IV:lem:middle-scalar-angular}
Let $n\ge10$, $h=2\pi/n$, $\theta=h/2$, $c=\cos\theta$,
$\sigma=\operatorname{sinc}\theta$ and $.353\le l\le q\le u\le1.459$,
where $l\le1\le u$. Put
$Q=q+2$, $F=q(q+1)^2$, $D_*=F+Q$. Use
$0\le s\le s_D$, $d\le z_d$, $z_t\le t_0$, $z\le Z$,
$M_D\le mhE_D$, $\lambda_0mh\le\kappa\le1/4$;
let $T_1,T_2$ bound the relative boundary norms of the two components.
Use the three cap bounds $(a,b,e)$ from the disk comparison. Allow a radial
error of $khF$ in units $h^4(q-1)^2$. If $R>0$ at $l,u$,
subtracting this error leaves an angular gain of at least
\begin{equation}\label{IV:eq:five-exact-angular}
\begin{aligned}
 G={}&sR/(1-\kappa)+dL_d+t_0\min(0,R_t)+h[(1/120-k)F-Q/90]\\
 &+h^3(aH_7+bH_5/240-eQ/12)-E_0-E_Q Q-E_*D_* .
\end{aligned}
\end{equation}
Here $R,R_t$ come from Lemma~\ref{IV:lem:exact-metric-angular},
$L_d=Q[1+h^2\operatorname{sinc}(hu)/(6\operatorname{sinc}h)]$,
and the affine error bounds the boundary and volume terms below,
including any negative contribution from the constant volume term. If this term
has a positive lower bound, it can be added to $G$.
For $C=l+q+u+2$ its divided difference satisfies
\begin{equation}\label{IV:eq:five-metric-curvature}
 [l,q,u]G/C\le-sM/(1-\kappa)+h^3
 \left(a\frac{[l,q,u]H_7}{C}+b/240\right)-h(k-1/120)-E_* ,
\end{equation}
where $M$ is the metric curvature in Lemma~\ref{IV:lem:exact-metric-angular}.
On a subinterval the same curvature bound holds with its three
divided-difference nodes, retaining the metric bounds from a window
containing those nodes and one.
\end{lemma}
\begin{proof}
Write $\rho=\lambda_0M_D/E_D$, $A_m=2T\omega-TH_A$.
The isotropic and mass terms give
$s_D[2T\omega-T(B_0^3/A_0)-\rho A_m]/(1-\rho)$.
This fraction decreases with $\rho$ since $TH_A<0<T(B_0^3/A_0)$,
and $R$ is concave. The cap term follows from
Proposition~\ref{IV:prop:expanded-cap-eighteen-nineteen}.
Since $s/(1-\kappa)\ge0$, concavity of $R+MD_*$ gives
$[l,q,u](sR/(1-\kappa))\le -sMC/(1-\kappa)$.
The curvature bound follows from concavity of $\min(0,R_t)$, $t_0\ge0$,
and $[l,q,u]D_*=[l,q,u]F=[l,q,u]H_5=C$.

Take lower bounds $I_0,\nu_0\ge0$ for the two moments of Proposition~\ref{IV:lem:projected-wall-pair}. On the reference chord put $p=\sin^2(\theta q)$, $p_1=\sin^2\theta$,
$Z_0=1-t^2$, $k(x)=g'(\sqrt{1-x})/x$, $j(x)=2K'(1-x)$.
The common factor in its two coefficients is
\[
 \omega(p-p_1)=\tfrac14h^2(q^5-q^3)C(q),\qquad
 C(q)=\frac{\operatorname{sinc}(hq)\operatorname{sinc}^2(\theta q)
 \operatorname{sinc}(\theta(q+1))\operatorname{sinc}(\theta(q-1))}
 {\operatorname{sinc}h\operatorname{sinc}^2\theta}.
\]
The remaining factors are $p_1Z_0(Z_0\overline{k'}-t^2\overline{j'})$
and $c vZ_0\overline{j'}$; bars average from $p_1 Z_0$ to $pZ_0$ and $v$ denotes chord height.
Differentiating the defining Bessel series gives
$k'(x)=\int_0^1 b[5K''(1-bx)+2(1-bx)K^{(3)}(1-bx)]\,db$,
$j'(x)=-2K''(1-x)$. Freeze $\overline{k'},\overline{j'}$ at $.355,-.308$. The normal and transverse polynomials are $.308Z_0+.047Z_0^2,-.308tZ_0$ ($Z_0=1-t^2$), with multiplier
$C_* = D_*C(q)+2[C(q)-C(1)]/(q-1)$. The sinc product $C(q)$ decreases on the aperture interval. Indeed the odd function $1/x-\cot x$ increases, so for $1/3\le q<1$ the factor $\operatorname{sinc}(hq)$ dominates the rising factor $\operatorname{sinc}(h(q-1)/2)$ in logarithmic derivative. For $q\ge1$ all numerator factors decrease. Hence $c_-D_*-2.208h^2\le C_*\le c_+D_*$ if $c_-\le C(u)\le C(l)\le c_+$, since $|C'|<1.104h^2$. Use $.70,1.11$ throughout, or the endpoints' tighter bounds. The lower multiplier is positive.

The two retained moments, the two centered lifts of Corollary~\ref{IV:cor:optimal-polynomial-chord-lift}, and the differentiated coefficient errors, respectively, bound the loss from the boundary integrals by
$(J_1-J_0)D_*+2.208h^2J_0/c_-+J_2$, where
\begin{equation}\label{IV:eq:angular-wall-charges}\begin{aligned}
 J_0&=c_-[(1-h^2/24)^2(.2303)I_0h^3/8
              +(1-h^2/8)(.1232)\nu_0h^2/2],\\
 J_1&=c_+\sqrt Z\sqrt{(.044h^3/8)^2+(.0156h^2/4)^2},\\
 J_2&=1.053\sqrt Z\{[c_+(.020\cdot.66+.012\cdot.18)D_*+.141h^2]T_1h^3/8\\
 &\hspace{43mm}+[c_+(.012)(.15)D_*+.15(.096)h^2]T_2h^2/2\}.
\end{aligned}\end{equation}
Indeed $p=\sin^2(\theta q)<.206$; on $[1-.206,1]$ the Bessel series gives $.1524<K''<.1599$, $-.0354<K'''<-.0336$, $|K^{(4)}|<.007$. Thus the averaged derivative errors are below $.020,.012$ and $|k''|<.115$, $|j''|<.082$. Differentiating the residuals times $C$ costs at most $2(1.104\cdot.020+1.11\cdot.115\cdot.748/2)h^2<.141h^2$ and $2(1.104\cdot.012+1.11\cdot.082\cdot.748/2)h^2<.096h^2$ because $|p'|<.748h^2$. The three angle norms are below $.66,.18,.15h$ by their squared polynomial integrals with $ds/dt<.518$; $1/c<1.053$.

For the varying volume loss use
\begin{equation}\label{IV:eq:angular-varying-volume}
 V=\tfrac{1.16}4\sqrt{mhZ}\,h^2[N_1(\eta Q+\epsilon)+N_2(\eta D_*+\epsilon)],\qquad
 N_1^2=2c^3\sin\theta/15,\quad N_2^2=2c\sin^3\theta/105 .
\end{equation}
Indeed write $K_q-K_1=(p-p_1)r^2(1-t^2)L$, where $L$ averages $-K'$ between the regular and transformed arguments. Thus $0<L<.68$, $0\le L'\le.05h^2q$ by $|K''|<.20$. After the factor $h^2r^2(1-t^2)/4$, the volume weights are $Q_i C_i L$, with $Q_i=q^3-q,q^5-q^3$ and, for $(j,k)=(3,1),(1,3)$,
\[
 C_i=A_0^j[\operatorname{sinc}(\theta q)/\sigma]^k
       \operatorname{sinc}(\theta(q+1))\operatorname{sinc}(\theta(q-1)).
\]
One has $C_i<(c\sigma)^{-3}<1.23$, $|(\log C_i)'|<.75$ by $|\tan x|\le1.27|x|$ for $|x|\le\theta u$, $|(\log\operatorname{sinc})'(x)|\le.373|x|$ for $|x|\le\theta(u+1)$. Since $Q_i(1)=0$, $Q_i'(1)=2$, this gives $|T(Q_i C_iL)|\le\eta TQ_i+\epsilon$ with $(\eta,\epsilon)=(.85,1.4)$ and $TQ_i=Q,D_*$. Integrating $r^7(1-t^2)^2$ against the two squared coordinates gives $N_i$. The contribution from the remaining constant volume term is $2(\lambda_0-\lambda)T(\omega-H_A)(M_D+\int U_2D_2)/h^4$.
\end{proof}

\begin{corollary}[One-sided chord remainder]\label{IV:cor:exact-two-curvature-moments}
Let $t\in(0,\pi)$, $S=\sin^2(t/2)$, $I=\operatorname{sinc}t$,
$\rho_0(s)=\cos(t/2)/\cos(ts)$ and
$D(e)=\int_{-1/2}^{1/2}\Phi((1+e)\rho_0(s))\,ds$. For
$0\le e\le\sqrt3-1$, set $m_k=S^k\int_0^1(1-u^2)^k\,du$,
$V_j=\sum_{i=0}^{j}\binom ji e^i\,2/[(i+1)(i+2)]$ and
\[
 M_1=1-V_2(1-m_1),\qquad
 M_2=1-2V_2(1-m_1)+V_4(1-2m_1+m_2).
\]
For every $a>0$ one has
\[
 -\frac2I\int_0^1(1-v)D''(ve)\,dv
 \ge 2-L+c_1 M_1+c_2 M_2-\frac{|c_3|m_3}{1-S/2}.
\]
With $B=2e+e^2$, $k=(1+B)S$, if $W(1-p)\ge L+p$ on
$-B\le p\le S$, the second derivative with respect to $w$ at $w=0$ in Lemma~\ref{IV:lem:exact-skew-chord} satisfies
\[
 -D_{ww}(e,0)/I\ge Q+Ak-D_1 B-D_2(8k^2/105-4kB/15+B^2/3)
             -\frac{4|c_3|S^3}{315(1-S/2)},
\]
where $Q=1/6+L/4$, $A=c_1/30+4/15$, $D_1=c_1/12+1/3$ and $D_2=|c_2|/4$.
\end{corollary}
\begin{proof}
Use the exact chord identity of Lemma~\ref{IV:lem:weighted-disk-curvature}
with $u$ uniform on $[0,1]$ and independently $v$ of density $2(1-v)$.
Set $p=S(1-u^2)$ and $x=1-(1+ve)^2(1-p)$.
Then $-2\le x\le p$, $\mathbb E(1+ve)^j=V_j$ and the two moments
$\mathbb Ex,\mathbb Ex^2$ are exactly $M_1,M_2$.
Average Lemma~\ref{IV:lem:all-side-chord-remainder} and use
$\mathbb E(x_+)^3\le m_3$.
For the skew identity use $p=k(1-u^2)-B\le S(1-u^2)$.
The uniform moments of $u^2p,(1-u^2)p,u^2p^2$ are
$2k/15-B/3,8k/15-2B/3,8k^2/105-4kB/15+B^2/3$;
the moment of $u^2(p_+)^3$ is at most $16S^3/315$.
\end{proof}

\subsection{Two signed projections of the boundary flux}\label{mid-wall}\label{mid-wall10}\label{mid-wall11}\label{mid-wall12}\label{mid-wall13}\label{mid-wall14}
\begin{lemma}[Second differences and derivatives of the boundary flux]\label{IV:lem:signed-wall-shape}
For the chord flux $F$ from Lemma~\ref{IV:thm:radial} let $W_e,W_w$ be upper bounds for
\[
 \frac2{qh}\int_{\rm chord} \partial_e^2 F(q,e,0)\cdot D\,dv_{\rm wall},\qquad
 \frac2{qh}\int_{\rm chord}\frac{F(q,e,w)+F(q,e,-w)-2F(q,e,0)}{w^2}\cdot D\,dv_{\rm wall},
\]
respectively, with the continuous value at $w=0$. For the parameter ranges of Section~\ref{mid-geometry}, one may take $W_e=.255$ and $W_w=.014$ for every $n\ge12$. For ten and eleven sides the bounds $W_e,W_w$ on each aperture band are given in the two tables of radial estimates in Section~\ref{mid-radial}.
\end{lemma}
\begin{proof}
Let $a=(1+e)^2$, $v=\sin^2(hq/2)$, $c=1-v$, $k=h\cot(h/2)$, $t\in[-1,1]$, $y=t^2$, $X=a(c+vy)$ and $X_m=a(1-v/2)$. Write $b=\sin(h/2)$, $c_h=\cos(h/2)$, $\kappa=\sinc(hq)/\sinc h$. With $K(X)=g(\sqrt X)/\sqrt X$, differentiation of the flux in Lemma~\ref{IV:thm:radial} gives the following two components at zero skew, after division by $q\kappa$ and multiplication of the transverse component by $h$,
\[
\begin{split}
 N(y,X)&=-K/2+a(-1+5v/2+cy/2)K'
      +a^2[vc+(3+c^2)y]K''+2a^3cyK''',\\
 kat B(y,X)&=kat[(1-2v)K'+a(2+v^2+vcy)K''+2va^2yK'''] .
\end{split}
\]
Let $N_m=N(1/3,X_m)$ and $B_m=B(1/2,X_m)$. The signed moments are $-\int D_1\,dv_{\rm wall}=Ih^3$ and $\int tD_2\,dv_{\rm wall}=\nu h^4/b$, both positive. Use $dv_{\rm wall}/d\theta\le1/c_h$ and let $u_i$ be the upper bounds for the $L^2(d\theta/h)$ norms of $t^2-1/3,t^2-1/2,t(t^2-1/2),t$, in that order. Their squares are bounded by $(\tan(h/2)/h)(8/45,7/30,11/210,2/3)$ by direct integration. Denote interval suprema of absolute values by $\|\cdot\|_\infty$.
If $f_4,g_4$ bound the two fourth skew derivatives of the reduced flux and $|w|\le w_0$, the two projection errors are at most
\[
\begin{split}
 N_*&=u_1\|N_y\|_\infty+u_2av\|N_X\|_\infty+f_4w_0^2/12,\\
 B_*&=2a\{u_3(\|B_y\|_\infty+av\|B_X\|_\infty)+u_4(B_m)_+\}+g_4w_0^2/12 .
\end{split}
\]
Indeed $X-X_m=av(t^2-1/2)$, $k\le2$, and the even Taylor difference contributes $w_0^2/12$. The transverse projection with negative coefficient can be dropped. By parity, Cauchy--Schwarz and the trace bounds $T_1,T_2$ for the regular sector,
\begin{equation}\label{IV:eq:wall-parameter}
 W_w\le 2\kappa\left[I_+h^2(-N_m)_++\frac{\sqrt Z}{c_h}
 (T_1h^2N_*+T_2h B_*)\right] .
\end{equation}
For the mean variation write $A=2K+10XK'+4X^2K''$, $B_1=A'$, $F=A+2XB_1$, $V=a(1-v)$ and $P_m=F(X_m)-av B_1(X_m)$. The two components are $A+2VB_1,2kavt B_1$. Their mean error bounds are $N_e=u_2av\|B_1+2VB_1'\|_\infty$ and $B_e=4u_3(av)^2\|B_1'\|_\infty$. The same projection gives
\[
 W_e\le2\kappa\left[I_+h^2(-P_m)_++\frac{\sqrt Z}{c_h}
 (T_1h^2N_e+T_2h B_e)\right]+4\kappa kav B_1(X_m)\nu h^2/b .
\]
Here $B_1(X_m)<0$ on the indicated intervals. For $n\ge12$ the signed term may be omitted. Positivity of $\nu$ there follows from \eqref{op-torque} with $s_D>.0096$, $h<.524$, $I<.044$ and the trace bounds below.

The following Bessel recurrence and endpoint evaluations apply at every side count. In the bounds on $|N_X|,|B_X|$ use $y=0,1$ since the derivatives are affine in $y$ at fixed $X$.

Put $r=1+\sup e$, $s=(1+\inf e)^2$, $v_*=\sup v$. The displaced squared argument is in $[s(1-v_*)-rw_0,r^2+rw_0+w_0^2/4]$. For an increment $z$, the nonnegative coefficient bounds for the factors $P,V,W$ of the boundary flux and for the increment of its argument are
\[
\begin{split}
 p&=r^2+w_0z/2+z^2/4,\quad D=(r+w_0/2)z+z^2/4,\\
 v_b&=r^2+rw_0/2+v_*w_0^2/4+(r+v_*w_0)z/2+v_*z^2/4,\\
 w_b&=v_*r^2+rw_0/2+w_0^2/4+(r+w_0)z/2+z^2/4 .
\end{split}
\]
With $J_i=\sum_{j=0}^4\|K^{(i+j)}\|_\infty D^j/j!$, take $f_4=24[z^4]p(J_0+2v_bJ_1)$ and $g_4=96[z^4]pw_bJ_1$.

For $n\ge12$ one pair of intervals bounds the flux on the entire chord. Put
$H=.524$, $u=1.403$, $E=.044$, $b_w=.202$, and take $0<h\le H$,
$|e|\le E$, $|w|\le b_wh^2$ from the geometric estimates.
The common spectral bounds give $Z=.015$, $T_1=.65$, $T_2=.50$ and $I_+=.044$.
For the last bound insert $d_->.560$, $m<.048$ and $\lambda_0<3.55$ in
\eqref{op-coordinate}: the resulting quadratic upper bound in $X=h^2$
has positive linear and quadratic coefficients and is less than $.044$ at $X=H^2$.
Set
\[
 a\in[(1-E)^2,(1+E)^2],\quad v\in[0,\sin^2(Hu/2)],\quad
 X\in[a_-(1-v_+),a_+],\quad X_m\in[a_-(1-v_+/2),a_+] .
\]
The Bessel equation $XK''+2K'+\lambda K/4=0$ reduces the two cancellations in the mean curvature to
\[
 \begin{split}
 B_1&=-\lambda(3K/2+XK'),\qquad B'_1=\lambda^2K/4-\lambda K'/2,\\
 F&=(2-4\lambda X)K+2X(1-\lambda X)K',\\
 F'&=\lambda[(\lambda X/2-9/2)K-4XK'],\\
 B_1+2VB'_1&=\lambda[(\lambda V/2-3/2)K-(X+V)K'] .
 \end{split}
\]
All lower and upper endpoints can now be evaluated without subdividing $a,v$.
Here are the intermediate bounds needed in the two projections:
\begin{equation}\label{eq:wall-endpoints}
\begin{array}{c|rrrr|rrrr}
 &(-P_m)_+&N_e&B_e&(-N_m)_+&N_*&B_*&f_4&g_4\\\hline
 \text{upper bound}&9.37&.256&.054&.187&.168&.0154&4.38&8.24
\end{array}
\end{equation}
To reproduce these entries let $S_j(x,r)=\sum_{i=0}^r(-\lambda x/4)^i/[i!(i+j+1)!]$.
For $(-1)^jK^{(j)}$ on $[l,r]$ use
$(\lambda/4)^j S_j(r,3)/S_0(1,4)$ below and
$(\lambda/4)^j S_j(l,4)/S_0(1,3)$ above, with $3.3899<\lambda<3.39$
inserted outwards. These bounds follow by alternation and monotonicity of the signed derivative.
In $-P_m$ one can drop the nonpositive $avB_1$; for an interval $[l,r]$
of $X_m$, use $F(X_m)\ge F(r)-(r-l)(\sup F')_+$.
The other entries follow by the displayed shape formulas, evaluating $N_X,B_X$ at $y=0,1$.
Use $w_0=b_wH^2$ in the fourth-derivative term and the endpoint
$\tan(H/2)/H$ in the four norm factors; the latter increases with $H$.
For example, the boundary contribution for the mean displacement is at most
\[
 \frac{2}{1-H^2/6}\left[.044H^2(9.37)+
 \frac{\sqrt{.015}}{1-H^2/8}(.65H^2\cdot.256+.50H\cdot.054)\right]<.255.
\]
Replacing the three shape entries by $.187,.168,.0154$ gives $W_w<.014$.
Here $H$ bounds only increasing positive factors; the four geometric
intervals include their disk-limit endpoints.

For $n=10,11$ use each aperture interval $[l,u]$ in the radial tables. Adjacent rows with the same boundary estimates can be treated together. Keep $v$ in $[\sin^2(hl/2),\sin^2(hu/2)]$. Use the two radius ranges $[(1+\inf e)^2,1]$ and $[1,(1+\sup e)^2]$ for $a$. The fourth-derivative estimate uses the full bound on $|w|$ for the aperture band. In the mean projection retain both signed terms: with $X_m=a(1-v/2)$ use $P_m=F(X_m)-avB_1(X_m)$ and the lower bound for $\nu$ from Section~\ref{mid-response}. One can bound $\kappa$ above by $\sinc(hl)/\sinc h$ and below by $\sinc(hu)/\sinc h$. Because $B_1<0$, the term containing $\nu$ is bounded above by
\[
 4\frac{\sinc(hu)}{\sinc h}\, h\cot(h/2)\,a_-v_- (\sup B_1)\,\nu_-\,h^2/\sin(h/2).
\]
For each radius range use the same Bessel differential identities and fourth-order alternating sums just displayed. For $F$ use the lower bound $F(r)-(r-l')(\sup F')_+$ on the $X_m$ interval $[l',r]$. Take the larger of the two boundary estimates.
For example, for ten sides on $1.30\le q\le1.35$ one has
$.15772<v<.16935$, $.94530<\kappa<.95400$ and $w_0<.06751$.
In increasing order of radius the two $X_m$ intervals lie in
$[.8035,.9212]$ and $[.9153,1.0429]$. The three mean shapes
$((-P_m)_+,N_e,B_e)$ are bounded above by
$(8.408,.1850,.07928)$ and $(8.779,.27642,.09718)$.
For the skew shapes $((-N_m)_+,N_*,B_*)$ the two triples are
$(.21145,.12837,.01109)$ and $(.12870,.16458,.06243)$.
The signed moment uses $B_1<-3.6292,-3.1607$, respectively.
These substitutions give $W_e<.298$ and $W_w<.0203$
on the entire band.
\end{proof}

\subsection{Second derivative of the volume term in the mean radius}\label{IV:lem:radial-smooth-tensor}
On the mean-radius line the volume is
\[2(\lambda-\lambda_0)q\kappa R^4\int K(R^2|B_0x|^2)B_0x\cdot B_0D,\]
with $B_0$ the unskewed map at $R=1$.
Two derivatives replace $R^4K$ by $R^2[(12-\lambda X)K+10XK']$ by the
Bessel equation. For $0\le X\le1.18$ the bracket lies between $-8$
and $19.08-13.6X$ by $K\le1.59-.46X$ and $-.675<K'<-.46$.
Its product with $\sqrt X$ has absolute value below $8.8$.

\begin{lemma}[Angular weights and curvature]\label{IV:lem:exact-metric-angular}\label{IV:lem:common-metric-curvature}
Let $0<l\le1\le u$, $0<h\le\pi/3$, $hu\le\arccos(19/272)$,
$\theta=h/2$, $A=\cos(\theta q)/\cos\theta$, $B=\sin(\theta q)/\sin\theta$,
$d=AB$, $\omega=dB^2$, $H_A=dA^2$, $Q=q+2$, $D_*=q(q+1)^2+Q$.
Write $Tf=(f(q)-f(1)-(q-1)f'(1))/(q-1)^2$, with continuous extension.
For $\lambda_0 M_D\le\kappa E_D$, $0\le\kappa\le1/4$, put
\[
 R=2T\omega-T(B^3/A)-\kappa\max(2T\omega-TH_A,T\omega,0),\quad
 R_d=T(B^3/A-d),\quad R_t=T(B^3/A-A^3/B).
\]
After applying $T$, the gain from the pure term, in units $h^4$, is at least
$Rz+R_d z_d+R_tz_t$, where $h^4(z,z_d,z_t)=(E_D,
\|\partial_1D_1\|^2+\|\partial_2D_2\|^2,\|\partial_2D_1\|^2)$.
On $[l,u]$, $R$ is concave, $R_t$ increases and
$R_d\ge Q[1+h^2\operatorname{sinc}(hu)/(6\operatorname{sinc}h)]$.

If $h\le\pi/5$ and $u\le3/2$, set $a=2(1-\kappa)$,
$r=\tan^2\theta$. Then $R=T(a\omega+\kappa H_A-B^3/A)$, both
$R+MD_*$ and $R_t$ are concave, and $0<T(\omega-H_A)\le V$, with
\[
 M=\frac{h^2\operatorname{sinc}(hu)}{60}
      [16(a-\kappa r)\cos(hu)-(a+\kappa r)],\qquad
 V=\frac{(u+2)(1+5r/3)}{\cos\theta\operatorname{sinc}^3\theta}.
\]
\end{lemma}
\begin{proof}
Subtracting the four gradient weights of the pure term
$d(1,A^2/B^2,B^2/A^2,1)-2\omega(1,1,1,1)$ from zero gives the three energy coefficients; the maximum in $R$ bounds both mass
coefficients. Throughout use
$(Tf)^{(j)}=\int_0^1 (1-v)v^j f^{(j+2)}(1+v(q-1))\,dv$.
Put $b=\sin\theta,c=\cos\theta,r=\tan^2\theta,x=\theta q$. Up to a positive factor
$B^3/A=\tan x-\sin(2x)/2$, whose Taylor coefficients from degree
three on are nonnegative. Indeed $t_j=[x^{2j+1}]\tan x$ and
$d_j=2^{2j}/(2j+1)!$ satisfy $t_2=d_2$ and
$(2j+3)t_{j+1}\ge2t_j$ by the tangent equation, whereas
$(2j+3)d_{j+1}=4d_j/(2j+2)$. Thus $T(B^3/A)\ge Q$, since its
coefficients in $q$ sum to one and $Tq^k\ge Q$ for $k\ge3$.
Also $-d''\ge h^2q\operatorname{sinc}(hu)/\operatorname{sinc}h$,
which yields the claim for $R_d$. Up to a positive factor
$A^3/B=\cot x-\sin(2x)/2$, whose third derivative is
$-2\csc^2x(1+3\cot^2x)+4\cos(2x)<0$ for $x<\pi/4$. Hence $R_t$ increases.

The sine weights are $8cb^3(\omega,H_A)=
(2\sin(2x)-\sin(4x),r[2\sin(2x)+\sin(4x)])$. For each
$(a',k)=(2(1-\kappa),\kappa),(2-\kappa,0),(2,0)$ they give
\[
 -(a'\omega+kH_A)^{(4)}=
 \frac{h^4\sin(hq)}{4cb^3}[16(a'-kr)\cos(hq)-(a'+kr)]\ge0
\]
because $r\le1/3$ and $\cos(hq)\ge19/272$. Subtracting
$B^3/A$ preserves the sign, so the minimum of the three Taylor
remainders is concave. In the narrower range $4\cos(hq)>1$ gives
$\omega''>0>H_A''$, selecting the first branch. Its negative fourth
derivative is at least $120M q$, by $b\le h/2$ and
$\sin(hq)\ge hq\operatorname{sinc}(hu)$. Since
$\int_0^1 v^2(1-v)[1+v(q-1)]\,dv=(3q+2)/60$, this gives
$-R''\ge M D_*''$. Similarly
$0<(\omega-H_A)''=h^2\sin(hq)[4(1+r)\cos(hq)-(1-r)]/(4cb^3)
\le2q(3+5r)/(c\operatorname{sinc}^3\theta)$, giving $V$.

For the remaining concavity put $t=\tan x$,
$L(t)=t^3(7+8t^2+3t^4)/(1+t^2)$. The fourth derivatives of
$\tan x-\sin(2x)/2$ and $\cot x-\sin(2x)/2$ are $8L(t),8L(1/t)$.
Their ratio is at most $7t^8/3\le\tan^4\theta$: indeed
$\tan^2(3\theta/2)/\tan\theta$ increases by the decreasing sinc function
and at $\pi/10$ is below $.800<(3/7)^{1/4}$. Restoring the two
normalizations gives $(B^3/A-A^3/B)^{(4)}\le0$, hence $R_t''\le0$.
\end{proof}

\begin{lemma}[A diagonal cubic on the regular sector]\label{IV:lem:common-diagonal-cubic}
Let $b=\sin(h/2)$, $c=\cos(h/2)$, $J=bc$, $\lambda_0\ge\lambda$,
$r^2=|x|^2$, and put $z_d=h^{-4}(\|\partial_1D_1\|^2+\|\partial_2D_2\|^2)$.
For $f(t)=t(1-t^2)$, $\Psi=(cf(x_1/c),bf(x_2/b))$ and
$B=\int(x\cdot\Psi)[1+(1-r^2)/2+.076(1-r^2)^2]$,
\[
 \begin{aligned}
 h^2\sqrt{8Jz_d/5}&\ge(\lambda_0-\lambda)B+4cb^3/15+.30(16cb^5/105)\\
 &\quad-\lambda_0[M_D/3+\sqrt{M_D}\sqrt{\|\Psi\|^2-2B/3+h(1-1/3.39)/18}].
 \end{aligned}
\]
\end{lemma}
\begin{proof}
For any nonconstant equivariant test $\Psi$ of diagonal gradient, put $N=\|\nabla\Psi\|$, $j=\nabla\Psi/N$ and choose any $H_U\ge\|U\|^2$. One may take the smaller of
$h(1-1/3.39)/2$ and $(\sin h-h/3.39)/2+3h^7/768$ for $n\ge10$:
the disk-sector integral equals $h(1-1/\lambda)/2$ and clipping the chord
subtracts the cap area minus $\int_{\rm cap}(1-g^2)\le3h^7/768$.
If $\alpha\ge0$, $0\le B\le\int U\cdot\Psi$, $M\ge M_D$, $\Lambda\ge\lambda_0$ and $0\le d_-\le(\lambda_0-\lambda)/h^2$, \eqref{op-weak-correction} gives
\begin{equation}\label{IV:eq:common-diagonal-test}
 \langle\nabla D/h^2,j\rangle\ge L_\alpha:=
 \frac{d_-h^2B-\int_{\rm chord}\partial_\nu U\cdot\Psi
 -\Lambda[\alpha M+\sqrt{M(\|\Psi\|^2-2\alpha B+\alpha^2 H_U)}]}{h^2N}.
\end{equation}
Indeed $\int D\cdot U=-M_D$, so pairing $D$ with $\Psi-\alpha U$ bounds the mass term. The stated cubic has radial ray traces and $\Psi_1=0$ on the chord. For $v=1-r^2$, the positive series for $K(1-v)$ and $K''(1)=(4-\lambda)/4$ give $\int U\cdot\Psi\ge B$ and $-2K'(1-v)\ge1+.30v$, hence the stated flux bound. In $(X,Y)=(x_1/c,x_2/b)$, $\int X^{2i}Y^{2j}=J/[(2j+1)(i+j+1)]$. Thus $N^2=8J/5$, $\|\Psi\|^2=J(c^2/12+29b^2/420)$, and $\sqrt{z_d}\ge L_{1/3}$ with the exact mass and spectral parameters.
\end{proof}

\section{Lower bounds for the angular deficit}\label{mid-angular}
\label{III:sec:ang10}\label{III:sec:ang11}\label{III:sec:ang12}\label{III:sec:ang13}\label{III:sec:ang14}\label{mid-cap}
\begin{proposition}[Angular gain at ten and eleven sides]\label{IV:thm:uniform-twenty-plus}\label{IV:thm:coupled-twentytwo}\label{IV:sec:twenty-angular}
Let $n=10,11$, $h=2\pi/n$, and use the parameter ranges of Section~\ref{mid-geometry}. The angular deficit $\Delta$ is bounded below by $k qh^5(q^2-1)^2$ with the two coefficients $k$ given in Section~\ref{mid-summation}. At $q=1$ the deficit vanishes; Proposition~\ref{prop:angular-from-twelve} treats every $n\ge12$.
\end{proposition}
\begin{proof}
Put $X=h^2$, $\theta=h/2$, $c=\cos\theta$, $b=\sin\theta$, $J=cb$ and $S=\sinc h$.
We apply the Taylor quotient of Lemma~\ref{IV:lem:exact-metric-angular}. Write $l,u$ for the aperture bounds and $F=q(q+1)^2$.
\paragraph{Two projections on the regular sector.}
In sector coordinates put $Y=h^{-2}(\partial_1D_1,\partial_2D_2)$,
viewed as a diagonal matrix field. Write $f(t)=t(1-t^2)$,
$\Psi=(cf(x_1/c),bf(x_2/b))$, $N=\sqrt{8J/5}$, and use the unit
fields $j=\nabla\Psi/N$, $j_0=\mathrm{Id}/\sqrt{2J}$. They are
orthogonal: $\Psi$ is tangential on the rays and its normal component
vanishes on the chord. For $p=1-|x|^2$ put
\[
 B_0=\int(x\cdot\Psi)(1+p/2+.076p^2),\qquad
 H=(\sin h-h/3.39)/2+3h^7/768,\qquad M_0=mZh^5 .
\]
Here $m,Z,\Lambda,d_-,d_+,s_-,\kappa,I_-,\nu_-$ are the
bounds for the regular sector in Section~\ref{mid-response}, with $\kappa$ the
upper bound on $\lambda_0M_D/E_D$.
Lemma~\ref{IV:lem:common-diagonal-cubic} with $\alpha=1/3$ and the
cap estimate $\|U\|^2\le H$ give $\langle Y,j\rangle\ge A_0$, where
\[
 A_0=\frac{d_-h^2 B_0+cb^3(4/15+.30 b^2\,16/105)
 -\Lambda[M_0/3+\sqrt{M_0(\|\Psi\|^2-2B_0/3+H/9)}]}{h^2N}.
\]
Indeed the disk mass is $h(1-1/\lambda)/2$, the removed cap area
is $(h-\sin h)/2$, and its lost mass relative to area is at most
$3h^7/768$. The chord flux identity gives
$-\langle Y,j_0\rangle=hI/\sqrt{2J}\ge B_1=hI_-/\sqrt{2J}$.
Hence $z_d=\|Y\|^2\ge D_0=A_0^2+B_1^2$ in each row, where
$A_0,B_1>0$. For ten and eleven sides this gives $D_0>.00499,.00476$,
respectively, for use in the mixed slope.

The retained volume term in \eqref{IV:eq:five-exact-angular} is
$a_q\langle w,Y\rangle$ after discarding its favorable mass, with
$a_q=2d_* T(\omega-H_A)$ and
$w_{22}(cr,v)=\frac12\int_{c^2r^2+v^2}^{r^2}K(t)\,dt$, $w_{11}=0$.
This is integration by parts on $|v|\le br$, where the primitive
vanishes at the two endpoints. The disk inequalities
$1+p/2\le K(|x|^2)\le1.59-.59|x|^2$ and
Lemma~\ref{IV:lem:two-chord-primitives} give
\[
 w_j=\frac{\int x_2\Psi_2(1+p/2)}N\le\langle w,j\rangle,\quad
 0\le\langle w,j_0\rangle\le w_0=
 \frac{\int x_2^2(1.59-.59|x|^2)}{\sqrt{2J}},\quad
 \|w\|\le W=h^{5/2}L_y .
\]
All integrands in these constants are polynomials;
$\int(x_1/c)^{2i}(x_2/b)^{2j}=J/[(2j+1)(i+j+1)]$ and
$\|\Psi\|^2=J(c^2/12+29b^2/420)$.

\paragraph{A single chord on each aperture half.}
Take $[v_0,v_1]=[l,1]$ or $[1,u]$, put $r=\tan^2\theta$ and
$D_L=1+X\sinc(hu)/(6S)$, and let
\[
 P(p)=\frac{h^2}{16cb^3}[-2(1-r)\sin(hp)+4(1+r)\sin(2hp)].
\]
Here $P=(\omega-H_A)''/2$ is concave since $16(1+r)\cos(hp)>1-r$.
In the integral Taylor quotient the point $1+t(q-1)$ has density
$2(1-t)\,dt$ and $\mathbb Et=1/3$. The chord and tangent bounds for $P$
therefore bound $a_q$ by $a_-$ and $a_+$, where
\[
 \begin{aligned}
 a_-&=\tfrac23d_-[2P(1)+\min(P(v_0),P(v_1))],\\
 a_+&=2d_+\big[P(1)+\tfrac13\max\{P'(1)(v_0-1),P'(1)(v_1-1)\}\big].
 \end{aligned}
\]
Put $L_*=(2+v_0)D_L$. In all the indicated intervals $a_->0$
and $2L_*B_1>a_+w_0$. Decompose $Y$ in $j,j_0$ and their
orthogonal complement. The two scalar quadratics increase on their
respective coefficient ranges $[A_0,\infty)$ and $[B_1,\infty)$,
after taking the minus sign of the identity coefficient into account.
Completing the square in the complement therefore gives
\[
 L_d\|Y\|^2+a_q\langle w,Y\rangle\ge L_dD_0+C_0,\qquad
 C_0=a_-w_j A_0-a_+w_0B_1-(a_+W)^2/(4L_*),
 \qquad L_d=(q+2)D_L .
\]
Use this bound in \eqref{IV:eq:five-exact-angular}: set $d=D_0$,
$t_0=Z-D_0$ and add $C_0$. To collect the boundary and volume terms in \eqref{IV:eq:angular-wall-charges}
and \eqref{IV:eq:angular-varying-volume}, write $J_2=J_{20}+J_{21}D_*$ and
$V_0=(1.16/4)\sqrt{mhZ}\,h^2$. Taking
$c_-=C(u)$, $c_+=C(l)$ and $(\eta,\epsilon)=(.85,1.4)$ gives
\[
\begin{aligned}
 E_0&=2.208h^2J_0/c_-+J_{20}+\epsilon V_0(N_1+N_2),\\
 E_Q&=\eta V_0N_1,\qquad E_*=J_1-J_0+J_{21}+\eta V_0N_2 .
\end{aligned}
\]
The disk cap triple is $(k_7,\ell_-/X,(\ell_+-X/3)/X^2)$, where
$\ell_-\le-L\le\ell_+$ are given by \eqref{common-cap-height}. In order of increasing $n$ one
can take
\[
 k_7=.00328,\ .00335
\]
by the three-moment formula of Proposition~\ref{IV:prop:expanded-cap-eighteen-nineteen}.
The mass identity and the cap coefficients used here involve only
quantities defined on the regular sector.

Write $F=q(q+1)^2$ and $Q=q+2$. Let $k$ be the proposed angular coefficient on the interval in question. After subtracting $khF$ from $\Delta/[h^4(q-1)^2]$, the remainder in units of $10^{-3}$ is at least
\[
 g(q)=aR+b_0\min(0,R_t)+cQ+fF+pH_7+r_5H_5+v .
\]
Set $e_*=(\ell_+-X/3)/X^2$. Before rounding, the seven coefficients are
\[
 \begin{aligned}
 (a,b_0,c)/1000
 &=\left(\frac{s_-}{1-\kappa},\ Z-D_0,\
 D_0D_L-\frac{h}{90}-\frac{h^3 e_*}{12}-E_Q-E_*\right),\\
 (f,p,r_5,v)/1000
 &=\left(h(1/120-k)-E_*,\ h^3k_7,\ \frac{h^3\ell_-}{240X},\ C_0-E_0\right).
 \end{aligned}
\]
Round $b_0$ up and every other entry down. The same fixed coefficients
then apply on the corresponding aperture interval. The last three columns below
give a lower bound at each endpoint and an upper bound on the
interpolation error bound $P_*$ defined below.
\[
\begin{array}{c|c|rrrrrrr|rrr}
 n&q&a&b_0&c&f&p&r_5&v&g(v_0)>&g(v_1)>&P_*<\\\hline
10&[.426,1]&14.6&11.1&-3.59&-5.79&.813&.361&-.866&9.28&.253&6.32\\
  &[1,1.413]&14.6&11.1&-3.59&-4.28&.813&.361&-.73&6.42&.22&5.81\\
11&[.463,1]&14.1&10.5&-2.67&-5.70&.624&.269&-.57&14.51&2.88&2.65\\
  &[1,1.391]&14.1&10.5&-2.67&-4.22&.624&.269&-.46&8.91&1.14&3.13
\end{array}
\]
For the endpoint substitutions, $R$ at $(l,1,u)$ exceeds
$(1.19,.984,.543)$ and $(1.40,1.31,1.02)$ in the two rows of side
counts. The corresponding lower bounds for $\min(0,R_t)$ are
$(-.271,0,0)$ and $(0,0,0)$. For instance at one take
$T(A^iB^j)=\theta^2[(j\cot\theta-i\tan\theta)^2-i\sec^2\theta-j\csc^2\theta]/2$;
at either outside endpoint subtract the tangent of slope
$\theta(j\cot\theta-i\tan\theta)$ before dividing by $(q-1)^2$.
The trigonometric arguments in these endpoint expressions are below one.
The alternating sine and cosine polynomials through terms 8 and 9 suffice.
For ten sides the first row gives
$g(1)\ge14.6(.984)-3(3.59)-4(5.79)+21(.813)+10(.361)-.866=.2534$.
At $q=1$ the table values refer to the continuous normalized quotient.
They give regular angular coefficients above $.01730,.01919$: add the narrow
endpoint value times $10^{-3}/(4h)$ to the coefficient $k$ on the narrow interval.

To check the whole interval, $R$ and $\min(0,R_t)$ are concave by
Lemma~\ref{IV:lem:exact-metric-angular}. Also $H_5-F$ is affine
and $H_7''\le12F''$, since $20q^3+24q^2-54q-40$ is convex and
negative at $0$ and $1.42$. Consequently the divided difference of $g$
at $v_0,q,v_1$ is bounded above by
$(12p+r_5+f)(v_0+q+v_1+2)$. For $t=(q-v_0)/(v_1-v_0)$ let
\[
 P_*=(v_1-v_0)^2(v_0+2v_1+2)(12p+r_5+f)_+ .
\]
Then $g(q)\ge(1-t)g(v_0)+t g(v_1)-t(1-t)P_*>0$:
in every row $g(v_0)>P_*$ and $g(v_1)>0$.

\end{proof}

\begin{proposition}[Angular deficit for twelve or more sides]\label{prop:angular-from-twelve}
Write $h=2\pi/n$, $l=.475$, $u=1.403$, $Q=q+2$ and $F=q(1+q)^2$.
For every $n\ge12$ and $l\le q\le u$,
\[
 \Delta(q)\ge h^4(q-1)^2\bigl(.0123h F+.002\bigr),
\]
strictly for $q\ne1$. In these units the leading positive term as $h\to0$
is bounded below by $.0105(.79)Q$, uniformly in $q$. The bound on the radial error
is of order $h$ and the remaining boundary and nonconstant volume errors are of
order $h^2$.
\end{proposition}

\begin{proof}
Set $H=.524$, $X=h^2$, $D_*=Q+F$ and
$H_j=\sum_{i=0}^{j-2}(j-1-i)q^i$. The regular-sector bounds of
Section~\ref{mid-response} imply
\[
 s_D>\alpha_D+h/200,\quad z<.015,\quad m<.049,\quad
 (T_1,T_2)<(.66,.50),\quad 2d_*<1.16,\quad
 \langle Y,j\rangle>\sqrt{.0072h}.
\]
The signed cubic estimate in Proposition~\ref{continuous-response} gives the last inequality.
We can take $\kappa=.092$ since $3.546(.049)H<.092$.
Only these bounds for the correction $D$ enter the proof.

\paragraph{The positive term.}
For $\theta=h/2$, $c=\cos\theta$, let $g(x)=\tan x-\sin(2x)/2$.
The positive coefficients from degree three in the proof of
Lemma~\ref{IV:lem:exact-metric-angular} imply that $g(x)\ge x^3$
and $g''(x)/(6x)$ increases. At $x=.368>Hu/2$ its value is below
$1.01$, using $g''(x)=2\sin x/\cos^3 x+2\sin(2x)$. Thus
$T(B^3/A)\le1.01Q$. On the same interval the sine formulas for the
metric weights give
\[
 T H_A\ge-\frac{5\tan^2(H/2)}{3\cos(H/2)\sinc^3(H/2)}Q>-.13Q,
 \qquad
 T\omega\ge (1+X/4)(Q-XH_5/4).
\]
Indeed $\omega''(p)\ge(6p-5Xp^3)/(c\sinc^3\theta)$, by the quadratic
lower bounds on cosine and sinc. Also $1/(c\sinc^3\theta)\ge1+X/4$
by the first two positive secant and reciprocal-sinc coefficients.
Applying $T$ gives the displayed bound. Consequently
\[
 R\ge .79Q+X(.45Q-.49H_5)>0 .
\]
For the sign, $H_5/Q=q^2+3-2/(q+2)$ increases and is below $4.4$
at $u$. The transverse term is favorable throughout:
$T(B^3/A)\ge Q$ by the metric lemma, whereas
\[
 T(A^3/B)\le\frac{\sinc\theta}{c^3}
            (1/q+2\theta^4Q/3)<1.12/q+.0036Q<Q .
\]
Here $(\cot x)''\le2/x^3$ follows by twice differentiating its
partial-fraction expansion at zero, in which every non-pole
coefficient is nonpositive. Apply $T$ to $A^3/B=(\sin\theta/c^3)
[\cot(\theta q)-\sin(2\theta q)/2]$ and use $T(1/q)=1/q$.
The metric lower bound above is greater than $.32Q$, since
$.79+(.45-.49\cdot4.4)H^2>.32$. Thus the contribution of
$sR/(1-\kappa)$ in \eqref{IV:eq:five-exact-angular} is at least
\[
 .0105[.79Q+X(.45Q-.49H_5)]+.0039h .
\]
Indeed $\alpha_D/(1-\kappa)>.0105$ and $.32(2+l)/200>.0039$.

\paragraph{Bounds for the loss terms.}
The signed cubic of Lemma~\ref{IV:lem:common-diagonal-cubic} points
in a favorable direction for the volume moment, since
$\langle w,j\rangle\ge0$. Completing the square on its orthogonal
complement costs at most $a_q^2L_y^2h^5/(4Q)$ while the diagonal
energy contributes $.0072hQ$. Here $a_q\le1.16V<1.16(4.09)$:
in the formula for $V$ in Lemma~\ref{IV:lem:exact-metric-angular}
replace $u,\theta$ by $1.403,H/2$. Lemma~\ref{IV:lem:two-chord-primitives}
gives $L_y<.044$ by $B>.57$. Completing the square therefore reduces the lower bound by less than $.0046h^5$.

For the other losses take $I_0=\nu_0=0$ in
\eqref{IV:eq:angular-wall-charges},
$N_1\le\sqrt{h/15}$ and $N_2\le h^{3/2}/\sqrt{420}$ in
\eqref{IV:eq:angular-varying-volume}. After dividing by $h^2$,
the boundary norm in the square-root term is a constant times
$\sqrt{(.044h/8)^2+(.0156/4)^2}$; the other boundary terms are
nonnegative multiples of $1,h,h^2,h^3$. The volume terms are
nonnegative multiples of $h,h^2$. All these upper bounds increase
with $h$. To check their endpoint let
$B_v=(1.16/4)\sqrt{.049H(.015)}\,H^2$,
$N=\sqrt{H/15}$, $N'=H\sqrt{H/420}$, $Z=.015$. In order of
the coefficients of $1,Q,D_*$ the full losses at $H$ are below
\[
\begin{aligned}
 &1.053\sqrt Z[.141(.66)H^5/8+.15(.096)(.50)H^4/2]
       +1.4 B_v(N+N')<.00056,\\
 &.85 B_vN<.00026,\\
 &1.11\sqrt Z\sqrt{(.044H^3/8)^2+(.0156H^2/4)^2}+.85 B_vN'\\
 &\qquad+1.053(1.11)\sqrt Z\big[(.020(.66)+.012(.18))(.66)H^3/8
       +.012(.15)(.50)H^2/2\big]<.00026 .
\end{aligned}
\]
Hence $(h/H)^2[.00056+.00026(Q+D_*)]$ bounds the sum of the three errors.

\paragraph{Two endpoint substitutions.}
The cap estimate in Proposition~\ref{IV:prop:expanded-cap-eighteen-nineteen}
gives $.0032$ since $H^2H_9(u)/H_7(u)<.657$ and
$.003584-.000467(.657)-.000114(.657)^2>.0032$.
Put $C=.0032H_7+.32H_5/240-.077Q/12$.
It is positive since $H_7\ge2Q$ and $H_5\ge Q$; its coefficients
of degrees two and above are positive. After subtracting $.0123hF+.002$,
the entire residual is bounded below by
\[
\begin{aligned}
 P(h,q)={}&.0105[.79Q+h^2(.45Q-.49H_5)]+.0072hQ-.0046h^5\\
 &+h[(1/120-.0123)F-Q/90]+h^3C+.0039h\\
 &-(h/H)^2[.00056+.00026(Q+D_*)]-.002 .
\end{aligned}
\]
The disk limit of the same estimate before subtracting $.002$ gives
the stated leading bound. Replace $h^3$ by its tangent $3t^2h-2t^3$
at $t=.46$. The resulting minorant $P_0$ is concave in $h$, since
$H_5\ge Q$. At $h=0$ it is concave in $q$, and
\[
 P_0(0,l)>.014,\qquad P_0(0,u)>.0002 .
\]
At $h=H$, differentiating the displayed polynomial at $q=l+v$ gives
\[
 \partial_q P_0(H,l+v)
 <-.0038-.0137v+.0013v^2+.0078v^3+.0023v^4<0
 \qquad(0\le v\le u-l<1).
\]
At the right endpoint $u$, $P_0(H,u)>.0016$.
At $q=1$ the tangent subtraction gives $\Delta(1)=0$.
\end{proof}

\section{The mixed aperture--radius slope}\label{mid-slope}
\label{mid-slope10}\label{mid-slope11}\label{mid-slope12}\label{mid-slope13}\label{mid-slope14}

\begin{proposition}[Mixed aperture--radius slope for ten or more sides]
\label{IV:prop:middle-aperture-radius}\label{IV:prop:all-large-radius-slope}
For every $n\ge10$ the slope $J=G_e(q,0,0)/[qh^3(q^2-1)]$,
defined at one by continuity, is positive throughout the admissible aperture interval.
For $n\ge12$, the explicit aperture-dependent bound has the
form $J<1/6+j_Qh^2$ on $q\le Q$, with $j_Q$ as given below.
For ten and eleven sides the proof gives separate bounds above and below unit aperture.
\end{proposition}
\begin{proof}
Use the aperture interval $[l,u]$ of Section~\ref{mid-geometry}
and its continuous bounds for $n\ge12$.
Put $X=h^2$, $y=X/4$, $b=\sin(h/2)$, $c=\cos(h/2)$, $s=\sinc h$.
Use $\Lambda,d_-,d_+,s_-,s_+,m,Z,T_1,T_2,I_-,I_+$ from the bounds for the
regular-sector correction in Section~\ref{mid-response}, with $z_d>d_0$ from its diagonal
test. At ten and eleven sides the two orthogonal projections in
Section~\ref{III:sec:ang11} give $d_0=.00499,.00476$, respectively.
For ten and eleven sides we obtain
\[
\begin{array}{c|cc}n&10&11\\\hline
 \sup_{q\le1}J<&.221&.207\\
 \sup_{q\ge1}J<&.250&.230\\
\end{array}
\qquad \inf_{l\le q\le u} J>.020.
\]
At ten sides the same upper estimate gives
$J<.231,.238,.242,.246,.249$ when the aperture has successive upper bounds
$1.15,1.25,1.30,1.35,1.40$. At eleven sides $J<.220$ for apertures up to $1.20$. Each bound follows by substituting the indicated upper endpoint into the same formula.

\paragraph{The secant slope of the disk term.}
Use the cap density bounds \eqref{common-cap-height} of Section~\ref{mid-response}.
In the polynomial of Lemma~\ref{IV:lem:convex-radius-disk} use the endpoint
bounds $d_2^\pm,d_3^+$ obtained from
\[
 d_2=[\lambda_0(\lambda-1)+\lambda^2-3\lambda+6]/4,\quad
 d_3=[\lambda_0(2\lambda-3)+3\lambda^2-14\lambda+24]/12,
 \qquad \lambda_-+d_-X<\lambda_0<\Lambda .
\]
To bound all apertures $q\le Q$, take $t=Q^2$, $Q\ge1$, $v=yt$,
$\mu=y(1+t)/2$, $\mu_2=y^2(1+t+t^2)/3$.
The derivative of the convex disk function increases, so one upper
estimate for its secant is
\[
\begin{split}
 D_+={}&(1+\ell_+)/6+4\mu(d_2^+/15-1/12-\ell_-/60)\\
 &+16\mu_2(\ell_+/1680+1/80-d_2^-/30+3d_3^+/140+.0049Xt+.0059X^2t^2).
\end{split}
\]
For the lower endpoint put $t_0=l^2$. Truncating the same derivative
from below gives
\[
 D_-=(1+\ell_-)/6+(d_2^-/15-1/12-\ell_+/60)X(1+t_0)/2
                         -4X^2(1+t_0+t_0^2)/90.
\]
Indeed the first two positive density orders contribute derivatives at least
$1/6-x/12$ and $x/15-x^2/30$, respectively, and the higher positive
orders have nonnegative derivative. Here $d_2<4$, $0<-L<.145$, $Xt<.89$. For the continuous range we may use $d_2^-=3.855$, together with
$.32X<-L<X/3+.077X^2$ from the regular cap estimate.
 
\paragraph{Boundary and volume terms.}
With $f_0=b^2/X$ and $B=f_0(.5+.172y-2c^2(.1525))$, bounds for
the mean and standard deviation of the normal coefficient are
\[
\begin{aligned}
 N_-&={4.9831-[7.56-(2/3)\lambda_-(1-4v/3)(1-v/3)]\mu\over2s}
       -2\lambda_+(f_0+2B/3),\\
 N_+&=4.984/(2s)-2\lambda_- f_0,\\
 \sigma&=\sqrt{4/45}\max\{.371,1.60\lambda_+ v/(2s),2\lambda_+ B\}.
\end{aligned}
\]
For Lemma~\ref{IV:lem:normal-wall-mean}, with $S=\sin^2\sqrt x$ and
$D=1-t_{\rm wall}^2$, $A=(1-S)K'+K/2$ satisfies $A(0)=0$ and
$A_S=-(1+D/2)K'-(1-S)DK''\ge1/2$. Thus $H_S=\lambda(K+DA)
\ge\lambda(1+SD)$, $S\ge x(1-v/3)$ and $(\sinc(2\sqrt x))^2\ge1-4x/3$.
Averaging $D$ gives $2/3$, yielding $N_-$. The two terms before subtraction
have nonnegative $D$-derivatives bounded, respectively, by the maximum of
the first two entries and the last entry in $\sigma/\sqrt{4/45}$, by
$A\le.80S$ and $A_S\le.80$. Now use $\operatorname{Var}D=4/45$.
 
In the transverse coefficient of Lemma~\ref{IV:lem:middle-radial-slope} put
$f(x)=\cos\sqrt x\,\sinc^3\sqrt x$ and
$A_x=(f+xf')K-xfS'DK'$. Define
\[
 \alpha=[(1-2v)(1-v/3)+(1-v)(1-2v/3)]/2 .
\]
The sine series gives $1-2x+1.12x^2\le f+xf'\le1-2x+1.2x^2$,
$1-x\le f\le1-x+.4x^2$ and $S'\ge1-2x/3$ for $0\le x\le v<.20$.
Differentiation using $.5\le-K'<.54$ and $0<K''<.172$ gives
$\alpha x\le\partial_D A_x\le(1.08+.172v)x$.
On the reference chord $.30b^2\le\partial_D(-J_*)\le.32b^2$.
Hence $T=2\lambda[-\overline A_x/f(y)-J_*]$ decreases in $D$,
since $\alpha(1+l^2)/2>.32$. Its Lipschitz constant and its value
at $D=2/5$ are bounded above by
\[
 L_T=2\lambda_+[(1.08+.172v)\mu/(1-y)-.30b^2],\quad
 T_+=2\lambda_+\{[v+.4y^2-1.12\mu_2-(2/5)\alpha\mu]_+/(1-y)+.128b^2\}.
\]
The latter follows by integrating the lower bound on $f+xf'$ and
then on $\partial_D A_x$. Its numerator increases up to $t$,
since $1-1.12y(1+2t)/3-\alpha/5>0$. The chord measure weighted by
$t_{\rm wall}^2$ has $\mathbb E D=2/5$, $\operatorname{Var}D=12/175$.
It follows that the retained signed transverse mean contributes at
most $ch\nu_+ T_+$, where the upper bound for $\nu$ in
Section~\ref{mid-response} gives
\[
 \nu_+={s_++d_+mZh^3+
 h[T_1\sqrt{Z/c}/40+hT_2\sqrt{Z}/(24\sqrt{12})
       -f_0(1+2(5/24)/3)I_-]\over c(1+b^2/6)} .
\]
For the lower slope it suffices to bound $|T|\le2\lambda_+
(v+.4y^2)/(1-y)$: $A_x\le1$ gives the lower inequality,
and $A_x\ge1-2x+\alpha xD$ with $\alpha(1+l^2)/2>.32$
gives the upper inequality.

Put $a=2.6\cos^2(hQ/2)$, $a_1=6.36$, $a_2=2a_1+4.06Xt/2$.
The absolute volume contribution is at most $V_1\sqrt z+V_2\sqrt{z_d}$, with
\[
 V_1={5d_+h^{7/2}c\sqrt{bcm/2}\over6sc^2}
 \sqrt{a_1^2-4a_1a/3+a^2/2},\qquad
 V_2={d_+hb^2\sqrt{2bc/45}\over1-y}
 \sqrt{a_2^2-3a_2a+12a^2/5}.
\]
Indeed the radial disk factor $W=4K(r)+2rK'(r)$ satisfies
$3\le W\le6.36-2.6r$ and $-4.06<W'<0$. With reference
coordinates $(cX_1,bY_1)$, $r_p'\le0$ because $Y_1^2\le X_1^2$.
The two derivative terms in the longitudinal weight have opposite signs,
so its magnitude divided by $p$ is bounded by $5X(a_1-aX_1^2)/(6sc^2)$.
This uses $5W/6>4.06/2$ and $|r_p'|\le pX/2$.
The transverse weight is at most $(a_2-2aX_1^2)/(1-y)$.
Cauchy--Schwarz on the first, and the fiber primitive on the second,
give the displayed $V_i$. The normalized radial weights for their
squares are $4X_1^3\,dX_1$ and $6X_1^5\,dX_1$.

\paragraph{Two scalar inequalities.}
Set $\rho=\Lambda mh$, $p_0=X/(6s)$ and
\[
\begin{aligned}
 P_+&=-2\lambda_-+
 {2[1+(1.2v^2-.8v^3)/(1-v)^2]+5\Lambda mh^3/(3sc^2)\over1-\rho},\\
 P_-&={1+1.5y/[(1-v/6)(1-y/2)^3]\over \min(1,q_{\min})^2},\\
 Y_+&={\sigma T_1X/\sqrt c+T_2hL_T\sqrt{12/175}/(\sqrt{12}c^2)+V_1
      \over\sqrt{1-\rho}},\\
 Y_-&=\sigma T_1X/\sqrt c+
 {2\lambda_+T_2h(v+.4y^2)\over\sqrt{12}c^2(1-y)}+V_1 .
\end{aligned}
\]
\paragraph{Subtracting the mean and first moment of the boundary terms for ten and eleven sides.}
Recall the normal and transverse coefficients $N$ and $T$ of Lemma~\ref{IV:lem:middle-radial-slope}. From $N$ subtract its mean on the chord. From the transverse load $tT$ subtract its projection onto $t$, namely $t$ times the mean of $T$ with density $3t^2/2$ on $[-1,1]$. The remaining transverse load has zero first moment. We then bound the remaining errors using the energy of each field
component separately.
Let $f$ be $L$-Lipschitz as a function of $D=1-t^2$ and write
$\bar f,\bar f_2$ for its means on $[-1,1]$ with densities $1/2,3t^2/2$.
For $w=f-\bar f$ or $w=t(f-\bar f_2)$, put $W(t)=\int_{-1}^tw$.
Minimizing Lemma~\ref{IV:lem:two-chord-primitives} over its exponent gives
\[
 \begin{gathered}
 4\inf_k N_k(w)^2=\sinc^2(h/2)
 \left[\sqrt{C(U+V\tan^2(h/2))}+\tfrac12 C\tan(h/2)\right],\\
 (U,V,C)=\left(\int w^2,\int t^2w^2,\int W^2\right).
 \end{gathered}
\]
On $[0,1]$ the first primitive vanishes at both ends. The second has
zero mean, since $\int_0^1W=-\int_0^1tw=0$. Wirtinger's inequality
and parity give $C\le U/\pi^2$ in either case. The variance identity
$2\operatorname{Var}(f)=\mathbb E(f(D)-f(D'))^2$ for independent copies
gives $U\le(8/45)L^2$ and $U\le(8/175)L^2$, respectively.
Indeed the two variances of $D$ are $4/45$ and $12/175$.
Using $V\le U$, $\sec(h/2)<1.052$ and $\tan(h/2)<1/3$ gives
dual norms $.126 L$ and $.064 L$, respectively.

For the normal coefficient $N$ one may take the Lipschitz constant $L_N=.34$ in both polygons.
To verify this, use $w(x)=\sinc(2\sqrt x)$ and $A=(1-S)K'+K/2$
as in Lemma~\ref{IV:lem:normal-wall-mean}. Before subtracting the reference-chord term, the derivative of the normal coefficient with respect to $D$ is the secant average of $\lambda[w'SA+w^2(A+SA_S)]/(2s)$.
Here $1/2\le A_S\le .80$, $S/2\le A\le .80S$, $-2/3\le w'\le0$,
$w^2\ge1-4v/3$ and $(1-v/3)x\le S\le x$.
Averaging $x$ on the secant gives the lower and upper bounds
\[
 \frac{\lambda_-(1-1.88v)(1-v/3)y(1+l^2)}{4s},\qquad
 \frac{1.6\lambda_+\mu}{2s}.
\]
The subtracted reference derivative is between
$2\lambda_- f_0(.5-2c^2\cdot.172)$ and $2\lambda_+ B$.
At $Q=u$ the first pair lies inside $[.059,.430]$ and the second inside
$[.302,.396]$ for both side counts. These intervals contain those at
smaller upper limits on the aperture as well, so $|\partial_D N|<.34$.
For the transverse coefficient keep $L_T$ above. Its mean in density $3t^2/2$
is at most $T_+$: in deriving this estimate we bounded the average of $A_x$
from below by an affine function of $D$, and this density has mean $D=2/5$.
The two centered errors pair with $\nabla D_1$ and $\nabla D_2$.
Cauchy--Schwarz in this pair of components thus permits the replacement
\[
 Y_+=\frac{\sqrt{(.126L_N h^2)^2+(.064cbL_T)^2}+V_1}{\sqrt{1-\rho}}.
\]
This improves the upper bound for the mixed slope on the aperture intervals above one at ten and eleven sides.

For the lower estimate use $q_{\min}=l$ on $q\le1$ and $q_{\min}=1$
on $q\ge1$. The leading boundary term is at most $-2\lambda h s_D$.
Only $\|\partial_1D_2\|^2$ can have adverse upper energy coefficient,
and this component is at most $h^4(z-z_d)$. In the opposite direction
both diagonal coefficients are at least $-p_0$ and
$2\lambda_-^2>4\Lambda/(1-y)$ absorbs the adverse pure mass.
Combining Lemma~\ref{IV:lem:middle-radial-slope} with the preceding
estimates yields
\[
\begin{aligned}
 J&\le D_++Y_+\sqrt{s_-}+P_+hs_--N_- I_-X
                 +ch\nu_+T_++V_2\sqrt{d_0}-2hd_0,\\
 J&\ge D_--2h(\lambda_++P_-)Z+2h(P_--p_0)d_0
                  -N_+I_+X-(Y_-+V_2)\sqrt Z .
\end{aligned}
\]
Here $Y_+<-2P_+h\sqrt{s_-}$ and $V_2<4h\sqrt{d_0}$, so the upper
bounds involving $s_D$ and $z_d$ take their maxima at the indicated
lower endpoints. The lower inequality uses $z_d\le z\le Z$ in its
adverse volume term. Evaluating these bounds at each upper aperture limit gives the table.

For the continuous range $0<h\le H= .524$ put $l=.475$ and $u=1.403$.
Use $\Lambda=3.546$, $d_+=.575$, $s_-=\alpha_D$, $s_+=.0136$,
$m=.048$, $Z=.015$, $(T_1,T_2)=(.65,\mathcal T_H)$,
$(I_-,I_+)=(.0225,.044)$ and $d_0=.0072h$. Here $\mathcal T_H$
always denotes the exact endpoint expression at $\pi/6$ in
Proposition~\ref{continuous-response}; $.524$ is the outward
endpoint used in the interval estimates. These are outward bounds in Section~\ref{mid-response}.
Put $L=.32X$, $U=X/3+.077X^2$ in place of $\ell_-,\ell_+$ and
$d_2^-=3.855$. We retain the negative boundary contribution of order $h$
when estimating $J$. This is what permits the disk secant and the
mixed slope to be compared on the whole interval up to twelve sides.

Here are explicit rules to compute the coefficients of one polynomial.
Write $f_0=(\sinc(h/2)/2)^2$ as before,
$P=P_+$, $Y=Y_+$, $v_2=V_2/h^{7/2}$, $p_*=2-2\lambda_-$, and set
\[
 \begin{aligned}
 A={}&1/18+(d_2^+/15-1/12)(1+t)/2,\\
 D={}&.077/6-.32(1+t)/120\\
 &+\big[U/1680+1/80-d_2^-/30+3d_3^+/140+.0049Xt+.0059(Xt)^2\big]
                 (1+t+t^2)/3 .
 \end{aligned}
\]
Then $D_+\le1/6+AX+DX^2$. If an overbar denotes an upper bound over
$[0,H]$, the four coefficients in
$J\le1/6-a h+b_2h^2+b_3h^3+b_4h^4$ may be taken as
\[
\begin{aligned}
 a&=(2\lambda_- -2)s_-,\\
 b_2&=\overline{A+(Y/h^2)\sqrt{s_-}+[(P-p_*)/h]s_- -N_-I_--.0144}-.019,\\
 b_3&=\overline{c\nu_+ T_+/h^2},\qquad
 b_4=\overline{D+v_2\sqrt{.0072}}.
\end{aligned}
\]
All endpoint substitutions here have finite values at zero.
For explicit evaluation put $h\in[0,H]$, $X\in[0,H^2]$ and use
$\cos(H/2)\le c\le1$, $\sinc(H)\le s\le1$ and
$\sinc(H/2)/2\le b/h\le1/2$.
In $Y/h^2$ divide $V_1$ by $h^2$ using $b/h$ and replace $L_T/h$
by $2\lambda_+h[(1.08+.172v)(1+t)/(8(1-y))-.30 f_0]$.
For the two other divisions use
\[
\begin{aligned}
 (P-p_*)/h
 &\le [2(1.2)h^3t^2/(16(1-v)^2)+5\Lambda m X/(3sc^2)+2\Lambda m]/(1-\rho),\\
 T_+/h^2
 &=2\lambda_+\{[t/4+.4 X/16-1.12X(1+t+t^2)/48-.4\alpha(1+t)/8]_+/(1-y)+.128f_0\}.
\end{aligned}
\]
The additional subtraction $.019$ keeps the positive first-order
correction to the disk limit. Indeed $P<-4.46$ and $Y/h<.121$
throughout. For the former use the displayed bound for $(P-p_*)/h$;
for the latter use the displayed formula for $Y_+/h^2$ with
$L_T/h$ and $V_1/h^2$ above. At $Q=1.403$ their positive majorants
give $(P-p_*)/h<.601$ and $Y/h^2<.220$, respectively; the same bounds hold at smaller upper limits. In particular for $s\ge\alpha_D$,
\[
 \partial_s(Phs+Y\sqrt s)
 \le h\left(-4.46+\frac{.121}{2\sqrt{.00969}}\right)<-3.8h .
\]
Since $s_D>\alpha_D+h/200$, integration improves the upper bound by $3.8h^2/200=.019h^2$.
This gives the following rational bounds on the coefficients. The last column
is chosen so that $J<1/6+j_Qh^2$ for $q\le Q$.
\[
\begin{array}{c|rrr|c}
 Q&b_2&b_3&b_4&j_Q\\\hline
 1&.219&.023&-.043&.132\\
 9/8&.244&.030&-.063&.155\\
 5/4&.273&.037&-.088&.180\\
 53/40&.291&.042&-.106&.197\\
 1.403&.312&.048&-.128&.214
\end{array}
\]
Indeed put $B=b_2-j_Q$. The four Bernstein coefficients of
$-a+Bh+b_3h^2+b_4h^3$ on $[0,H]$ are
$-a$, $-a+BH/3$, $-a+2BH/3+b_3H^2/3$ and
$-a+BH+b_3H^2+b_4H^3$; each is negative. For the lower sign,
inserting the same intervals in the displayed formula for $J$ gives
$J>.006$ on $q\le1$ and $J>.061$ on $q\ge1$. In the latter check use
$Q=1.403$, $q_{\min}=1$. The endpoint directions for the two energies
follow from $Y/h<-2P\sqrt{s_-}$ and $v_2X<4\sqrt{.0072}$.
\end{proof}

\section{The radial cost and the case of equality}\label{mid-radial}

\label{mid-radial10}\label{mid-radial11}\label{mid-radial12}\label{mid-radial13}\label{mid-radial14}
\begin{proposition}[Radial cost of the projected trial]
\label{IV:prop:large-radius-profile}\label{IV:prop:profile-fifteen-sixteen}
Let $n\ge10$. Center the fan at a point where the projected trial field has mean zero. If its angles and endpoint radii satisfy the bounds of Section~\ref{mid-geometry}, then every sector satisfies
\[
 G(q,e,w)-G(q,0,0)\le k qh^5(q^2-1)^2 .
\]
At ten sides take $k=.015$ below unit aperture and $.0127$ above. At eleven sides take $.013$ below and $.012$ above.
At every $n\ge12$ take $k=.0123$
on both halves.
At $q=1$ equality requires $e=w=0$.
\end{proposition}
\begin{proof}
Write $d=h^2$, $v=q^2-1$ and $x=e/d$. Let $l,u,E,b,\alpha$ be the
bounds at this center, so $l<q<u$, $|e|<E$, $|w|/d<b$,
and $|w|/d<\alpha q$ below one. We apply \eqref{mid-radial-identity}
throughout.

\paragraph{Error bounds for the affine correction.}
Let $Z,m,\Lambda,d_+$ be the upper bounds on $E_D/h^4$, $M_D/(h E_D)$,
$\lambda_0$ and $(\lambda_0-\lambda)/d$ from Section~\ref{mid-response}.
Write $k=1/\operatorname{sinc}h$, $c=\cos(h/2)$,
$C=1/\operatorname{sinc}(h/2)$, $m_0=(h/2)/\tan(h/2)$,
$m_1=\tan(hu/2)/(u\tan(h/2))$, $R=1+E$ and $B_0=2d_+\sqrt{Zm/2}$.
Let $W_e,W_w$ be the upper bounds on the two signed boundary terms of Section~\ref{mid-wall}.
For the mean displacement and the skew displacement, respectively, use the following bounds. The purely affine skew term is estimated separately. We take
\begin{align*}
 C_e(L_0)&=W_e+2Zh^3 k\max\{1,(m_1u)^2,(m_0L_0)^{-2}\}
            +8.8 B_0 R\max(c^{-1},uC)h^4 k,\\
 C_w&=W_w+B_0 h^3\left[ahk/2+1.59 f a h+2fC
              +2.92 p a h/4+1.59 p C\right]+5b^2d^2/12 ,
\end{align*}
on a band with $q\ge L_0$, where $a=R\max(c^{-1},uC)+bhC$,
$p=R^2 k$, $f=bdk/2$. Indeed the singular-value and determinant
bounds are $a,pq$; the four skew derivatives
$(|\dot T x|,h\|\dot T\|,|\dot{\det T}|/q,|\ddot{\det T}|/q)$
are at most $(1/2,C,f,k/2)$. These substitutions in the twice
differentiated volume identity of Lemma~\ref{IV:thm:radial} give
the bracket. The mean volume estimate follows from
Section~\ref{IV:lem:radial-smooth-tensor}, since $R^2<1.18$;
the fourth-order disk remainder follows from
Lemma~\ref{IV:lem:affine-disk-skew}, since $\sin(hq)>.22$ and
$|w|/(2(1-E))<.05$ for ten and eleven sides. On
$0<h\le.524$ the small-angle hypothesis of the same affine estimate follows
from $.202h^2/(1-.044)<.475h\sinc(.475\cdot.524)\le\sin(hq)$.

For the remaining affine cost put $t=q^2$, $Y=d/4$, $D=Y/c$,
$z=bh/(1-E)$, $A_*=\cos(hu/2)$, $B_*=\operatorname{sinc}(hu/2)$,
and set
\begin{equation}\label{mid-affine-majorant}
\begin{aligned}
 p_0&=\max\{k,(c^3 B_*t)^{-1},C^3t/A_*\},\qquad p_4=C^4/(A_*B_*t),\\
 p_2&=2Yk+\max\{2C^2/(m_0t)+2m_1D^2t,\,
 C^2(m_1+\sqrt{m_1^2+(m_0t)^{-2}})+2D^2/m_0+m_1D^2t\},\\
 B&=p_2+Yp_0,\qquad
 P(t)=2Zh[B+z^2(p_4+YB)/(1-Yz^2)]
        +2\Lambda mZdYkR^2\max(c^{-2},C^2t).
\end{aligned}
\end{equation}
To see this bound, substitute the aspect-ratio inequalities
$m_0q\le\tan(hq/2)/\tan(h/2)\le m_1q$ in the two affine gradient
blocks of Section~\ref{IV:sec:background}. Their off-diagonal
stretch entries have magnitudes at most $C,Dq$ and product at most
$Ykq$. Bounding the constant, quadratic and quartic coefficients gives
$p_0,p_2,p_4$ above. With $z_0=w/(h(1+e))$, the even quartic difference divided by $w^2/2$ is
\[
 2h^{-2}[P_2+Y P_0+z_0^2(P_4+Y(P_2+YP_0))/(1-Yz_0^2)] ,
\]
where $P_i/(qE_D)\le p_i$. This gives the first term of $P$. Dropping the favorable mass
square gives the second. With $w=0$ the same constant block gives
the pure mean term in $C_e(L_0)$. The radical is the norm of a constant
and a positive multiple of $1/t$. Thus $P$ is convex and its maximum
on each closed interval occurs at an endpoint.

\paragraph{The three disk moments.}
Let $I_q=\operatorname{sinc}(hq)$, $S=\sin^2(hq/2)$ and
$S_*=\sin^2(hu/2)<.19$. For $t$ uniform on $[0,1]$, put $p=S(1-t^2)$.
Its first three moments $m_1,m_2,m_3$ are $2S/3,8S^2/15,16S^3/35$.
In this paragraph $m_1$ denotes this moment, not the ratio $\tan(hu/2)/(u\tan(h/2))$ defined above. Let $c_*d\le-L$, with $c_*=.3498,.3472$ at ten and eleven sides and $.32$ thereafter
from the cap estimates. Those estimates give $-L<.363d$ in the full range.
For $\lambda_*=3.3899$ and a lower bound $d_-$ on $(\lambda_0-\lambda)/h^2$ as in
Section~\ref{mid-response}, set
\[
 c_1^-=(\lambda_*+d_-d)(\lambda_*-1)+\lambda_*^2-3\lambda_*+3 .
\]
By the cap recurrence $4d_2-3=\lambda_0(\lambda-1)+\lambda^2-3\lambda+3
\ge c_1^->12.42$. Set
\[
 \begin{aligned}
 H(q)&=I_q[2+c_*d+c_1^-m_1-7.64m_2-10.47m_3/(1-S_*/2)],\\
 D(q)&=I_q[1/6-.363d/4+.64 S],\\
 M(q)&=4(m_1-m_2)/3-(2+E)^2E/6,\\
 \gamma(q)&=I_q[13(2/3+E/6)(1-m_1)-\min(7.36M(q),7.64M(q))],\\
 \gamma_w&=17(2+E)/12+2E(2+E)^2/3 .
 \end{aligned}
\]
The disk terms for mean and skew displacement are bounded below by $H(q)-\gamma(q)e_+$ and
$D(q)-\gamma_w e_+$, respectively, before subtracting the fourth-order skew term
already included in $C_w$. For the first bound use the chord moments of
Corollary~\ref{IV:cor:exact-two-curvature-moments} and apply
Corollary~\ref{IV:cor:disk-annulus-product} and its coefficient estimates:
the first coefficient is between $c_1^-$ and $13$, the second lies between
$-7.64$ and $-7.36$, and the tail from degree three is bounded below by
$-10.47p^3/(1-S_*/2)$ on $0\le p\le S_*$. For a negative argument the
alternating tail starts positive. On moving the mean by $e\ge0$ the new
argument in the Taylor average is $p'=p-[(1+ze)^2-1](1-p)$, with $z$
independent of $p$ and having density $2(1-z)$ on $[0,1]$. Consequently
\[
 \mathbb E p'\ge m_1-e(2/3+E/6)(1-m_1),\qquad
 \mathbb E (p')^2\le m_2-eM(q).
\]
Here $\mathbb E[(1+ze)^2-1]\ge2e/3$ and the second moment of this
increment is at most $(2+E)^2e^2/6$. Since $p'\le p$, the same
positive-argument tail bound applies. This proves the stated $\gamma$,
including both signs of $M$. The bracket in $\gamma$ is positive and
decreases with $S$, since $m_1-m_2$ increases up to $S_*$.
Both $H,D$ increase. Indeed for $W=dq^2$, $dS/dW=I_q/4$ and
$-1/6\le dI_q/dW\le0$; the bracket in $H$ is less than $4$ and has
derivative with respect to $S$ greater than $6$, while $I_q>.87$.

For the skew term use the chord identity in
Lemma~\ref{IV:lem:exact-real-cap}. With $K=(1+e)^2S$, the positive
aperture terms contribute at least
\[
 (c_1^-/30+4/15)K-16K^2/105
       -4(10.47)S^3/[315(1-S/2)]\ge .64 S .
\]
Indeed it suffices that $c_1^-/30+4/15-16(1+E)^2S_*/105
-4(10.47)S_*^2/[315(1-S_*/2)]>.64$, which follows at the largest allowed $E,S_*$ in each range. The remaining linear error is at most $\gamma_w e$ by the
same chord identity and $2e\le(1+e)^2-1\le(2+E)e$.
For $e\le0$ both $-\Phi''$ and the disk density increase when the radius
decreases along this chord. The bounds therefore hold with $e_+=0$.

\paragraph{Positive and negative displacements on a narrow aperture interval.}
Let $d=h^2$, $s=|e|/d$, $t=|q^2-1|$. On a band $[a,r]$ let $B$
bound $|w|/d$, let $j$ bound the positive mixed slope, and put
\[
 G=H(a)-C_e(a),\qquad
 \delta=D(a)-C_w-\max\{P(a^2),P(r^2)\},\qquad
 g=d\gamma(a),\quad A=d\gamma_w .
\]
These endpoint bounds follow from the monotonicity of the disk
terms and convexity of $P$. If $X$ bounds the positive displacement
$s$, the radial Taylor bound before division by $t^2$ is
\[
 j s t\,\mathbf1_{e(q-1)\ge0}-
 \tfrac12\bigl(G-gs\mathbf1_{e\ge0}\bigr)s^2
 +\tfrac12 B^2\bigl(A s\mathbf1_{e\ge0}-\delta\bigr)_+ .
\]
The normalization here is $qh^5$. For $q\le1$, the mixed term can
increase the cost only for $e\le0$, when both displacement corrections
vanish. Completing the square gives at most
\[
 \frac{j^2}{2G}t^2+\tfrac12B^2(-\delta)_+ .
\]
For $e\ge0$ the mixed term can be dropped. Take $K=6/5$, increased to $13/10$ on the final band at each side count, and $A_*=5/4$.
On the aperture intervals below one, $G-gX>K$ and $A<A_*$: at the lowest aperture,
$H(l)-C_e(l)-\gamma(l)\sup e$ exceeds $1.23,1.45$ at ten and eleven
sides, respectively. The same substitution on the two final bands gives $G-gX>1.51$. Applying
$B^2A_*s/2-Ks^2/2\le B^4A_*^2/(8K)$ on the active skew branch
gives an upper bound valid on both branches,
\[
 \tfrac12 B^2\left(\frac{B^2A_*^2}{4K}-\delta\right)_+ .
\]
In each row of the table divide the constant terms above by
$(1-r^2)^2\le t^2$ to get the normalized coefficients $k_-,k_+$.
In the final band for each polygon both constant terms are zero.
The successive bands begin at $.426,.463$ and use $j=.221,.207$,
respectively. Each displayed coefficient is an upper bound.
\[
\begin{array}{c|c|rrr|rr|rr}
 n&r&G&\delta&B&10^3 k_-&10^3 k_+&W_e&W_w\\\hline
10&.48&1.81&-.158&.087&14.6&1.1&.341&.011\\
  &.52&1.85&-.101&.113&14.5&1.3&.339&.0111\\
  &.60&1.88&-.069&.150&14.9&2.1&.339&.0113\\
  &.70&1.95&-.022&.227&14.8&3.9&.337&.0123\\
  &1&2.04&.016&.227&12.0&0&.340&.0151\\\hline
11&.55&1.89&-.064&.155&13.0&1.8&.271&.008\\
  &.60&1.95&-.008&.169&11.3&.61&.268&.008\\
  &1&1.99&.015&.209&10.8&0&.273&.0102
\end{array}
\]
For example, for ten sides on $[.52,.60]$ the negative-displacement
coefficient is at most $.221^2/(2\cdot1.88)+.150^2(.069)/[2(1-.60^2)^2]
<.0149$. At $q=1$ the final row has $\delta>B^2A_*^2/(4K)>0$.
For $e\ge0$ the radial upper bound is strictly negative on the inactive skew branch for positive $s$,
and on the active branch the completed square above is strictly
negative. For $e<0$ use $G>0$. Thus equality gives $e=0$, and strictness in the skew estimate then gives $w=0$.

\paragraph{A single endpoint on each wide band at ten and eleven sides.}
Here $t=q^2-1\ge0$. Negative $x$ gives no positive mixed term and
both curvature terms have positive lower bounds. For $s=x\ge0$ let $X$ be its
upper bound. With $G,\delta,g,A,j,B$ as above, the two possible
upper bounds on the radial increase divided by $qh^5$ have the form
\[
 F_i(t,s)=jts-\tfrac12(G-gs)s^2+\tfrac12 i B^2(As-\delta),
 \qquad i=0,1 .
\]
For every row below, $G>3gX$, $\delta>0$ and
$27g\delta^2>B^2A^3$. The following calculation replaces all cubic
sign tests on these wide intervals. The function $F_i$ is strictly
concave in $s$. At an interior maximum of $F_i/t^2$, differentiating
its optimized value in $t$ has the sign of
\[
 t\,\partial_t F_i-2F_i
   =\tfrac12\{gs^3+iB^2(2\delta-As)\}\ge0 .
\]
The last inequality follows by minimizing the cubic on $s\ge0$, using
$27g\delta^2>B^2A^3$. The maximizing displacement therefore gives an
increasing normalized cost until it reaches $X$. For $[t_0,t_1]=[q_0^2-1,q_1^2-1]$ put
$D=iB^2$ and $t_c=(GX-3gX^2/2-DA/2)/j$.
When $t_1\le t_c$, the maximum occurs at $t=t_1$ and equals
\[
 \frac{GS_0^2/2-gS_0^3-D\delta/2}{t_1^2},\qquad
 S_0=\frac{2jt_1+DA}{G+\sqrt{G^2-6g(jt_1+DA/2)}} .
\]
When $t_1>t_c$ put $C=(G-gX)X^2+D(\delta-AX)>0$. On the remaining
branch the normalized cost is $jX/t-C/(2t^2)$. Its maximum is at
$C/(jX)$ clipped to $[\max(t_0,t_c),t_1]$. The table gives the upper
endpoint of each successive band, its rounded inputs and an upper bound
$10^3 k$ for both $i$. The two values of $A$ are $1.23,1.00$.
\[
\begin{array}{c|c|rrrrrr|r|rr}
 n&q_1&j&G&g&X&\delta&B&10^3 k&W_e&W_w\\\hline
10&1.15 &.231&2.339&2.96&.163&.078&.227&12.0&.322&.0157\\
  &1.25 &.238&2.49&2.78&.163&.092&.227&12.3&.311&.0176\\
  &1.30 &.242&2.59&2.65&.163&.101&.192&12.2&.302&.0186\\
  &1.35 &.246&2.645&2.58&.163&.104&.171&12.7&.298&.0203\\
  &1.37 &.249&2.69&2.51&.135&.111&.137&12.7&.289&.0186\\
  &1.385&.249&2.71&2.49&.125&.112&.137&12.7&.286&.019\\
  &1.40 &.249&2.73&2.47&.125&.114&.137&12.7&.286&.019\\
  &1.413&.250&2.74&2.45&.109&.116&.073&12.7&.280&.0181\\\hline
11&1.20 &.220&2.31&2.50&.151&.101&.209&11.1&.256&.0109\\
  &1.391&.230&2.49&2.32&.151&.118&.209&11.9&.247&.0124
\end{array}
\]
For the decagon on $[1.30,1.35]$, Section~\ref{mid-geometry} gives
$B=.171$, and the mixed-slope bound in Section~\ref{mid-slope} gives
$j=.246$. With $W_e=.298$ and $W_w=.0203$, the two endpoint formulas give
$G>2.64581$ and $\delta>.104926$. Using the rounded row above for both
$i=0,1$ gives a radial cost below $.01268915$.
On $[1.385,1.40]$ both switch points exceed $1.07$, so $t_1=.96$ applies twice. Direct substitution gives both costs below $.0127$. The table proves the required bounds on
the closed wide intervals, with the regular sector treated above.

\paragraph{Curvature bounds for twelve or more sides.}
Set $H=.524$, $D=H^2$, $l=.475$, $u=1.403$, $E=.160D$ and $b=.202$.
On the full range $0<h\le H$, Section~\ref{mid-geometry} gives
$|e|\le.160h^2$ and $|w|\le h^2\min(.202,.259q)$ for $q\le1$.
The bounds for the regular sector are $\Lambda=3.546$, $d_+=.575$, $Z=.015$, $m=.048$,
and Section~\ref{mid-wall} gives $W_e<.255$, $W_w<.014$.
Use the preceding bounds on the affine terms, increasing $h$ to $H$
in each positive factor: the cosine and sinc decrease, and
$\tan(uz)/(u\tan z)$ increases since $\sinc(2uz)<\sinc(2z)$.
Thus the nonaffine terms for the skew displacement satisfy $C_w<.025$. We use $.025$ below.

Let $W=h^2q^2$; then $W/4-W^2/48\le S\le W/4<.155$.
In the three-moment estimate use $c_1^->12.42$ and replace the two
adverse coefficients by $8,11.2$. Set
\[
 B_1=12.42/6,\qquad B_2=B_1/12+4/15,\qquad
 B_3=11.2/[140(1-.155/2)].
\]
The disk term for the mean displacement, evaluated at zero displacement, is at least
$(1-W/6)(2+.32d+B_1W-B_2W^2-B_3W^3)$.
Take $A_1=B_1-(2+.32D)/6$, $A_2=B_2+B_1/6$,
$A_3=B_3-B_2/6>0$. For a band $[a,b']$, a lower bound for the
coefficient of $d$ above two is
\[
 H_d=.32+a^2(A_1-A_2Da^2-A_3D^2a^4).
\]
Indeed $A_1-2A_2 W-3A_3W^2>0$ at $W=Du^2$, hence throughout.
For the skew displacement, retain the positive term $.64S$. It improves the bound most
at the largest apertures. With $P$
from \eqref{mid-affine-majorant} evaluated at $H$, use
\[
 \begin{aligned}
 \delta_0&=\sinc(Hb')(1/6-.363D/4)-C_w
                     -\max\{P(a^2),P((b')^2)\},\\
 \delta_1&=.64\sinc(Hb')(a^2/4-Da^4/48).
 \end{aligned}
\]
Here $S\ge d(a^2/4-Da^4/48)$ because $W/4-W^2/48$ increases.
The following table rounds the mean error bound upwards and
$H_d,\delta_0,\delta_1$ downwards. Use the preceding slope bound
with $Q=\max(1,b')$ for $c_J$.
\[
\begin{array}{c|c|rrrr|c}
 a&b'&C_e&H_d&\delta_0&\delta_1&c_J\\\hline
 .475&.55&.299&.697&-.041&.035&.132\\
 .55&.625&.294&.821&-.003&.047&.132\\
 .625&1&.290&.959&.019&.059&.132\\
 1&1.125&.287&1.825&.064&.147&.155\\
 1.125&1.25&.287&2.152&.064&.182&.180\\
 1.25&1.325&.287&2.481&.064&.222&.197\\
 1.325&1.403&.287&2.674&.063&.245&.214
\end{array}
\]
\paragraph{Square completion and two derivative tests.}
Put $D=.524^2$, $S=.160$, $r=.111$, $r'=.123$, $k=.0123$.
On each consecutive aperture band $[a,b]$ in the table write
$t=|q^2-1|$, $s=|e|/h^2\le S$, $d=h^2$,
$G=2-C_e+H_d d$, $\delta=\delta_0+\delta_1d$,
$j=1/6+c_Jd$ and $A=3.06d$.
Write $\beta=\min(.202,.259b)$ on the narrow side and $\beta=.202$
on the wide side. The displacement coefficient can be taken to be $9$
because $13(2/3+E/6)+7.64(2+E)^2E/6<8.996$ for $E=.160D$.
For fixed $s,t$ and the sign of $e$, the lower bounds on the surplus
with and without the signed skew term are affine in
$d$. It suffices to check $d=0,D$.

On a narrow band with $e\le0$, the only terms to bound are $jst$ and
the possible constant skew error $\beta^2(-\delta)_+/2$. If $b<1$ put
$k'=k-\beta^2(-\delta)_+/[2(1-b^2)^2]$, and on the band ending at one
put $k'=k$. The six endpoint substitutions give
$2Gk'-j^2>.0008$, so $jst\le k't^2+Gs^2/2$ bounds the mixed term.
For $e\ge0$ the mixed term adds no error. In the first two bands even
$2k(1-b^2)^2-\beta^2(AS-\delta)_+>.0058$. In the band ending at one,
$G-9dS-(A\beta)^2/(4\delta)>1$, and
$(As-\delta)_+\le (As)^2/(4\delta)$.

On a wide band only $e\ge0$ can give a positive increment. At $d=0$ use
$2kG-j^2>.014$. At $d=D$ all four rows satisfy
$\delta>Ar'$ and $2k(G-9Dr)-j^2>.0026$. For example the last
row, at the outward endpoint $D=.524^2$, gives
\[
 2(.0123)[2-.287+2.674D-9D(.111)]
       -(1/6+.214D)^2>.00263 .
\]
Equivalently $2k(G-9Dr)>1.05j^2$ on this row. Consequently the skew
error vanishes on $0\le s\le r$, where Young's inequality suffices.
For $s\ge r$ the derivative of the residual before subtracting the skew error is
$Gs-13.5Ds^2-jt$ and increases with $s$ since $G>27DS$. With
$B(s)=Gs-13.5Ds^2-j(b^2-1)$, substitution at the two radii gives
\[
 B(r)>.0076,\qquad B(r')-\beta^2A/2>.0094 .
\]
The first bound applies while the skew error is zero. When it is positive,
$s\ge r'$ and its derivative is at most $\beta^2A/2$, covered by
the second bound. Thus the residual is increasing on the entire
remaining interval. At $q=1$ strictness of the mean-displacement estimate forces
$e=0$ in equality and then $\delta>0$ forces $w=0$.
\end{proof}

\paragraph{The quadratic lower bound at the regular sector.}
Let $a_n=\lim_{q\to1}\Delta(q)/[qh^5(q^2-1)^2]$, and put $d=h^2$,
$x=e/d$, $y=w/d$, $v=q^2-1$. At the regular sector the quadratic
part of the surplus, divided by $h^5$, is
\[
 a_n v^2-J(1)vx+\tfrac12K_e(1,0,0)x^2+\tfrac12K_w(1,0,0)y^2 .
\]
The absence of terms linear in $y$ follows from reflection and the
Taylor formula in \eqref{mid-radial-identity}. The regular angular
endpoint test and the narrow-side limit of the mixed slope give
$(a_{10},a_{11})>(.01730,.01919)$ and
$(J_{10}(1),J_{11}(1))<(.221,.207)$, both slopes positive.
At $q=1$, $e=w=0$, the values of $H(1)-C_e(1)$ on the wide side at ten and eleven sides are
greater than $2.307,2.247$ and the skew curvature coefficients are positive.
Minimizing the quadratic form over $x,y$ at fixed $v$ thus leaves
more than $.0067v^2$ at ten sides and $.0096v^2$ at eleven:
\[
 .01730-\frac{.221^2}{2(2.307)}>.0067,\qquad
 .01919-\frac{.207^2}{2(2.247)}>.0096 .
\]
These are local coefficients of this trial at the regular sector.

\section{Summation on the projected fan}\label{mid-summation}\label{IV:sec:qualitative-tail}

At ten and eleven sides the intervals below and above unit aperture have lower bounds $a$ on the angular gain and upper bounds $b$ on the radial error, in units of $qh^5(q^2-1)^2$. The endpoints enclose the admissible aperture interval from Section~\ref{mid-geometry}. Section~\ref{mid-angular} proves the angular bounds and Section~\ref{mid-radial} proves the radial bounds.
\[
\begin{array}{c|l|l|l|c}
 n&\text{consecutive endpoints}&a&b&\min(a-b)\\\hline
10&.426,1,1.413&.0172;.0148&.0150;.0127&.0021\\
11&.463,1,1.391&.0181;.0155&.0130;.0120&.0035\\
\end{array}
\]
For every $n\ge12$, Proposition~\ref{IV:prop:entry-twentyfour} gives
\begin{equation}\label{IV:eq:large-local-window}
 .475<q_i<1.403,\qquad |e_i|<.160h^2,\qquad |w_i|<.202h^2,\qquad
 |w_i|<.259q_i h^2\quad(q_i\le1).
\end{equation}

\begin{proposition}[Projected comparison from ten sides]\label{prop:common-projected-comparison}\label{IV:thm:qualitative-large-polygon}
Let $n\ge10$, $h=2\pi/n$ and $\lambda_0=\mu_1(R_n)$. Give the convex
$n$-gon $P$ the area of the circumradius-one regular polygon and let $Q$
be an interior zero-mean center of the projected trial $V$. Suppose
its fan satisfies the geometric bounds of Section~\ref{mid-geometry}
when $n=10,11$, including the joint radius and skew bounds there,
and \eqref{IV:eq:large-local-window} when $n\ge12$. Then
\begin{equation}\label{eq:common-projected-comparison}
 \int_P\left(|\nabla V|^2-\lambda_0|V|^2\right)
 \le -\frac{h^4}{1000}\sum_{i=1}^n(q_i-1)^2 .
\end{equation}
If all $q_i=1$, the inequality is strict unless all endpoint radii are
one. Every competitive polygon with at most $n$ sides admits such a fan,
and hence is similar to $R_n$.
\end{proposition}
\begin{proof}
Write $C=G(q,e,w)-G(q,0,0)$ for the radial cost and $\Delta$ for
the angular deficit. For $n=10,11$ the two coefficients $a,b$ in
the table give, on their indicated bands,
\[
 \frac{\Delta-C}{h^4(q-1)^2}\ge (a-b)\,h q(1+q)^2
 \qquad(q\ne1).
\]
On $q<1$ the table gives $a-b\ge 1/500$. Each lower aperture
endpoint $l_n$ in its second column satisfies $12l_n(1+l_n)^2>n$.
Since $2\pi>6$ and $q(1+q)^2$ increases, $h q(1+q)^2>1/2$.
On $q>1$ use $a-b\ge1/1000$ and $h q(1+q)^2\ge 4h>1$.
Thus the angular gain minus the radial error in each sector is at least $h^4(q-1)^2/1000$. For every $n\ge12$,
Proposition~\ref{prop:angular-from-twelve} gives the same lower bound
after subtracting $.0123qh^5(q^2-1)^2$. Proposition~\ref{IV:prop:large-radius-profile}
bounds the radial error by the subtracted term, including the twelve-sided endpoint. In both ranges the radial estimates at $q=1$ give
$C<0$ unless $e=w=0$, and then $C=\Delta=0$.
Sum the local bounds using Lemma~\ref{IV:lem:projected-fan-identities}.
The common area cancels the constant area term, $\sum_i(q_i-1)=0$
cancels the angular tangent, and the trial mean vanishes by construction.
This proves \eqref{eq:common-projected-comparison} and its strictness
when all apertures are regular but an endpoint radius differs from one.

Theorem~\ref{IV:thm:corner-gap} forces a competitive polygon to
have exactly $n$ strictly convex vertices. Section~\ref{mid-geometry}
and Proposition~\ref{IV:prop:entry-twentyfour} provide the
required centers in their respective ranges. A negative residual
contradicts competitiveness by the zero-mean Rayleigh inequality,
applied to the components of $V$ and summed. The only remaining fan has
unit radii and successive angles $2\pi/n$, and is therefore regular.
\end{proof}

\section{Completion of the extremal comparison}
\begin{proof}[Completion of the proof of Theorem~\ref{thm:main}]
The triangle bound and the quadrilateral comparison are established in Sections~\ref{I:sec:lowertri}--\ref{I:sec:all-real-quads}. At five and six sides Section~\ref{paper-small-sum} compares the gain from a single narrow sector with the error in the Rayleigh quotient caused by using the polynomial as the reference field for the entire fan. Section~\ref{paper-seven-sum} applies the same comparison for seven through nine sides. Proposition~\ref{prop:common-projected-comparison} now applies to every $n\ge10$. The bounds for polygons with fewer than $n$ strictly convex vertices in Part~\ref{part:small} and Theorem~\ref{IV:thm:corner-gap} cover polygons with fewer sides and those with flat or reflex corners in the respective ranges. These arguments show that equality requires the polygon to be similar to the regular polygon. Conversely, scale invariance gives equality for every polygon similar to the regular polygon.
\end{proof}
\clearpage\part{Quantitative stability}
\section{Global stability in relative symmetric difference}\label{sec:stability}
Write $D_n(P)=F(R_n)-F(P)$. The proof of
Theorem~\ref{thm:stability} begins by showing that a small deficit implies the required geometric
bounds. For $n\ge5$ we bound the deficit below both for polygons with too
few convex corners and for convex $n$-gons that violate those bounds.

\begin{proposition}[The deficit outside the geometric window]
\label{prop:quantitative-window-entry}
Let $n\ge5$ and let $P$ be a simple polygon with at most $n$ sides. Normalize its area to that of $R_n$. Write $D_n(P)=F(R_n)-F(P)$
and set
\[
 \varepsilon_n=
 \begin{cases}
  8/10^4,&n=5,\\
  1/10^5,&n=6,\\
  2/(25n^5),&7\le n\le9,\\
  4/(5n^5),&n\ge10 .
 \end{cases}
\]
If $D_n(P)\le\varepsilon_n$, then $P$ is a strictly convex $n$-gon. Its area-centroid fan satisfies the bounds
of the five- and six-sided geometric argument or of
Section~\ref{seven-window-proof}, respectively, for $5\le n\le9$.
For $n\ge10$ it has an interior zero-mean projected fan satisfying all
the bounds in Proposition~\ref{prop:common-projected-comparison}.
In particular, failure of this conclusion costs strictly more than
$\varepsilon_n$.
\end{proposition}

\begin{proof}
If a strictly convex corner is missing, the deficit exceeds $10.3526-56\pi/17$
at five sides, by Lemma~\ref{lem:degree-six-pentagon-floor} and the
four-corner bound. At six sides the deficit exceeds $.0008$ and from seven
sides it exceeds $5/n^4$, by Theorem~\ref{IV:thm:corner-gap}.
These gaps exceed the respective $\varepsilon_n$. We can therefore
work with a strictly convex $n$-gon.

Write $\lambda=(j'_{1,1})^2$ and, for any polygon $K$, put
\[
 d(K)=(\lambda-1)(\pi\lambda/F(K)-1),\qquad H(K)=d(K)/\lambda .
\]
Given an upper bound $d^*$ on $d(P)$, the quantity
$\mu_*=\lambda/(1+d^*/(\lambda-1))$ is a lower bound for $\mu_1(P)$
at area $\pi$. This is exactly the hypothesis in
Proposition~\ref{IV:prop:harmonic-budget-65}, with
$H=(\lambda-1)/\mu_*-(1-1/\lambda)=d^*/\lambda$. The harmonic estimates below therefore apply whenever $d(P)\le d^*$.

At five sides the aperture and skew proof uses $F(P)>10.3518$,
which follows from $F(R_5)>10.3526$ and $D_5(P)\le.0008$.
At six sides the three harmonic pairs and the transfer to the centroid
use $H\le H_0=.009764$. The eigenvalue lower bound for the regular hexagon in
Section~\ref{I:sec:five-six-regular-data} gives $H(R_6)<.009763$,
as calculated in the six-sided geometric argument. For $0\le H\le H_0$,
\[
 \left|\frac{d}{dH}\frac{\pi\lambda}{1+\lambda H/(\lambda-1)}\right|
 =\frac{\pi\lambda^2}{(\lambda-1)(1+\lambda H/(\lambda-1))^2}>10.
\]
Consequently $H(P)>H_0$ implies a deficit greater than $10^{-5}$. When $H(P)\le H_0$,
choose $d^*=\lambda H_0$ in the preceding harmonic bound. The lower
bound for the multiplier of the three harmonic pairs was evaluated at
$H_0$.

For $n\ge7$ put $B=4\zeta(3)$ and
$d^*=B/n^3+U(n)/n^5$ with $U(n)$ as in
\eqref{IV:eq:regular-continuous-band}. The Schur estimates have the
following unused margin:
\[
 d(R_n)< d^*-\sigma_n/n^5,\qquad
 \sigma_n=
 \begin{cases}
  1/50,&7\le n\le9,\\
  1/5,&n\ge10 .
 \end{cases}
\]
For seven through nine subtract the displayed bounds
$92.29,85.16,80.72$ from $92.5,85.25,80.75$.
For $1/13\le x\le1/10$ compare $59.93+197x$ with the chord joining
$(1/13,74.2)$ to $(1/10,79.4)$. The difference is affine, with endpoint
values greater than $.88$ and $.22$. At $n=14$ use
$59.93+197/14-(68+80.6/14)>.24$. At $n\ge15$ use
$74-(68+80.6/15)>.62$. These are precisely the convex majorants in
Lemma~\ref{IV:lem:uniform-regular-remainder}.

Since $d^*<.021$ and $3.3899<\lambda<3.39$,
\[
 F(R_n)-\frac{\pi\lambda}{1+d^*/(\lambda-1)}
 =\frac{\pi\lambda(\lambda-1)(d^*-d(R_n))}
 {(\lambda-1+d(R_n))(\lambda-1+d^*)}
 >4\sigma_n/n^5.
\]
Thus $D_n(P)\le4\sigma_n/n^5$ implies $d(P)<d^*$.
For $7\le n\le11$ the geometric proofs use
$H_+=(4.809+U(n)/n^2)/(3.3899n^3)\ge d^*/\lambda$.
For $n\ge12$ Proposition~\ref{IV:prop:entry-twentyfour} uses
$H=(B+\max(74,59.93+197/12)/n^2)/(\lambda n^3)\ge d^*/\lambda$.
In either case choose the corresponding $\mu_*$.
We identify the spectral assumption used to obtain the geometric bounds
for each side count. Whenever we use an eigenvalue of the polygon $P$, we scale its area to $\pi$. Bounds on the trial on $R_n$ depend only on $n$.
\begin{enumerate}
\item At five sides put $d^*=(\lambda-1)(\pi\lambda/10.3518-1)$.
In Section~\ref{sec:pent-entry} the condition $d(P)<d^*$ is precisely
$F(P)>10.3518$. The covariance tests give the radius ratio
$K=70\pi/(9F(P))$ as an upper bound. Replacing $F(P)$ by $10.3518$
and applying the sine rule gives the upper bounds on aperture and skew,
as well as the bound on skew in terms of aperture below $x=.65$.
\item At six sides take $d^*=\lambda H_0$. Equality at $H_0$ is
also excluded by the strict derivative estimate above. In
Section~\ref{sec:hex-entry}, use the eigenvalue bound for this $d^*$ in
Proposition~\ref{IV:prop:harmonic-budget-65}. The first two harmonic
bounds give the multiplier lower bound $3.20$ at $H_0$. The three additional
pairs give $T<.03957$ and $C<116.3$. The support and centroid estimates
then bound $|G-O|/r_*<.080$. The six-sided transfer uses this last
bound, $C<117$ and the aperture--skew inequality at $O$ to move both
bounds to $G$. The regular hexagon's eigenvalue lower bound enters only in the
derivation of the threshold $H_0$. The remaining estimates on $P$ use
the bounds just listed.
\item At seven through nine sides Section~\ref{V:sec:geo-seven-nine}
uses $H_+$ in the two leading Fourier coefficients and in $P(D)$.
The two coefficients multiplying their squares are positive for
$D>.23$, so this substitution is an upper bound. Before the additional tests, the coefficients of $w_i$, $i\ge3$, are at least $.2448-H_+$. This gives $T<.025$ and $\mu_* /(1+T)>3.19$.
Section~\ref{V:sec:geo-bubble} then bounds the tail. The resulting upper bounds on the conformal
excess are $109,96.2,88.6$. Equation~\eqref{V:eq:finite-angle-budget}
bounds the aperture and skew at $O$. The support and displacement
bounds and the differential transfer in Section~\ref{seven-transport}
give the angular and skew bounds at $G$. Once the coefficients of $w_i$ have been bounded,
this argument uses only the bounds on the conformal excess and the map coefficients.
\item At ten and eleven sides the estimate of the conformal excess in
Section~\ref{mid-geometry} uses the corresponding $H_+$ in the same
two coefficients and to prove $\mu_* /(1+T)>3.35$.
This gives $C=83.1,80$ as upper bounds. The path estimate and
ellipse~\eqref{paper-joint-ellipse} use $C$ and the fixed side count.
The support and moment estimates give a ball containing a center for which
the full trial has mean zero, by \eqref{IV:eq:moment-center-ball}. Applying the same ellipse
with the displacement of that center gives \eqref{paper-geometric-window},
including the one-sided radius and skew bounds in the paragraphs
following it. To locate this center we combine the geometric bounds for $P$ with the bounds on the correction $D$ on $R_n$ from Section~\ref{mid-response}.
\item For twelve or more sides, the proof of
Proposition~\ref{IV:prop:entry-twentyfour} uses the stated bound on $H$ first to prove $\mu_* /(1+T)>3.37$ and then, after the additional harmonic tests, to bound the coefficient of each $w_i$, $i\ge3$, below by $1/\lambda-x^3V>.290$. Its scaled bounds
$T_3<5.45$, $C_*<79.28$ permit $C=79.3$. The path estimate of Lemma~\ref{IV:lem:middle-conformal-entry} and the estimate for the selected edge
use this upper bound, and the choice of a center where the trial field has mean zero uses $m_*<.048$,
$z<.015$ on the regular polygon. The bound for the center displacement
is inserted into the tilt and aperture estimates in the same proof.
These are the geometric bounds for the projected trial, including $|w|/(qh^2)<.259$
below unit aperture.
\end{enumerate}
In each sector comparison, the extremal argument used $D_n(P)\le0$
only to obtain the spectral bound verified above. Strict convexity, the other geometric input, was established
separately at the start of this proof.
\end{proof}

\subsection{The centroid and affine trials}
\begin{lemma}[Quadratic control on the central centroid fans]
\label{lem:central-centroid-coercivity}
Let $5\le n\le9$. Normalize $P$ and $R_n$ to the same area, and use
the first eigenfield of $R_n$ on the complete centroid fan of a convex
$n$-gon. Assume all its sectors satisfy the bounds for thick sectors in
Section~\ref{I:sec:five-six-thick} at five and six sides, or the
central aperture and skew bounds of Sections~\ref{seven-central}
and~\ref{seven-shear} at seven through nine sides. Then
\[
 F(R_n)-F(P)\ \ge\ \frac{F(R_n)}{n}
  \sum_i\left[\frac{(x_i-1)^2}{20000}+\frac{z_i^2}{1000}\right].
\]
In particular these inequalities hold throughout a neighborhood of
the regular polygon in the space of its ordered vertices.
\end{lemma}
\begin{proof}
Set $d=x-1$, $u=z^2$ and $D=T-H_1$ on one sector. We prove
$D\ge d^2/20000+u/1000$.

At five and six sides, use the upper bound on $H_1$ obtained by
increasing $\tau$ in Section~\ref{I:sec:five-six-thick}; write $D_*$ for the
resulting lower bound on $D$. It is concave in $u$. At $u=0$ the
zero-skew estimate there gives $D_*\ge .0155d^2/M$,
where $M=1-v+vx^2\le9$. At the upper edge
$u=z_*^2$ its lower bounds are $.00939$ and $.000184$, respectively.
These exceed the required lower bound even at the largest $|d|$ on each
rectangle, since
\[
 2^2/20000+.41^2/1000<.00939,\qquad
 1.32^2/20000+.302^2/1000<.000184.
\]
Interpolation between the two edges proves the sector inequality.

At seven through nine sides write $H_0=H_1(x,0)$, and use $m_0,e_0$
for the normalized mass and energy there. The estimate in
Section~\ref{seven-central} gives $T-H_0\ge .00012d^2$ throughout
all three central windows. To bound the additional deficit due to skew, use
the ratio in Section~\ref{seven-shear}. Let $A(u)$ be its numerator,
$U_s$ the largest allowed $u$, and take
$\gamma=.256,.256,.52$ in order. The proved lower bound
$m_1/m_0\ge\gamma$ bounds the denominator from below by
$D_\gamma(u)=(1-u)^{3/2}\sqrt{1+\gamma u}$. Rationalizing the
square root as in that section gives
\[
  D_\gamma(u)-A(u)\ge\frac{u}{1+e_0}
  \left[(1+e_0)(1-U_s)^{3/2}
  \frac{\gamma}{1+\sqrt{1+\gamma U_s}}-\mathcal F\right].
\]
The bracket is greater than $.23-.02$, $.23-.14$, and
$.494-.482$, respectively. Since $A(u)>0$, division now gives
\[
 H_0-H_1(x,u)\ge H_0\!\left(1-\sqrt{H_1(x,u)/H_0}\right)
 \ge \frac{x(1+e_0)}{4m_0}\,\frac{.012u}{D_\gamma(u)}
 \ge \frac{u}{1000}.
\]
Here $x\ge.67$, $m_0<1.205$ in all three rows,
$1+e_0\ge1$ and $D_\gamma(u)\le\sqrt{1+.52(.16)^2}<1.01$.
Adding the lower bound at zero skew proves the claim. Summing it in the
centroid variational inequality of Section~\ref{I:sec:fan} with
$R=\mu_1(R_n)$ gives the stated deficit. The central rectangles
contain $(x,z)=(1,0)$ in their relative interior, and the centroid,
radii and apertures vary continuously with ordered vertices near
$R_n$.
\end{proof}

\begin{lemma}[Distance from a regular polygon in centroid coordinates]
\label{lem:centroid-to-area-distance}
Let $5\le n\le9$, and let the centroid fan of a convex $n$-gon have
parameters $x_i,z_i$ as in Section~\ref{I:sec:fan}. Set
$D^2=\sum_i[(x_i-1)^2+20z_i^2]$. Then
\[
                       \mathcal A_n(P)\le 21\sqrt n\,D .
\]
\end{lemma}
\begin{proof}
If $\sqrt nD\ge1/10$, use $\mathcal A_n\le2$.
Otherwise translate to the centroid, and scale and rotate so that the first
vertex is $(1,0)$. Let $v_j=r_j e^{i\theta_j}$ be the remaining vertices
in positive cyclic order and $u_j=e^{2\pi i(j-1)/n}$ the regular ones.
Since $r_{j+1}/r_j=(1+z_j)/(1-z_j)$ and $|z_j|<1/10$,
$|\log r_j|\le3\sum_i|z_i|\le(3/\sqrt{20})\sqrt nD<.672\sqrt nD$
when $D>0$. Since $2k=2\tan(\pi/n)\le2\tan(\pi/5)<1.46$,
\[
 |\theta_j-2\pi(j-1)/n|
 \le2\sum_i|\arctan(kx_i)-\arctan k|\le1.46\sqrt nD .
\]
Using $|e^a-1|\le |a|e^{|a|}$ and $.672e^{.0672}+1.46<2.18$, we get
$\max_j|v_j-u_j|\le2.18\sqrt nD=:\delta$. Taking convex hulls bounds
the difference of the two support functions by $\delta$. The support of
the circumradius-one regular polygon $R$ is at least
$c=\cos(\pi/n)>4/5$ in every unit direction. Hence
$(1-t)R\subset P\subset(1+t)R$ with $t=\delta/c<2.75\sqrt nD$ when
$D>0$. The same inclusions hold with $P$ replaced by the centered homothet
of $R$ having area $|P|$. Their symmetric difference therefore has
relative area at most $4t/(1-t)^2\le21\sqrt nD$, since $t<.275$ and $11/(1-.275)^2<21$.
The case $D=0$ is the regular polygon.
\end{proof}

\begin{proposition}[Stability at five through nine sides]
\label{prop:centroid-global-stability}
Every simple polygon $P$ with at most $n$ sides, $5\le n\le9$, satisfies
\[
 F(R_n)-F(P)\ \ge\ \frac{1}{10^6 n^2}\,\mathcal A_n(P)^2 .
\]
\end{proposition}
\begin{proof}
In the geometric window of Proposition~\ref{prop:quantitative-window-entry},
suppose first that all centroid sectors are thick or central. Apply
Lemmas~\ref{lem:central-centroid-coercivity} and~\ref{lem:centroid-to-area-distance}.
They give the coefficient $F(R_n)/(20000\cdot21^2n^2)$, which exceeds
the stated one since $F(R_n)\ge F(R_4)=\pi^2$ by Theorem~\ref{thm:main}.

If one sector is narrow, use the polynomial as the reference field for the entire fan
of Sections~\ref{paper-small-sum} and~\ref{paper-seven-sum}. With their
notation $\mu_1(P)\le\bar R_n(1-b_n/n)$ in circumradius-one reference
units. Thus the deficit in $F$ is bounded below by
\[
 \frac{|R_n|\bar R_n}{n}
 \left[b_n-n(1-L_n/\bar R_n)\right] .
\]
Each of the five bracketed margins in the two cited tables exceeds
$1/5000$. The prefactor exceeds one because
$|R_n|\bar R_n\ge F(R_n)\ge\pi^2>9$.
Outside the geometric window the lower bounds of
Proposition~\ref{prop:quantitative-window-entry} apply. Both this case
and the narrow-sector case give the asserted coefficient using
$\mathcal A_n\le2$.
\end{proof}

\begin{proposition}[Explicit stability at three and four sides]
\label{prop:tri-square-stability}
For every admissible polygon $P$,
\[
 F(R_3)-F(P)\ge\frac{F(R_3)}{144}\mathcal A_3(P)^2
 \quad(n=3),\qquad
 F(R_4)-F(P)\ge\frac1{60000}\mathcal A_4(P)^2
 \quad(n=4).
\]
\end{proposition}
\begin{proof}
An area-matched triangle is an affine image of $R_3$. After isometries
its linear part is positive definite with eigenvalues $s,1/s$, $s\ge1$.
The three unit side directions $v_j$ of the equilateral triangle satisfy
$\sum_jv_jv_j^\top=3I/2$. Thus for this positive matrix $T$,
\[
 \frac{|\partial P|}{|\partial R_3|}
 =\frac13\sum_j|Tv_j|\ge\frac13\sum_jv_j^\top Tv_j
 =\frac{s+s^{-1}}2 .
\]
By the triangle perimeter bound in Section~\ref{I:sec:lowertri}, the
relative spectral deficit $1-F(P)/F(R_3)$ is at least
$w^2$, where $w=(s^2-1)/(s^2+1)$. Give $R_3$ circumradius one. Each
point of $R_3$ moves at most $s-1$ under $T$. Its concentric incircle
has radius $1/2$, so $TR_3\subset(1+2(s-1))R_3$. Equal areas now
give $\mathcal A_3\le\min(2,8s(s-1))\le12w$: for $1<s\le5/4$ use
$8s(s^2+1)/(s+1)\le 8(5/9)(41/16)<12$, and above $5/4$ use
$w\ge9/41$. The undistorted triangle gives equality.

For the square set $p=\pi^2$ and $\delta=p-F(P)$. If $\delta\ge.01$,
the bound follows from $\mathcal A_4\le2$. Otherwise Section~\ref{I:sec:lowertri}
and the bound $F^2<97$ for $b+c\ge.85$ in Section~\ref{I:sec:all-real-quads}
reduce the problem to a strictly convex quadrilateral in the coordinates
of the latter section, with $r=b^2+c^2<S=.7225$. The form
$B_\delta=B+\delta\operatorname{Cov}$ is positive semidefinite, by the
Rayleigh principle at $F(P)$. The linear tests give $M\succeq(p-\delta)\Sigma_Q$.
Since $p-\delta>9.85$, the proof of \eqref{I:eq:quad-t-affine} still gives
$t\le .201+.405r<.494$. The trace of the covariance matrix on the
sines $U,V$ is $1-I^2-J_0^2\le1$. Thus the least-eigenvalue test on
these sines now gives
\[
 (C-t/4)t\le 16r/p^2+\delta/p,\qquad t\le .8r+.5\delta,
\]
using $C-t/4>.207$.

We show $r\le200\delta$. If $r>0$ and $\delta\le r/200$, then
$t<.803r$ and
$t^2/r\le.201(.803)+.405(.494)<.362$.
Write $P_\delta=[bB_\delta(U,H)-cB_\delta(V,H)]/r$ with $H=xy$.
By the covariance Cauchy--Schwarz inequality,
$|P_\delta-P_0|\le\delta\sqrt{\operatorname{Var}(H)/r}$.
The displayed formula for $\operatorname{Var}(H)$ in
Section~\ref{I:sec:all-real-quads} gives $\operatorname{Var}(H)<.017$.
In \eqref{I:eq:quad-mixed-lower} the coefficient of the total loss in
$rw$ is at most $.01+.21(.803)+.21(.803)S/8<.21$, using
$w=2bc/r\le1$ and $rw<S/2$. Since $r(1+w)<S$, this yields
\[
 P_\delta>.82-.21S-.005\sqrt{.017 S}>2/3 .
\]
On the other hand,
\[
 \begin{aligned}
 B_\delta(H,H)&\le .099+.494^2(2+S)/24
                   +S(.494/12-.039)+.005S(.017)<.129,\\
 P_\delta^2&\le B_\delta(H,H)
       \frac{B_\delta(U,U)+B_\delta(V,V)}{r}
       <.129\bigl(1.625+.005+5(.362)\bigr)<4/9 .
 \end{aligned}
\]
The second line is the Schur test \eqref{I:eq:quad-schur} for $B_\delta$.
This contradiction proves the claim; $r=0$ satisfies it directly.

In area-one coordinates, align the square with the orthogonal polar
factor of the linear map $T$ that takes the normalized quadrilateral $Q$ to $P$. Its largest singular value is
$\sqrt{\sqrt{1+t^2}+t}\le1+t$. Its smallest is the reciprocal, so
the distance to that polar factor in operator norm is at most $t$.
In the affine model the vertex shifts have length $\sqrt r/2$.
Since $t<.494$, the four vertices of $P$ are therefore at distance at most
$d=t+\sqrt r$ from the corresponding vertices of the aligned square. By convexity $P$
lies in the square's parallel set of width $d$, whose area is $1+4d+\pi d^2$.
Thus $\mathcal A_4\le\min(2,8d+2\pi d^2)\le10d$ by splitting at $d=1/4$.
Using $t^2\le t/2$ and the two parameter bounds above gives
\[
 \mathcal A_4^2\le200(r+t^2)\le200(1.4r+.25\delta)\le56050\delta.
\]
\end{proof}

\subsection{The projected trial}
\begin{proposition}[Quantitative control of both endpoint radii]
\label{prop:coercive-physical-fan}
Let $n\ge10$, $h=2\pi/n$, and normalize $|P|=|R_n|=A_0$ with
circumradius one for $R_n$. Suppose $P$ admits a fan, centered at a point where the projected trial field has mean zero,
that satisfies the geometric bounds of Proposition~\ref{prop:common-projected-comparison}.
Then
\begin{equation}\label{eq:coercive-physical-fan}
 F(R_n)-F(P)\ \ge\ \frac1{25}\sum_i
 \left[\frac{h^4}{2000}(q_i-1)^2+
              \frac{h}{10000}(e_i^2+w_i^2)\right].
\end{equation}
\end{proposition}
\begin{proof}
Write $d=h^2$, $s=|e|/d$, $y=w/d$ and $t=|q^2-1|$ for one sector.
In all cases the geometric bounds give $s<.163$, $|y|<.227$ and $h<.629$. When
$|q-1|\ge.1$, half the lower bound for the angular gain left after subtracting the radial error in the proof of
Proposition~\ref{prop:common-projected-comparison} exceeds the required radial
term, since
\[
 \frac{.629(.163^2+.227^2)}{10000}<\frac{(.1)^2}{2000}.
\]
It suffices to retain a small part of each curvature on $.9\le q\le1.1$.
Put $\eta=1/2500$. In \eqref{mid-radial-identity} replace both
$K_e,K_w$ by $K_e-\eta,K_w-\eta$. We show that the resulting radial term
still satisfies the upper bound in Proposition~\ref{IV:prop:large-radius-profile}.
This leaves $\eta qh(e^2+w^2)/2$ in addition to the original lower bound on the net sector gain.

With the notation $G,\delta,g,A,B,j,X$ of the radial estimate, the
two modified curvature coefficients are at least
$G-\eta-gs\mathbf1_{e>0}$ and $\delta-\eta-As\mathbf1_{e>0}$.
For $q\le1$, use the row for the last aperture interval below one. If $e\le0$, only
the mixed term $jst$ remains. It is bounded above by $kt^2+(G-\eta)s^2/2$.
Indeed $j^2/[2(G-\eta)]<.0121,.0108$ at ten and eleven sides.
For twelve or more sides, the endpoint estimates $2Gk-j^2>.0008$ remain
positive after deducting $2k\eta$. The coefficient of the squared skew is positive in
all three ranges. If $e\ge0$, the mixed term is nonpositive. Here
$G-gX>1.3$ by the estimates for narrow sectors (for twelve or more sides use
$2-.290-9(.524)^2(.160)>1.3$).
The inequality $A s y^2\le (1.3-\eta)s^2+
(AB)^2 y^2/[4(1.3-\eta)]$ bounds the term involving both radial displacement and skew, since
\[
 \delta-\eta>\frac{(AB)^2}{4(1.3-\eta)}.
\]
For this substitution the triples $(\delta,A,B)$, rounded down in
the first position and up in the other two, are
$(.016,1.25,.227)$, $(.015,1.25,.209)$ and $(.019,.841,.202)$.

For $1\le q\le1.1$ use the first row for wide sectors, so $t\le.21$.
Both modified curvature coefficients are positive when $e\le0$, and the mixed
term has favorable sign. When $e\ge0$, put $r=.062$, $L=G-\eta$
and $b_0=\delta-\eta$. For $s\le r$ the two sufficient comparisons
are $b_0>Ar$ and $j^2/[2(L-gr)]<k$, by square completion.
For $s\ge r$, use $(As-b_0)_+\le A^2s^2/(4b_0)$.
The upper bound for the radial term is then nonpositive if
$[L-gX-(AB)^2/(4b_0)]r>2j(.21)$.
The following bounds are obtained by substitution in these three
expressions, in that order:
\[
\begin{array}{c|rrr}
 &b_0-Ar>&j^2/[2(L-gr)]<&[L-gX-(AB)^2/(4b_0)]r-.42j>\\\hline
 10&.0013&.0124&.0024\\
 11&.038&.0113&.020\\
 n\ge12,\ d=0&.063&.0082&.036\\
 n\ge12,\ d=.524^2&.051&.0112&.018
\end{array}
\]
The last two rows allow $j=1/6+.173d$ by the first wide row of the
mixed-slope estimate. They use $G=1.713+1.825d$, $\delta=.064+.147d$,
$g=9d$, $A=3.06d$, $B=.202$ and $X=.160$.
For fixed $s,y,t$ and sign of $e$, the modified Taylor majorant
is affine in $d$, so the same upper bound holds for every intermediate $d$. Since $q\ge.9$, the reserved
coefficient $q\eta/2$ exceeds $1/10000$.

Summing the lower bounds over the sectors as in
Lemma~\ref{IV:lem:projected-fan-identities} bounds the negative
trial residual from below. Lemma~\ref{IV:thm:radial} gives
$|U|\le1$. Orthogonal projection gives $nM_D\le\int_{R_n}|U|^2\le A_0$.
Adding the norms of the diagonal and off-diagonal blocks of the sector
matrix gives
\[
 \|T\|\le1.064\max\{1/\cos(h/2),1.413/\sinc(h/2)\}
              +.227h/\sinc(h/2)<2 .
\]
Thus $|\det T|\le4$. Changing variables on each sector gives
$\|V\|_{L^2(P)}\le\sqrt{A_0}+4\sqrt{nM_D}\le5\sqrt{A_0}$.
The Rayleigh inequality for the zero-mean components
now gives \eqref{eq:coercive-physical-fan}.
\end{proof}

\begin{theorem}[Global quadratic stability from ten sides]
\label{thm:projected-global-stability}
For every $n\ge10$ and every simple polygon $P$ with at most $n$ sides,
\[
 F(R_n)-F(P)\ \ge\ c_n\,\mathcal A_n(P)^2,\qquad
 c_n=\min\left\{\frac{1}{5n^5},\frac{h}{3750000}\right\},\quad h=\frac{2\pi}{n}.
\]
In particular $c_n\ge1/(60n^5)$.
\end{theorem}
\begin{proof}
One always has $\mathcal A_n(P)\le2$. Thus a deficit exceeding
$4/(5n^5)$ suffices. Otherwise Proposition~\ref{prop:quantitative-window-entry}
provides a fan for which the trial has mean zero and the geometric bounds hold. Normalize the regular
circumradius to one. Put $S=\sum_i(e_i^2+w_i^2)$ and
$T=\sum_i(q_i-1)^2$. Proposition~\ref{prop:coercive-physical-fan} gives the lower bound
$hS/250000+h^4T/50000$ for the spectral deficit.

Place the regular polygon at the fan center and align its first vertex
with the first vertex ray of $P$. If $\theta_k$ is the difference of the
two directions after $k$ edges, then
\[
 |\theta_k|=h\left|\sum_{i=1}^n
 (\mathbf 1_{i\le k}-k/n)(q_i-1)\right|
 \le h\sqrt{k(1-k/n)}\,\sqrt T\le h\sqrt{nT}/2 .
\]
Here $\sum_i(q_i-1)=0$. Each endpoint radius differs from one by
$e_i\pm w_i/2$, of absolute value at most $(\sqrt5/2)\sqrt S$.
The piecewise linear angular map from the rays of $R_n$ to the rays through the vertices of $P$
has displacement at most $h\sqrt{nT}/2$ everywhere by interpolation.
For aperture $a$ and fractional angle $u\in[0,1]$ the unit-endpoint
chord has radius $f(a,u)=\cos(a/2)/\cos((u-1/2)a)$.
For endpoints $r_-,r_+$ the reciprocal radius is
\[
 \frac{\sin((1-u)a)/r_-+\sin(ua)/r_+}{\sin a}.
\]
Thus the radial function of $P$ differs from $f(a,u)$ by at most
$(\sqrt5/2)\sqrt S$, since $f\le1$. On $a\le1.413h$,
$|\partial_a f|\le\tan(a/2)<.8h$ (use $\tan x\le x/\cos x$).
The regular polygon's radial graph is globally Lipschitz in angle
with constant $\tan(h/2)<.53h$. First change the aperture to $h$
at fixed $u$, then undo the angular displacement. Since $n\ge10$,
the difference of the two radii on the same ray is at most
\[
 (\sqrt5/2)\sqrt S+
 (.8+.265\sqrt n)h^2\sqrt T
 \le(\sqrt5/2)\sqrt S+.52h^2\sqrt{nT}.
\]
The radial function of $P$ is at most $1+.064+.227h^2/2<1.12$.
Polar integration of the absolute difference of the squared radii,
with $|R_n|=\pi\sinc h>.9\pi$, therefore bounds $\mathcal A_n$
by $2.4$ times the last display. Squaring gives
$\mathcal A_n^2\le15S+3.2nh^4T$. The lower bounds on the deficit in terms of radial and angular displacement give the
local coefficient $\min(h/3750000,1/(160000n))=h/3750000$.
The complementary case gives $1/(5n^5)$ as above. Since $n\ge10$
and $2\pi>6.25$, the stated uniform bound follows.
\end{proof}

\begin{remark}[Dependence on the number of sides]
For all fans satisfying the geometric bounds, the stability coefficient obtained here is of order
$n^{-1}$. A comparison using only distances between corresponding
vertices loses the small angular derivative of the regular radial
graph. The order $n^{-5}$ in the global estimate enters through the
deficit threshold that ensures the geometric bounds hold and the estimate
$\mathcal A_n\le2$ when they fail. The optimal dependence on $n$ remains undetermined.
\end{remark}

\begin{proof}[Proof of Theorem~\ref{thm:stability}]
Proposition~\ref{prop:tri-square-stability} covers triangles and
quadrilaterals. Proposition~\ref{prop:centroid-global-stability} covers
five through nine sides. Theorem~\ref{thm:projected-global-stability}
gives the asserted constant for each $n\ge10$.
\end{proof}

\subsection{Optimality of the power of the asymmetry}

\begin{proposition}\label{prop:optimal-power}
For each fixed $n\ge4$ and each $p<2$, no positive constant works in
$D_n(P)\ge c\,\mathcal A_n(P)^p$ even in a neighborhood of $R_n$.
For triangles the optimal power is one.
\end{proposition}
\begin{proof}
Write $v_i$ for the vertices of $R_n$. Move them to $v_i+tV_i$,
then dilate about the origin to restore the area. For small $|t|$ the
polygons are strictly convex and all the vertex paths are analytic.
Restrict the $2n$ velocities by the single independent linear
constraint of vanishing area derivative (dilation has nonzero derivative).
The restoring dilation then has zero first derivative.

Triangulate from the origin and pull the mass and Dirichlet forms
$m_t,a_t$ back by the piecewise affine vertex map, keeping the origin
fixed. On each triangle its Jacobian matrix $J_t$ is analytic and
invertible near zero. The coefficients are $\det J_t$ and
$(\det J_t)J_t^{-1}J_t^{-\top}$, respectively. In particular the
compressed first variation
\[
 V\longmapsto (a'_0-\mu_1(R_n)m'_0)|_{E\times E}
\]
is linear, with $E$ the two-dimensional real first eigenspace
(Lemma~\ref{IV:lem:regular-character}). Its target has dimension
three. Thus its kernel in the area tangent space has dimension at least
$2n-4>3$ for $n\ge4$. Choose a velocity in this kernel outside the
three-dimensional space of infinitesimal rigid motions.

Use $a_0+m_0$ as the fixed Hilbert inner product on $H^1$.
If $B_t,M_t$ represent $a_t+m_t,m_t$, the compact self-adjoint family
$B_t^{-1/2}M_tB_t^{-1/2}$ is analytic near zero by the binomial series
($B_0=I$). Its isolated cluster at $(1+\mu_1)^{-1}$ has an analytic
two-dimensional spectral projection by the resolvent formula.
Project a basis and orthonormalize to obtain an analytic basis, then
transport it back by $B_t^{-1/2}$. In this basis the restriction of
$a_t-\mu_1m_t$ has zero constant and linear terms: basis derivatives
pair to zero with $E$ by the weak eigenvalue equation. The restricted
mass is uniformly positive, so both eigenvalues move by $O(t^2)$.
The spectral gaps on either side persist. Consequently $D_n(P_t)=O(t^2)$.

It remains to bound the symmetric difference below along the chosen perturbation. For a vertex
displacement $W$ let $L_i(W)$ be its affine interpolation on edge $i$
and put $N(W)=\sum_i\int_i|L_i(W)\cdot\nu_i|\,ds$.
If the integral on every edge is zero, the interpolated normal velocity vanishes along each edge. Each vertex velocity is then tangent to both incident edges and must vanish. Thus the formula defines a norm. For two polygons whose vertices differ from those
of $R_n$ by at most $\eta$, their symmetric difference equals $N$ of the
difference of their displacements, up to $O_n(\eta^2)$. Indeed discard
disks of radius $C_n\eta$ about each vertex. In the remaining edge strips
each displaced line is a normal graph; its height is the interpolated
normal displacement up to $O_n(\eta^2)$, since the tangential shift
and the slope are both $O_n(\eta)$. The discarded pieces and the
omitted end strips have area $O_n(\eta^2)$.

A best aligned regular copy exists by compactness of rotations and
boundedness of translations for copies with positive overlap. As
$t\to0$ these alignments converge to the identity modulo the finite
symmetries (along subsequences). The norm formula applied to the
regular copy and $R_n$ bounds its translation and rotation angle by
$O(|t|)$: their linearized rigid vertex displacement is injective and
their area difference is $O(|t|)$ by optimality and the triangle
inequality. Dividing the strip formula by $|t|$ now bounds the lower
limit of the optimal relative difference from below by
$\inf_{a,\omega}N(V_i-a-\omega v_i^\perp)/|R_n|>0$. This proves the claim for $n\ge4$. The same strip estimate applies to an area-preserving affine path from an equilateral triangle whose initial vertex velocities do not come from an infinitesimal rigid motion. Along such a path the asymmetry is bounded above and below by positive multiples of the absolute value of the path parameter.

For $n=3$ take the area-preserving affine map from the equilateral
triangle, with singular values $s,1/s$, $s\ge1$. Its energy matrix
has eigenvalues $m\pm t$, $t=(s^2-s^{-2})/2$, $m=\sqrt{1+t^2}$.
For the first mode $u_1$ of the unit-side triangle the tensor
$\int\nabla u_1\nabla u_1^\top$ is diagonal. The difference of its
diagonals divided by their sum is $\beta=243/(32\pi^2)$ by the
integrals in \cite[Section~4.1]{LS}. Rotating this trial by multiples
of $2\pi/3$ gives three traceless tensor directions at equal angular
spacing. At least one has normalized scalar projection at most $-1/2$
onto the metric distortion. These transported trials have zero mean
because the Jacobian is constant. Hence $F(P)/F(R_3)\le m-\beta t/2$.
For $t\le1/4$, $m-1\le t^2/2$ gives $D_3/F(R_3)\ge t/4$.
In the notation of Proposition~\ref{prop:tri-square-stability},
$w=t/(m+1)$ and $\mathcal A_3\le12w$. For $t\ge1/4$, the same
proposition's perimeter estimate gives $D_3/F(R_3)\ge w^2>1/100$.
Together these bounds imply $D_3(P)\ge F(R_3)\mathcal A_3(P)/200$.
Finally the affine metric is bounded below by $(m-t)I$, so $D_3=O(t)$
along a small nonzero pure stretch. The strip argument gives
$\mathcal A_3\asymp t$ there, excluding any power smaller than one.
\end{proof}

\clearpage\appendix\section{Numerical assertions and margins}\label{sec:verification}
The two tables follow the proof from the bounds for the regular polygon to the
strict inequality for the competitor. For quantitative stability, the
deficit thresholds are in Proposition~\ref{prop:quantitative-window-entry}
and the constants in Theorem~\ref{thm:stability}. A finite decimal denotes the exact
rational number written. A tuple is ordered by increasing side count
unless its two entries are identified as aperture halves. The last column
gives the margin in the indicated comparison. For $a>b>0$ a percentage
uses $(a-b)/a$ unless the table specifies a different denominator. A comparison
with zero is reported in its displayed units. A dash marks an input bound for which no single margin is given.

\begingroup\small
\setlength{\tabcolsep}{4pt}
\begin{longtable}{@{}>{\raggedright\arraybackslash}p{.135\linewidth}>{\raggedright\arraybackslash}p{.263\linewidth}>{\raggedright\arraybackslash}p{.342\linewidth}>{\raggedright\arraybackslash}p{.18\linewidth}@{}}
\caption{Finite sides and missing convex vertices. The lower eigenvalue bounds $L_n$ and upper Rayleigh quotient bounds $\bar R_n$ use circumradius one. The bounds for $F=|P|\mu_1(P)$ are scale invariant.}\label{tab:centroid-assertions}\\
\toprule
Case & Assertion entering the proof & How to obtain the bound & Margin\\\midrule\endfirsthead
\toprule Case & Assertion entering the proof & How to obtain the bound & Margin\\\midrule\endhead
\bottomrule\endfoot
Triangles and quadrilaterals &
$F(R_3)=4\pi^2/(3\sqrt3)$, $F(R_4)=\pi^2$. A reflex quadrilateral has $F(P)<\sqrt{89}$. &
Use the perimeter bound and Heron's formula in Section~\ref{I:sec:lowertri}. For a convex quadrilateral with positive squared offset $r$, equation~\eqref{I:eq:quad-schur} contradicts $.129(1.625+5\cdot.361)<(2/3)^2$. & $4/9-.129(1.625+5\cdot.361)>.00197$.\\\addlinespace
Missing corner, $n=5$ &
$F(P)<56\pi/17$ and $F(R_5)>10.3526$. &
Integrate the three pieces of the circle-mean inequality with the primitive $U_j$ in Section~\ref{sec:five-six-missing}. Lemma~\ref{lem:degree-six-pentagon-floor} gives the lower bound for the regular pentagon. & $10.3526-56\pi/17>.0038$.\\\addlinespace
Missing corner, $n=6$ &
$d_P>2/57$, $d_R<(4.809+98/36)/216$; hence $F(P)<10.4959$ and $F(R_6)>10.4967$. &
Compare Fourier order one, then use a quadratic discriminant for all orders $k\ge2$. Apply the Schur bound in Lemma~\ref{small-schur-regular} and solve for the normalized eigenvalue as in Theorem~\ref{IV:thm:corner-gap}. & $>.0008$ in $F$ (over $.0076\%$ of the lower bound).\\\addlinespace
Missing corner, $n\ge7$ &
$m^3J_m>(4.832,4.804,4.792)$ for $m=6,7,8$, respectively. The gap in $F$ exceeds $5/n^4$. &
Use the five retained terms and the tail in Lemma~\ref{IV:lem:positive-turn-kernel}, the upper bound \eqref{IV:eq:regular-continuous-band} for $d_R$, and Theorem~\ref{IV:thm:corner-gap}. & $>5/n^4$ in $F$.\\\addlinespace
Regular data, $n=5,6$ &
$(L_5,L_6)=(4.35415,4.04315)$; the two trial quotients are below $4.355245$, $4.043369$. &
Apply the area and chord monomial pairings in Section~\ref{I:sec:five-six-regular-data}; the six-sided field here has degree nineteen. & --- (reference bounds).\\\addlinespace
Thick centroid sectors, $n=5,6$ &
$H_1\le T$ for $x\ge .65,.77$ in the respective windows. At the lower aperture and maximal skew the lower bounds for the deficit exceed $.00939,.000184$. &
Section~\ref{I:sec:five-six-thick} reduces the estimate to zero and maximal skew. At maximal skew substitute $B_-,C_+,v_-$ and bound the arctangent by its alternating series. & At $n=6$, $>.000184$ in $T-H_1$ ($.0229\%$ of $.802899$).\\\addlinespace
Narrow centroid sectors, $n=5,6$ &
A single narrow sector contributes at least $b_n=(.005,.0017)$ to the deficit. After subtracting the spectral error for the entire fan, the margins exceed $.003742,.001375$. &
Section~\ref{paper-small-sum} compares $b_n$ with $n(1-L_n/\bar R_n)$ using each field's own tangent. & $>.003742$, $.001375$ after subtracting the spectral error.\\\addlinespace
Reference data, $n=7,8,9$ &
$U_n=(92.5,85.25,80.75)$; $L_n=(3.86037,3.74649,3.66960)$; $\bar R_n=(112/29,3.74722,3.67006)$. &
Integrate the twelve monomial pairings per side count in Section~\ref{seven-moments}, then insert the constants $U_n$ and the lower bound $G_n$ for the second eigenvalue in rotation character one, as in Section~\ref{seven-spectral}. & Unused $U_n$: $>.21,.09,.03$ with the strict upper endpoints in the proof.\\\addlinespace
Centroid window, $n=7,8,9$ &
The upper bounds on $x$ are $1.78,1.64,1.59$; those on $|z|/x$ are $.230,.129,.085$. &
In Section~\ref{seven-window-proof} apply the logarithmic bound on the aperture, then transfer the fan to the centroid using the differential bounds in Section~\ref{seven-transport}. & --- (simultaneous bounds on the coordinates).\\\addlinespace
Central centroid sectors, $n=7,8,9$ &
The first-field transverse fractions exceed $.07285,.05491,.04295$. At zero skew on $.67\le x\le1$ the quotient $s=(T-H_1)/(x-1)^2$, continued at $x=1$, exceeds $.00263,.00118,.00012$. &
Section~\ref{seven-central} retains positive Fourier terms, then checks $p_{s_0}(t,E)$ at four corners. At $n=9$ the integer parts of $10^9p_{s_0}$ are $1581,1566,643,628$. The right half and the shear terms are treated there and in Section~\ref{seven-shear}. & At $n=9$, $.00046-.33(.00072)-.33^2(.0008)=.00013528>.00012$ for $s$.\\\addlinespace
Narrow centroid sectors, $n=7,8,9$ &
For $x<.67$ one sector contributes at least $b_n=(.00368,.00219,.00138)$ to the deficit. After subtracting the spectral error for the entire fan, the margins exceed $.000600,.000631,.000251$. &
Use both skew endpoints and the minimum of the cubic in Section~\ref{sec:seven-narrow}. Subtract $n(1-L_n/\bar R_n)$ using the exact fractions in Section~\ref{paper-seven-sum}. & $>.000600$, $.000631$, $.000251$ after subtracting the spectral error.\\
\addlinespace
Quadrilateral stability &
$P_\delta>.6677>2/3$ when $0<r$ and $\delta\le r/200$. &
In Proposition~\ref{prop:tri-square-stability} put $S=.7225$ in
$.82-.21S-.005\sqrt{.017S}$. The Schur upper bound on $P_\delta^2$ is
$.129(1.625+.005+5\cdot.362)<4/9$. &
$>.00103$ ($>.154\%$ of $.6677$).\\
\end{longtable}
\endgroup

\begingroup\small
\setlength{\tabcolsep}{4pt}
\begin{longtable}{@{}>{\raggedright\arraybackslash}p{.135\linewidth}>{\raggedright\arraybackslash}p{.263\linewidth}>{\raggedright\arraybackslash}p{.342\linewidth}>{\raggedright\arraybackslash}p{.18\linewidth}@{}}
\caption{The projected trial. Here $h=2\pi/n$, $q$ is the normalized aperture, $e,w$ are the mean and skew displacements, and $D_q=qh^5(q^2-1)^2$. In a pair indexed by aperture, the order is $q<1$, $q>1$.}\label{tab:projected-assertions}\\
\toprule Case & Assertion entering the proof & How to obtain the bound & Margin\\\midrule\endfirsthead
\toprule Case & Assertion entering the proof & How to obtain the bound & Margin\\\midrule\endhead
\bottomrule\endfoot
Sector angles and radii &
$.426<q<1.413$ at ten sides and $.463<q<1.391$ at eleven. For $n\ge12$, $.475<q<1.403$, $|e|<.160h^2$, $|w|<.202h^2$. &
Use the tangent bounds for the joint ellipse in Section~\ref{mid-geometry} and Proposition~\ref{IV:prop:entry-twentyfour}. Equation~\eqref{paper-geometric-window} and the estimates following it give the bounds on the skew for each aperture range. & For $n\ge12$, $79.3-79.28=.02$ in $C$.\\\addlinespace
Bounds for the correction $D$, $n\ge12$ &
$s_D>\zeta(3)/(4\pi^3)+h/200$, $T_2\le\mathcal T_H<1/2$, $m_*<.048$, $\rho<.09$, $z<.015$, $T_1<.65$. &
Insert the upper bound~\eqref{regular-common-ceiling} in Proposition~\ref{continuous-response}. Keep the exact endpoint expression $\mathcal T_H$ at $H=\pi/6$; the squared cubic projection exceeds $.0072h$. & $1/2-\mathcal T_H>.00399$ (over $.798\%$ of the upper bound).\\\addlinespace
Boundary flux and mixed slope &
For $n\ge12$, $W_e<.255$, $W_w<.014$. On the two aperture halves $J$ is at most $(.221,.250)$ at ten sides and $(.207,.230)$ at eleven. &
Use the Bessel recurrence in Lemma~\ref{IV:lem:signed-wall-shape} at the interval ends. Proposition~\ref{IV:prop:middle-aperture-radius} subtracts the two orthogonal projections on the boundary before bounding $J$. & --- (curvature and slope inputs).\\\addlinespace
Ten sides &
$\Delta/D_q\ge(.0172,.0148)$ and $C/D_q\le(.0150,.0127)$ on the respective halves. The smaller surplus is $.0021$. &
Use a chord on each half in Proposition~\ref{IV:thm:uniform-twenty-plus}. The narrow-band square completions and wide-band endpoint estimate are in Proposition~\ref{IV:prop:large-radius-profile}. & $\ge .0021$ ($14.1\%$ of $.0148$ on the wide half).\\\addlinespace
Eleven sides &
$\Delta/D_q\ge(.0181,.0155)$ and $C/D_q\le(.0130,.0120)$. The smaller surplus is $.0035$. &
Use the same two propositions with the eleven-sided window. Section~\ref{mid-summation} restores the factor $D_q$ before summing. & $\ge .0035$ ($22.5\%$ of $.0155$ on the wide half).\\\addlinespace
All $n\ge12$ &
$C\le .0123D_q$ and $\Delta\ge .0123D_q+.002h^4(q-1)^2$. &
Combine Propositions~\ref{prop:angular-from-twelve} and~\ref{IV:prop:large-radius-profile} on the same window. At $q=1$ use the strict radial inequality in Proposition~\ref{prop:common-projected-comparison}. & $\ge .002$ in units of $h^4(q-1)^2$, $q\ne1$.\\
\addlinespace
Radial estimate, $n=10$ &
$\delta-1/2500=.0156$ in the narrow row ending at $q=1$. &
In Proposition~\ref{prop:coercive-physical-fan} the coefficient needed to bound the mixed term is
$(1.25\cdot.227)^2/[4(1.3-1/2500)]<.01549$. &
$>.00011$ ($>.70\%$ of $.0156$).\\
\end{longtable}
\endgroup

\section*{AI use declaration}
The author used AI tools to suggest improvements and to help polish the exposition. The author is responsible for all content.

\end{document}